\documentclass[12pt]{amsart}
\usepackage{amsmath,amsfonts,amssymb,amscd}
\usepackage{latexsym}

\usepackage[
hyperindex=true,pagebackref=true,bookmarks=true,
colorlinks=true,linkcolor=blue,citecolor=red]
{hyperref}
\def\~{{\rm --}} 

\title [Jones polynomials of torus knots via DAHA]
{Jones polynomials of torus knots via DAHA} 
\author[Ivan Cherednik]{Ivan Cherednik $^\dag$}
\date{August 4, 2012}

\thanks{$^\dag$  \date\ \ \ Partially supported by NSF grant
DMS--1101535 and the Fulbright Program}

\address[I. Cherednik]{Department of Mathematics, UNC
Chapel Hill, North Carolina 27599, USA\\
chered@email.unc.edu}

 \def\bysame{{\bf --- }}
 \def\~{{\bf --}}
\newcommand{\comment}[1]{}
\renewcommand{\tilde}{\widetilde}
\renewcommand{\hat}{\widehat}

\newcommand{\HOM}{\hbox{\fontfamily{pcr}\fontseries{m}\fontshape{sl}
\selectfont H\!O\!M\!F\,}}

\renewcommand{\tilde}{\widetilde}
\renewcommand{\hat}{\widehat}

\newcommand{\Z}{{\mathbb Z}}
\newcommand{\Q}{{\mathbb Q}}
\newcommand{\N}{{\mathbb N}}
\newcommand{\C}{{\mathbb C}}
\newcommand{\R}{{\mathbb R}}

\def\HH{\mbox{${\mathcal H}$\kern-5.2pt${\mathcal H}$}}

\newtheorem{theorem}{Theorem}[section]

\newtheorem{proposition}[theorem]{Proposition}
\newtheorem{definition}[theorem]{Definition}
\newtheorem{lemma}[theorem]{Lemma}

\newtheorem{theorem }{Theorem}[section]
\newtheorem{maintheorem }[theorem]{Main Theorem}
\newtheorem{proposition }[theorem]{Proposition}
\newtheorem{definition }[theorem]{Definition}
\newtheorem{lemma }[theorem]{Lemma}
\newtheorem{corollary }[theorem]{Corollary}
\newtheorem{notation }[theorem]{Notation}
\newtheorem{remark }[theorem]{Remark}
\newtheorem{example }[theorem]{Example}

\newtheorem{ maintheorem }[theorem]{Main Theorem}
\newtheorem{ theorem}{Theorem}[section]
\newtheorem{ proposition}[theorem]{Proposition}
\newtheorem{ definition}[theorem]{Definition}
\newtheorem{ lemma}[theorem]{Lemma}
\newtheorem{ corollary}[theorem]{Corollary}
\newtheorem{ notation}[theorem]{Notation}
\newtheorem{ remark}[theorem]{Remark}
\newtheorem{ example}[theorem]{Example}

 \newcommand{\rmk}{{\bf Comment.\ }}

\def\for{\  \hbox{ for } \ }
\def\iif{ \ \hbox{ if } \ }

\def\where{\  \hbox{ where } \ }
\def\and{\  \hbox{ and } \ }
\def\and{\  \hbox{ and } \ }

\def\equal{\stackrel{\,\mathbf{def}}{= \kern-3pt =}}

\def\la{\lambda}

\def\om{\omega}

\def\al{\alpha}

\def\ga{\gamma}
\def\ep{\epsilon}

\def\de{\delta}

\def\ka{\kappa}

\def\si{\sigma}

\def\ze{\zeta}

\def\vph{\varphi}

\def\vth{{\vartheta}}

\def\tal{\tilde{\alpha}}

\def\tga{\tilde{\gamma}}
\def\tGa{\tilde{\Gamma}}

\def\tw{\widetilde w}
\def\tW{\widetilde W}

\def\tz{\tilde z}
\def\tb{\tilde b}

\def\tR{\tilde R}

\def\hw{\widehat{w}}
\def\hW{\widehat{W}}

\def\hv{\hat{v}}

\def\S{\mathbf{S}}

\def\0{\mathbf{0}}

\def\f{\mathcal{F}}
\def\çF{\mathcal{F}}

\def\k{\mathcal{K}}

\def\h{\mathcal{H}}

\def\v{\mathcal{V}}

\def\x{\mathcal{X}}

\def\j{\mathcal{J}}

\def\lan{\langle}

\def\ran{\rangle}

 \def\dim{{\hbox{\rm dim}}_{\mathbb C}\,}
\def\lng{\hbox{\rm{\tiny lng}}}
\def\sht{\hbox{\rm{\tiny sht}}}

\newcommand{\tr}{\operatorname{tr}}

\newcommand{\sq}{\phantom{1}\hfill$\qed$}

\newcommand{\phk}{\phi^{(k)}}

\newcommand{\sgn}{\mbox{sgn}}

\def\HH{\mathfrak{H}}

\def\HH{\hbox{${\mathcal H}$\kern-5.2pt${\mathcal H}$}}

\font\smm=msbm10 at 12pt 
\def\symbol#1{\hbox{\smm #1}}
\def\lsmash{{\symbol n}}

\def\#{\sharp}

\begin{document}
\par
{\centering 
Dedicated to  Yuri Ivanovich Manin \\
on the occasion of his 75th birthday
\medskip
\par} 

\maketitle

\renewcommand{\baselinestretch}{1.2}

{
\renewcommand{\baselinestretch}{2.} 
\tableofcontents
\renewcommand{\baselinestretch}{2.} 
}
\renewcommand{\baselinestretch}{1.2} 

\vfill\eject

\renewcommand{\natural}{\wr}

\setcounter{section}{-1}
\setcounter{equation}{0}
\section{Introduction}

This work is mainly inspired by paper \cite{AS}, where
a construction was presented for a   
$q,t$\~version of the Jones polynomials of torus knots
and the corresponding {\em super-polynomials} in terms 
of the {\em generalized Verlinde algebras}. The latter algebras
are symmetric parts of perfect DAHA modules at roots of unity 
$q$ with $t=q^k$ for proper integral $k$; see \cite{C101}, Section 0.4.
The approach of \cite{AS} is based on the relation of the Jones 
polynomials  to the ``usual" Verlinde algebra, i.e., that 
defined for $t=q$ (describing integrable Kac-Moody representations).

An important but not quite rigorous component of this construction 
is the switch from the roots of unity, inevitable in the 
Verlinde algebras, to generic $q$ from the 
Jones and HOMFLYPT polynomials. Also, the universality of the 
approach of \cite{AS} remained unclear; only $A_n$ and mainly
$\om_1$  were considered there (and the simplest knots).

\medskip
We suggest a conjectural DAHA interpretation 
of the Quantum Group - Jones polynomials of torus knots, 
as well as their extension to the DAHA super 
and hyper-polynomials including a link to the Khovanov-Rozansky 
polynomials in the stable case and for exceptionally small $N$ 
(for some knots). The focus is on the 3-parametric
super-polynomials 
\cite{GSV}, \cite{DGR}, \cite{KhR2}, \cite{AS}, \cite{GS},
\cite{ORS}.
 
Our approach is based on the technique that does not involve 
the roots of unity and the Verlinde algebras. The construction of 
the DAHA-Jones polynomials holds for arbitrary 
(reduced, twisted) root systems and any weights. Their 
stabilizations are certain DAHA counterparts of 
super-polynomials for the root systems of type $A,D$, as well as 
the hyper-polynomials of type $B$ and $C$ (with an extra
parameter, the fourth).

The super-polynomials of torus knots have deep  relations to
the refined BPS theory, the Khovanov-Rozansky homology and 
the Hilbert schemes of $\C^2$ and singular curves.
In the DAHA theory, it opens a new and challenging direction
related to the DAHA-coinvariants and elliptic Hall functions
(see \cite{C102}). 
\smallskip 

\subsection{Our approach}
Our interpretation of the QG-Jones colored polynomials of torus knots 
is based directly on the PBW theorem for DAHA. Given a torus knot,
we apply the corresponding element of the projective $PSL_2(\Z)$
to the Macdonald polynomial representing the color 
and then take the DAHA {\em evaluation coinvariant}. The latter
is calculated using a Shapovalov-type machinery; 
see \cite{C102} for a general theory of Shapovalov functionals 
in DAHA. 

The connection with the colored Jones polynomials from
Conjecture \ref{MAINCONJ} is proved for $A_n$ via
the technique of the DAHA shift operator; we simply compare
our formula with (5.4) from \cite{St} and that from
\cite{LZ} (Theorem 5.1). The Jones polynomials 
of torus knots are well known for $A_1$; a complete 
proof is provided in this case. See Section \ref{sec:escqg} 
concerning the non-$A$ types. 

The theory of the
shift operators results in a uniform and relatively 
straightforward expressions for any root 
systems, so Conjecture \ref{MAINCONJ} would follow from
the counterparts of the formulas from \cite{St,LZ,CC} 
for arbitrary root systems. For instance formula (5.6) from
\cite{St} seems sufficient to verify Conjecture \ref{MAINCONJ} 
for type $D$ (see also \cite{CC}), 
but we did not check all details.

\smallskip
Proposition \ref{SH-JONES}
provides explicit formulas $k=1$ ($t=q$). Our method and formulas
resemble some steps in papers \cite{AS}, \cite{GSV}, \cite{LM} 
(formula (3.29) there) 
and other papers involving integrable modules of Kac-Moody algebras 
(or quantum groups at roots of unity). However, it is an
analogy rather than an exact connection at this moment. 
Our approach is also related to the BPS theory and matrix models.  

The Verlinde algebras are clearly in the intersection of
the Chern-Simons theory and DAHA; such algebras and their 
generalizations are well understood in the DAHA theory.  
For instance, the projective action of $PSL_2(\Z)$ on the Verlinde 
algebras is a corollary of its action on DAHA, which is directly 
linked to the topology of the elliptic configuration space. 
We note that the Pieri rules, important in the theory of the polynomial
DAHA representation and related to this paper,
were justified in \cite{C3} using the passage through the
roots of unity.
\medskip

Our conjectures on the stabilization, duality and the special 
values of the DAHA super and hyper-polynomials entirely
belong to the DAHA theory and do not depend
on the coincidence with
the QG-Jones, HOMFLYPT and Kauffman polynomials (upon proper
reductions of the parameters).
Conjecture \ref{CONJKHR} relates the DAHA super- polynomials
to the Khovanov- Rozansky polynomials
(including small $N$ for some knots).

\subsection{Perspectives}
It is worth mentioning that the generalized Verlinde algebras 
have $q$\~ deformations from roots of unity to any $q$
(the case $|q|=1$ is of importance), 
though this line is not directly related to this paper.
All structures are preserved under such deformations
but the integrality and positivity of the $N_{ij}^k$ numbers. 
See, e.g.,  \cite{C101}, Section 2.10.5  and remark after
Theorem 2.9.9.  

This provides a link from the $q,t$\~theory to the
{\em rational} DAHA, already known to be connected with 
the super-polynomials and Khovanov-Rozansky homology of the
torus knots due to Gorsky, Oblomkov, 
Rasmussen, Shende, based on prior works on Hilbert schemes by 
Haiman, Gordon, Stafford, Nakajima and others.  
We see some connections of these two approaches, but so far 
there are no formal claims or conjectures in this direction. 
The match at level of formulas is impressive. 

The standard relation of the $q,t$\~DAHA to
the rational one is upon the limiting procedure $q\to 1$. 
It is compatible with our formula, but does not 
connect it with the aforementioned theory so far.  
The rational case is managed in the $A_1$\~case in detail
in this paper. 
\smallskip

There is likeness between the DAHA rational limiting 
procedure and that in the hyperbolic volume conjecture. 
However we do not set $q=\exp(2\pi i/N)$ before sending 
$N\to \infty$, where $N-1$ is the number of colors.
See \cite{MM}, Conjecture 5.1; the normalization must be 
$J(\hbox{unknot})=1$. Our limit is non-trivial in 
contrast to the hyperbolic volume, which is zero for torus knots. 

Also, we can take $|q|<1$ and try to replace the Macdonald 
polynomials in the construction of our $q,t$\~invariants by 
the {\em global
$q,t$\~hypergeometric function} introduced in \cite{C5},
which is in the focus of the
DAHA theory of the last years. It readily adds an analytic
dimension to the theory and presumably can provide 
generating functions and recurrence relations for the DAHA-Jones 
polynomials and super-polynomials. 

\smallskip
This paper indicates that there are deep unknown relations  
of DAHA to the theory of singularities and related
physics (the refined BPS states and 
matrix models), topology (knot Floer homology, 
Khovanov-Rozansky homology) and Hilbert schemes.  
Concerning possible applications to the 
QG-Jones polynomials, we provide 
in Section \ref{sec:escqg} some DAHA-generated
formulas for exceptional root systems in a hope that they
can be directly verified against the standard formulas.


Our approach has a clear topological meaning,
matching some topological ingredients of the 
Khovanov-Rozansky theory \cite{KhR1, KhR2} for such knots. 
Generally, Hecke algebras are closely related to the 
topology of braids and links; in the case of DAHA, it is for
the elliptic configuration spaces. However no exact
connections of the topological foundations of both theories
are established at the moment.


\subsection{Acknowledgements}
The author would like to thank E.~Gorsky, S.~Gukov, 
M.~Kho\-va\-nov, J.~Rasmussen and L.~Ro\-zan\-sky for very 
useful discussions and sharing with me information about their 
ongoing research. 
I thank Daniel Orr for his participation in creating the 
software for DAHA super-polynomials. I want to thank
Sergei Gukov and Tetsuji Miwa (Caltech and the Kyoto
University) for hospitality. My special thanks to
David Kazhdan for his attention to this work and the Hebrew 
University for the invitation. I would like to thank
the referee of this paper for various useful suggestions.

\setcounter{equation}{0}
\section{Double Hecke algebras}
We will begin with the basic DAHA definitions in the
twisted case.
Let $R=\{\al\}   \subset \R^n$ be a root system of type
$A,B,...,F,G$
with respect to a euclidean form $(z,z')$ on $\R^n
\ni z,z'$,
$W$ the {\em Weyl group} 
generated by the reflections $s_\al$,
$R_{+}$ the set of positive  roots ($R_-=-R_+$)
corresponding to fixed simple 
roots $\al_1,...,\al_n.$

\subsection{Affine root systems}
The root lattice and the weight lattice are:
\begin{align}
& Q=\oplus^n_{i=1}\Z \al_i \subset P=\oplus^n_{i=1}\Z \om_i,
\notag
\end{align}
where $\{\om_i\}$ are fundamental weights:
$ (\om_i,\al_j^\vee)=\de_{ij}$ for the
coroots $\al^\vee=2\al/(\al,\al).$
Replacing $\Z$ by $\Z_{+}=\{m\in\Z, m\ge 0\}$ we obtain
$Q_+, P_+.$
Here and further see  \cite{Bo}, \cite{Hu} and \cite{C101}. 

The form will be normalized
by the condition  $(\al,\al)=2$ for the
{\em short} roots in this paper. 
Thus,

\centerline{
$\nu_\al\equal (\al,\al)/2$ can be either $1,$ or $\{1,2\},$ or
$\{1,3\}.$ }
\noindent

The vectors $\ \tal=[\al,\nu_\al j] \in
\R^n\times \R \subset \R^{n+1}$
for $\al \in R, j \in \Z $ form the
{\em affine root system}
$\tR \supset R$ ($z\in \R^n$ are identified with $ [z,0]$).
We add $\al_0 \equal [-\vth,1]$ to the simple
roots for the {\em maximal short root} $\vth\in R_+$.
It is also the {\em maximal positive
coroot} because of the choice of normalization. 
The corresponding set
$\tR_+$ of positive roots equals 
$R_+\cup \{[\al,\nu_\al j],\ \al\in R, \ j > 0\}$.

We complete the Dynkin diagram of $R$
by $\al_0$ (by $-\vth$, to be more
exact); it is called {\em affine Dynkin diagram}
$\tGa$.  The number of laces
between $\al_i$ and $\al_j$ will be denoted by $m_{ij}$.
One can obtain $\tGa$ from the usual
{\em completed\,} Dynkin diagram from \cite{Bo} for 
the {\em dual system}
$R^\vee$ by reversing all arrows there.
 
The set of the indices of the images of $\al_0$ by all
the automorphisms of $\tGa$ will be denoted by $O$
($O=\{0\} \for E_8,F_4,G_2$). Let $O'\equal\{r\in O, r\neq 0\}$.
The elements $\om_r$ for $r\in O'$ are the so-called 
{\em minuscule
weights}: $(\om_r,\al^\vee)\le 1$ for
$\al \in R_+$.

Given $\tal=[\al,\nu_\al j]\in \tR,  \ b \in P$, let
\begin{align}
&s_{\tal}(\tz)\ =\  \tz-(z,\al^\vee)\tal,\
\ b'(\tz)\ =\ [z,\ze-(z,b)]
\label{ondon}
\end{align}
for $\tz=[z,\ze] \in \R^{n+1}$.

\subsection{Affine Weyl groups}
The 
{\em affine Weyl group} $\tW$ is generated by all $s_{\tal}$
(we write $\tW = \lan s_{\tal}, \tal\in \tR_+\ran)$. One can take
the simple reflections $s_i=s_{\al_i}\ (0 \le i \le n)$
as its
generators and introduce the corresponding notion of the
length. This group is
the semidirect product $W\lsmash Q'$ of
its subgroups $W=$ $\lan s_\al,
\al \in R_+\ran$ and $Q'=\{a', a\in Q\}$, where
\begin{align}
& \al'=\ s_{\al}s_{[\al,\,\nu_{\al}]}=\
s_{[-\al,\,\nu_\al]}s_{\al}\for
\al\in R.
\label{ondtwo}
\end{align}

The {\em extended Weyl group} $ \hW$ generated by $W\and P'$
(instead of $Q'$) is isomorphic to $W\lsmash P'$:
\begin{align}
&(wb')([z,\ze])\ =\ [w(z),\ze-(z,b)] \for w\in W, b\in B.
\label{ondthr}
\end{align}
From now on,  $b$ and $b',$ $P$ and $P'$ will be identified.

Given $b\in P_+$, let $w^b_0$ be the longest element
in the subgroup $W_0^{b}\subset W$ of the elements
preserving $b$. This subgroup is generated by simple
reflections. We set
\begin{align}
&u_{b} = w_0w^b_0  \in  W,\ \pi_{b} =
b( u_{b})^{-1}
\ \in \ \hW, \  u_i= u_{\om_i},\pi_i=\pi_{\om_i},
\label{xwo}
\end{align}
where $w_0$ is the longest element in $W,$
$1\le i\le n.$

The elements $\pi_r\equal\pi_{\om_r}, r \in O'$ and
$\pi_0=\hbox{id}$ leave $\tGa$ invariant
and form a group denoted by $\Pi$,
 which is isomorphic to $P/Q$ by the natural
projection $\{\om_r \mapsto \pi_r\}$. As to $\{ u_r\}$,
they preserve the set $\{-\vth,\al_i, i>0\}$.
The relations $\pi_r(\al_0)= \al_r= ( u_r)^{-1}(-\vth)$
distinguish the
indices $r \in O'$. Moreover,
\begin{align}
& \hW  = \Pi \lsmash \tW, \where
  \pi_rs_i\pi_r^{-1}  =  s_j \iif \pi_r(\al_i)=\al_j,\
 0\le j\le n.
\end{align}

Setting
$\hw = \pi_r\tw \in \hW,\ \pi_r\in \Pi, \tw\in \tW,$
the length $l(\hw)$
is by definition the length of the reduced decomposition
$\tw= $ $s_{i_l}...s_{i_2} s_{i_1} $
in terms of the simple reflections
$s_i, 0\le i\le n.$ 
Alternatively, 
\begin{align}
&l(\hw)=|\la(\hw)| \for \la(\hw)\equal\tR_+\cap \hw^{-1}(-\tR_+).
\label{xlambda}
\end{align}
\smallskip

For an arbitrary
weight $b\in P$, there exists $w\in W$ such that
$w(b)\in P_+$ and $b_+\equal w(b)$ is unique such; 
$b_+=c_+$ simply means that $b,c$
belong to the same $W$\~orbit. Let
\begin{align}
&b \prec c, c \succ b \for b, c\in P_+ \iif 0\neq c-b \in Q_+.
\label{succ}
\end{align}
\medskip

\subsection{Main definition}
By  $m,$ we denote in this section the least natural number
such that  $(P,P)=(1/m)\Z.$  Thus
$m=2 \for D_{2k},\ m=1 \for B_{2k} \and C_{k},$
otherwise $m=|\Pi|$.

The double affine Hecke algebra depends
on the parameters
$q, t_\nu,\, \nu\in \{\nu_\al\}.$ It will be defined
over the ring
$\Q[q^{\pm 1/2m},t_\nu^{\pm 1/2}]$
formed by
polynomials in terms of $q^{\pm 1/m}$ and
$\{t_\nu \}.$ The coefficient of Macdonald polynomials
will belong to
$$
\Q_{q,t}'\equal \Q(q^{\pm 1/2m},t_\nu^{\pm 1/2})
$$
(actually $q^{\pm 1/m},t_\nu$ are sufficient).

It will be convenient to use the parameters
$\{k_\nu\}$ together with  $\{t_\nu \},$ setting
\begin{align}
&   t_{\tal} =t_{\al}=t_{\nu_\al}=q_\al^{k_\nu} ,\ \, 
q_{\tal}=q^{\nu_\al}, \ \, t_i = t_{\al_i},  q_i=q_{\al_i},
\notag\\
&\where \tal=[\al,\nu_\al j] \in \tR,\ 0\le i\le n.
\label{taljx}
\end{align}

We will sometimes use $\sht,\lng$ instead of $\nu$:
\begin{align}\label{qtnuk}
\rho_k\equal \frac{1}{2}\sum_{\al>0} k_\al \al=
k_{\sht}\rho_{\sht}+k_{\lng}\rho_{\lng},\ \,
\rho_\nu=\frac{1}{2}\sum_{\nu_\al=\nu} k_\al \al.
\end{align}

For pairwise commutative $X_1,\ldots,X_n,$
\begin{align}
& X_{\tb}\ =\ \prod_{i=1}^nX_i^{l_i} q^{ j}
\iif \tb=[b,j],\ \hw(X_{\tb})\ =\ X_{\hw(\tb)},
\label{Xdex}\\
&\hbox{where\ } b=\sum_{i=1}^n l_i \om_i\in P,\ j \in
\frac{1}{ m}\Z,\ \hw\in \hW.
\notag \end{align}
For instance, $X_0\equal X_{\al_0}=qX_\vth^{-1}$.
\medskip

We note that $\pi_r^{-1}$ is $\pi_{r^*}$ and
$u_r^{-1}$ is $u_{r^*}$
for $r^*\in O\ ,$  $u_r=\pi_r^{-1}\om_r.$
The reflection $^*$ is
induced by the standard involution (sometimes
trivial) of the nonaffine Dynkin diagram.

\begin{definition}
The double affine Hecke algebra $\HH\ $
is generated over $\Q[q^{\pm 1/m},t_\nu]$ by
the elements $\{ T_i,\ 0\le i\le n\}$,
pairwise commutative $\{X_b, \ b\in P\}$ satisfying
(\ref{Xdex}),
and the group $\Pi,$ where the following relations are imposed:

(o)\ \  $ (T_i-t_i^{1/2})(T_i+t_i^{-1/2})\ =\
0,\ 0\ \le\ i\ \le\ n$;

(i)\ \ \ $ T_iT_jT_i...\ =\ T_jT_iT_j...,\ m_{ij}$
factors on each side;

(ii)\ \   $ \pi_rT_i\pi_r^{-1}\ =\ T_j \iif
\pi_r(\al_i)=\al_j$;

(iii)\  $T_iX_b \ =\ X_b X_{\al_i}^{-1} T_i^{-1} \iif
(b,\al^\vee_i)=1,\
0 \le i\le  n$;

(iv)\ $T_iX_b\ =\ X_b T_i\ $ if $\ (b,\al^\vee_i)=0
\for 0 \le i\le  n$;

(v)\ \ $\pi_rX_b \pi_r^{-1}\ =\ X_{\pi_r(b)}\ =\
X_{ u^{-1}_r(b)}
 q^{(\om_{r^*},b)},\  r\in O'$.
\label{double}
\end{definition}

Given $\tw \in \tW, r\in O,\ $ the product
\begin{align}
&T_{\pi_r\tw}\equal \pi_r\prod_{k=1}^l T_{i_k},\where
\tw=\prod_{k=1}^l s_{i_k},
l=l(\tw),
\label{Twx}
\end{align}
does not depend on the choice of the reduced decomposition
(because $T_i$ satisfy the same ``braid'' relations
as $s_i$ do).
Moreover,
\begin{align}
&T_{\hv}T_{\hw}\ =\ T_{\hv\hw}\  \hbox{ whenever}\
 l(\hv\hw)=l(\hv)+l(\hw) \for
\hv,\hw \in \hW. \label{TTx}
\end{align}
In particular, we arrive at the pairwise
commutative elements: 
\begin{align}
& Y_{b} = 
\prod_{i=1}^nY_i^{l_i} \iif
b=\sum_{i=1}^n l_i\om_i\in P,\ 
Y_i\equal T_{\om_i},b\in P.
\label{Ybx}
\end{align}
\smallskip

\subsection{Polynomial representation}
The {\em Demazure-Lusztig operators}
are as follows:
\begin{align}
&T_i\  = \  t_i s_i\ +\
(t_i-1)(X_{\al_i}-1)^{-1}(s_i-1),
\ 0\le i\le n;
\label{Demazx}
\end{align}
they obviously preserve $\Q[q,t_\nu][X_b]$.
We note that only the formula for $T_0$ involves $q$:
\begin{align}
&T_0\  = \ t_0s_0\ +\ (t_0-1)
(X_0 -1)^{-1}(s_0-1),\hbox{\ where\ }\notag\\
&X_0=qX_\vth^{-1},\
s_0(X_b)\ =\ X_bX_{\vth}^{-(b,\vth)}
 q^{(b,\vth)},\
\al_0=[-\vth,1].
\end{align}

The map sending $ T_j$ to the corresponding operator from
(\ref{Demazx}), $X_b$ to the operator of multiplication by $X_b$
(see (\ref{Xdex})) and 
$\pi_r\mapsto \pi_r$ induces a
$ \Q_{ q,t}'$\~linear
homomorphism from $\HH\, $ to the algebra of
linear endomorphisms
of $\Q_{ q,t}'[X]$.
This $\HH\,$-module is faithful
and remains faithful when  $q,t$ take
any complex values assuming that
$q\neq 0$ is not a root of unity.
It will be called the
{\em polynomial representation};
the notation is
$$
\v\equal \Q_{q,t}'[X_b]\ =\ \Q_{q,t}'[X_b, b\in P].
$$

The images of the $Y_b$ are called the
{\em difference-trigonometric Dunkl operators}.

The polynomial representation
is the $\HH\,$\~module induced from the one-dimensional
representation $T_i\mapsto t_i,\,$ 
$Y_b\mapsto q^{2(\rho_k,b)}$
of the affine Hecke subalgebra 
$\h_Y=\lan T_i,Y_b\ran.$ Here we extend the ring of
constants to $\Q_{q,t}'$ in the definition of $\HH\,$.
\medskip

The {\em nonsymmetric Macdonald polynomials} are the eigenvectors
of the $Y$\~operators in $\v$. Their $t$\~symmetrizations 
are the {\em symmetric Macdonald polynomials}; see \cite{C4},
\cite{C101} and references therein.
 
The following direct definition is due to Macdonald \cite{M2}
for arbitrary root systems and Kadell
for the classical root systems. For $b\in P_+$,
\begin{align}\label{macdsymf}
&P_b-\sum_{b'\in W(b)}\, X_{b'}\, 
\in\, \oplus_{\,c_+<b\,}\,\Q_{q,t}' X_c,\ \,
\hbox{CT} \bigl(P_b X_{c}\,\de(X;q,t)\bigr) = 0
\\
&\hbox{for\ \ }\de(X;q,t)\equal\prod_{\al \in R_+}
\prod_{j=0}^\infty \frac{(1-X_\al q_\al^{j})
(1-X_\al^{-1}q_\al^{j})
}{
(1-X_\al t_\al q_\al^{j})
(1-X_\al^{-1}t_\al^{}q_\al^{j})}\,,
\label{desymf}
\end{align}
where $CT$ is the constant term (the coefficient of $X_0$); 
$\de$ is considered
a Laurent series of $X_b$ with 
the coefficients expanded in terms of
positive powers of $q$.

We note that when $k_{\al}=1=k_{\nu}$ for all $\al,\nu$, 
then $t_\al=q_\al$ and  $\de$ becomes the
standard discriminant
$\prod_{\al \in R_+}(1-X_\al)(1-X_\al^{-1})$.
Thus the {\em symmetric}
Macdonald polynomials do not depend on $q$ 
and become the standard finite-dimensional characters 
in this case. 

It is important for any aspects of the theory of $P$\~polynomials 
that they are eigenfunctions of the Macdonald-Ruijsenaars operators, 
which are symmetric ($W$\~invariant) polynomials in terms of $Y$
upon the restriction to the $W$\~invariant part of $\v$.

\subsection{Automorphisms}\label{sect:Aut}
The following map can be uniquely extended to
an automorphism of $\HH\,$ where 
proper fractional powers of $q$ are added
(see \cite{C15},\cite{C4}):
\begin{align}
& \tau_+:\  X_b \mapsto X_b, \ T_i\mapsto T_i\, (i>0),\
\ Y_r \mapsto X_rY_r q^{-\frac{(\om_r,\om_r)}{2}}\,,
\notag\\
\label{tauplus}
& \tau_+:\ T_0\mapsto  q^{-1}\,X_\vth T_0^{-1},\
\pi_r \mapsto q^{-\frac{(\om_r,\om_r)}{2}}X_r\pi_r\
(r\in O'),\\
& \label{taumin}
\tau_-:\ Y_b \mapsto \,Y_b, \ T_i\mapsto T_i\, (i\ge 0),\
\ X_r \mapsto Y_r X_r q^\frac{(\om_r,\om_r)}{ 2},\\
&\tau_-(X_{\vth})= 
q T_0 X_\vth^{-1} T_{s_{\vth}}^{-1};\ \
\si\equal \tau_+\tau_-^{-1}\tau_+\, =\,
\tau_-^{-1}\tau_+\tau_-^{-1}.
\label{taux}
\end{align}
These automorphisms fix $\ t_\nu,\ q$
and their fractional powers, as well as the
following {\em anti-involution}:
\begin{align}
&\phi:\ 
X_b\mapsto Y_b^{-1},\, Y_b\mapsto X_b^{-1},\
T_i\mapsto T_i\ (1\le i\le n),\label{starphi}\\
&\phi(\tau_+)\equal \phi\circ \tau_+\circ \phi\ =\ \tau_-\,,\ 
\ \, \phi(\tau_-)\ =\ \tau_+.\notag 
\end{align} 
\medskip

This anti-involution is the key in proving
the Macdonald duality, evaluation and norm
conjectures; in this paper, we will need only
the evaluation formula for $P_b\, (b\in P_+)$:
\begin{align}\label{pebeb}
P_{b}(q^{-\rho_k})\ =&\ q^{(\rho_k,b)}
\prod_{\al>0}\prod_{j=0}^{(\al^\vee,b)-1}
\Bigl(
\frac{
1- q_\al^{j}t_\al X_\al(q^{\rho_k})
 }{
1- q_\al^{j}X_\al(q^{\rho_k})
}
\Bigr).
\end{align}
\medskip
It is a corollary of the following construction.

\subsection{The evaluation map}
Following \cite{C3,C4}, we set for $f,g\in \v$,
\begin{align}\label{symfg}
&\{f,g\}\equal \{L_{\imath(f)}(g(X))\}\ =\
\{L_{\imath(f)}(g(X))\}(q^{-\rho_k}),\\
X_b(q^{-\rho_k})&=q^{-(b,\rho_k)},\ \,
\imath(X_b)=X_{-b}=X_b^{-1},\ \,
\imath(z)=z \for
z\in \Q_{q,t}\, ,
\notag 
\end{align}
where $L_f\equal f(Y)$.

This pairing is symmetric and corresponds to the
anti-involution $\phi$ in $\HH\,$, i.e.,
$\{Hf,g\}=\{f,\phi(H)g\}$ for $H\in \HH\,$.
Indeed, it can be represented as $\{L_f(g(X))\}$
for the following $\phi$\~invariant 
{\em evaluation functional} $\{\cdot\}$ on $\HH\,$.

We use the PBW theorem to express any $H\in \HH$ in the form 
\,$\sum_{a,w,b} c_{a,w,b}\, X_a T_{w} Y_b$\, for $w\in W$,
$a,b\in P$ (this presentation is unique). Then we substitute:
\begin{align}\label{evfunct}
X_a \ \mapsto\  q^{-(\rho_k,a)},\ 
Y_b \ \mapsto\  q^{(\rho_k,b)},\ 
T_i \ \mapsto\  t_i^{1/2}. 
\end{align}
The resulting functional $\HH\ni H\mapsto \{H\}$
acts via the projection $H\mapsto H(1)$ of $\HH\,$
onto $\v$, namely,
 $\{\,H\,\}=H(1)(q^{-\rho_k})=\{H(1),1\}$.
\smallskip

Generalizing (\ref{symfg}) one can consider here any 
character $\chi$ of the {\em non-affine} Hecke algebra 
generated by $\{ T_1,T_2,\ldots,T_n\}$ and also the character 
(an algebra homomorphism) 
$\ze: \C[X_a,a\in P]\to \C$. Then we replace: 
\begin{align}\label{evfunchi}
T_{w}\mapsto \chi(T_w),\ X_a \ \mapsto\ze(X_a) ,\ 
Y_b \ \mapsto\  \ze(X_b^{-1}). 
\end{align}
The property $\chi(T_uT_w)=\chi(T_wT_u)$ for
$u,w\in W$ readily results in the
$\phi$\~invariance of such functional.

There is also a possibility of using the
other two major DAHA functionals (linear maps to $\C$)
instead of (\ref{symfg}). 
For $|q|<1$, we set 
\begin{align}\label{ctfunct}
&\lan f\ran_0=\hbox{CT}\bigl(f\,\de(X;q,t)\bigr) 
\for f\in \C[X_b,\,b\in P],\,\de \hbox{\ from
(\ref{desymf})\,},
\\
\label{ctgafunct}
&\lan f\ran_1=\hbox{CT}\bigl(f\,\vth(X;q)\de(X;q,t)\bigr),\ 
\vth(X;q)\equal\sum_{b\in P}\,X_{b}\,q^{b^{2}/2}.
\end{align}
Correspondingly, 
$$
\lan f \ran'_0\,=\,
\frac{\lan f \ran_0}{\lan 1\ran_0},\ \ 
\lan f \ran'_1\,=\,\frac{\lan f \ran_1}{\lan 1\ran_1}\,.
$$
See \cite{C101} and \cite{C102} (Theorem 2.15 and, 
especially, Section ``Polynomial case" after it and 
formula (2.40) there). Let us mention here
the following important relation:
\begin{align}\label{evaga}
&\lan \tau_-(H)(1)\ran_1'\ =\ \{\,H\,\} \for H\in \HH\,.
\end{align}
See, e.g., Proposition 3.3.4 from \cite{C101}.

\setcounter{equation}{0}
\section{Super-polynomials via DAHA}
There are many sources devoted to the torus knots, 
including the following site: http://katlas.org/wiki/.
Concerning Jones polynomials and the
Quantum Groups generalizations (referred to as QG-Jones
polynomials) see \cite{Jon},\cite{Resh},\cite{Ros}.
The HOMFLY-PT polynomials were introduced in \cite{FYH}
and \cite{PT}. 

The works
\cite{RJ},\cite{Mo},\cite{Hi}, \cite{LZ} and \cite{St}
are basically sufficient in this section. See also
\cite{GSV},\cite{DGR} and \cite{AS}, \cite{GS} about 
the ongoing theory of super-polynomials of torus knots.  

\subsection{Jones polynomials}
Given a torus knot $K_{r,s}$ in $S^3$ of type $\{r,s\}$,
let $\tga_K=\tga_{r,s}$ be a product 
$\ \ \ldots (\tau_-)^w\,(\tau_+)^v\,(\tau_-)^u$ for 
$u,v,w,\ldots \in \Z$
such that the first column of the corresponding element
$\ga_K\in PSL_2(\Z)$ is $(r,s)^{tr}$ ($tr$ is the transposition). 
Here we send 
$$
\tau_+\mapsto \binom{11}{ 01},\  \tau_-\mapsto
\binom{10}{  11}.
$$
We allow $r,s$ to be arbitrary
relatively prime integers, including zero and 
negative numbers. We will use
the symmetric Macdonald polynomials $P_b\, (b\in P_+)\,$
and their evaluations $P_b(q^{-\rho_k})=\{\, P_b\,\}$.
\medskip

For a polynomial $F$ in terms of positive and negative 
fractional powers of $q$ and $t$ ($t_\nu$ to be exact), 
the {\em tilde-normalization}
$\tilde{F}$ will be the result of its division by the lowest 
$q,t$\~monomial (assuming that it is well defined).
In the following conjecture, the lowest monomials will always 
exist (it can be checked)
and the corresponding $\tilde{F}$ will contain only integral 
non-negative powers of $q,t_\nu$ upon this normalization.

\medskip
\newtheorem{conjecture}[theorem]{Conjecture}
\begin{conjecture} [DAHA-Jones polynomials]\label{MAINCONJ} 
Given a knot $K=K_{r,s}$ for $r,s\in \Z$, the root system $R$
and a weight $b\in P_+\,$, let
\begin{align}\label{jones-d}
J\!D_{r,s}^{R}(b\,;\,q,t)=
J\!D_{r,s}(b\,;\,q,t)\equal \{\,\tga_{r,s}(P_b)/P_b(q^{-\rho_k})\,\}.
\end{align}
Then $\tilde{J\!D}_{r,s}(b\,;\,q,t)$, the corresponding
DAHA-Jones polynomial, 
does not depend on the particular choice of $\ga_{r,s}\,$ 
representing $K$; it is a polynomial in terms of $q,t_\nu$.
For $t\mapsto q$, more exactly, upon the substitution $k_\nu=1$
in (\ref{qtnuk}), one has:
$$
\tilde{J\!D}_{r,s}(b\,;\,q,t\mapsto q)\ =\ \tilde{\j}_{r,s}(b\,;\,q)
$$ 
for the Jones polynomial $\j_{r,s}(b\,;\,q)$
of $K_{r,s}\,$ defined for the quantum 
group associated with the roots system $\tR$ and its 
irreducible representation with the highest weight $b$.
The normalization here is by the condition  
$\j(\hbox{unknot})=1$.
\end{conjecture}

The polynomiality of the (reduced) $J\!D$\~invariants
can be deduced from the theory 
of the polynomial DAHA-modules. The zeros of the evaluations
of Macdonald polynomials are directly connected with the values of
parameters $q,t$ where $\v$ becomes reducible; see \cite{C103}.

The coincidence with Quantum Groups - Jones polynomials
(referred to as QG-Jones polynomials)
was established for $A_n$ for any weights using the formulas from
\cite{GMV}, \cite{LM} and, especially, the HOMFLYPT 
formula (5.4) from \cite{St} and 
\cite{LZ} (Theorem 5.1). The following proposition provides the
necessary tools. We need the simplest case of the theory of 
difference shift operators
from Theorem 2.4 of \cite{C3} or other author's papers.

Let us assume that $k=\{k_\nu=0,1\}$ and define
\begin{align}\label{shiftx}
\x_k  = \prod_{\al\in R_+}\prod_{j_\al=0}^{k_\al-1}
\bigl((q_\al^{j_\al} X_{\al})^{1/2}-
(q_\al^{j_\al} X_{\al})^{-1/2})\bigr).
\end{align}
Recall that $q_\al=q^{\nu_{\al}}$,
where $\nu_\al\equal (\al,\al)/2$.

We put $H=H^{(k)}$ for elements $H$ from the algebra
$\HH^{(k)}$ with the structural parameters
$q$ and $t_\nu=q_\nu^{k_\nu}$. The image 
of $H=H^{(k)}$ in the
corresponding polynomial representation $\v^{(k)}$
will be denoted by $\hat{H}^{(k)}$
or by $(H^{(k)})^{\wedge}$. It is an operator in terms
of $X_b, \hat{w}\in \hat{W}$ acting there by multiplication; 
the dependence of $X$ can be via rational functions. 
Correspondingly, $\tau_{\pm}^{(k)}$ will be the automorphisms of the
double Hecke algebra for $t_\nu=q_\nu^{k_\nu}$. 

For instance, the operator 
$\hat{H}^{(0)}$ is obtained from $H$ by
replacing every $Y_b$ by the corresponding difference
operators $b^{-1}$ and $T_w (w\in W)$ replaced by $w$. Also,
\begin{align}\label{tauzero}
&\tau_{+}^\circ (Y_b)=q^{-(b,b)/2}X_b Y_b,\ \,
\tau_{-}^\circ (X_b)=q^{(b,b)/2}Y_b X_b,\\
&\hbox{where\ \ }\tau_{\pm}^\circ\, \equal\, \tau_{\pm}^{(0)} \for
b\in P,\, 0=\{k_\nu=0\}.\notag
\end{align}

\begin{lemma}\label{SHIFTAU}
Let $H^{(k)}$ be an algebraic expression with constant coefficients
from $\C$ in terms of the standard $W$\~invariant polynomials 
(say, the symmetrizations of monomials) with respect to
$\{Y_b\}$ and those with respect to $\{X_b\}$. 
Then $\hat{H}^{(k)}$ can be
considered as an operator acting in 
the subsubspace $(\v^{(k)})^W$ of $W$\~invariant elements 
in $\v^{(k)}$, denoted by $\hat{H}^{(k)}_{sym}$. 
One has:
\begin{align}\label{xh0x}
&\hat{H}^{(k)}_{sym}\ =\ \x_k^{-1}\, \hat{H}_{sym}^{(0)}\, \x_k\,,\\
&\bigr(\tau_{\pm}^{(k)}(H^{(k)})\bigl)^{\wedge}_{sym}\ =\ 
\x_k^{-1}\, \bigl( \tau_{\pm}^\circ(
\hat{H}^{(0)})\bigr)^{\wedge}_{sym}\, \x_k.\label{xhtaux}
\end{align}
\sq 
\end{lemma}
\medskip

We are going to make $\tau_{\pm}^\circ$ inner automorphisms
in a proper $\HH^{(0)}$\~module and extend them to the action 
of $SL(2,\Z)$ there.
The necessary extension of $\v^{(0)}$ is
the linear span $\tilde{\v}$ of products $X_\la q^{M x^2/2}$
for complex (sufficiently general) $\la, M$,
where we formally set 
$$
X_\la=q^{(\la,x)},\ 
x^2=(x,x)=\sum_{i=1}^n (x,\om_i)(x,\al_i^\vee).
$$
We naturally extend the action of $W$ and $\hat{W}$ in this
space through its action on $\{x_\la\}$; then 
$q^{M x^2/2}$ will be symmetric ($W$\~invariant).

The main formulas we need are as follows:

\begin{align}\label{sl2Mq}
&\si^\circ(X_\la\,q^{-M x^2/2})\,=\,
\frac{q^{\frac{\la^2}{2M}}}{M^{1/2}}\,
X_{\la/M}\,q^{+x^2/(2M)} \for M\neq 0,\\
&(\tau_+^\circ)^N\,(X_\la\,q^{-M x^2/2})\ =
\ X_\la\,q^{(N-M)x^2/2}\for N\in \Z
\and \notag\\ 
&(\tau_-^\circ)^N\,(X_\la\,q^{-M x^2/2})\ =\ 
((\si^\circ)^{-1}\,(\tau_+^\circ)^{-N}\,
\si^\circ)(X_\la\,e^{-M x^2/2})
\notag\\
&=\ \frac{1}{(1-MN)^{1/2}}\,
q^{\frac{\la^2\,N}{2(1-MN)}}\
X_{\la/(1-MN)}\, q^{-x^2\frac{M}{2(1-MN)}}.\notag
\end{align}

Assuming that $M\mapsto M^{1/2}$ is a homomorphism, 
these formulas can be readily extended to $SL(2,\Z)$:

\begin{align}\label{sl2zq}
&\ga^\circ\,(X_\la\,q^{z x^2/2})\,=\,
\frac{1}{(cz+d)^{1/2}}\ q^{-\frac{\la^2\, c}{2(cz+d)}}\
X_{\frac{\la}{cz+d}}\ q^{\frac{az+b}{cz+d}\,x^2/2} \for\\
&\ga= \left(
  \begin{array}{cc}
    a & b \\
    c & d \\
  \end{array}
\right)\in SL(2,\Z),\ \la \in \C 
\hbox{\ and for generic\, } z\in \C\,. \notag
\end{align}
We will call it the {\em free $SL(2,\Z)$\~action}.
The mere fact of its existence makes it possible
to deduce  relation (\ref{xhtaux}) from
(\ref{xh0x}). Indeed, 
$\tau_+$ and $\tau_-$ become conjugations 
by $W$\~invariant series in terms of $\{X_b\}$ 
and $\{Y_b\}$ upon the
proper completion of $\HH$ . For instance,
we can set 
$$
Y_b^{(0)}=q^{y_b}=e^{-\partial_b} \for
y_b=-\frac{\partial_b}{\log q},\ 
\partial_b(x_a)=(b,a),
$$
where the notation $ x_a=(x,a)$ is used
($a\in \C^n$).
Then $\tau_-^{(0)}$ acts in $\tilde{\v}$ as
$\tau_-^\circ=q^{-y^2/2}=e^{-\frac{\partial^2}{2\log q}}$;
here $y^2$ and $\partial^2$
are defined for the standard inner product
$(\,,\,)$ as for $x$. 

Then $\tga^{(k)}$ acts in the space of
$W$\~invariant elements of $\tilde{\v}$ and
coincides there with $\x_k^{-1}\,\ga^\circ\,\x_k$ for
$\ga^\circ$ from  (\ref{sl2zq}). We 
come to the following proposition.

\begin{proposition}\label{SH-JONES}
We continue to assume that $k_\nu=0,1$.
For $\ga=\ga_{r,s}$,
\begin{align}\label{jones-shift}
&J\!D_{r,s}(b\,;\,q,t)\ =\ \{\, R_{r,s}^{k,b}\,\}_k\ =\ 
R_{r,s}^{k,b}(q^{-\rho_k}),\where\\ 
&R_{r,s}^{k,b}\ \equal\ \x_k^{-1}\,
\ga^\circ \,\bigl((P^{(k)}_b/P^{(k)}_b(q^{-\rho_k}))\,
(\ga^\circ)^{-1}(\x_k)\,\bigr).\notag\\
\notag
\end{align}
\sq
\end{proposition}

The proposition makes the calculation of 
the DAHA-Jones invariants straightforward 
for any root systems when $\{k_\nu=0,1\}$.

The Macdonald polynomial $P^{(k)}_b$, generally, can be obtained,
for instance, using the shift operators; see formula (2.18)
from Theorem 2.4 of \cite{C3}. When $k_\nu=1$, what we really need
here, they are simply the characters of
finite-dimensional irreducible representations of
the simple Lie algebra associated with the root system $R$;
correspondingly,
$
\x_1=\sum_{w\in W} (-1)^{\hbox{\tiny sgn}(w)}X_{w(\rho)}.
$

\smallskip
The output of (\ref{jones-shift}) can be directly 
identified with formula (5.4) from \cite{St} and 
Theorem 5.1 from \cite{LZ}. For instance (cf. \cite{St},
Appendix A), 
$$
(\ga_{r,s}^\circ)^{-1} (\x_1)=
\frac{1}{r^{1/2}}\, q^{\frac{s\rho^2}{2r}-\frac{u x^2}{2r}}
\sum_{w\in W} (-1)^{\hbox{\tiny sgn}(w)}
X_{\frac{w(\rho)}{r}},\ 
\ga_{r,s}=\left(
\begin{array}{cc}
r & u \\
s & v \\
\end{array}
\right).
$$
{\em Thus, the coincidence
part of Conjecture \ref{MAINCONJ} is verified for
$A_n$.} It is possible to extend this calculation to obtain
the formula there for the Quantum Groups - HOMFLYPT polynomials too. 

We will not discuss here the details, but the calculation
is really similar to that from
Proposition \ref{MAINCONJ1} below in the case of $A_1$,
which provided in full. We note that formula (5.6) from
\cite{St} seems sufficient to check Conjecture \ref{MAINCONJ} 
for type $D$ (the case of Kauffman polynomials).

\smallskip
The coincidence of the $J\!D$\~polynomials for $t=q$ ($k_\nu=1$)
with the {\em reduced} (divided by the corresponding $q$\~dimension) 
and {\em tilde-normalized} QG-Jones invariants was extensively
checked using the program \ QuantumKnotInvariants[$q^{-1/2}$]\  
at http://katlas.org/wiki/. It was done for $A_n$ 
(many weights) and for $B_{2-4}$, $C_{3-4}$, $D_{4-5}, G_2$ 
(mainly minuscule and quasi-minuscule weights); 
see also Section \ref{sec:escqg}.

\subsection{Three super-conjectures}

The {\em tilde-normalization\,} from the previous
section will be applied in this section to the coefficient
of $a^0$ (the value at $a=0$). 
The coefficients of $a^l\,$ 
in the expressions $\tilde{F}$ in the conjectures below
may contain negative 
integral powers of $t$ (only for $l>0$); the powers of 
$q,a$ will be always integral non-negative. Actually, the
existence of the tilde-normalization is not proven in
full at the moment, so it is part of the conjectures below. 

\begin{conjecture} [Stabilization]\label{HOMFLY}
In the $A_n$\~case, we naturally interpret $b\in P_+$ 
as weights for greater systems $A_{m}\, (m\ge n)\,$. 
There exists a DAHA 
super-polynomial $H\!D_{r,s}(b\,;\,q,t,a)$ in terms of 
non-negative integral powers of $\,a,q\,$ and integral 
powers of $\,t\,,$ possibly negative, such that for 
all $m\ge n$,
\begin{align}\label{super-polyn}
&H\!D_{r,s}(b\,;\,q,t,a\mapsto -t^{m+1})=
\tilde{J\!D}^{m+1}_{r,s}(b\,;\,q,t)\equal
\tilde{J\!D}^{A_{m}}_{r,s}(b\,;\,q,t),\\
&H\!D_{r,s}\,(b\,;\,q,\,t\mapsto q,\,a\mapsto -a)\ \, =\ \, 
\tilde{\HOM}_{r,s}(b\,;\,q,a),\ \where\notag
\end{align}
$\HOM_{r,s}(b\,;\,q,a)$ is the HOMFLYPT polynomial for the
weight $b\in P_+$ normalized by the condition
$\HOM(unknot)=1$. Here 
$H\!D_{r,s}(b\,;\,q,t,a=0)$ is a polynomial in terms 
of $q,t$ with the constant term $1$. The polynomial
$H\!D_{r,s}(b\,;\,q,t,a)$ depends only on the Young diagram
representing $b$ (not on $n$) and its 
$\,a,q,t$\~coefficients are integers, possibly negative. 
They are non-negative for
$b=i\,\om_j$, corresponding to the rectangles $j\times i$.  
\end{conjecture}

\comment{
Sometimes the greatest power of $a$ in  
$H\!D_{r,s}(b\,;\,v,w,a)$ 
is a pure monomial in the form $a^A q^B t^C$, where
$$
B=\ \hbox{max\,deg}_{\,q}\,(H\!D),\ \  
C=\ \hbox{min\,deg}_{\,t}\,(H\!D).
$$
Thus, $q^{-B}t^{-C}H\!D_{r,s}$ is a polynomial 
in terms of non-negative powers of $a, q^{-1},t$.
This provides the most universal
normalization of $H\!D_{r,s}$, though there are some
advantages of the tilde-normalization too. Let
\begin{align}\label{duanorm}
&H\!D_{r,s}^\dag(b\,;\,q,t,a)\equal 
q^{-B} t^{-C} H\!D_{r,s}(b\,;\,q,t,a)\\
&=\ a^A+ \sum_{m=0}^{A-1}\sum_{i=0}^B\sum_{j=0}^{C}\, 
c_{ij}^m\, a^m q^{-i}\, t^j.
\notag
\end{align}
The coefficients $c_{ij}^m$ here are integers by construction.
}
\smallskip

\begin{conjecture} [Duality]\label{DUALIT}
Continuing the pervious conjecture, for an arbitrary weight 
$b\in P_+$ (equivalently, for an arbitrary Young diagram),
the $a$\~constant term  $H\!D_{r,s}(b\,;\,q,t,a=0)$, 
a polynomial in terms of non-negative powers of $q,t$, 
contains the greatest $q,t$\~monomial $q^A t^B$ (any $q^i t^j$
there are with $i\le A, j\le B$); its 
the coefficient is $1$. Then
\begin{align}\label{duaeq}
t^A\,q^B\,H\!D_{r,s}(b\,;\,q\mapsto t^{-1},
t\mapsto q^{-1},a)\ =
\ H\!D_{r,s}(b^{tr};q,t,a)\,,
\end{align}
where $b^{tr}$ corresponds to the transpose of the
Young diagram $\mu_b$ associated with $b$ (for sufficiently
large initial $n$).
\end{conjecture}

\begin{conjecture}[Evaluations]\label{CONJEVAL}
The evaluation $H\!D_{r,s}(b;q,t=1,a)$
is the product of the
$H\!D_{r,s}(b_i;q,t=1,a)$ calculated
for $b_i=m_i\om_1$, where $m_i$ is the number of boxes in the
$i$-th row of the Young diagram $\mu_b$. 
Moreover, 
$H\!D_{r,s}(b_i;q=1,t=1,a)$=$\bigl(H\!D_{r,s}(\om_1;q=1,t=1,a)
\bigr)^{\hbox{\small ord}}$ for the order of $\mu_b$. Also,
\,deg${}_{\,a\,} (H\!D_{r,s})=\,ord\times (s-1)\,$, where $0<s<r$.
\end{conjecture} 

The formula for the top
$a$\~degree obviously coincides with that for the
evaluation at $t=1=q$ provided the positivity of all 
$q,t$\~coefficients in the $a$\~leading term of 
$H\!D_{r,s}(\om_1;q,t,a)$, for instance, for the
rectangles $\mu_b$ if the end of Conjecture \ref{HOMFLY}
holds true. The maximal degree of $a$ in 
$H\!D_{r,s}(\om_1;q=1,t=1,a)$ can be calculated and it 
agrees with that obtained via the Gorsky presentation 
of the super-polynomials for $\om_1$; see \cite{Gor}.
\smallskip

\subsection{Toward KhR polynomials}
Continuing with the $A$\~case, let us introduce the
{\em restricted DAHA super-polynomials}  
$H\!D_{r,s}^{[n]}$
for any $n\in \N$ but no smaller than the number of rows in the
Young diagram representing the weight $b$. We 
follow the construction from Conjecture 3.1 \cite{DGR} (see
also examples there). I am thankful to Sergei
Gukov and Jacob Rasmussen for clarifying discussions on 
\cite{DGR,Ras}.

Let $H\!D_{r,s}(b\,;\,q,t,a)=
\sum_{d,i,j}\, c_{d,\,i,\,j}\,a^d t^i q^j\,$. 
For $d=1$, we find all
$c_{d,\,i,\,j}$ such that 
$\,c_{d,\,i,\,j}\,c_{d-1,\,i+n,\,j}>0\,$ and 
diminish $c_{d,\,i,\,j}$ and $c_{d-1,\,i+n,\,j}$ (both) by 
\begin{align}\label{redprocedure}
&\sgn (c_{d-1,\,i,\,j}\,)\,
\min \{|\,c_{d,\,i,\,j}\,|\,,\, |\,c_{d-1,\,i+n,\,j}\,|\}. 
\end{align}
After examining all monomials of $a$\~degree $d=1$, we 
go to $d=2$ and so on till the greatest power of $a$.

The resulting polynomial will be denoted by
$H\!D_{r,s}^{[n]}(b\,;\,q,t,a)$.
By construction,
$H\!D_{r,s}^{[n]}(b\,;\,q,t,a=-t^n)$
$=\tilde{J\!D}_{r,s}^n(b\,;\,q,t)$. The number of
different monomials in the latter is no greater
than in the former.
Moreover, negative coefficients can appear
in the restricted super-polynomial only from the 
corresponding negative coefficients in the original one.

If at least one reduction occurs in this
process, we call $n$ {\em exceptional}; 
sufficiently large $n$ are then {\em regular}. 
By definition, {\em DAHA-KhR polynomials} 
$K\!R_{r,s}^{n}(b\,;\,q,t)$ are obtained
from $H\!D_{r,s}^{[n]}(b\,;\,q,t,a)$ upon
the substitution
$a\mapsto t^n \sqrt{t/q}$. It is reasonable
to switch here to the standard
parameters (used, e.g., in \cite{DGR} and
http://katlas.org/wiki/)\,:
\begin{align}\label{qtareli}
t=q^2_{st},\  q=(q_{st}t_{st})^2,\  a=a_{st}^2 t_{st}.
\end{align}
Then the substitution to the KhR polynomials
becomes $a_{st}=q_{st}^n$.

Here $n$ is no smaller than the number $p$
of rows of $\mu_b$. If $n=p$, 
then we set
$\tilde{J\!D}^n_{r,s}(b\,;\,q,t)\equal$ 
$H\!D_{r,s}(b\,;\,q,t,a=-t^n)$.


\medskip
This definition is directly related to  
Conjectures 1.7 and 3.1 from \cite{DGR}. The later 
includes the existence of a family of anti-commutative 
differentials $\{d_n\}$ on the triply-graded stable
Khovanov- Rozansky space 
$\h=\oplus\h_{i,j,k}$ such that the Poincar\'e polynomial
of the homology of $d_n$ is the Khovanov- Rozansky
polynomial for $SL(n)$. The grading here is associated with 
the parameters $\{\,a_{st},q_{st},t_{st}\,\}$; \ 
$d_n$ acts from $\h_{i,j,k}\to \h_{i-2,j+2n,\,k-1}$\,, by
zero if one of these spaces is zero.

Using \cite{KhR1,KhR2}, Jacob Rasmussen \cite{Ras}
established the existence of such differentials and
conjectured the degeneracy of the 
corresponding spectral sequence, which guarantees 
the desired relation to the Khovanov- Rozansky polynomials
with the reservation that the differentials can be very
difficult to calculate.

The construction from (\ref{redprocedure}) is directly related
to the following additional assumption on the images and 
kernels of $d_n$:
\begin{align}\label{nonzerod}
&\dim (\hbox{Im}(\h_{i+2,\,j-2n,\,k+1}))\,=\,\\
\hbox{Min}\,\bigl(\,
\dim &(\h_{i,j,\,k}\,/\,\hbox{Ker}(\h_{i,j,\,k}))\,,\,
\dim (\h_{i+2,\,j-2n,\,k+1})\,\bigr)
\notag
\end{align}
for all indices $\{i,j,k\}$. If such conditions hold, 
the extraction of the 
Khovanov- Rozansky polynomial for this $n$ from the
super-polynomial becomes a combinatorial problem (as
in the examples from \cite{DGR}).

\begin{conjecture}\label{CONJKHR}
Let us assume that all coefficients of
$H\!D_{r,s}(b\,;\,q,t,a)$ are positive,
the triply-graded space $\h$ is well defined for the torus 
knot $\{r,s\}$ and the weight $b$\,, and that
(\ref{nonzerod}) holds for $d_n$ including the degeneracy
of the corresponding spectral sequence from \cite{Ras}. Then :

(i) the dimensions of  $H\!D_{r,s}^{[n]}(b\,;\,q,t,a)$
and $\tilde{J\!D}^n_{r,s}(b\,;\,q,t)$, defined as the
sums of absolute values of all coefficients, coincide;

(ii) the polynomial
$K\!R_{r,s}^{n}(b\,;\,q,t)$ coincides with the
corresponding reduced Khovanov-Rozansky polynomial upon the
tilde-normalization (the standard parameters are used);

(iii) the dimension of $\tilde{J\!D}^n_{r,s}(b\,;\,q,t)$ 
coincides with $\dim\h$ and the
$q_{st}$\~grading of $\h$ (associated with $j$)  
results from this polynomial up to a general
shift in $j$.
\end{conjecture}

Following
\cite{DGR}, a similar relation
to the knot Floer homology of torus knots can be expected, 
but we will not define the differential $d_0$ and
discuss this line here.

Claim $(i)$ is entirely algebraic. We have no
counterexamples so far, but do not conjecture
that it is always true (provided the positivity of
$H\!D_{r,s}(b\,;\,q,t,a)$).
For instance, $(i)$ gives that there will be no cancelation
of different monomials in $H\!D_{r,s}^{[n]}(b\,;\,q,t,a)$
upon the substitution $a\mapsto -t^n$ to
$\tilde{J\!D}^n_{r,s}(b\,;\,q,t)$,
though some monomials of the same sign may result in one 
monomial in $\tilde{J\!D}^n_{r,s}(b\,;\,q,t)$ with the coefficient
that is the corresponding sum.

Claim $(ii)$ holds for $n=2$ for the torus knots:
$\{3,2\}, \{5,2\}, \{4,3\},$ $ \{7,3\}, \{13,3\}, \{5,4\}, \{9,4\}, 
\{6,5\}, \{8,5\}, \{9,5\}, \{7,6\}, \{11,6\}.$
The corresponding reduced Khovanov polynomials were calculated
using the program KhReduced from 
http://katlas.org/wiki/. The simplest counterexample
we found is $\{12,7\}$, though $(iii)$ still holds for this knot.
Thus (\ref{nonzerod}) must be not valid for $\{12,7\}$; 
see Section \ref{sec:Khr-p} below. We note that if the order 
in the procedure from (\ref{redprocedure}) is changed to
the opposite one (from the top $a$\~degree), then $\{9,4\}$
and $\{8,5\}$ will not satisfy $(ii)$, though $(iii)$ will
hold for them.

The passage from $\tilde{J\!D}^n_{r,s}(b\,;\,q,t)$ to
$K\!R_{r,s}^{n}(b\,;\,q,t)$ in $(iii)$
is the multiplication of each monomial in 
$\tilde{J\!D}^{n}$ by a proper non-negative
power of $t_{st}$ (the $a_{st}^2$\~degree),
which must be odd if the corresponding coefficient is
negative. If the coefficient of the monomial is greater
than one, it must be split into the sum of monic ones.
After this, the minus-signs must be dropped and similar
terms collected when necessary. 

Knowing these powers of 
$t_{st}$ is the key here for the
total recovery of $K\!R_{r,s}^{n}$; the restricted 
super-polynomials
$H\!D_{r,s}^{[n]}(b\,;\,q,t,a)$ provide them modulo
the conjecture.
\smallskip

I thank Evgeny Gorsky, Sergei Gukov, Mikhail Khovanov, 
Jacob Rasmussen and Lev Rozansky for
multiple talks on these issues; the results of
E.~G. and J.~R. indicate that  (\ref{nonzerod})
is not generally true.
\smallskip 

We note that $(i)$ fails for
the torus knot $\{4,3\}$ and $b=\om_1+\om_2$ (the 3-hook);
the consecutive cancelations of the monomials 
for the {\em neighboring} $a$\~powers only are
simply not sufficient in this case.  However Conjecture 
\ref{CONJKHR}, $(i)$ holds
for  the torus knot $\{3,2\}$ and such $b$. 

There are $3$ exceptional (unstable)
DAHA-KhR polynomials in the latter case.
The formulas for the first two are as follows ({\em in terms
of the standard parameters $q,t$\,}):  
\begin{align*}
&K\!R_{3,2}^{2}(\om_1+\om_2\,;\,q,t)=
1+q^4 t^2-q^6 t^2,\  \ K\!R_{3,2}^{3}(\om_1+\om_2\,;\,q,t)\ =
\end{align*}
\renewcommand{\baselinestretch}{0.5} 
{\small
\(
1+2 q^4 t^2-q^6 t^2+q^6 t^3+2 q^8 t^4-q^{10} t^4+2 q^{10} t^5
-q^{12} t^5+q^{10} t^6+q^{12} t^6-q^{14} t^6+3 q^{14} t^7
-q^{16} t^7+2 q^{16} t^8+2 q^{18} t^9-q^{20} t^9
+q^{20} t^{10}-q^{22} t^{10}+q^{22} t^{11}+q^{24} t^{12}.
\)
}
\renewcommand{\baselinestretch}{1.2} 
\smallskip 

We note that 
$$K\!R_{3,2}^{2}(m\om_1\,;\,q,t)=1=
K\!R_{3,2}^{2}(m\om_2\,;\,q,t) \for m\ge 1.
$$

\subsection{Discussion}
\subsubsection{Confirmations}
The $J\!D(t\mapsto q)$\~polynomials coincide with the QG-Jones 
polynomials for $A_n$ (any weights); it was also checked for 
$k_\nu=1$ and quite a few root systems of type $B,C,D,G_2$. 
The existence ($a$\~stabilization) of the super-DAHA-polynomials 
($A$\~type) is confirmed now; it is actually similar to the 
approach from \cite{SV}; see Lemma 4.3 and formula (4.1) there.
The super-DAHA-polynomials were explicitly computed in the 
following cases:

\smallskip
\noindent
$(i)\ \, \ b=\om_1$ with many torus knots, ``up to" 
$\{12,7\}-\{11,10\}$, mainly
for $A_n$ and for some other classical root 
systems;

\noindent
$(ii) \ \, b=2\om_1, 3\om_1,\ b=\om_2,\om_3$ for the knots 
$\{3,2\}$,\, $\{5,2\}$,\, $\{4,3\},\ldots,\,$  $\{9,4\}$\,, 
only $b=\om_1$ and $b=2\om_1$ were computed for the latter knot);

\noindent
$(iii)$
self-dual $b=2\om_2, \om_1+\om_2$ for the knots 
$\{3,2\},\{5,2\},\{4,3\}$ and $b=2\om_3, 3\om_2$ 
(the $2\times 3$ rectangle) for the knot $\{3,2\}$.

\smallskip
In these cases, the exceptional DAHA-KhR polynomials are all
known and obey Conjecture \ref{CONJKHR}, $(i)$.
For $b=\om_1+\om_2$, negative coefficients appear in
the $a,q,t$\~expansions for both 
knots (see the formulas below). Such negative coefficients
for this and any sufficiently general $b$ seem inevitable 
if Conjecture \ref{CONJEVAL}, is supposed to hold,
more precisely, if these polynomials are supposed to 
 
\smallskip

$(i)\ \,$  be $a$\~polynomials and also polynomials
in terms of $q,t^{\pm1}$,

$(ii)\ $ extend the corresponding HOMFLYPT polynomials 
$\tilde{\HOM}$,

$(iii)$ have the expected evaluations at $\,q=1$ and
$\,t=\, 1\ (=q_{st}^2)\,$

\smallskip
\noindent
for $q_{st}$ from (\ref{qtareli}). 
One can add to this list the duality 
$b\mapsto b^{\tr}$\,; the duality and evaluation
conjectures are actually reflected in $(iii)$. 
Both conjectures hold in all known cases of 
DAHA super-polynomials (including $b=\om_1+\om_2$).
The $a,q,t$\~coefficients are all positive in the examples of
rectangle Young diagrams we calculated, including the 
rectangle $2\times 3$ in the case of trefoil.
\smallskip

The numerical calculations for $A_n$ are based on the 
programs written with participation of Daniel Orr. 

\subsubsection{Physics aspects}
{\em We expect that our polynomials coincide with the
physics $A_n$\~super-polynomials} $\h_{r,s}(b\,;\,q,t,a)$ 
defined for symmetric and wedge powers of the fundamental 
representation. At the moment, no explicit super-formulas are posted
in the literature beyond such $b$ (they are mainly for $b=\om_1$).
See, e.g., \cite{GSV},\cite{DGR}, \cite{GIKV}, \cite{AS} and \cite{GS} 
for definitions and references. The match is solid for $\om_1$,
as well as the coincidence with the posted formulas for 
$b= j\om_1 (j>0)$ (for the trefoil and in some 
other simple torus knots).

A ``universal" explicit formula for the colored
super-poly\-no\-mials we know is that for the trefoil 
and $j\om_1$ (all $j>0$) due to Gukov and 
Sto$\breve{\hbox{s}}$i\'{c}. 
We are thankful for sending it to us; 
it matches the examples $j=1,2,3$ we considered. 
See also \cite{DMMSS}, \cite{GS} and \cite{GIKV}. 
We currently cannot connect the connection of formulas 
(67) and (69) from the latter reference with our calculations 
but it seems not impossible (we thank Sergei Gukov for the 
reference).

Our approach generally provides  ``universal" formulas for any 
given knots and reasonable families of weights (say, for 1-row 
or 2-row diagrams), but
they are expected to be involved even for simple knots. For the
trefoil, the calculations are not too difficult for 2-row
diagrams. 
\smallskip

The approach via the refined Chern-Simons
theory and the BPS states are the major but not the only 
possible physics methods; for instance, the matrix
models can be used here too. The applications
of the refined BPS theory for knots are 
mainly based on the localization formula with some links to 
the Hilbert schemes. Its foundations are actually 
in the classical theory of singularities associated to
simple Lie algebras and their representations. We hope that
the latter can be linked to the DAHA approach.
\smallskip

We note that quite a few instances of super-polynomials were 
actually obtained combinatorially (without much theory) by imposing 
and using their rich (expected) symmetries. 
It is not always possible to distinguish in the literature
on super-polynomials, especially in physics papers, which formulas 
are calculated rigorously and which were found/guessed 
combinatorially on the basis of certain evaluations 
(cf. $(1,ii,iii)$ above). Unless the $a$\~degrees of 
the super-polynomials are {\em a priori} known,
using only finitely many evaluations, which is common
in physics papers, does not fix them
uniquely. Our approach is based on evaluations at
$a=-t^{m}$ for {\em all} sufficiently big $m$, so
it determines these polynomials rigorously.   

The existence 
of the duality $b\leftrightarrow b^{tr}$ was predicted by 
Gukov and Sto$\breve{\hbox{s}}$i\'{c} in \cite{GS}. 
It was supported by explicit calculations with certain
symmetric and wedge powers (for the trefoil and some  
simple knots). The physical motivation was provided in
Section 5.3 there.

In more detail, Conjecture 1.1 in \cite{GS} concerns the 
duality between the super-poly\-no\-mials for the symmetric 
and wedge powers. 
Their expectations for arbitrary Young diagrams are 
summarized in (1.12) and (5.18) there,
including the existence of the triple grading, hence, 
the polynomiality of the $a,q,t$\~expansions of 
super-polynomials and the positivity of the corresponding
(integral) coefficients.

Importantly, the positivity does not always hold for the DAHA 
super-polynomials; see our Conjecture \ref{DUALIT}. 
The simplest example is for the trefoil and the hook 
with 3 boxes; it is considered below in detail. Thus,
the DAHA colored super-polynomials generally are {\em not} 
those predicted in \cite{GS} beyond the symmetric and wedge
powers. Correspondingly, either the physics motivation 
of the duality and other symmetries of super-polynomials
fails for arbitrary Young diagrams or our theory is not 
generally connected with the BPS states. 

\medskip
\subsubsection{Mathematics origins}
Mathematically, the super-polynomials are related to 
the Haiman theory  of the Macdonald polynomials via the 
Hilbert schemes of $\C^2$ (started with Garsia and developed 
since then with/by quite a few researches). The latest breakthrough 
project due to Gorsky, Oblomkov, Rasmussen and Schende 
\cite{ORS,Gor} connects the super-polynomials with the 
Hilbert schemes of 
singular curves and rational DAHA. 
{\em A coincidence of our super-polynomials (for torus knots)
with the corresponding (stable) Poincar\'e polynomials} 
is expected, when the latter are well-defined. 
We mention that the Hilbert schemes appear in the refined BPS 
theory, sometimes in equivalent forms.

I am thankful to Evgeny Gorsky for multiple talks on their and 
his own ongoing theories, prior and related directions. 
Gorsky's approach (see \cite{Gor} and Appendix of \cite{ORS}) 
is based on the $0$\~dimensio\-nal formal deformation
of the torus singularity $x^r=y^s$. 
It provides a uniform interpretation of the
super-polynomials for $b=\om_1$ in terms of the 
combinatorics of perfect rational DAHA modules of
type $A_{n-1}$ with $k=-m/n$ for the torus knot 
$K_{m,n}$. The decomposition of such modules 
in terms of the action of the spherical (symmetric) DAHA 
subalgebra is the key; it results in combinatorial
formulas for super-polynomials for any torus knots ($b=\om_1$).
This approach has a clear potential 
of reaching $b=j\om_1$ via more advanced theory of 
rational DAHA modules.

The key motivation of the super-polynomials is their
expected coincidence with the Khovanov-Rozansky polynomials for 
sufficiently large matrix dimensions $N$
(see Conjectures 1.7, 3.1  from \cite{DGR} and \cite{KhR2}).
It is connected with the latest conjectural extension of
the Khovanov-Rozansky theory via the categorification
of arbitrary irreducible representations of quantum groups; 
see \cite{We} and references therein. The polynomiality 
of the resulting Poincar\'e series is not established 
(and not expected) for generic representations, which may match 
the appearance of the negative coefficients of super-polynomials
in our approach.

I thank Mikhail Khovanov
and Lev Rozansky for multiple discussions of their theory.
Let us also mention an important 
physics-mathematics interpretation of the super-polynomials via
the classical theory of singularities; I am thankful to
Sergei Gukov for explaining it to me.
  
\medskip
\subsubsection{Parameters}
The parameters $\,a,q,t\,$ that 
are mostly used in the recent mathematics and physics works 
on super-polynomials, we call them {\em \,standard\,}, 
are different from our ones. The passage is as follows:
\begin{align}\label{qtarel}
q^2_{st}=t,\  t_{st}=\sqrt{q/t},\  a_{st}^2=a\sqrt{t/q}.
\end{align}
For example, $t^2_{st}=\pm 1$  
corresponds to our $t=\pm q$. We have already used
the inverse of this transformation in (\ref{qtareli}).

Such parameters are used, for instance,
in the formulas from \cite{AS}, those calculated/posted 
by Gorsky (coinciding with our ones in all checked cases), 
those from \cite{Sh} (coinciding with ours apart from 
$K_{8,5}$), in \cite{DGR} and in recent \cite{GS}; see
also \cite{DMMSS}. For instance, formula (23) 
from \cite{DGR} in the tilde normalization results
from our 
$$
K\!R^{2}_{4,3}(\om_1\,;\, q,t)= 1+ qt+qt^2+qt^2\sqrt{q/t} +
q^2t^3\sqrt{q/t}.
$$ 

Note that the HOMFLYPT polynomial 
$\HOM_{r,s}(b\,;\,q,a)$ 
becomes the Jones polynomial $\j_{r,s}(b\,;\,q)$ as $a=q^{n+1}$
in our notations. There are other normalizations used 
in the vast mathematical and physical literature on the
knot invariants, more specifically, on the Khovanov and 
Rozansky theory and super-polynomials. 

Changing the parameters is not always plain and square.
For instance, the (conjectural) relations of DAHA 
super-polynomials 
to the classical HOMFLY polynomials in terms of our 
$q,t$ in Conjecture \ref{HOMFLY} is not exactly a 
reformulation of the known implications. It is combined 
with the (expected) fact that the parity of $a_{st}^2$ 
and $t_{st}$ is the same in the monomials of the known 
formulas for the super-polynomials. Another example
is the duality conjecture, where the choice of the 
parameters reflects the filtrations of the corresponding
homology theory, so has a clear topological-geometric
meaning. 

Generally, the standard parameters from
(\ref{qtarel}) are adjusted to the  
natural filtrations 
in the Khovanov-Rozansky theory and in other topological 
and geometric constructions. The  DAHA-parameters  
(algebraically, quite natural) are sufficient and
convenient to 
formulate the conjectures above. We note that
the duality conjecture may suggest a connection 
with the symmetry $q\leftrightarrow t$
of the Macdonald $GL$\~stable polynomials $\tilde{H}_b$,
though a different transposition, $q\leftrightarrow t^{-1}$
($a\mapsto a$), is exactly what is needed for the 
super-polynomials.

\smallskip
\subsubsection{Using other coinvariants}
We note that  $r,s$ are allowed to 
be negative in Conjecture \ref{MAINCONJ}. 
The corresponding tilde-Jones polynomial $\tilde{\j}$
coincides with that for $|r|,|s|$. See (\ref{jones-formula}).
It is directly connected with the $\phi$\~invariance of
the coinvariant we use. This invariance also provides the
compatibility with the relations from the projective
$PSL_2(\Z)$; we will leave a systematic discussion
of this and such matters for the sequel of this paper.
Almost always we use  $r>s>0$ ($(r,s)=1$)
in this paper.
\smallskip

It is important that the functional  $\{\cdot\}$ can be replaced 
by $\lan\cdot\ran_1'$ from (\ref{ctgafunct}). It 
results in changing $\tga$ by $\tau_-^{-1}\tga$ inside $\{\cdot\}$ 
in (\ref{jones-d}), i.e., in changing the corresponding
torus knot. The reason for this is 
the difference theory of Macdonald-Mehta integrals from
\cite{C101}; see (\ref{evaga}) above.

Using this functional is expected to establish a link to the 
physics research on super-polynomials based on matrix models. 
We will not discuss this approach here.
Algebraically, these two functionals are the only non-equivalent 
DAHA {\em coinvariants} of level $1$.  
The functional (\ref{ctfunct}) is the only DAHA-coinvariant
of level zero; using it, presumably,
corresponds to the theory of super-polynomials in
$S^2\times S^1$ (cf. \cite{AS}).

The vector space of general DAHA coinvariants is isomorphic to
the corresponding Looijenga space \cite{C102},
spanned by the Kac-Moody characters of the corresponding level.
This is directly related to the theory of elliptic
Hall functions ({\em ibid.}). 

Topologically, we glue two solid tori twisting their 
boundaries by $\ga,\ga'\in PSL_2(\Z)$. When the product
$\ga'\ga$ equals $\tau_-$ or $\si$ (in $PSL_2(\Z)$)
we obtain $S^3$; the image of the circumference 
at the boundary of the first solid torus 
associated with the $X$\~direction
is $K_{r,s}$ corresponding to $\ga$.
These two cases somehow match using $\{\cdot\}$ and 
$\lan\cdot\ran'_1$. The resulting space will be a 
{\em lens space} for arbitrary $\ga, \ga'$.
The DAHA coinvariants of levels 
greater than $1$ in the sense of \cite{C102} are expected
to serve the Jones-type invariants of torus knots in these spaces
and the corresponding super-polynomials.
I am thankful to Lev Rozansky for discussions of 
these and related topological matters.

\medskip

\subsection{The rank one case}\label{sect:rank-one}
Let $\al=\al_1$, $s=s_1$ and $\om=\om_1$, the fundamental weight;
then $\alpha=\al_1=2\omega$ and
$\rho=\om$. The extended affine Weyl group
$\widehat{W}=<s,\om>$
can be presented as the free group generated by
the involutions  $s$ and $\pi=\om s$.
We will denote the weights $b\om$ ($b\in \Z$)
simply by $b$.

The double affine Hecke algebra
$\HH$ is generated by $Y=Y_{\om_1}=\pi T, T=T_1, X=X_{\om_1}$
subject to the quadratic relation $(T-t^{1/2})(T+t^{-1/2})=0$
and the cross-relations:
\begin{align}\label{dahaone}
&TXT=X^{-1},\ T^{-1}YT^{-1}=Y^{-1},\ Y^{-1}X^{-1}YXT^2q^{1/2}=1.
\end{align}
Setting $\pi\equal YT^{-1}$, the second relation 
becomes $\pi^2=1$.
The field of definition will be $\Q(q^{1/4},t^{1/2})$
although $\Z[q^{\pm 1/4},t^{\pm 1/2}]$ is sufficient
for many constructions.
Here $q^{\pm 1/4}$ is needed in the automorphisms $\tau_{\pm}$:
\begin{align}\label{tau+def}
&\tau_+(X)=X,\ \ \tau_+(T)=T,\ \ \tau_+(Y)=q^{-1/4}XY,\ \\
\label{tau-def}
&\tau_-(Y)=Y,\ \tau_-(T)=T,\ \tau_-(X)=q^{1/4}YX,\ 
\tau_-(\pi)=\pi.
\end{align}

For $r,s\ge 0$, it is not difficult to check that 
$q^{b^2rs/4}t^{b(r+s-1)/2}J\!D_{r,s}(b)$
is a polynomial in terms of the non-negative 
integral powers of $q,t$ with the constant term 
$1$, which addresses a natural question concerning
the exact normalization needed in Conjecture \ref{MAINCONJ}.

In this section we will involve the
non-reduced DAHA-Jones invariants 
\begin{align}\label{nonrjone}
&J\!D^\#_{r,s}(b\,;\,q,t) \equal \{\,\tga\bigl(P^{(k)}_b)\, \}_k\,,
\for t \equal q^k,\, \ga=\ga_{r,s},
\end{align}
where we show explicitly the dependence of the $P$\~polynomial
and the evaluation coinvariant on the parameter $k$.
Generally, they are inconvenient because of
non-trivial $t$\~denominators. However, they
become $q$\~polynomials as $t=q$ and somewhat simplify 
the considerations. 
\smallskip

The formula for the colored Jones polynomials
of the torus knot $K_{r,s}$ for the representation
of weight $b\ge 0$ and dimension $b+1$
is well known; see \cite{Mo} and also \cite{Hi}. For $r,s\ge 0$,
\begin{align}\label{jones-formula}
&\j^{\#}_{r,s}(b\,;\,q)=\frac{q^{\frac{rsb(b+2)}{2}}}
{q^{\frac{1}{2}}-q^{-\frac{1}{2}}}
\sum_{p=-b/2}^{b/2}\bigr(q^{rsp^2+(r+s)p+\frac{1}{2}}-
q^{rsp^2+(r-s)p-\frac{1}{2}}\bigl),
\end{align}
where the $p$\~summation step is $1$  and
we set 
$$
\j^{\#}_{r,s}(b\,;\;q)\equal
\frac{q^{(b+1)/2}-q^{-(b+1)/2}}{q^{1/2}-q^{-1/2}}
\,\j_{r,s}(b),
$$
i.e., normalize $\j_{unknot}^{\#}(b)$ to be the $q$\~dimension
for $b=b\om_1$. The greatest positive power of $q$ in 
(\ref{jones-formula}) is  obviously 
$(\frac{3}{4}b+1)rsb+(r+s)b$.

In this normalization, Conjecture \ref{MAINCONJ} 
becomes as follows:
\begin{align}\label{conjaone}
\j^{\#}_{r,s}(b\,;\, q\mapsto q^{-1})=
q^{-\frac{rsb(b+2)}{2}}\,J\!D_{r,s}^\#(b\,;\,q,t\mapsto q).
\end{align}
It is equivalent to the following
proposition.

\begin{proposition}\label{MAINCONJ1}
In the notation from (\ref{nonrjone}),
\begin{align}\label{jonesf-1}
&J\!D_{r,s}^\#(b\,;\,q,q)= 
\sum_{p=-b/2}^{b/2}
\frac{q^{-rsp^2/4}}{q^{-1/2}-q^{1/2}}
\bigr(q^{(r+s)p-1/2}-
q^{(r-s)p+1/2}\bigl).
\end{align}
In particular for $b=1$,
\begin{align}\label{jones1b1}
&J\!D_{r,s}^\#(1\,;\,q,q)= 
\frac{q^{-rs/4-(r+s)/2}}{1-q}
(1-q^{r+1}-q^{s+1}+q^{r+s}),\\
&\hbox{and\ \ \ \, }
\tilde{J\!D}_{r,s}(b=1\,;\,q,q)\ \,=\ \,
\frac{1-q^{r+1}-q^{s+1}+q^{r+s}}{1-q^2}.\notag
\end{align}
\end{proposition}
{\it Proof.}
We formally set $X_\la=q^{\la x}$ for $\la\in \C$,
which extends $X_b=X_{b\om_1}=X^b$ for $b$ in $\Z$
identified with the lattice $P=\Z\om_1$.
The {\em free action} of $SL(2,\Z)$ 
on $X_\la\,q^{M x^2}$ from (\ref{sl2Mq})
reads as follows:

\begin{align}\label{hankelfunq}
&\si^\circ(X_\la\,q^{-Mx^2})\,=\,
\frac{q^{\frac{\la^2}{4M}}}{M^{1/2}}\,
X_{\la/M}\,q^{+x^2/M} \for M\neq 0,\\
&(\tau_+^\circ)^N\,(X_\la\,q^{-Mx^2})\ =\ X_\la\,q^{(N-M)x^2}\,
\and \notag\\ 
&(\tau_-^\circ)^N\,(X_\la\,q^{-Mx^2})\ =\ 
((\si^\circ)^{-1}\,(\tau_+^\circ)^{-N}\,
\si^\circ)(X_\la\,e^{-Mx^2})
\notag\\
&=\ \frac{1}{(1-MN)^{1/2}}\,
q^{\frac{\la^2\,N}{4(1-MN)}}\
X_{\la/(1-MN)}\, q^{-x^2\frac{M}{1-MN}}.\notag
\end{align}

Here $N\in \Z$; the parameter $M$ is a sufficiently
general complex number, for instance $MN\neq 1$ 
in the last formula. We also assume that $(MM')^{1/2}$
$= M^{1/2}(M')^{1/2}$ for $M,M'\in \C$. The general $SL(2,\Z)$\~action
from (\ref{sl2zq}) is adjusted to $A_1$ as follows:

\begin{align}\label{gamfunq}
&\ga^\circ\,(X_\la\,q^{z x^2})\,=\,
\frac{1}{(cz+d)^{1/2}}\ q^{-\frac{\la^2\, c}{4(cz+d)}}\
X_{\frac{\la}{cz+d}}\ q^{\frac{az+b}{cz+d}\,x^2} \for\\
&\ga= \left(
  \begin{array}{cc}
    a & b \\
    c & d \\
  \end{array}
\right)\in SL(2,\Z),\ \la \in \C 
\hbox{\ and for generic\, } z\in \C\,. \notag
\end{align}
We emphasize that this action does not depend on $t$.

We will need the following formula:
\begin{align}\label{garsmfunq}
&\ga^\circ\,\bigl(X_\la\,
((\ga^{-1})^\circ(X_\mu))\,\bigr)
=\ q^{-ac\,\la^2/4-c\la\mu/2}\, X_{a\la+\mu}\,.
\end{align}


Using this action and Proposition \ref{SH-JONES},
\begin{align}\label{jonesf-1r}
&J\!D_{r,s}^\#(b\,;\,q,t\mapsto q)\ =\ R_{r,s}^b(X\mapsto q^{-1/2})
\for\\ 
&R_{r,s}^b(X)\ \equal\ (X-X^{-1})^{-1}\ga^\circ \,\bigl(P_b^{(1)}\,
(\ga^\circ)^{-1}(X-X^{-1})\,\bigr),\hbox{ where}\notag\\
&P^{(1)}_b=\frac{X^{b+1}-X^{-b-1}}{X-X^{-1}}=
\sum_{p=-b/2}^{b/2}X^{2p}\ (\hbox{ step }=\ 1) \for b\in \Z_+.
\notag
\end{align}
Now applying (\ref{garsmfunq}) for $\mu=\pm 1$ and
$\la=b,b-2,\cdots,-b$,\, we readily arrive at 
formula (\ref{jonesf-1}).
\sq

\comment{
The summation step is $1$ in (\ref{jones-formula}).
The following recursion from \cite{Hi} is expected to
be used to check Conjecture \ref{MAINCONJ} in this
case:
\begin{align}
\j^\#_{r,s}(b)=&q^{\frac{(rs-r-s)b}{2}}
\frac{1-q^{rb+1}-q^{sb+1}+q^{(r+s)b}}{1-q}+
q^{\frac{rsb}{2}}\j^\#_{r,s}(b-2),\notag\\
\label{jones-rec}
&\j^\#_{r,s}(0)=1,\ 
\j^\#_{r,s}(1)=q^{\frac{rs-r-s}{2}}
\frac{1-q^{r+1}-q^{s+1}+q^{r+s}}{1-q}.
\end{align}
\medskip

For $r,s\ge 0$, it is not difficult to check that 
$q^{b^2rs/4}t^{b(r+s-1)/2}J\!D_{r,s}(b)$
is a polynomial in terms of the non-negative 
integral powers of $q,t$ with the constant term 
$1$, which addresses a natural question concerning
the exact normalization needed in Conjecture \ref{MAINCONJ}.

For $b=1$ and $r,s\ge 0$ (for the sake
of definiteness), this conjecture becomes:
\begin{align*}
&J\!D_{r,s}(1\,;\,t\mapsto q)=q^{-\frac{rs-2(r+s-1)}{4}}\,
\frac{1-q^{r+1}-q^{s+1}+q^{r+s}}{1-q}=
q^{-3\frac{rs}{4}}\,\j_{r,s}(1).
\end{align*}
}
\medskip

Let us go back to arbitrary $t$
and give the formula for the torus knot $K_{9,4}$
for $b=1$. One has:
\begin{align}\label{jd9-4}
&\tilde{J\!D}_{9,4}(1\,;\,q,t)=1 + (1 - t)\bigl(qt  + q^2 t(1 + t)\\
&+ q^3 t (1 + t) + 2q^4 t^2 + q^5 t^2 (1 + t - t^2)\notag\\
&+ q^6 t^2 (1 + 2t) + q^7 t^3 (1 - t^2) + q^8 t^4  
\bigr).\notag
\end{align}
This torus knot corresponds to $\tga=\tau_+^2\tau_-^4$.
The substitution $t\mapsto q$ results in many reductions here,
including the cancelation of the (potential) leading term $q^{12}$:
\begin{align}\label{j9-4}
\tilde{J\!D}_{9,4}(1\,;\,t\mapsto q)\ &=\ 
\tilde{\j}_{9,4}(1)\ =\ (1-q^5-q^{10}+q^{13})/(1-q^2)\notag\\
&=\ 1+q^2+q^4-q^5+q^6-q^7+q^8-q^9-q^{11}.
\end{align}

For the torus knot $K_{8,5}$, the corresponding product
of $\tau$\~matrices is $\tau_+\tau_-\tau_+\tau_-^2$ and
\begin{align}\label{jd8-5}
&\tilde{J\!D}_{8,5}(1)=
1+q t+q^2 t+q^3 t+q^4 t-q t^2\\
+q^4 t^2+2 q^5 t^2&+2 q^6 t^2-q^2 t^3-q^3 t^3-2 q^4 t^3-2 q^5 t^3
+2 q^7 t^3+q^8 t^3\notag\\
-q^5 t^4-3 q^6 t^4&-
3 q^7 t^4+q^9 t^4+q^5 t^5+q^6 t^5-2 q^8 t^5-q^9 t^5+q^7 t^6+q^8 t^6.
\notag
\end{align}

Let us provide the formula for $\tilde{J\!D}_{9,4}(2\,;q,t)$
defined (see above) as 
$$
\{\,\tau_+^2\tau_-^4\bigl(P^{(k)}_2/P^{(k)}_2(t^{1/2})\bigr)\,\}_k,
$$
divided then by a proper (fractional) power of $q,t=q^k$ to make 
it a series in terms of positive powers of $q,t$ starting
with $1$. It equals

\renewcommand{\baselinestretch}{0.5} 
{\small
\noindent
\( 
1+q^2 t+q^3 t+q^4 t+q^5 t+q^6 t+q^7 t-q^2 t^2-q^3 t^2+q^6 t^2+
q^7 t^2+3 q^8 t^2+3 q^9 t^2+3 q^{10} t^2+2 q^{11} t^2+2 q^{12} t^2+
q^{13} t^2+q^{14} t^2-q^4 t^3-2 q^5 t^3-2 q^6 t^3-3 q^7 t^3
-3 q^8 t^3-3 q^9 t^3-q^{10} t^3+q^{11} t^3+3 q^{12} t^3+4 q^{13} t^3
+5 q^{14} t^3+6 q^{15} t^3+4 q^{16} t^3+3 q^{17} t^3+2 q^{18} t^3
+q^{19} t^3+q^5 t^4-q^8 t^4-2 q^9 t^4-5 q^{10} t^4-7 q^{11} t^4
-8 q^{12} t^4-8 q^{13} t^4-7 q^{14} t^4-5 q^{15} t^4+3 q^{17} t^4
+6 q^{18} t^4+6 q^{19} t^4+7 q^{20} t^4+4 q^{21} t^4+3 q^{22} t^4
+q^{23} t^4+q^{24} t^4+q^7 t^5+q^8 t^5+2 q^9 t^5+3 q^{10} t^5
+3 q^{11} t^5+2 q^{12} t^5-q^{13} t^5-4 q^{14} t^5-8 q^{15} t^5
-12 q^{16} t^5-16 q^{17} t^5-14 q^{18} t^5-11 q^{19} t^5
-6 q^{20} t^5
+4 q^{22} t^5+6 q^{23} t^5+5 q^{24} t^5+4 q^{25} t^5+2 q^{26} t^5
+q^{27} t^5+q^{11} t^6+2 q^{12} t^6+5 q^{13} t^6+6 q^{14} t^6
+8 q^{15} t^6+8 q^{16} t^6+8 q^{17} t^6+q^{18} t^6-5 q^{19} t^6
-12 q^{20} t^6-16 q^{21} t^6-18 q^{22} t^6-14 q^{23} t^6
-7 q^{24} t^6
-2 q^{25} t^6+3 q^{26} t^6+3 q^{27} t^6+3 q^{28} t^6+q^{29} t^6
+q^{30} t^6-q^{12} t^7-q^{13} t^7-q^{14} t^7+q^{16} t^7+4 q^{17} t^7
+8 q^{18} t^7+13 q^{19} t^7+15 q^{20} t^7+13 q^{21} t^7
+9 q^{22} t^7-8 q^{24} t^7-14 q^{25} t^7-13 q^{26} t^7-9 q^{27} t^7
-3 q^{28} t^7+q^{29} t^7+q^{30} t^7+q^{31} t^7-q^{15} t^8-q^{16} t^8
-2 q^{17} t^8-3 q^{18} t^8-4 q^{19} t^8-3 q^{20} t^8+2 q^{21} t^8
+6 q^{22} t^8+12 q^{23} t^8+15 q^{24} t^8+16 q^{25} t^8+7 q^{26} t^8
+q^{27} t^8-7 q^{28} t^8-8 q^{29} t^8-4 q^{30} t^8-q^{31} t^8
+q^{32} t^8
-q^{20} t^9-3 q^{21} t^9-5 q^{22} t^9-5 q^{23} t^9-6 q^{24} t^9
-2 q^{25} t^9+4 q^{26} t^9+9 q^{27} t^9+11 q^{28} t^9+7 q^{29} t^9
+q^{30} t^9-4 q^{31} t^9-2 q^{32} t^9+q^{22} t^{10}-2 q^{25} t^{10}
-4 q^{26} t^{10}-6 q^{27} t^{10}-4 q^{28} t^{10}+3 q^{30} t^{10}
+6 q^{31} t^{10}+q^{32} t^{10}-q^{33} t^{10}-q^{34} t^{10}
+q^{26} t^{11}
+q^{27} t^{11}-q^{29} t^{11}-3 q^{30} t^{11}-2 q^{31} t^{11}
+2 q^{33} t^{11}+q^{34} t^{11}-q^{35} t^{11}+q^{30} t^{12}
-q^{33} t^{12}+q^{35} t^{12}.
\)
}
\renewcommand{\baselinestretch}{1.2} 

\smallskip
Here the Macdonald polynomial $P_2$ and its evaluation are 
$$
P_2=X^2+X^{-2}+\frac{(1-t)(1+q)}{1-qt},\
P_2(t^{1/2})=\frac{(1-qt^2)(1+t)}{t(1-qt)}.
$$ 
The corresponding reduced tilde-normalized 
Jones polynomial $\tilde{\j}_{9,4}(2)$ (in 
our notations) is
$$
(1-q^9-q^{19}+q^{26}+q^{49}-q^{50})/(1-q^3).
$$ 
\medskip

\subsection{The rational limit}
We will conclude this section with the calculation
of the rational limit of $J\!D_{r,s}(b)$ in the 
case of $A_1$. It can be dealt directly within the
rational DAHA theory. The construction is
general, but the rank one case is especially
simple due to fact that the corresponding eigenfunction
function is a function of the product of the argument
and the spectral parameter.

Generally,
the corresponding limiting procedure requires sending 
$b\to\infty$ together with making $q\mapsto 1$ subject
to $t=q^k$. It results in an important construction
of the Bessel functions from the Macdonald (Jack-
Heckman- Opdam) polynomials. 
We will not discuss here this limit and remind the reader
the definition of rational DAHA; only the action of
$\tau_{\pm}$ on the Bessel functions is really needed. 

The limit can be expected
to catch the leading term of $\j_{r,s}(b)$,  
which is $q^{b^2\,rs/4}$; it really does! We mainly
follow Section 2.4 of \cite{C101}. 
Only the case of $A_1$ will be considered.

We set $\phk_\lambda (x)\equal
\phk (\lambda x)$ for the Bessel-type function
$$
\phk (t)=\sum_{n=0}^\infty \frac{t^{2n}\Gamma (k+1/2)}{n!
\Gamma (k+n+1/2)},
\; \C\ni k\not\in -1/2-\Z_+.
$$
See formula (2.1.15) from \cite{C101}. It is a unique 
{\em even} solution of the eigenvalue problem:
\begin{align*}
L\phi_\la\ =\ 4\la^2 \phi_\la,\  \phi(0)=1\, \for\, 
L\equal\frac{d^2}{dx^2}+\frac{2k}{x}\frac{d}{dx}\,.
\end{align*}

Then, 
\begin{align}\label{tauoper} 
\tau_+&(f(x))\,=\, e^{x^2}f(x),\ \ \tau_+(A)\ =\ e^{x^2}\circ 
A\, \circ e^{-x^2},\\
\tau_-&(f(x))=e^{-L/4}(f(x)),\ 
\tau_-(A)=e^{-\frac{L}{4}}\circ 
A\, \circ e^{\frac{L}{4}}\label{Ltaum}
\end{align}
for functions $f$ and operators $A$.

The {\em Hankel transform} is defined as
$\si\equal\tau_+\tau_-^{-1}\tau_+.$ The calculations
below mainly use this algebraic definition; the
original (classical) definition is as follows:
\begin{equation}\label{fhankel}
\si(f)(\mu)\equal \frac{1}{\Gamma (k+\frac{1}{2})}
\int_{-\infty}^{\infty}
f(x)\phk_\mu (x)\,|x|^{\,2k}dx, \ k\not\in -\frac{1}{2}-\Z_+.
\end{equation}

The rational variant of the $J\!D_{r,s}$ depends
on $\la\in\C$, which replaces $b\in \Z_+$ in the 
$q$\~theory:
\begin{align}\label{ratjones}
J\!D^{rat}_{r,s}(\la)\,\equal\,&\Bigl(\tga_{r,s}\bigl(
\,(\phk_\la)^{op}\,\bigr)
\Bigr)(1)\mid_{\,x\mapsto 0}\notag\\
=\,&\Bigl(\tga_{r,s}\bigl(\phk_\lambda(x) 
\tga_{r,s}^{-1}(1)\bigr)\Bigr)\mid_{\,x\mapsto 0}.
\end{align}
We will see below that it does not depend on $k\in \C$ due to  
our special choice of the input function $\phk$ in
(\ref{ratjones}).
In the first line, $(\phk_\la)^{op}$ is the operator
of multiplication by  $\phk_\lambda(x)$;
then we conjugate it by $\tga$, 
apply the result to the function $1$ and, finally, 
evaluate the output (a function) at $x=0$. 

The following formulas are needed to perform
the calculation of $J\!D^{rat}$ (see Theorem 2.4.1
from \cite{C101} and formula (\ref{hankelfunq}) above):
\begin{align}\label{hankelfun}
&\si(\phk_\la(x)e^{-Mx^2})\,=\,\frac{e^{\la^2/M}}{M^{k+1/2}}\,
\phk_{\la/M}(x)\,e^{+x^2/M} \for M\neq 0,\\
&\tau_+^N\,(\phk_\la(x)\,e^{-Mx^2})\ =\ \phk_\la(x)\,e^{(N-M)x^2}\,
\and \notag\\ 
&\tau_-^N\,(\phk_\la(x)\,e^{-Mx^2})\ =\ 
(\si^{-1}\,\tau_+^{-N}\,\si)(\phk_\la(x)\,e^{-Mx^2})\\
\notag
&=\ \frac{1}{(1-MN)^{k+1/2}}\,
e^{-\la^2 \frac{N}{1-MN}}\
\phk_{\la/(1-MN)}(x)\, e^{-x^2\frac{M}{1-MN}}.\notag
\end{align}

Here $N\in \Z$; the parameter $M$ is allowed to
be a complex number, provided that $(MM')^{k+1/2}= 
M^{k+1/2}(M')^{k+1/2}$ and $M\neq 1/N$ in the
last formula.
We disregard the convergence matters
here; for instance, $\Re M>0$ must hold if we take  
(\ref{fhankel}) as the definition of $\si$ and 
$\Re (1/M)>N$ is necessary for $\tau_-^N$. The only fact we
really need is that these formulas, treated formally, can be 
extended naturally to all $\tga$.  

The corresponding formula is
{\em directly} connected with the action of $SL_2(\Z)$ on
the modular functions of weight $k+1/2$ when $\la=0$.
For any 
$\ga= \left(
  \begin{array}{cc}
    a & b \\
    c & d \\
  \end{array}
\right)\in SL(2,\Z)$, complex $k\not\in -1/2-\Z_+$,
arbitrary $\la \in \C$ and sufficiently general $z\in \C$,
the formula reads as follows:

\begin{align}\label{gamfun}
&\ga\,(\phk_\la(x)\,e^{z x^2})\,=\,
\frac{1}{(cz+d)^{k+1/2}}\ e^{-\frac{c}{cz+d}\,\la^2}\
\phk_{\frac{\la}{cz+d}}(x)\ e^{\frac{az+b}{cz+d}\,x^2}.
\end{align}
Cf. (\ref{gamfunq}). Similar to (\ref{garsmfunq}),

\begin{align}\label{garsmfun}
&\ga\,\bigl(\phk_\la(x)\,
(\ga^{-1}(1))\,\bigr)
=\ e^{-ac\,\la^2}\, \phk_{a\la}\,.
\end{align}


Concerning the meaning of (\ref{gamfun}), one can
present it as a set of entirely algebraic (and quite
nontrivial) coefficient-wise identities upon the 
expansion in terms of the non-negative powers of 
$z$ and $\la$. We note that $\ga$ can be assumed 
from $SL_2(\C)$ in the rational theory.

These considerations readily result in the 
following proposition.

\begin{proposition}\label{JONESRAT}
Let $\tga_{r,s}$ correspond to a 
matrix $\ga\in PSL_2(\Z)$
with the first column $(r,s)^{tr}$, representing 
the torus knot $K_{r,s}$. 
In the setting above,
\begin{align}\label{jonratf}
&\Bigl(\tga_{r,s}\bigl(\,(\phk_\la)^{op}\,\bigr)
\Bigr)(1)\ =\ e^{-\la^2\, rs}\,\phk_{r\la},\notag\\
&\ \ J\!D^{rat}_{r,s}(\la)\ \, =\, \ e^{-\la^2\, rs} \for \la\in \C.
\end{align}
\sq
\end{proposition}

This rational DAHA limit, namely sending $b$ to $\infty$,
is similar to that from the so-called
{\em volume conjecture}. It states that the limit $N\equal
b+1\to \infty$
of 
\begin{align*}
&\frac{\log\,|\j\bigl(b\,;\,q\mapsto\exp(2\pi i/N)\bigr)|}
{N}\ =\\
&\frac{1}{N}\,\log\,\mid\frac{(q^{1/2}-q^{-1/2})\,
\j^\#\bigl(b\,;\,q\mapsto\exp(2\pi i/N)\bigr)}
{(q^{N/2}-q^{-N/2})} \mid
\end{align*}
coincides up to a global constant with
the hyperbolic volume of $S^3\setminus K$,
which is zero for the torus knots.
This limit is known to be $0$ too for the torus knots 
due to the special choice of $q$, so the coincidence holds. 
See \cite{MM}, Conjecture 5.1 and discussion there. 

We note that this limit vanishes for any root systems
if the DAHA -formula (\ref{jones-d}) for $t_\nu=q_\nu$ 
is used. This fact and its extensions to arbitrary
$t$ requires the theory of the polynomial representation
at roots of unity, which we do not discuss in this paper.

\setcounter{equation}{0}
\section{Numerical examples}
\subsection{Uncolored super-formulas}
We want to demonstrate what can be expected and provide
material for experiments of the advanced readers.
Let us begin with the formula for the DAHA super-polynomial 
$$
H\!D_{8,5}(b=\om_1;q,t,a)\ =
$$

\renewcommand{\baselinestretch}{0.5} 
{\small
\noindent
\(
a^0 \bigl( 1+q t+q^2 t+q^3 t+q^4 t+q^2 t^2+q^3 t^2+2 q^4 t^2
+2 q^5 t^2+2 q^6 t^2+q^3 t^3+q^4 t^3+2 q^5 t^3+3 q^6 t^3
+3 q^7 t^3+q^8 t^3+q^4 t^4+q^5
t^4+2 q^6 t^4+3 q^7 t^4+4 q^8 t^4+2 q^9 t^4+q^5 t^5+q^6 t^5
+2 q^7 t^5+3 q^8 t^5+4 q^9 t^5+2 q^{10} t^5+q^6 t^6+q^7 t^6
+2 q^8 t^6+3 q^9 t^6+4 q^{10}
t^6+q^{11} t^6+q^7 t^7+q^8 t^7+2 q^9 t^7+3 q^{10} t^7
+3 q^{11} t^7+q^8 t^8+q^9 t^8+2 q^{10} t^8+3 q^{11} t^8
+2 q^{12} t^8+q^9 t^9+q^{10} t^9+2 q^{11}
t^9+2 q^{12} t^9+q^{10} t^{10}+q^{11} t^{10}+2 q^{12} t^{10}
+q^{13} t^{10}+q^{11} t^{11}+q^{12} t^{11}+q^{13} t^{11}
+q^{12} t^{12}+q^{13} t^{12}+q^{13}
t^{13}+q^{14} t^{14}\bigr)
\)
\vfill
\noindent
\(
+a^1 \bigl(q+q^2+q^3+q^4+q^2 t+2 q^3 t+3 q^4 t+4 q^5 t+3 q^6 t
+q^7 t+q^3 t^2+2 q^4 t^2+4 q^5 t^2+6 q^6 t^2+7
q^7 t^2+4 q^8 t^2+q^9 t^2+q^4 t^3+2 q^5 t^3+4 q^6 t^3+7 q^7 t^3
+9 q^8 t^3+7 q^9 t^3+2 q^{10} t^3+q^5 t^4+2 q^6 t^4+4 q^7 t^4
+7 q^8 t^4+10 q^9 t^4+8
q^{10} t^4+2 q^{11} t^4+q^6 t^5+2 q^7 t^5+4 q^8 t^5+7 q^9 t^5
+10 q^{10} t^5+7 q^{11} t^5+q^{12} t^5+q^7 t^6+2 q^8 t^6
+4 q^9 t^6+7 q^{10} t^6+9 q^{11}
t^6+4 q^{12} t^6+q^8 t^7+2 q^9 t^7+4 q^{10} t^7+7 q^{11} t^7
+7 q^{12} t^7+q^{13} t^7+q^9 t^8+2 q^{10} t^8+4 q^{11} t^8
+6 q^{12} t^8+3 q^{13} t^8+q^{10}
t^9+2 q^{11} t^9+4 q^{12} t^9+4 q^{13} t^9+q^{11} t^{10}
+2 q^{12} t^{10}+3 q^{13} t^{10}+q^{14} t^{10}+q^{12} t^{11}
+2 q^{13} t^{11}+q^{14} t^{11}+q^{13}
t^{12}+q^{14} t^{12}+q^{14} t^{13}\bigr)
\)
\vfill
\noindent
\(+a^2 \bigl(q^3+q^4+2 q^5+q^6+q^7+q^4 t+2 q^5 t+4 q^6 t+5 q^7 t
+4 q^8 t+2 q^9 t+q^5 t^2+2 q^6 t^2+5 q^7 t^2+7 q^8 t^2+8 q^9 t^2
+4 q^{10}
t^2+q^{11} t^2+q^6 t^3+2 q^7 t^3+5 q^8 t^3+8 q^9 t^3+10 q^{10} t^3
+5 q^{11} t^3+q^{12} t^3+q^7 t^4+2 q^8 t^4+5 q^9 t^4+8 q^{10} t^4
+10 q^{11} t^4+4
q^{12} t^4+q^8 t^5+2 q^9 t^5+5 q^{10} t^5+8 q^{11} t^5
+8 q^{12} t^5+2 q^{13} t^5+q^9 t^6+2 q^{10} t^6+5 q^{11} t^6
+7 q^{12} t^6+4 q^{13} t^6+q^{10}
t^7+2 q^{11} t^7+5 q^{12} t^7+5 q^{13} t^7+q^{14} t^7
+q^{11} t^8+2 q^{12} t^8+4 q^{13} t^8+q^{14} t^8+q^{12} t^9
+2 q^{13} t^9+2 q^{14} t^9+q^{13}
t^{10}+q^{14} t^{10}+q^{14} t^{11}\bigr)
\)
\vfill
\noindent
\(+a^3 \bigl(q^6+q^7+q^8+q^9+q^7 t+2 q^8
t+3 q^9 t+3 q^{10} t+q^{11} t+q^8 t^2+2 q^9 t^2+4 q^{10} t^2
+4 q^{11} t^2+2 q^{12} t^2+q^9 t^3+2 q^{10} t^3+4 q^{11} t^3
+4 q^{12} t^3+q^{13} t^3+q^{10}
t^4+2 q^{11} t^4+4 q^{12} t^4+3 q^{13} t^4+q^{11} t^5
+2 q^{12} t^5+3 q^{13} t^5+q^{14} t^5+q^{12} t^6+2 q^{13} t^6
+q^{14} t^6+q^{13} t^7+q^{14} t^7+q^{14}
t^8\bigr)+
\)
\vfill
\noindent
\(
a^4 \bigl(q^{10}+q^{11} t+q^{12} t+q^{12} t^2+q^{13} t^2
+q^{13} t^3+q^{14} t^4\bigr).
\)
}
\renewcommand{\baselinestretch}{1.2}

\medskip
Its special case under $a=-t^2$ results in
(\ref{jd8-5}).
Switching to the ``standard notations" from
(\ref{qtarel}),
this polynomial coincides with the super-polynomial
$\h_{8,5}(b=\om_1)$ calculated 
(for the first time) by Gorsky based on 
the approach due to Gorsky, Oblomkov, 
Rasmussen, Schende \cite{ORS,Gor}. Their construction 
automatically provides the positivity 
(categorification) of the resulting $q,t,a$\~coefficients, 
but it is developed at the moment only for $b=\om_1$.

Our formula can be used for any torus knots and weights, 
but the positivity can be seen only {\em a posteriori}, at 
level of concrete formulas.
In spite of the fact that both constructions are based on DAHA,
the exact connection is not established so far. Importantly, 
our (conjectural)
formula opens a road to the theory 
of super-polynomials for the most general classical 
root systems $C^\vee C_n$.

Relatively direct (combinatorial) methods of finding this
super- polynomial seem insufficient and the connection 
to the Khovanov-Rozansky homology is non-trivial for $K_{8,5}$. 
It is one of the first uncolored cases where 
the number of $q,t$\~monomials in $H\!D^{[n]}_{r,s}$ can be greater
than that in  $J\!D^{n}_{r,s}$. Actually, it occurs 
only when $N=3$ for $K_{8,5}$ (the dimensions always coincide 
as claimed in Conjecture \ref{CONJKHR}). 

\medskip

Furthermore, let us also provide the $\{9,5\}$ formula: 
$$
H\!D_{9,5}(b=\om_1;q,t,a)\ =
$$

\renewcommand{\baselinestretch}{0.5} 
{\small
\noindent
\(
a^0\bigl(
1+q t+q^2 t+q^3 t+q^4 t+q^2 t^2+q^3 t^2+2 q^4 t^2+2 q^5 t^2+2 q^6 t^2
+q^7 t^2+q^3 t^3+q^4 t^3+2 q^5 t^3+3 q^6 t^3+3 q^7 t^3+2 q^8 t^3
+q^9 t^3+q^4 t^4+q^5 t^4+2 q^6 t^4+3 q^7 t^4+4 q^8 t^4+3 q^9 t^4
+2 q^{10} t^4+q^5 t^5+q^6 t^5+2 q^7 t^5+3 q^8 t^5+4 q^9 t^5
+4 q^{10} t^5
+2 q^{11} t^5+q^6 t^6+q^7 t^6+2 q^8 t^6+3 q^9 t^6+4 q^{10} t^6
+4 q^{11} t^6+2 q^{12} t^6+q^7 t^7+q^8 t^7+2 q^9 t^7+3 q^{10} t^7
+4 q^{11} t^7+3 q^{12} t^7+q^{13} t^7+q^8 t^8+q^9 t^8+2 q^{10} t^8
+3 q^{11} t^8+4 q^{12} t^8+2 q^{13} t^8+q^9 t^9+q^{10} t^9
+2 q^{11} t^9
+3 q^{12} t^9+3 q^{13} t^9+q^{14} t^9+q^{10} t^{10}+q^{11} t^{10}
+2 q^{12} t^{10}+3 q^{13} t^{10}+2 q^{14} t^{10}+q^{11} t^{11}
+q^{12} t^{11}+2 q^{13} t^{11}+2 q^{14} t^{11}+q^{12} t^{12}
+q^{13} t^{12}+2 q^{14} t^{12}+q^{15} t^{12}+q^{13} t^{13}
+q^{14} t^{13}+q^{15} t^{13}+q^{14} t^{14}+q^{15} t^{14}
+q^{15} t^{15}
+q^{16} t^{16}\bigr)
\)
\vfill
\noindent
\(
+a^1 \bigl(q+q^2+q^3+q^4+q^2 t
+2 q^3 t+3 q^4 t+4 q^5 t+3 q^6 t+2 q^7 t+q^3 t^2+2 q^4 t^2+4 q^5 t^2
+6 q^6 t^2+7 q^7 t^2+6 q^8 t^2+3 q^9 t^2+q^{10} t^2+q^4 t^3+2 q^5 t^3
+4 q^6 t^3+7 q^7 t^3+9 q^8 t^3+9 q^9 t^3+6 q^{10} t^3+2 q^{11} t^3
+q^5 t^4+2 q^6 t^4+4 q^7 t^4+7 q^8 t^4+10 q^9 t^4+11 q^{10} t^4
+8 q^{11} t^4+3 q^{12} t^4+q^6 t^5+2 q^7 t^5+4 q^8 t^5+7 q^9 t^5
+10 q^{10} t^5+12 q^{11} t^5+8 q^{12} t^5+2 q^{13} t^5+q^7 t^6
+2 q^8 t^6
+4 q^9 t^6+7 q^{10} t^6+10 q^{11} t^6+11 q^{12} t^6+6 q^{13} t^6
+q^{14} t^6+q^8 t^7+2 q^9 t^7+4 q^{10} t^7+7 q^{11} t^7+10 q^{12} t^7
+9 q^{13} t^7+3 q^{14} t^7+q^9 t^8+2 q^{10} t^8+4 q^{11} t^8
+7 q^{12} t^8
+9 q^{13} t^8+6 q^{14} t^8+q^{10} t^9+2 q^{11} t^9+4 q^{12} t^9
+7 q^{13} t^9+7 q^{14} t^9+2 q^{15} t^9+q^{11} t^{10}
+2 q^{12} t^{10}
+4 q^{13} t^{10}+6 q^{14} t^{10}+3 q^{15} t^{10}+q^{12} t^{11}
+2 q^{13} t^{11}+4 q^{14} t^{11}+4 q^{15} t^{11}+q^{13} t^{12}
+2 q^{14} t^{12}+3 q^{15} t^{12}+q^{16} t^{12}+q^{14} t^{13}
+2 q^{15} t^{13}+q^{16} t^{13}+q^{15} t^{14}+q^{16} t^{14}
+q^{16} t^{15}\bigr)
\)
\vfill
\noindent
\(
+a^2 \bigl(q^3+q^4+2 q^5+q^6+q^7+q^4 t+2 q^5 t+4 q^6 t+5 q^7 t
+5 q^8 t
+3 q^9 t+q^{10} t+q^5 t^2+2 q^6 t^2+5 q^7 t^2+7 q^8 t^2+9 q^9 t^2
+7 q^{10} t^2+4 q^{11} t^2+q^{12} t^2+q^6 t^3+2 q^7 t^3+5 q^8 t^3
+8 q^9 t^3+11 q^{10} t^3+10 q^{11} t^3+6 q^{12} t^3+q^{13} t^3
+q^7 t^4
+2 q^8 t^4+5 q^9 t^4+8 q^{10} t^4+12 q^{11} t^4+11 q^{12} t^4
+6 q^{13} t^4+q^{14} t^4+q^8 t^5+2 q^9 t^5+5 q^{10} t^5+8 q^{11} t^5
+12 q^{12} t^5+10 q^{13} t^5+4 q^{14} t^5+q^9 t^6+2 q^{10} t^6
+5 q^{11} t^6+8 q^{12} t^6+11 q^{13} t^6+7 q^{14} t^6+q^{15} t^6
+q^{10} t^7+2 q^{11} t^7+5 q^{12} t^7+8 q^{13} t^7+9 q^{14} t^7
+3 q^{15} t^7+q^{11} t^8+2 q^{12} t^8+5 q^{13} t^8+7 q^{14} t^8
+5 q^{15} t^8+q^{12} t^9+2 q^{13} t^9+5 q^{14} t^9+5 q^{15} t^9
+q^{16} t^9+q^{13} t^{10}+2 q^{14} t^{10}+4 q^{15} t^{10}
+q^{16} t^{10}
+q^{14} t^{11}+2 q^{15} t^{11}+2 q^{16} t^{11}+q^{15} t^{12}
+q^{16} t^{12}+q^{16} t^{13}\bigr)
\)
\vfill
\noindent
\(
+a^3 \bigl(q^6+q^7+q^8+q^9+q^7 t+2 q^8 t+3 q^9 t+4 q^{10} t
+2 q^{11} t
+q^{12} t+q^8 t^2+2 q^9 t^2+4 q^{10} t^2+5 q^{11} t^2+4 q^{12} t^2
+2 q^{13} t^2+q^9 t^3+2 q^{10} t^3+4 q^{11} t^3+6 q^{12} t^3
+5 q^{13} t^3
+2 q^{14} t^3+q^{10} t^4+2 q^{11} t^4+4 q^{12} t^4+6 q^{13} t^4
+4 q^{14} t^4+q^{15} t^4+q^{11} t^5+2 q^{12} t^5+4 q^{13} t^5
+5 q^{14} t^5+2 q^{15} t^5+q^{12} t^6+2 q^{13} t^6+4 q^{14} t^6
+4 q^{15} t^6+q^{13} t^7+2 q^{14} t^7+3 q^{15} t^7+q^{16} t^7
+q^{14} t^8+2 q^{15} t^8+q^{16} t^8+q^{15} t^9+q^{16} t^9
+q^{16} t^{10}\bigr)
\)
\vfill
\noindent
\(
+a^4 \bigl(q^{10}+q^{11} t+q^{12} t+q^{13} t
+q^{12} t^2+q^{13} t^2+q^{14} t^2+q^{13} t^3+q^{14} t^3+q^{15} t^3
+q^{14} t^4+q^{15} t^4+q^{15} t^5+q^{16} t^6\bigr).
\)
}
\renewcommand{\baselinestretch}{1.2}
\medskip

It coincides with that obtained by Gorsky
subject to the substitutions from (\ref{qtarel}),
as well as the following example of the super-polynomial:

$$
H\!D_{11,6}(b=\om_1;q,t,a)\ =
$$

\renewcommand{\baselinestretch}{0.5} 
{\small
\noindent
\(
a^0 (1+q t+q^2 t+q^3 t+q^4 t+q^5 t+q^2 t^2+q^3 t^2+2 q^4 t^2
+2 q^5 t^2+3 q^6 t^2+2 q^7 t^2+2 q^8 t^2+q^9 t^2+q^3 t^3+q^4 t^3
+2 q^5 t^3+3 q^6 t^3+4 q^7 t^3+4 q^8 t^3+5 q^9 t^3+3 q^{10} t^3
+2 q^{11} t^3+q^{12} t^3+q^4 t^4+q^5 t^4+2 q^6 t^4+3 q^7 t^4
+5 q^8 t^4+5 q^9 t^4+7 q^{10} t^4+6 q^{11} t^4+5 q^{12} t^4
+3 q^{13} t^4+q^{14} t^4+q^5 t^5+q^6 t^5+2 q^7 t^5+3 q^8 t^5
+5 q^9 t^5+6 q^{10} t^5+8 q^{11} t^5+8 q^{12} t^5+8 q^{13} t^5
+5 q^{14} t^5+3 q^{15} t^5+q^6 t^6+q^7 t^6+2 q^8 t^6+3 q^9 t^6
+5 q^{10} t^6+6 q^{11} t^6+9 q^{12} t^6+9 q^{13} t^6+10 q^{14} t^6
+8 q^{15} t^6+4 q^{16} t^6+q^{17} t^6+q^7 t^7+q^8 t^7+2 q^9 t^7
+3 q^{10} t^7+5 q^{11} t^7+6 q^{12} t^7+9 q^{13} t^7+10 q^{14} t^7
+11 q^{15} t^7+9 q^{16} t^7+5 q^{17} t^7+q^{18} t^7+q^8 t^8
+q^9 t^8+2 q^{10} t^8+3 q^{11} t^8+5 q^{12} t^8+6 q^{13} t^8
+9 q^{14} t^8+10 q^{15} t^8+12 q^{16} t^8+9 q^{17} t^8
+5 q^{18} t^8+q^{19} t^8+q^9 t^9+q^{10} t^9+2 q^{11} t^9
+3 q^{12} t^9+5 q^{13} t^9+6 q^{14} t^9+9 q^{15} t^9+10 q^{16} t^9
+12 q^{17} t^9+9 q^{18} t^9+4 q^{19} t^9+q^{10} t^{10}
+q^{11} t^{10}+2 q^{12} t^{10}+3 q^{13} t^{10}+5 q^{14} t^{10}
+6 q^{15} t^{10}+9 q^{16} t^{10}+10 q^{17} t^{10}+11 q^{18} t^{10}
+8 q^{19} t^{10}+3 q^{20} t^{10}+q^{11} t^{11}+q^{12} t^{11}
+2 q^{13} t^{11}+3 q^{14} t^{11}+5 q^{15} t^{11}+6 q^{16} t^{11}
+9 q^{17} t^{11}+10 q^{18} t^{11}+10 q^{19} t^{11}+5 q^{20} t^{11}
+q^{21} t^{11}+q^{12} t^{12}+q^{13} t^{12}+2 q^{14} t^{12}
+3 q^{15} t^{12}+5 q^{16} t^{12}+6 q^{17} t^{12}+9 q^{18} t^{12}
+9 q^{19} t^{12}+8 q^{20} t^{12}+3 q^{21} t^{12}+q^{13} t^{13}
+q^{14} t^{13}+2 q^{15} t^{13}+3 q^{16} t^{13}+5 q^{17} t^{13}
+6 q^{18} t^{13}+9 q^{19} t^{13}+8 q^{20} t^{13}+5 q^{21} t^{13}
+q^{22} t^{13}+q^{14} t^{14}+q^{15} t^{14}+2 q^{16} t^{14}
+3 q^{17} t^{14}+5 q^{18} t^{14}+6 q^{19} t^{14}+8 q^{20} t^{14}
+6 q^{21} t^{14}+2 q^{22} t^{14}+q^{15} t^{15}+q^{16} t^{15}
+2 q^{17} t^{15}+3 q^{18} t^{15}+5 q^{19} t^{15}+6 q^{20} t^{15}
+7 q^{21} t^{15}+3 q^{22} t^{15}+q^{16} t^{16}+q^{17} t^{16}
+2 q^{18} t^{16}+3 q^{19} t^{16}+5 q^{20} t^{16}+5 q^{21} t^{16}
+5 q^{22} t^{16}+q^{23} t^{16}+q^{17} t^{17}+q^{18} t^{17}
+2 q^{19} t^{17}+3 q^{20} t^{17}+5 q^{21} t^{17}+4 q^{22} t^{17}
+2 q^{23} t^{17}+q^{18} t^{18}+q^{19} t^{18}+2 q^{20} t^{18}
+3 q^{21} t^{18}+4 q^{22} t^{18}+2 q^{23} t^{18}+q^{19} t^{19}
+q^{20} t^{19}+2 q^{21} t^{19}+3 q^{22} t^{19}+3 q^{23} t^{19}
+q^{20} t^{20}+q^{21} t^{20}+2 q^{22} t^{20}+2 q^{23} t^{20}
+q^{24} t^{20}+q^{21} t^{21}+q^{22} t^{21}+2 q^{23} t^{21}
+q^{24} t^{21}+q^{22} t^{22}+q^{23} t^{22}+q^{24} t^{22}
+q^{23} t^{23}+q^{24} t^{23}+q^{24} t^{24}+q^{25} t^{25})
\)
\vfill
\noindent
\(
+a^1 (q+q^2+q^3+q^4+q^5+q^2 t+2 q^3 t+3 q^4 t+4 q^5 t+5 q^6 t
+4 q^7 t+3 q^8 t+2 q^9 t+q^3 t^2+2 q^4 t^2+4 q^5 t^2+6 q^6 t^2
+9 q^7 t^2+10 q^8 t^2+11 q^9 t^2+9 q^{10} t^2+6 q^{11} t^2
+3 q^{12} t^2+q^{13} t^2+q^4 t^3+2 q^5 t^3+4 q^6 t^3+7 q^7 t^3
+11 q^8 t^3+14 q^9 t^3+18 q^{10} t^3+18 q^{11} t^3+15 q^{12} t^3
+11 q^{13} t^3+5 q^{14} t^3+2 q^{15} t^3+q^5 t^4+2 q^6 t^4
+4 q^7 t^4+7 q^8 t^4+12 q^9 t^4+16 q^{10} t^4+22 q^{11} t^4
+25 q^{12} t^4+25 q^{13} t^4+20 q^{14} t^4+13 q^{15} t^4
+5 q^{16} t^4+q^{17} t^4+q^6 t^5+2 q^7 t^5+4 q^8 t^5+7 q^9 t^5
+12 q^{10} t^5+17 q^{11} t^5+24 q^{12} t^5+29 q^{13} t^5
+32 q^{14} t^5+29 q^{15} t^5+20 q^{16} t^5+9 q^{17} t^5
+2 q^{18} t^5+q^7 t^6+2 q^8 t^6+4 q^9 t^6+7 q^{10} t^6
+12 q^{11} t^6+17 q^{12} t^6+25 q^{13} t^6+31 q^{14} t^6
+36 q^{15} t^6+35 q^{16} t^6+25 q^{17} t^6+12 q^{18} t^6
+3 q^{19} t^6+q^8 t^7+2 q^9 t^7+4 q^{10} t^7+7 q^{11} t^7
+12 q^{12} t^7+17 q^{13} t^7+25 q^{14} t^7+32 q^{15} t^7
+38 q^{16} t^7+37 q^{17} t^7+27 q^{18} t^7+12 q^{19} t^7
+2 q^{20} t^7+q^9 t^8+2 q^{10} t^8+4 q^{11} t^8+7 q^{12} t^8
+12 q^{13} t^8+17 q^{14} t^8+25 q^{15} t^8+32 q^{16} t^8
+39 q^{17} t^8+37 q^{18} t^8+25 q^{19} t^8+9 q^{20} t^8
+q^{21} t^8+q^{10} t^9+2 q^{11} t^9+4 q^{12} t^9+7 q^{13} t^9
+12 q^{14} t^9+17 q^{15} t^9+25 q^{16} t^9+32 q^{17} t^9
+38 q^{18} t^9+35 q^{19} t^9+20 q^{20} t^9+5 q^{21} t^9
+q^{11} t^{10}+2 q^{12} t^{10}+4 q^{13} t^{10}+7 q^{14} t^{10}
+12 q^{15} t^{10}+17 q^{16} t^{10}+25 q^{17} t^{10}
+32 q^{18} t^{10}+36 q^{19} t^{10}+29 q^{20} t^{10}
+13 q^{21} t^{10}+2 q^{22} t^{10}+q^{12} t^{11}
+2 q^{13} t^{11}+4 q^{14} t^{11}+7 q^{15} t^{11}
+12 q^{16} t^{11}+17 q^{17} t^{11}+25 q^{18} t^{11}
+31 q^{19} t^{11}+32 q^{20} t^{11}+20 q^{21} t^{11}
+5 q^{22} t^{11}+q^{13} t^{12}+2 q^{14} t^{12}
+4 q^{15} t^{12}+7 q^{16} t^{12}+12 q^{17} t^{12}
+17 q^{18} t^{12}+25 q^{19} t^{12}+29 q^{20} t^{12}
+25 q^{21} t^{12}+11 q^{22} t^{12}+q^{23} t^{12}
+q^{14} t^{13}+2 q^{15} t^{13}+4 q^{16} t^{13}
+7 q^{17} t^{13}+12 q^{18} t^{13}+17 q^{19} t^{13}
+24 q^{20} t^{13}+25 q^{21} t^{13}
+15 q^{22} t^{13}+3 q^{23} t^{13}+q^{15} t^{14}+2 q^{16} t^{14}
+4 q^{17} t^{14}+7 q^{18} t^{14}+12 q^{19} t^{14}+17 q^{20} t^{14}
+22 q^{21} t^{14}+18 q^{22} t^{14}+6 q^{23} t^{14}+q^{16} t^{15}
+2 q^{17} t^{15}+4 q^{18} t^{15}+7 q^{19} t^{15}+12 q^{20} t^{15}
+16 q^{21} t^{15}+18 q^{22} t^{15}+9 q^{23} t^{15}+q^{17} t^{16}
+2 q^{18} t^{16}+4 q^{19} t^{16}+7 q^{20} t^{16}+12 q^{21} t^{16}
+14 q^{22} t^{16}+11 q^{23} t^{16}+2 q^{24} t^{16}+q^{18} t^{17}
+2 q^{19} t^{17}+4 q^{20} t^{17}+7 q^{21} t^{17}+11 q^{22} t^{17}
+10 q^{23} t^{17}+3 q^{24} t^{17}+q^{19} t^{18}+2 q^{20} t^{18}
+4 q^{21} t^{18}+7 q^{22} t^{18}+9 q^{23} t^{18}+4 q^{24} t^{18}
+q^{20} t^{19}+2 q^{21} t^{19}+4 q^{22} t^{19}+6 q^{23} t^{19}
+5 q^{24} t^{19}+q^{21} t^{20}+2 q^{22} t^{20}+4 q^{23} t^{20}
+4 q^{24} t^{20}+q^{25} t^{20}+q^{22} t^{21}+2 q^{23} t^{21}
+3 q^{24} t^{21}+q^{25} t^{21}+q^{23} t^{22}+2 q^{24} t^{22}
+q^{25} t^{22}+q^{24} t^{23}+q^{25} t^{23}+q^{25} t^{24})
\)
\smallskip

\noindent
\(
+a^2 
(q^3+q^4+2 q^5+2 q^6+2 q^7+q^8+q^9+q^4 t+2 q^5 t+4 q^6 t+6 q^7 t
+8 q^8 t+8 q^9 t+8 q^{10} t+5 q^{11} t+3 q^{12} t+q^{13} t+q^5 t^2
+2 q^6 t^2+5 q^7 t^2+8 q^8 t^2+13 q^9 t^2+16 q^{10} t^2
+19 q^{11} t^2+17 q^{12} t^2+14 q^{13} t^2+8 q^{14} t^2
+4 q^{15} t^2+q^{16} t^2+q^6 t^3+2 q^7 t^3+5 q^8 t^3+9 q^9 t^3
+15 q^{10} t^3+21 q^{11} t^3+28 q^{12} t^3+29 q^{13} t^3
+28 q^{14} t^3+20 q^{15} t^3+12 q^{16} t^3+4 q^{17} t^3
+q^{18} t^3+q^7 t^4+2 q^8 t^4+5 q^9 t^4+9 q^{10} t^4+16 q^{11} t^4
+23 q^{12} t^4+33 q^{13} t^4+38 q^{14} t^4+40 q^{15} t^4
+33 q^{16} t^4+21 q^{17} t^4+9 q^{18} t^4+2 q^{19} t^4
+q^8 t^5+2 q^9 t^5+5 q^{10} t^5+9 q^{11} t^5+16 q^{12} t^5
+24 q^{13} t^5+35 q^{14} t^5+43 q^{15} t^5+49 q^{16} t^5
+42 q^{17} t^5+29 q^{18} t^5+12 q^{19} t^5+3 q^{20} t^5
+q^9 t^6+2 q^{10} t^6+5 q^{11} t^6+9 q^{12} t^6+16 q^{13} t^6
+24 q^{14} t^6+36 q^{15} t^6+45 q^{16} t^6+53 q^{17} t^6
+47 q^{18} t^6+31 q^{19} t^6+12 q^{20} t^6+2 q^{21} t^6
+q^{10} t^7+2 q^{11} t^7+5 q^{12} t^7+9 q^{13} t^7+16 q^{14} t^7
+24 q^{15} t^7+36 q^{16} t^7+46 q^{17} t^7+54 q^{18} t^7
+47 q^{19} t^7+29 q^{20} t^7+9 q^{21} t^7+q^{22} t^7+q^{11} t^8
+2 q^{12} t^8+5 q^{13} t^8+9 q^{14} t^8+16 q^{15} t^8
+24 q^{16} t^8+36 q^{17} t^8+46 q^{18} t^8+53 q^{19} t^8
+42 q^{20} t^8+21 q^{21} t^8+4 q^{22} t^8+q^{12} t^9+2 q^{13} t^9
+5 q^{14} t^9+9 q^{15} t^9+16 q^{16} t^9+24 q^{17} t^9
+36 q^{18} t^9+45 q^{19} t^9+49 q^{20} t^9+33 q^{21} t^9
+12 q^{22} t^9+q^{23} t^9+q^{13} t^{10}+2 q^{14} t^{10}
+5 q^{15} t^{10}+9 q^{16} t^{10}+16 q^{17} t^{10}
+24 q^{18} t^{10}+36 q^{19} t^{10}+43 q^{20} t^{10}
+40 q^{21} t^{10}+20 q^{22} t^{10}+4 q^{23} t^{10}
+q^{14} t^{11}+2 q^{15} t^{11}+5 q^{16} t^{11}+9 q^{17} t^{11}
+16 q^{18} t^{11}+24 q^{19} t^{11}+35 q^{20} t^{11}
+38 q^{21} t^{11}+28 q^{22} t^{11}+8 q^{23} t^{11}
+q^{15} t^{12}+2 q^{16} t^{12}+5 q^{17} t^{12}+9 q^{18} t^{12}
+16 q^{19} t^{12}+24 q^{20} t^{12}+33 q^{21} t^{12}
+29 q^{22} t^{12}+14 q^{23} t^{12}+q^{24} t^{12}+q^{16} t^{13}
+2 q^{17} t^{13}+5 q^{18} t^{13}+9 q^{19} t^{13}
+16 q^{20} t^{13}+23 q^{21} t^{13}+28 q^{22} t^{13}
+17 q^{23} t^{13}+3 q^{24} t^{13}+q^{17} t^{14}
+2 q^{18} t^{14}+5 q^{19} t^{14}+9 q^{20} t^{14}
+16 q^{21} t^{14}+21 q^{22} t^{14}+19 q^{23} t^{14}
+5 q^{24} t^{14}+q^{18} t^{15}+2 q^{19} t^{15}
+5 q^{20} t^{15}+9 q^{21} t^{15}+15 q^{22} t^{15}
+16 q^{23} t^{15}+8 q^{24} t^{15}+q^{19} t^{16}
+2 q^{20} t^{16}+5 q^{21} t^{16}+9 q^{22} t^{16}
+13 q^{23} t^{16}+8 q^{24} t^{16}+q^{25} t^{16}
+q^{20} t^{17}+2 q^{21} t^{17}+5 q^{22} t^{17}
+8 q^{23} t^{17}+8 q^{24} t^{17}+q^{25} t^{17}
+q^{21} t^{18}+2 q^{22} t^{18}+5 q^{23} t^{18}
+6 q^{24} t^{18}+2 q^{25} t^{18}+q^{22} t^{19}
+2 q^{23} t^{19}+4 q^{24} t^{19}+2 q^{25} t^{19}
+q^{23} t^{20}+2 q^{24} t^{20}+2 q^{25} t^{20}
+q^{24} t^{21}+q^{25} t^{21}+q^{25} t^{22})
\)
\vfill

\noindent
\(
+a^3 (q^6+q^7+2 q^8+2 q^9+2 q^{10}+q^{11}+q^{12}+q^7 t+2 q^8 t
+4 q^9 t+6 q^{10} t+8 q^{11} t+8 q^{12} t+7 q^{13} t+5 q^{14} t
+2 q^{15} t+q^{16} t+q^8 t^2+2 q^9 t^2+5 q^{10} t^2+8 q^{11} t^2
+13 q^{12} t^2+15 q^{13} t^2+17 q^{14} t^2+14 q^{15} t^2
+10 q^{16} t^2+5 q^{17} t^2+2 q^{18} t^2+q^9 t^3+2 q^{10} t^3
+5 q^{11} t^3+9 q^{12} t^3+15 q^{13} t^3+20 q^{14} t^3+25 q^{15} t^3
+24 q^{16} t^3+19 q^{17} t^3+12 q^{18} t^3+4 q^{19} t^3+q^{20} t^3
+q^{10} t^4+2 q^{11} t^4+5 q^{12} t^4+9 q^{13} t^4+16 q^{14} t^4
+22 q^{15} t^4+30 q^{16} t^4+31 q^{17} t^4+27 q^{18} t^4
+17 q^{19} t^4+7 q^{20} t^4+q^{21} t^4+q^{11} t^5+2 q^{12} t^5
+5 q^{13} t^5+9 q^{14} t^5+16 q^{15} t^5+23 q^{16} t^5
+32 q^{17} t^5+35 q^{18} t^5+31 q^{19} t^5+19 q^{20} t^5
+7 q^{21} t^5+q^{22} t^5+q^{12} t^6+2 q^{13} t^6+5 q^{14} t^6
+9 q^{15} t^6+16 q^{16} t^6+23 q^{17} t^6+33 q^{18} t^6+36 q^{19} t^6
+31 q^{20} t^6+17 q^{21} t^6+4 q^{22} t^6+q^{13} t^7+2 q^{14} t^7
+5 q^{15} t^7+9 q^{16} t^7+16 q^{17} t^7+23 q^{18} t^7
+33 q^{19} t^7+35 q^{20} t^7+27 q^{21} t^7+12 q^{22} t^7
+2 q^{23} t^7+q^{14} t^8+2 q^{15} t^8+5 q^{16} t^8+9 q^{17} t^8
+16 q^{18} t^8+23 q^{19} t^8+32 q^{20} t^8+31 q^{21} t^8
+19 q^{22} t^8+5 q^{23} t^8+q^{15} t^9+2 q^{16} t^9+5 q^{17} t^9
+9 q^{18} t^9+16 q^{19} t^9+23 q^{20} t^9+30 q^{21} t^9+24 q^{22} t^9
+10 q^{23} t^9+q^{24} t^9+q^{16} t^{10}+2 q^{17} t^{10}
+5 q^{18} t^{10}+9 q^{19} t^{10}+16 q^{20} t^{10}+22 q^{21} t^{10}
+25 q^{22} t^{10}+14 q^{23} t^{10}+2 q^{24} t^{10}+q^{17} t^{11}
+2 q^{18} t^{11}+5 q^{19} t^{11}+9 q^{20} t^{11}+16 q^{21} t^{11}
+20 q^{22} t^{11}+17 q^{23} t^{11}+5 q^{24} t^{11}+q^{18} t^{12}
+2 q^{19} t^{12}+5 q^{20} t^{12}+9 q^{21} t^{12}+15 q^{22} t^{12}
+15 q^{23} t^{12}+7 q^{24} t^{12}+q^{19} t^{13}+2 q^{20} t^{13}
+5 q^{21} t^{13}+9 q^{22} t^{13}+13 q^{23} t^{13}+8 q^{24} t^{13}
+q^{25} t^{13}+q^{20} t^{14}+2 q^{21} t^{14}+5 q^{22} t^{14}
+8 q^{23} t^{14}+8 q^{24} t^{14}+q^{25} t^{14}+q^{21} t^{15}
+2 q^{22} t^{15}+5 q^{23} t^{15}+6 q^{24} t^{15}+2 q^{25} t^{15}
+q^{22} t^{16}+2 q^{23} t^{16}+4 q^{24} t^{16}+2 q^{25} t^{16}
+q^{23} t^{17}+2 q^{24} t^{17}+2 q^{25} t^{17}+q^{24} t^{18}
+q^{25} t^{18}+q^{25} t^{19})
\)
\smallskip

\noindent
\(
+a^4 (q^{10}+q^{11}+q^{12}+q^{13}+q^{14}+q^{11} t+2 q^{12} t
+3 q^{13} t+4 q^{14} t+5 q^{15} t+3 q^{16} t+2 q^{17} t+q^{18} t
+q^{12} t^2+2 q^{13} t^2+4 q^{14} t^2+6 q^{15} t^2+8 q^{16} t^2
+7 q^{17} t^2+6 q^{18} t^2+3 q^{19} t^2+q^{20} t^2+q^{13} t^3
+2 q^{14} t^3+4 q^{15} t^3+7 q^{16} t^3+10 q^{17} t^3+10 q^{18} t^3
+10 q^{19} t^3+5 q^{20} t^3+2 q^{21} t^3+q^{14} t^4+2 q^{15} t^4
+4 q^{16} t^4+7 q^{17} t^4+11 q^{18} t^4+12 q^{19} t^4+11 q^{20} t^4
+7 q^{21} t^4+2 q^{22} t^4+q^{15} t^5+2 q^{16} t^5+4 q^{17} t^5
+7 q^{18} t^5+11 q^{19} t^5+12 q^{20} t^5+11 q^{21} t^5+5 q^{22} t^5
+q^{23} t^5+q^{16} t^6+2 q^{17} t^6+4 q^{18} t^6+7 q^{19} t^6
+11 q^{20} t^6+12 q^{21} t^6+10 q^{22} t^6+3 q^{23} t^6+q^{17} t^7
+2 q^{18} t^7+4 q^{19} t^7+7 q^{20} t^7+11 q^{21} t^7+10 q^{22} t^7
+6 q^{23} t^7+q^{24} t^7+q^{18} t^8+2 q^{19} t^8+4 q^{20} t^8
+7 q^{21} t^8+10 q^{22} t^8+7 q^{23} t^8+2 q^{24} t^8+q^{19} t^9
+2 q^{20} t^9+4 q^{21} t^9+7 q^{22} t^9+8 q^{23} t^9+3 q^{24} t^9
+q^{20} t^{10}+2 q^{21} t^{10}+4 q^{22} t^{10}+6 q^{23} t^{10}
+5 q^{24} t^{10}+q^{21} t^{11}+2 q^{22} t^{11}+4 q^{23} t^{11}
+4 q^{24} t^{11}+q^{25} t^{11}+q^{22} t^{12}+2 q^{23} t^{12}
+3 q^{24} t^{12}+q^{25} t^{12}+q^{23} t^{13}+2 q^{24} t^{13}
+q^{25} t^{13}+q^{24} t^{14}+q^{25} t^{14}+q^{25} t^{15})
\)

\smallskip
\noindent
\(
+a^5 (q^{15}+q^{16} t+q^{17} t+q^{18} t+q^{19} t+q^{17} t^2
+q^{18} t^2+2 q^{19} t^2+q^{20} t^2+q^{21} t^2+q^{18} t^3
+q^{19} t^3+2 q^{20} t^3+2 q^{21} t^3+q^{22} t^3+q^{19} t^4
+q^{20} t^4+2 q^{21} t^4+2 q^{22} t^4+q^{23} t^4+q^{20} t^5
+q^{21} t^5+2 q^{22} t^5+q^{23} t^5+q^{21} t^6+q^{22} t^6
+2 q^{23} t^6+q^{24} t^6+q^{22} t^7+q^{23} t^7+q^{24} t^7
+q^{23} t^8+q^{24} t^8+q^{24} t^9+q^{25} t^{10}).
\)
}
\renewcommand{\baselinestretch}{1.2} 

\smallskip
We think that these examples confirm beyond
a reasonable doubt the coincidence of the non-colored
DAHA super-polynomials for torus knots originated 
in topology, physics and geometry.
The rich expected symmetries of super-polynomials
can be an explanation of their universality, but we
hope that more fundamental roots will be found.
The relation to the classical theory of singularities
associated with simple Lie algebras and their irreducible
representations via the BPS theory seems of significant 
importance here. 

Also, the special role of the torus knots in this theory
requires careful analysis, though here the refined
Chern-Simons theory and the BPS theory already provide some
explanations.
\smallskip

We would like to 
mention here another approach to 
super-poly\-no\-mials from \cite{DMMSS} and 
other works of these authors, but cannot comment much 
on the methods they use. 
The formulas they obtain (some are involved) match 
or are expected to match our ones and those obtained by Gorsky
(there are some direct links);
Gorsky's computer program is by far the main source of explicit 
non-colored formulas.

The case 
$b=\om_1$ is the simplest and standard in almost all
mathematical and physical papers on super-polynomials.
All non-colored polynomials we and others constructed
satisfy Conjecture \ref{DUALIT} ($\om_1^{tr}=\om_1$).
The confirmations of this conjecture in the colored
case are significantly more limited at the moment; 
our construction is universal but the formulas are getting 
involved when the weights grow. 
  
Our method seems the only universal one
now that can be applied to construct $q,t$\~deformations
of Jones polynomials of arbitrary torus knots for all root systems
and any weights. For $A_n$, it also provides
super-polynomials, though our ones can have negative
coefficients for generic weights, in contrast to
the physics expectations.
\medskip

\subsection{Columns and rows}
We will provide the formulas for the diagrams with 2-6 boxes 
for the knots $\{3,2\}$, $\{5,2\}$ and $\{4,3\}$.
Due to the duality
(which holds in all calculated examples, including these ones), 
we will mainly skip here $b=j\om_1$, focusing on $b=\om_j$.
Let us begin with the 2-box formulas:
\smallskip

\begin{align}\label{super3-2+2-0}
&H\!D_{3,2}(b=\om_2;\, q,t,a)\ =\\
&1+q t+q t^2+q^2 t^4+a \bigl(q+\frac{q}{t}+q^2 t+q^2 t^2\bigr)
+\frac{a^2 q^2}{t},\notag\\
\label{super3-2+1-1}
&H\!D_{3,2}(b=2\om_1;\, q,t,a)\ =\\
&1+q^2 t+q^3 t+q^4 t^2 +a \bigl(q^2+q^3+q^4 t+q^5 t\bigr)+a^2 q^5.
\notag
\end{align}
\smallskip

The last formula coincides with that calculated
by Gukov and Sto$\breve{\hbox{s}}$i\'{c} and with 
the prediction from Section 7.7 of \cite{AS} upon the substitution 
from (\ref{qtarel}) and up to a (general) proportionality
coefficient; see also formulas (169) and (176) 
from \cite{DMMSS}.

This latter polynomial becomes the former one
upon the duality substitution  $t\mapsto q^{-1}$,
$q\mapsto t^{-1}$ followed by multiplication by
$q^2 t^4$.

\medskip

{\em The case of $b=\om_3$,\, knots$\ =\{3,2\}, \{4,3\}$.}
These formulas coincide via the duality
with the formulas for symmetric powers 
posted in \cite{GS} and other works:

$$
H\!D_{3,2}(b=\om_3;\, q,t,a)\ =
$$

\renewcommand{\baselinestretch}{0.5} 
{\small
\noindent
\(
1+q t+q t^2+q t^3+q^2 t^4+q^2 t^5+q^2 t^6
+q^3 t^9+a \bigl(q+q^2+\frac{q}{t^2}+\frac{q}{t}+2 q^2 t+2 q^2 t^2
+q^2 t^3+q^3 t^4+q^3 t^5+q^3 t^6\bigr)
\)

\noindent
\(
+a^2 \bigl(q^3+\frac{q^2}{t^3}+\frac{q^2}{t^2}+
\frac{q^2}{t}+q^3 t+q^3 t^2\bigr)+\frac{a^3 q^3}{t^3}
\)\,,

}
\renewcommand{\baselinestretch}{1.2} 

$$
H\!D_{4,3}(b=\om_3;\, q,t,a)\ =
$$
\renewcommand{\baselinestretch}{0.5} 
{\small

\noindent
\(
1+\frac{a^6 q^9}{t^6}+q t+q^2 t+q t^2+2 q^2 t^2+q t^3+2 q^2 t^3
+2 q^3 t^3+2 q^2 t^4+3 q^3 t^4+q^4 t^4+q^2 t^5+4 q^3 t^5+2 q^4 t^5
+q^2 t^6+3 q^3 t^6+4 q^4 t^6+2 q^3 t^7+5 q^4 t^7+2 q^5 t^7+q^3 t^8
+5 q^4 t^8+3 q^5 t^8+q^3 t^9+3 q^4 t^9+5 q^5 t^9+q^6 t^9+2 q^4 t^{10}
+5 q^5 t^{10}+q^6 t^{10}+q^4 t^{11}+5 q^5 t^{11}+2 q^6 t^{11}
+q^4 t^{12}+3 q^5 t^{12}+4 q^6 t^{12}+2 q^5 t^{13}+4 q^6 t^{13}
+q^7 t^{13}+q^5 t^{14}+4 q^6 t^{14}+q^7 t^{14}+q^5 t^{15}
+3 q^6 t^{15}+2 q^7 t^{15}+2 q^6 t^{16}+2 q^7 t^{16}+q^6 t^{17}
+3 q^7 t^{17}+q^6 t^{18}+2 q^7 t^{18}+2 q^7 t^{19}+q^8 t^{19}
+q^7 t^{20}+q^8 t^{20}+q^7 t^{21}+q^8 t^{21}+q^8 t^{22}
+q^8 t^{23}+q^8 t^{24}+q^9 t^{27}
\)

\noindent
\(
+a^1 \bigl(q+3 q^2+4 q^3+q^4+\frac{q}{t^2}+\frac{q^2}{t^2}
+\frac{q}{t}+\frac{2 q^2}{t}+\frac{q^3}{t}+3 q^2 t+7 q^3 t
+3 q^4 t+2 q^2 t^2+8 q^3 t^2+7 q^4 t^2+q^5 t^2+q^2 t^3+6 q^3 t^3
+11 q^4 t^3+4 q^5 t^3+4 q^3 t^4+13 q^4 t^4+9 q^5 t^4+q^6 t^4
+2 q^3 t^5+11 q^4 t^5+13 q^5 t^5+2 q^6 t^5+q^3 t^6+7 q^4 t^6
+16 q^5 t^6+5 q^6 t^6+4 q^4 t^7+15 q^5 t^7+9 q^6 t^7+q^7 t^7
+2 q^4 t^8+12 q^5 t^8+13 q^6 t^8+2 q^7 t^8+q^4 t^9+7 q^5 t^9
+15 q^6 t^9+4 q^7 t^9+4 q^5 t^{10}+14 q^6 t^{10}+6 q^7 t^{10}
+2 q^5 t^{11}+11 q^6 t^{11}+8 q^7 t^{11}+q^5 t^{12}+7 q^6 t^{12}
+10 q^7 t^{12}+q^8 t^{12}+4 q^6 t^{13}+10 q^7 t^{13}+2 q^8 t^{13}
+2 q^6 t^{14}+9 q^7 t^{14}+3 q^8 t^{14}+q^6 t^{15}+6 q^7 t^{15}
+4 q^8 t^{15}+4 q^7 t^{16}+5 q^8 t^{16}+2 q^7 t^{17}+5 q^8 t^{17}
+q^7 t^{18}+4 q^8 t^{18}+3 q^8 t^{19}+q^9 t^{19}+2 q^8 t^{20}
+q^9 t^{20}+q^8 t^{21}+q^9 t^{21}+q^9 t^{22}+q^9 t^{23}
+q^9 t^{24}\bigr)
\)

\noindent
\(
+a^2 \bigl(4 q^3+10 q^4+8 q^5+q^6+\frac{q^3}{t^4}+\frac{q^2}{t^3}
+\frac{2 q^3}{t^3}+\frac{q^4}{t^3}+\frac{q^2}{t^2}+\frac{4 q^3}{t^2}
+\frac{3 q^4}{t^2}+\frac{q^5}{t^2}+\frac{q^2}{t}+\frac{4 q^3}{t}
+\frac{7 q^4}{t}+\frac{3 q^5}{t}+2 q^3 t+11 q^4 t+12 q^5 t+3 q^6 t
+q^3 t^2+8 q^4 t^2+17 q^5 t^2+7 q^6 t^2+q^7 t^2+5 q^4 t^3
+17 q^5 t^3+13 q^6 t^3+2 q^7 t^3+2 q^4 t^4+15 q^5 t^4+18 q^6 t^4
+5 q^7 t^4+q^4 t^5+9 q^5 t^5+21 q^6 t^5+7 q^7 t^5+5 q^5 t^6
+19 q^6 t^6+12 q^7 t^6+q^8 t^6+2 q^5 t^7+15 q^6 t^7+14 q^7 t^7
+2 q^8 t^7+q^5 t^8+9 q^6 t^8+17 q^7 t^8+4 q^8 t^8+5 q^6 t^9
+15 q^7 t^9+6 q^8 t^9+2 q^6 t^{10}+13 q^7 t^{10}+8 q^8 t^{10}
+q^6 t^{11}+8 q^7 t^{11}+9 q^8 t^{11}+5 q^7 t^{12}+9 q^8 t^{12}
+q^9 t^{12}+2 q^7 t^{13}+8 q^8 t^{13}+q^9 t^{13}+q^7 t^{14}
+6 q^8 t^{14}+2 q^9 t^{14}+4 q^8 t^{15}+2 q^9 t^{15}+2 q^8 t^{16}
+3 q^9 t^{16}+q^8 t^{17}+2 q^9 t^{17}+2 q^9 t^{18}+q^9 t^{19}
+q^9 t^{20}\bigr)
\)

\noindent
\(
+a^3 \bigl(q^4+7 q^5+12 q^6+4 q^7+\frac{q^4}{t^5}
+\frac{q^5}{t^5}+\frac{2 q^4}{t^4}+\frac{2 q^5}{t^4}+\frac{q^3}{t^3}
+\frac{3 q^4}{t^3}+\frac{4 q^5}{t^3}+\frac{2 q^6}{t^3}
+\frac{3 q^4}{t^2}+\frac{6 q^5}{t^2}+\frac{4 q^6}{t^2}
+\frac{q^7}{t^2}+\frac{2 q^4}{t}+\frac{8 q^5}{t}+\frac{8 q^6}{t}
+\frac{2 q^7}{t}+5 q^5 t+13 q^6 t+7 q^7 t+q^8 t+2 q^5 t^2
+12 q^6 t^2+10 q^7 t^2+2 q^8 t^2+q^5 t^3+9 q^6 t^3+13 q^7 t^3
+4 q^8 t^3+5 q^6 t^4+13 q^7 t^4+6 q^8 t^4+2 q^6 t^5+12 q^7 t^5
+8 q^8 t^5+q^6 t^6+8 q^7 t^6+9 q^8 t^6+q^9 t^6+5 q^7 t^7+9 q^8 t^7
+q^9 t^7+2 q^7 t^8+8 q^8 t^8+2 q^9 t^8+q^7 t^9+6 q^8 t^9+3 q^9 t^9
+4 q^8 t^{10}+3 q^9 t^{10}+2 q^8 t^{11}+3 q^9 t^{11}+q^8 t^{12}
+3 q^9 t^{12}+2 q^9 t^{13}+q^9 t^{14}+q^9 t^{15}\bigr)
\)

\noindent
\(
+a^4 \bigl(q^6+5 q^7+4 q^8+\frac{q^6}{t^6}+\frac{q^5}{t^5}
+\frac{2 q^6}{t^5}+\frac{q^7}{t^5}+\frac{q^5}{t^4}+\frac{3 q^6}{t^4}
+\frac{q^7}{t^4}+\frac{q^5}{t^3}+\frac{3 q^6}{t^3}+\frac{3 q^7}{t^3}
+\frac{q^8}{t^3}+\frac{3 q^6}{t^2}+\frac{4 q^7}{t^2}+\frac{2 q^8}{t^2}
+\frac{2 q^6}{t}+\frac{6 q^7}{t}+\frac{3 q^8}{t}+4 q^7 t+5 q^8 t+q^9 t
+2 q^7 t^2+5 q^8 t^2+q^9 t^2+q^7 t^3+4 q^8 t^3+2 q^9 t^3+3 q^8 t^4
+2 q^9 t^4+2 q^8 t^5+3 q^9 t^5+q^8 t^6+2 q^9 t^6+2 q^9 t^7+q^9 t^8
+q^9 t^9\bigr)
\)

\noindent
\(
+a^5 \bigl(q^9+\frac{q^7}{t^6}
+\frac{q^8}{t^6}+\frac{q^7}{t^5}+\frac{q^8}{t^5}+\frac{q^7}{t^4}
+\frac{q^8}{t^4}+\frac{q^8}{t^3}+\frac{q^9}{t^3}+\frac{q^8}{t^2}
+\frac{q^9}{t^2}+\frac{q^8}{t}+\frac{q^9}{t}+q^9 t+q^9 t^2\bigr)
\)\,.

}
\renewcommand{\baselinestretch}{1.2} 
\medskip

\subsection{Young diagram 2x2}
$$
H\!D_{3,2}(b=2\om_2;\,q,t,a)\ =
$$

\renewcommand{\baselinestretch}{0.5} 
{\small

\noindent
\(
1+\frac{a^4 q^{10}}{t^2}+q^2 t+q^3 t+q^2 t^2+q^3 t^2+q^4 t^2
+q^4 t^3+q^5 t^3+q^4 t^4+q^5 t^4+q^6 t^4+q^6 t^5+q^7 t^5
+q^6 t^6+q^7 t^6+q^8 t^8+a^3 \bigl(q^9+q^{10}+\frac{q^7}{t^2}
+\frac{q^8}{t^2}+\frac{q^7}{t}+\frac{q^8}{t}+q^9 t
+q^{10} t\bigr)
\)

\noindent
\(
+a \bigl(q^2+q^3+q^4+q^5
+\frac{q^2}{t}+\frac{q^3}{t}+2 q^4 t+3 q^5 t+q^6 t+q^4 t^2
+2 q^5 t^2+2 q^6 t^2+q^7 t^2+2 q^6 t^3+3 q^7 t^3+q^8 t^3
+q^6 t^4+2 q^7 t^4+q^8 t^4+q^8 t^5+q^9 t^5+q^8 t^6
+q^9 t^6\bigr)
\)

\noindent
\(
+a^2 \bigl(q^5+q^6+2 q^7+q^8+\frac{q^5}{t^2}
+\frac{q^4}{t}+\frac{2 q^5}{t}+\frac{q^6}{t}+q^6 t+3 q^7 t
+2 q^8 t+q^7 t^2+q^8 t^2+q^9 t^2+q^8 t^3+2 q^9 t^3
+q^{10} t^3+q^9 t^4\bigr)\,.
\)

}
\renewcommand{\baselinestretch}{1.2} 
\smallskip
Here all coefficients here are positive, as well
as in the following formulas, which 
confirms the positivity part of 
Conjecture \ref{HOMFLY}.
We note that 
$$
H\!D_{3,2}(b=2\om_2;\,q,t=1,a)\ =(H\!D_{3,2}(b=2\om_1;\,q,t=1,a)^2.
$$
The sum of all coefficients is $3^4$. See Conjecture \ref{CONJEVAL}.

\medskip

{\em The case of $b=2\om_2$ for knot\ $\{5,2\}$.}
Let us give another example for the self-dual 
$2\times 2$\~square corresponding to $2\om_2$.
It demonstrate the duality (which always holds)
and the positivity of all coefficients, which
seems (experimentally) a special property
of rectangles.
 
$$
H\!D_{5,2}(b=2\om_2;\,q,t,a)\ =
$$

\renewcommand{\baselinestretch}{0.5} 
{\small

\noindent
\(
a^0\bigl(1+q^2 t+q^3 t+q^2 t^2+q^3 t^2+q^4 t^2+q^5 t^2+q^6 t^2
+q^4 t^3+2 q^5 t^3+2 q^6 t^3+q^7 t^3+q^4 t^4+q^5 t^4+2 q^6 t^4
+2 q^7 t^4+2 q^8 t^4+q^6 t^5+2 q^7 t^5+3 q^8 t^5+2 q^9 t^5+q^6 t^6
+q^7 t^6+2 q^8 t^6+3 q^9 t^6+2 q^{10} t^6+q^8 t^7+2 q^9 t^7
+3 q^{10} t^7+2 q^{11} t^7+q^8 t^8+q^9 t^8+2 q^{10} t^8+3 q^{11} t^8
+2 q^{12} t^8+q^{10} t^9+2 q^{11} t^9+2 q^{12} t^9+q^{13} t^9
+q^{10} t^{10}+q^{11} t^{10}+2 q^{12} t^{10}+2 q^{13} t^{10}
+q^{14} t^{10}+q^{12} t^{11}+2 q^{13} t^{11}+q^{14} t^{11}
+q^{12} t^{12}+q^{13} t^{12}+q^{14} t^{12}+q^{14} t^{13}
+q^{15} t^{13}+q^{14} t^{14}+q^{15} t^{14}+q^{16} t^{16}\bigr)
\)

\noindent
\(
+a^1 \bigl(q^2+q^3+q^4+2 q^5+q^6+\frac{q^2}{t}+\frac{q^3}{t}
+2 q^4 t+4 q^5 t+3 q^6 t+2 q^7 t+q^8 t+q^4 t^2+2 q^5 t^2
+3 q^6 t^2+5 q^7 t^2+5 q^8 t^2+2 q^9 t^2+2 q^6 t^3+5 q^7 t^3
+7 q^8 t^3+6 q^9 t^3+2 q^{10} t^3+q^6 t^4+2 q^7 t^4+4 q^8 t^4
+7 q^9 t^4+7 q^{10} t^4+3 q^{11} t^4+2 q^8 t^5+5 q^9 t^5
+8 q^{10} t^5+8 q^{11} t^5+3 q^{12} t^5+q^8 t^6+2 q^9 t^6
+4 q^{10} t^6+8 q^{11} t^6+7 q^{12} t^6+2 q^{13} t^6+2 q^{10} t^7
+5 q^{11} t^7+7 q^{12} t^7+6 q^{13} t^7+2 q^{14} t^7+q^{10} t^8
+2 q^{11} t^8+4 q^{12} t^8+7 q^{13} t^8+5 q^{14} t^8+q^{15} t^8
+2 q^{12} t^9+5 q^{13} t^9+5 q^{14} t^9+2 q^{15} t^9+q^{12} t^{10}
+2 q^{13} t^{10}+3 q^{14} t^{10}+3 q^{15} t^{10}+q^{16} t^{10}
+2 q^{14} t^{11}+4 q^{15} t^{11}+2 q^{16} t^{11}+q^{14} t^{12}
+2 q^{15} t^{12}+q^{16} t^{12}+q^{16} t^{13}+q^{17} t^{13}
+q^{16} t^{14}+q^{17} t^{14}\bigr)
\)

\noindent
\(
+a^2 \bigl(q^5+q^6+4 q^7+4 q^8+2 q^9+\frac{q^5}{t^2}+\frac{q^4}{t}
+\frac{2 q^5}{t}+\frac{q^6}{t}+\frac{q^7}{t}+\frac{q^8}{t}+q^6 t
+4 q^7 t+5 q^8 t+5 q^9 t+4 q^{10} t+q^{11} t+q^7 t^2+2 q^8 t^2
+6 q^9 t^2+8 q^{10} t^2+6 q^{11} t^2+q^{12} t^2+q^8 t^3+4 q^9 t^3
+7 q^{10} t^3+9 q^{11} t^3+6 q^{12} t^3+q^{13} t^3+q^9 t^4
+2 q^{10} t^4+7 q^{11} t^4+9 q^{12} t^4+6 q^{13} t^4+q^{14} t^4
+q^{10} t^5+4 q^{11} t^5+7 q^{12} t^5+9 q^{13} t^5+6 q^{14} t^5
+q^{15} t^5+q^{11} t^6+2 q^{12} t^6+7 q^{13} t^6+8 q^{14} t^6
+4 q^{15} t^6+q^{12} t^7+4 q^{13} t^7+6 q^{14} t^7+5 q^{15} t^7
+2 q^{16} t^7+q^{13} t^8+2 q^{14} t^8+5 q^{15} t^8+4 q^{16} t^8
+q^{17} t^8+q^{14} t^9+4 q^{15} t^9+4 q^{16} t^9+q^{17} t^9
+q^{15} t^{10}+q^{16} t^{10}+q^{17} t^{10}+q^{16} t^{11}
+2 q^{17} t^{11}+q^{18} t^{11}+q^{17} t^{12}\bigr)
\)

\noindent
\(
+a^3 \bigl(2 q^9+4 q^{10}+3 q^{11}+q^{12}+\frac{q^7}{t^2}
+\frac{q^8}{t^2}+\frac{q^7}{t}+\frac{q^8}{t}+\frac{q^9}{t}
+\frac{2 q^{10}}{t}+\frac{q^{11}}{t}+q^9 t+2 q^{10} t+3 q^{11} t
+4 q^{12} t+2 q^{13} t+2 q^{11} t^2+5 q^{12} t^2+5 q^{13} t^2
+2 q^{14} t^2+q^{11} t^3+2 q^{12} t^3+4 q^{13} t^3+5 q^{14} t^3
+2 q^{15} t^3+2 q^{13} t^4+5 q^{14} t^4+4 q^{15} t^4+q^{16} t^4
+q^{13} t^5+2 q^{14} t^5+3 q^{15} t^5+3 q^{16} t^5+q^{17} t^5
+2 q^{15} t^6+4 q^{16} t^6+2 q^{17} t^6+q^{15} t^7+2 q^{16} t^7
+q^{17} t^7+q^{17} t^8+q^{18} t^8+q^{17} t^9+q^{18} t^9\bigr)
\)

\noindent
\(
+a^4 \bigl(q^{12}+q^{13}+q^{14}+\frac{q^{10}}{t^2}+\frac{q^{12}}{t}
+\frac{q^{13}}{t}+q^{14} t+q^{15} t+q^{14} t^2+q^{15} t^2
+q^{16} t^2+q^{16} t^3+q^{17} t^3+q^{16} t^4+q^{17} t^4
+q^{18} t^6\bigr).
\)
}
\renewcommand{\baselinestretch}{1.2} 
\medskip

The sum of all coefficients here is $5^4$; moreover,
$$
H\!D_{5,2}(b=2\om_2;\,q,t=1,a)\ =(H\!D_{5,2}(b=2\om_1;\,q,t=1,a)^2.
$$
Correspondingly, $H\!D_K(b=\om_1+\om_2;\,q,t=1,a)$
$$
=\,H\!D_K(b=2\om_1;\,q,t=1,a)\, H\!D_K(b=\om_1;\,q,t=1,a) 
$$
for $K=\{3,2\},\{4,3\}$ in the next subsection, which
confirms Conjecture \ref{CONJEVAL}. 

The leading term will be $a^8$ in next example of $\{4,3\}$, 
which agrees  
with the formula $(\hbox{min}(4,3)-1)*(4\, \hbox{boxes})$
from this conjecture.

\medskip

{\em The case of $b=2\om_2$ for the knot\ $\{4,3\}$}.
In contrast to the previous example, the leading
$a$\~term here has a pure 
monomial $q,t$\~coefficient. Notice the positivity of all
coefficients, which seems the most unexpected part
of Conjecture \ref{HOMFLY}. The dimension is $11^4$;
the other claims from Conjecture \ref{CONJEVAL} hold as well.

$$
H\!D_{4,3}(b=2\om_2;\,q,t,a)\ =
$$

\renewcommand{\baselinestretch}{0.5}
{\small 
\noindent
\(
1+\frac{a^8 q^{28}}{t^4}+q^2 t+q^3 t+q^4 t+q^5 t+q^2 t^2+q^3 t^2+2 
q^4 t^2+2 q^5 t^2+2 q^6 t^2+q^7 t^2+q^8 t^2+q^4 t^3+2 q^5 t^3+4 q^6 
t^3+4 q^7 t^3+4 q^8 t^3+3 q^9 t^3+q^4 t^4+q^5 t^4+3 q^6 t^4+4 q^7 
t^4+7 q^8 t^4+6 q^9 t^4+5 q^{10} t^4+2 q^{11} t^4+q^6 t^5+2 q^7 
t^5+5 q^8 t^5+7 q^9 t^5+10 q^{10} t^5+9 q^{11} t^5+3 q^{12} 
t^5+q^{13} t^5+q^6 t^6+q^7 t^6+3 q^8 t^6+5 q^9 t^6+9 q^{10} t^6+11 
q^{11} t^6+11 q^{12} t^6+5 q^{13} t^6+q^{14} t^6+q^8 t^7+2 q^9 
t^7+5 q^{10} t^7+8 q^{11} t^7+13 q^{12} t^7+14 q^{13} t^7+7 q^{14} 
t^7+2 q^{15} t^7+q^8 t^8+q^9 t^8+3 q^{10} t^8+5 q^{11} t^8+10 
q^{12} t^8+14 q^{13} t^8+15 q^{14} t^8+7 q^{15} t^8+2 q^{16} 
t^8+q^{10} t^9+2 q^{11} t^9+5 q^{12} t^9+8 q^{13} t^9+15 q^{14} 
t^9+16 q^{15} t^9+7 q^{16} t^9+2 q^{17} t^9+q^{10} t^{10}+q^{11} 
t^{10}+3 q^{12} t^{10}+5 q^{13} t^{10}+10 q^{14} t^{10}+15 q^{15} 
t^{10}+15 q^{16} t^{10}+7 q^{17} t^{10}+q^{18} t^{10}+q^{12} 
t^{11}+2 q^{13} t^{11}+5 q^{14} t^{11}+8 q^{15} t^{11}+14 q^{16} 
t^{11}+14 q^{17} t^{11}+5 q^{18} t^{11}+q^{19} t^{11}+q^{12} 
t^{12}+q^{13} t^{12}+3 q^{14} t^{12}+5 q^{15} t^{12}+10 q^{16} 
t^{12}+13 q^{17} t^{12}+11 q^{18} t^{12}+3 q^{19} t^{12}+q^{14} 
t^{13}+2 q^{15} t^{13}+5 q^{16} t^{13}+8 q^{17} t^{13}+11 q^{18} 
t^{13}+9 q^{19} t^{13}+2 q^{20} t^{13}+q^{14} t^{14}+q^{15} 
t^{14}+3 q^{16} t^{14}+5 q^{17} t^{14}+9 q^{18} t^{14}+10 q^{19} 
t^{14}+5 q^{20} t^{14}+q^{16} t^{15}+2 q^{17} t^{15}+5 q^{18} 
t^{15}+7 q^{19} t^{15}+6 q^{20} t^{15}+3 q^{21} t^{15}+q^{16} 
t^{16}+q^{17} t^{16}+3 q^{18} t^{16}+5 q^{19} t^{16}+7 q^{20} 
t^{16}+4 q^{21} t^{16}+q^{22} t^{16}+q^{18} t^{17}+2 q^{19} 
t^{17}+4 q^{20} t^{17}+4 q^{21} t^{17}+q^{22} t^{17}+q^{18} 
t^{18}+q^{19} t^{18}+3 q^{20} t^{18}+4 q^{21} t^{18}+2 q^{22} 
t^{18}+q^{20} t^{19}+2 q^{21} t^{19}+2 q^{22} t^{19}+q^{23} 
t^{19}+q^{20} t^{20}+q^{21} t^{20}+2 q^{22} t^{20}+q^{23} 
t^{20}+q^{22} t^{21}+q^{23} t^{21}+q^{22} t^{22}+q^{23} 
t^{22}+q^{24} t^{24}+
\)

\vfill
\noindent
\(
+a \bigl(q^2+q^3+2 q^4+3 q^5+3 q^6+3 
q^7+2 
q^8+q^9+\frac{q^2}{t}+\frac{q^3}{t}+\frac{q^4}{t}+\frac{q^5}{t}+2 
q^4 t+4 q^5 t+7 q^6 t+9 q^7 t+9 q^8 t+7 q^9 t+3 q^{10} t+q^{11} 
t+q^4 t^2+2 q^5 t^2+5 q^6 t^2+9 q^7 t^2+14 q^8 t^2+16 q^9 t^2+15 
q^{10} t^2+11 q^{11} t^2+4 q^{12} t^2+q^{13} t^2+2 q^6 t^3+5 q^7 
t^3+12 q^8 t^3+19 q^9 t^3+27 q^{10} t^3+29 q^{11} t^3+19 q^{12} 
t^3+9 q^{13} t^3+2 q^{14} t^3+q^6 t^4+2 q^7 t^4+6 q^8 t^4+12 q^9 
t^4+23 q^{10} t^4+33 q^{11} t^4+37 q^{12} t^4+31 q^{13} t^4+15 
q^{14} t^4+5 q^{15} t^4+q^{16} t^4+2 q^8 t^5+5 q^9 t^5+13 q^{10} 
t^5+24 q^{11} t^5+40 q^{12} t^5+51 q^{13} t^5+41 q^{14} t^5+22 
q^{15} t^5+7 q^{16} t^5+q^{17} t^5+q^8 t^6+2 q^9 t^6+6 q^{10} 
t^6+13 q^{11} t^6+26 q^{12} t^6+43 q^{13} t^6+55 q^{14} t^6+48 
q^{15} t^6+25 q^{16} t^6+8 q^{17} t^6+q^{18} t^6+2 q^{10} t^7+5 
q^{11} t^7+13 q^{12} t^7+25 q^{13} t^7+46 q^{14} t^7+63 q^{15} 
t^7+53 q^{16} t^7+28 q^{17} t^7+8 q^{18} t^7+q^{19} t^7+q^{10} 
t^8+2 q^{11} t^8+6 q^{12} t^8+13 q^{13} t^8+27 q^{14} t^8+47 q^{15} 
t^8+62 q^{16} t^8+53 q^{17} t^8+25 q^{18} t^8+7 q^{19} t^8+q^{20} 
t^8+2 q^{12} t^9+5 q^{13} t^9+13 q^{14} t^9+25 q^{15} t^9+47 q^{16} 
t^9+63 q^{17} t^9+48 q^{18} t^9+22 q^{19} t^9+5 q^{20} t^9+q^{12} 
t^{10}+2 q^{13} t^{10}+6 q^{14} t^{10}+13 q^{15} t^{10}+27 q^{16} 
t^{10}+46 q^{17} t^{10}+55 q^{18} t^{10}+41 q^{19} t^{10}+15 q^{20} 
t^{10}+2 q^{21} t^{10}+2 q^{14} t^{11}+5 q^{15} t^{11}+13 q^{16} 
t^{11}+25 q^{17} t^{11}+43 q^{18} t^{11}+51 q^{19} t^{11}+31 q^{20} 
t^{11}+9 q^{21} t^{11}+q^{22} t^{11}+q^{14} t^{12}+2 q^{15} 
t^{12}+6 q^{16} t^{12}+13 q^{17} t^{12}+26 q^{18} t^{12}+40 q^{19} 
t^{12}+37 q^{20} t^{12}+19 q^{21} t^{12}+4 q^{22} t^{12}+2 q^{16} 
t^{13}+5 q^{17} t^{13}+13 q^{18} t^{13}+24 q^{19} t^{13}+33 q^{20} 
t^{13}+29 q^{21} t^{13}+11 q^{22} t^{13}+q^{23} t^{13}+q^{16} 
t^{14}+2 q^{17} t^{14}+6 q^{18} t^{14}+13 q^{19} t^{14}+23 q^{20} 
t^{14}+27 q^{21} t^{14}+15 q^{22} t^{14}+3 q^{23} t^{14}+2 q^{18} 
t^{15}+5 q^{19} t^{15}+12 q^{20} t^{15}+19 q^{21} t^{15}+16 q^{22} 
t^{15}+7 q^{23} t^{15}+q^{24} t^{15}+q^{18} t^{16}+2 q^{19} 
t^{16}+6 q^{20} t^{16}+12 q^{21} t^{16}+14 q^{22} t^{16}+9 q^{23} 
t^{16}+2 q^{24} t^{16}+2 q^{20} t^{17}+5 q^{21} t^{17}+9 q^{22} 
t^{17}+9 q^{23} t^{17}+3 q^{24} t^{17}+q^{20} t^{18}+2 q^{21} 
t^{18}+5 q^{22} t^{18}+7 q^{23} t^{18}+3 q^{24} t^{18}+2 q^{22} 
t^{19}+4 q^{23} t^{19}+3 q^{24} t^{19}+q^{25} t^{19}+q^{22} 
t^{20}+2 q^{23} t^{20}+2 q^{24} t^{20}+q^{25} t^{20}+q^{24} 
t^{21}+q^{25} t^{21}+q^{24} t^{22}+q^{25} t^{22}\bigr) 
\)

\vfill
\noindent
\(
+a^2 \bigl(q^5+2 q^6+6 q^7+9 q^8+14 
q^9+15 q^{10}+15 q^{11}+9 q^{12}+4 
q^{13}+q^{14}+\frac{q^5}{t^2}+\frac{q^6}{t^2}+\frac{2 
q^7}{t^2}+\frac{q^8}{t^2}+\frac{q^9}{t^2}+\frac{q^4}{t}+\frac{2 
q^5}{t}+\frac{3 q^6}{t}+\frac{5 q^7}{t}+\frac{5 q^8}{t}+\frac{5 
q^9}{t}+\frac{4 q^{10}}{t}+\frac{2 q^{11}}{t}+\frac{q^{12}}{t}+q^6 
t+4 q^7 t+9 q^8 t+17 q^9 t+25 q^{10} t+31 q^{11} t+28 q^{12} t+19 
q^{13} t+10 q^{14} t+3 q^{15} t+q^{16} t+q^7 t^2+3 q^8 t^2+10 q^9 
t^2+20 q^{10} t^2+35 q^{11} t^2+45 q^{12} t^2+51 q^{13} t^2+37 
q^{14} t^2+21 q^{15} t^2+8 q^{16} t^2+2 q^{17} t^2+q^8 t^3+4 q^9 
t^3+11 q^{10} t^3+25 q^{11} t^3+44 q^{12} t^3+66 q^{13} t^3+70 
q^{14} t^3+54 q^{15} t^3+31 q^{16} t^3+11 q^{17} t^3+3 q^{18} 
t^3+q^9 t^4+3 q^{10} t^4+11 q^{11} t^4+24 q^{12} t^4+49 q^{13} 
t^4+73 q^{14} t^4+89 q^{15} t^4+71 q^{16} t^4+41 q^{17} t^4+16 
q^{18} t^4+4 q^{19} t^4+q^{20} t^4+q^{10} t^5+4 q^{11} t^5+11 
q^{12} t^5+27 q^{13} t^5+52 q^{14} t^5+88 q^{15} t^5+103 q^{16} 
t^5+82 q^{17} t^5+46 q^{18} t^5+16 q^{19} t^5+4 q^{20} t^5+q^{11} 
t^6+3 q^{12} t^6+11 q^{13} t^6+25 q^{14} t^6+54 q^{15} t^6+88 
q^{16} t^6+109 q^{17} t^6+85 q^{18} t^6+46 q^{19} t^6+16 q^{20} 
t^6+3 q^{21} t^6+q^{12} t^7+4 q^{13} t^7+11 q^{14} t^7+27 q^{15} 
t^7+55 q^{16} t^7+95 q^{17} t^7+109 q^{18} t^7+82 q^{19} t^7+41 
q^{20} t^7+11 q^{21} t^7+2 q^{22} t^7+q^{13} t^8+3 q^{14} t^8+11 
q^{15} t^8+25 q^{16} t^8+55 q^{17} t^8+88 q^{18} t^8+103 q^{19} 
t^8+71 q^{20} t^8+31 q^{21} t^8+8 q^{22} t^8+q^{23} t^8+q^{14} 
t^9+4 q^{15} t^9+11 q^{16} t^9+27 q^{17} t^9+54 q^{18} t^9+88 
q^{19} t^9+89 q^{20} t^9+54 q^{21} t^9+21 q^{22} t^9+3 q^{23} 
t^9+q^{15} t^{10}+3 q^{16} t^{10}+11 q^{17} t^{10}+25 q^{18} 
t^{10}+52 q^{19} t^{10}+73 q^{20} t^{10}+70 q^{21} t^{10}+37 q^{22} 
t^{10}+10 q^{23} t^{10}+q^{24} t^{10}+q^{16} t^{11}+4 q^{17} 
t^{11}+11 q^{18} t^{11}+27 q^{19} t^{11}+49 q^{20} t^{11}+66 q^{21} 
t^{11}+51 q^{22} t^{11}+19 q^{23} t^{11}+4 q^{24} t^{11}+q^{17} 
t^{12}+3 q^{18} t^{12}+11 q^{19} t^{12}+24 q^{20} t^{12}+44 q^{21} 
t^{12}+45 q^{22} t^{12}+28 q^{23} t^{12}+9 q^{24} t^{12}+q^{25} 
t^{12}+q^{18} t^{13}+4 q^{19} t^{13}+11 q^{20} t^{13}+25 q^{21} 
t^{13}+35 q^{22} t^{13}+31 q^{23} t^{13}+15 q^{24} t^{13}+2 q^{25} 
t^{13}+q^{19} t^{14}+3 q^{20} t^{14}+11 q^{21} t^{14}+20 q^{22} 
t^{14}+25 q^{23} t^{14}+15 q^{24} t^{14}+4 q^{25} t^{14}+q^{20} 
t^{15}+4 q^{21} t^{15}+10 q^{22} t^{15}+17 q^{23} t^{15}+14 q^{24} 
t^{15}+5 q^{25} t^{15}+q^{26} t^{15}+q^{21} t^{16}+3 q^{22} 
t^{16}+9 q^{23} t^{16}+9 q^{24} t^{16}+5 q^{25} t^{16}+q^{26} 
t^{16}+q^{22} t^{17}+4 q^{23} t^{17}+6 q^{24} t^{17}+5 q^{25} 
t^{17}+2 q^{26} t^{17}+q^{23} t^{18}+2 q^{24} t^{18}+3 q^{25} 
t^{18}+q^{26} t^{18}+q^{24} t^{19}+2 q^{25} t^{19}+q^{26} 
t^{19}+q^{25} t^{20}\bigr)
\)

\vfill
\noindent
\(
+a^3 \bigl(3 q^9+7 q^{10}+15 q^{11}+24 q^{12}+32 
q^{13}+33 q^{14}+25 q^{15}+15 q^{16}+6 q^{17}+2 
q^{18}+\frac{q^9}{t^3}+\frac{q^{10}}{t^3}+\frac{q^{11}}{t^3}+
\frac{q^{12}}{t^3}+\frac{q^7}{t^2}+\frac{2 
q^8}{t^2}+\frac{4 q^9}{t^2}+\frac{5 q^{10}}{t^2}+\frac{6 
q^{11}}{t^2}+\frac{5 q^{12}}{t^2}+\frac{3 q^{13}}{t^2}+\frac{2 
q^{14}}{t^2}+\frac{q^7}{t}+\frac{2 q^8}{t}+\frac{5 q^9}{t}+\frac{8 
q^{10}}{t}+\frac{12 q^{11}}{t}+\frac{14 q^{12}}{t}+\frac{14 
q^{13}}{t}+\frac{12 q^{14}}{t}+\frac{6 q^{15}}{t}+\frac{3 
q^{16}}{t}+\frac{q^{17}}{t}+q^9 t+3 q^{10} t+10 q^{11} t+20 q^{12} 
t+36 q^{13} t+48 q^{14} t+50 q^{15} t+42 q^{16} t+25 q^{17} t+12 
q^{18} t+4 q^{19} t+q^{20} t+3 q^{11} t^2+9 q^{12} t^2+24 q^{13} 
t^2+44 q^{14} t^2+67 q^{15} t^2+76 q^{16} t^2+62 q^{17} t^2+39 
q^{18} t^2+17 q^{19} t^2+6 q^{20} t^2+q^{21} t^2+q^{11} t^3+3 
q^{12} t^3+11 q^{13} t^3+25 q^{14} t^3+53 q^{15} t^3+80 q^{16} 
t^3+90 q^{17} t^3+76 q^{18} t^3+44 q^{19} t^3+20 q^{20} t^3+6 
q^{21} t^3+q^{22} t^3+3 q^{13} t^4+9 q^{14} t^4+26 q^{15} t^4+54 
q^{16} t^4+90 q^{17} t^4+105 q^{18} t^4+85 q^{19} t^4+51 q^{20} 
t^4+20 q^{21} t^4+6 q^{22} t^4+q^{23} t^4+q^{13} t^5+3 q^{14} 
t^5+11 q^{15} t^5+26 q^{16} t^5+59 q^{17} t^5+95 q^{18} t^5+107 
q^{19} t^5+85 q^{20} t^5+44 q^{21} t^5+17 q^{22} t^5+4 q^{23} t^5+3 
q^{15} t^6+9 q^{16} t^6+26 q^{17} t^6+57 q^{18} t^6+95 q^{19} 
t^6+105 q^{20} t^6+76 q^{21} t^6+39 q^{22} t^6+12 q^{23} t^6+2 
q^{24} t^6+q^{15} t^7+3 q^{16} t^7+11 q^{17} t^7+26 q^{18} t^7+59 
q^{19} t^7+90 q^{20} t^7+90 q^{21} t^7+62 q^{22} t^7+25 q^{23} 
t^7+6 q^{24} t^7+q^{25} t^7+3 q^{17} t^8+9 q^{18} t^8+26 q^{19} 
t^8+54 q^{20} t^8+80 q^{21} t^8+76 q^{22} t^8+42 q^{23} t^8+15 
q^{24} t^8+3 q^{25} t^8+q^{17} t^9+3 q^{18} t^9+11 q^{19} t^9+26 
q^{20} t^9+53 q^{21} t^9+67 q^{22} t^9+50 q^{23} t^9+25 q^{24} 
t^9+6 q^{25} t^9+3 q^{19} t^{10}+9 q^{20} t^{10}+25 q^{21} 
t^{10}+44 q^{22} t^{10}+48 q^{23} t^{10}+33 q^{24} t^{10}+12 q^{25} 
t^{10}+2 q^{26} t^{10}+q^{19} t^{11}+3 q^{20} t^{11}+11 q^{21} 
t^{11}+24 q^{22} t^{11}+36 q^{23} t^{11}+32 q^{24} t^{11}+14 q^{25} 
t^{11}+3 q^{26} t^{11}+3 q^{21} t^{12}+9 q^{22} t^{12}+20 q^{23} 
t^{12}+24 q^{24} t^{12}+14 q^{25} t^{12}+5 q^{26} t^{12}+q^{27} 
t^{12}+q^{21} t^{13}+3 q^{22} t^{13}+10 q^{23} t^{13}+15 q^{24} 
t^{13}+12 q^{25} t^{13}+6 q^{26} t^{13}+q^{27} t^{13}+3 q^{23} 
t^{14}+7 q^{24} t^{14}+8 q^{25} t^{14}+5 q^{26} t^{14}+q^{27} 
t^{14}+q^{23} t^{15}+3 q^{24} t^{15}+5 q^{25} t^{15}+4 q^{26} 
t^{15}+q^{27} t^{15}+2 q^{25} t^{16}+2 q^{26} t^{16}+q^{25} 
t^{17}+q^{26} t^{17}\bigr)
\)

\vfill
\noindent
\(
+a^4 \bigl(q^{12}+4 q^{13}+12 q^{14}+21 q^{15}+33 q^{16}+36 q^{17}
+33 q^{18}+21 q^{19}+11 q^{20}+4 
q^{21}+q^{22}+\frac{q^{14}}{t^4}+\frac{q^{11}}{t^3}+\frac{2 
q^{12}}{t^3}+\frac{3 q^{13}}{t^3}+\frac{4 q^{14}}{t^3}+\frac{3 
q^{15}}{t^3}+\frac{2 
q^{16}}{t^3}+\frac{q^{17}}{t^3}+\frac{q^{10}}{t^2}+\frac{2 
q^{11}}{t^2}+\frac{5 q^{12}}{t^2}+\frac{7 q^{13}}{t^2}+\frac{10 
q^{14}}{t^2}+\frac{9 q^{15}}{t^2}+\frac{9 q^{16}}{t^2}+\frac{5 
q^{17}}{t^2}+\frac{3 
q^{18}}{t^2}+\frac{q^{19}}{t^2}+\frac{q^{11}}{t}+\frac{3 
q^{12}}{t}+\frac{7 q^{13}}{t}+\frac{13 q^{14}}{t}+\frac{18 
q^{15}}{t}+\frac{22 q^{16}}{t}+\frac{20 q^{17}}{t}+\frac{14 
q^{18}}{t}+\frac{8 q^{19}}{t}+\frac{3 
q^{20}}{t}+\frac{q^{21}}{t}+q^{13} t+4 q^{14} t+12 q^{15} t+26 
q^{16} t+40 q^{17} t+49 q^{18} t+44 q^{19} t+30 q^{20} t+16 q^{21} 
t+6 q^{22} t+2 q^{23} t+q^{14} t^2+4 q^{15} t^2+14 q^{16} t^2+30 
q^{17} t^2+53 q^{18} t^2+61 q^{19} t^2+56 q^{20} t^2+34 q^{21} 
t^2+17 q^{22} t^2+6 q^{23} t^2+q^{24} t^2+q^{15} t^3+4 q^{16} 
t^3+13 q^{17} t^3+32 q^{18} t^3+54 q^{19} t^3+66 q^{20} t^3+55 
q^{21} t^3+34 q^{22} t^3+16 q^{23} t^3+4 q^{24} t^3+q^{25} 
t^3+q^{16} t^4+4 q^{17} t^4+14 q^{18} t^4+33 q^{19} t^4+59 q^{20} 
t^4+66 q^{21} t^4+56 q^{22} t^4+30 q^{23} t^4+11 q^{24} t^4+3 
q^{25} t^4+q^{17} t^5+4 q^{18} t^5+13 q^{19} t^5+33 q^{20} t^5+54 
q^{21} t^5+61 q^{22} t^5+44 q^{23} t^5+21 q^{24} t^5+8 q^{25} 
t^5+q^{26} t^5+q^{18} t^6+4 q^{19} t^6+14 q^{20} t^6+32 q^{21} 
t^6+53 q^{22} t^6+49 q^{23} t^6+33 q^{24} t^6+14 q^{25} t^6+3 
q^{26} t^6+q^{19} t^7+4 q^{20} t^7+13 q^{21} t^7+30 q^{22} t^7+40 
q^{23} t^7+36 q^{24} t^7+20 q^{25} t^7+5 q^{26} t^7+q^{27} 
t^7+q^{20} t^8+4 q^{21} t^8+14 q^{22} t^8+26 q^{23} t^8+33 q^{24} 
t^8+22 q^{25} t^8+9 q^{26} t^8+2 q^{27} t^8+q^{21} t^9+4 q^{22} 
t^9+12 q^{23} t^9+21 q^{24} t^9+18 q^{25} t^9+9 q^{26} t^9+3 q^{27} 
t^9+q^{22} t^{10}+4 q^{23} t^{10}+12 q^{24} t^{10}+13 q^{25} 
t^{10}+10 q^{26} t^{10}+4 q^{27} t^{10}+q^{28} t^{10}+q^{23} 
t^{11}+4 q^{24} t^{11}+7 q^{25} t^{11}+7 q^{26} t^{11}+3 q^{27} 
t^{11}+q^{24} t^{12}+3 q^{25} t^{12}+5 q^{26} t^{12}+2 q^{27} 
t^{12}+q^{25} t^{13}+2 q^{26} t^{13}+q^{27} t^{13}+q^{26} 
t^{14}\bigr)
\)

\noindent
\(
+a^5 \bigl(q^{16}+3 q^{17}+9 
q^{18}+16 q^{19}+22 q^{20}+21 q^{21}+15 q^{22}+9 q^{23}+3 
q^{24}+q^{25}+\frac{q^{16}}{t^4}+\frac{q^{17}}{t^4}
+\frac{q^{18}}{t^4}+\frac{q^{19}}{t^4}+\frac{q^{14}}{t^3}+\frac{2 
q^{15}}{t^3}+\frac{4 q^{16}}{t^3}+\frac{5 q^{17}}{t^3}+\frac{5 
q^{18}}{t^3}+\frac{4 q^{19}}{t^3}+\frac{2 
q^{20}}{t^3}+\frac{q^{21}}{t^3}+\frac{q^{14}}{t^2}+\frac{2 
q^{15}}{t^2}+\frac{5 q^{16}}{t^2}+\frac{7 q^{17}}{t^2}+\frac{9 
q^{18}}{t^2}+\frac{9 q^{19}}{t^2}+\frac{7 q^{20}}{t^2}+\frac{5 
q^{21}}{t^2}+\frac{2 q^{22}}{t^2}+\frac{q^{23}}{t^2}+\frac{3 
q^{16}}{t}+\frac{6 q^{17}}{t}+\frac{12 q^{18}}{t}+\frac{17 
q^{19}}{t}+\frac{17 q^{20}}{t}+\frac{14 q^{21}}{t}+\frac{8 
q^{22}}{t}+\frac{4 q^{23}}{t}+\frac{q^{24}}{t}+3 q^{18} t+8 q^{19} 
t+18 q^{20} t+25 q^{21} t+25 q^{22} t+19 q^{23} t+9 q^{24} t+4 
q^{25} t+q^{26} t+q^{18} t^2+3 q^{19} t^2+10 q^{20} t^2+20 q^{21} 
t^2+28 q^{22} t^2+25 q^{23} t^2+15 q^{24} t^2+8 q^{25} t^2+2 q^{26} 
t^2+3 q^{20} t^3+9 q^{21} t^3+20 q^{22} t^3+25 q^{23} t^3+21 q^{24} 
t^3+14 q^{25} t^3+5 q^{26} t^3+q^{27} t^3+q^{20} t^4+3 q^{21} 
t^4+10 q^{22} t^4+18 q^{23} t^4+22 q^{24} t^4+17 q^{25} t^4+7 
q^{26} t^4+2 q^{27} t^4+3 q^{22} t^5+8 q^{23} t^5+16 q^{24} t^5+17 
q^{25} t^5+9 q^{26} t^5+4 q^{27} t^5+q^{28} t^5+q^{22} t^6+3 q^{23} 
t^6+9 q^{24} t^6+12 q^{25} t^6+9 q^{26} t^6+5 q^{27} t^6+q^{28} 
t^6+3 q^{24} t^7+6 q^{25} t^7+7 q^{26} t^7+5 q^{27} t^7+q^{28} 
t^7+q^{24} t^8+3 q^{25} t^8+5 q^{26} t^8+4 q^{27} t^8+q^{28} t^8+2 
q^{26} t^9+2 q^{27} t^9+q^{26} t^{10}+q^{27} t^{10}\bigr)
\)

\noindent
\(
+a^6 \bigl(q^{21}+3 q^{22}+6 q^{23}+6 q^{24}+6 
q^{25}+3 q^{26}+q^{27}+\frac{q^{19}}{t^4}+\frac{q^{20}}{t^4}+
\frac{2 q^{21}}{t^4}+\frac{q^{22}}{t^4}+\frac{q^{23}}{t^4}+
\frac{q^{18}}{t^3}+\frac{2 
q^{19}}{t^3}+\frac{3 q^{20}}{t^3}+\frac{4 q^{21}}{t^3}+\frac{3 
q^{22}}{t^3}+\frac{2 
q^{23}}{t^3}+\frac{q^{24}}{t^3}+\frac{q^{19}}{t^2}+\frac{2 
q^{20}}{t^2}+\frac{4 q^{21}}{t^2}+\frac{4 q^{22}}{t^2}+\frac{5 
q^{23}}{t^2}+\frac{3 q^{24}}{t^2}+\frac{2 
q^{25}}{t^2}+\frac{q^{26}}{t^2}+\frac{q^{20}}{t}+\frac{3 
q^{21}}{t}+\frac{6 q^{22}}{t}+\frac{8 q^{23}}{t}+\frac{6 
q^{24}}{t}+\frac{4 q^{25}}{t}+\frac{2 q^{26}}{t}+q^{22} t+3 q^{23} 
t+6 q^{24} t+8 q^{25} t+5 q^{26} t+2 q^{27} t+q^{28} t+q^{23} t^2+3 
q^{24} t^2+6 q^{25} t^2+4 q^{26} t^2+3 q^{27} t^2+q^{28} t^2+q^{24} 
t^3+3 q^{25} t^3+4 q^{26} t^3+4 q^{27} t^3+2 q^{28} t^3+q^{25} 
t^4+2 q^{26} t^4+3 q^{27} t^4+q^{28} t^4+q^{26} t^5+2 q^{27} 
t^5+q^{28} t^5+q^{27} t^6\bigr)
\)

\noindent
\(
+a^7 \bigl(q^{27}+q^{28}+\frac{q^{23}}{t^4}+\frac{q^{24}}{t^4}+
\frac{q^{25}}{t^4}+\frac{q^{26}}{t^4}+\frac{q^{23}}{t^3}+
\frac{q^{24}}{t^3}+\frac{q^{25}}{t^3}+\frac{q^{26}}{t^3}+
\frac{q^{25}}{t^2}+\frac{q^{26}}{t^2}+\frac{q^{27}}{t^2}+
\frac{q^{28}}{t^2}+\frac{q^{25}}{t}+\frac{q^{26}}{t}+
\frac{q^{27}}{t}+\frac{q^{28}}{t}+q^{27} 
t+q^{28} t\bigr).
\)
}
\renewcommand{\baselinestretch}{1.2} 

\medskip

\subsection{The 3-hook}\label{sec:3hook}
We will consider now $b=\om_1+\om_2$, beginning with 
the knot$\ =\{3,2\}$.
There can be a potential relation to formula (69) 
from (\cite{GIKV}), but we cannot establish it so far.
One has:
$$
H\!D_{3,2}(b=\om_1+\om_2;\,q,t,a)\ =
$$
\renewcommand{\baselinestretch}{0.5} 
{\small
\noindent
\(
1+\frac{a^3 q^6}{t}+2 q t-q t^2+2 q^2 t^2+q^3 t^2-q^2 t^3
+2 q^3 t^3-q^3 t^4+2 q^4 t^4+q^5 t^5
\)

\noindent
\(
+a \bigl(2 q^2+q^3+\frac{q}{t}-q^2 t+3 q^3 t+q^4 t-q^3 t^2
+3 q^4 t^2+q^5 t^2-q^4 t^3+2 q^5 t^3+q^6 t^4\bigr)
\)

\noindent
\(
+a^2 \bigl(q^4
+q^5+\frac{q^3}{t}+\frac{q^4}{t}-q^4 t+q^5 t+q^6 t+q^6 t^2\bigr)
\)\,.

}
\renewcommand{\baselinestretch}{1.2} 
\smallskip

This is the simplest case where negative coefficients
appear. They seem inevitable due to the following reasoning.
As it was mentioned above, the super-polynomials  can be
found for quite a few relatively simple cases from the
symmetries and special values. We will assume that it
goes through the corresponding HOMFLYPT polynomial and
has the desired evaluations at $q=1$ and $t=1$.
Thus, let us suppose that there
exists $\h_{r,s}(b\,;\,q,t,a)$ with all positive coefficients
such that:
\smallskip

$(i)\  \,$  
$\h_{r,s}(b\,;\,q,q,a)$= 
${H\!D}_{r,s}(b\,;\,q,q,a)$ (so it extends HOMFLYPT${}_{r,s}$),

$(ii)\,$  
$\h_{r,s}(b\,;\,q,1,a)$= 
${H\!D}_{r,s}(b\,;\,q,1,a)$ \ (provided by Conjecture \ref{CONJEVAL}),

$(iii)$ 
$\h_{r,s}(b\,;\,1,t,a)$= 
${H\!D}_{r,s}(b^{tr};\,t^{-1},1,a)$ \, 
(given via Conjecture \ref{DUALIT}).
\smallskip

The evaluations at $q=1$ and $t=1$ are conjectured
to have only non-negative coefficients. In this case,
$b=b^{tr}$, so it suffices to know only
\begin{align*}
&{H\!D}_{3,2}(\om_1+\om_2\,;\,1,t,a)\\
&=\frac{(1+a+t) \bigl(a^2+a (1+t+t^2+t^3)+
t (1+t+t^2+t^4)\bigr)}{t}.
\end{align*}

Let us re-normalize ${H\!D}_{3,2}$ by
making the leading monomial  $a^3$ with the coefficient $1$:
 $$
(t/q^6)\,{H\!D}_{3,2}(b=\om_1+\om_2;\,q,t,a)\ =
$$
$$
a^3+\frac{t}{q^6}+\frac{2 t^2}{q^5}-\frac{t^3}{q^5}
+\frac{2 t^3}{q^4}+\frac{t^3}{q^3}-\frac{t^4}{q^4}
+\frac{2 t^4}{q^3}-\frac{t^5}{q^3}+\frac{2 t^5}{q^2}
+\frac{t^6}{q}
$$
$$
+a \bigl(\frac{1}{q^5}
+\frac{2 t}{q^4}+\frac{t}{q^3}-\frac{t^2}{q^4}
+\frac{3 t^2}{q^3}+\frac{t^2}{q^2}-\frac{t^3}{q^3}
+\frac{3 t^3}{q^2}+\frac{t^3}{q}-\frac{t^4}{q^2}
+\frac{2 t^4}{q}+t^5\bigr)
$$
$$
+a^2 \bigl(\frac{1}{q^3}+\frac{1}{q^2}
+\frac{t}{q^2}+\frac{t}{q}+t^2-\frac{t^2}{q^2}
+\frac{t^2}{q}+t^3\bigr).
$$

It is self-dual with respect to $q\leftrightarrow t^{-1}$,
$a\mapsto a$.
Let us assume that there exists a polynomial
$\h_{3,2}(b\,;\,q,t,a)$ satisfying $(i,ii,iii)$
and with non-negative integer coefficients.
Let us focus on the coefficient $C_2$ of $a^2$ in
$(t/q^6)\h_{3,2}$.

We know the $q$\~exponents of all monomials in $C_2$,
the $t$\~exponents of its monomials, and the $z$\~ monomials of
$C_2(q\mapsto z,t\mapsto z)$. They are:
\begin{align*}
\{\frac{1}{q^3},\, \frac{1}{q^2},&\,\ \frac{1}{q},
\,\frac{1}{q}\,,1\,,1\,\},\ \ 
\{1\,,1\,,t\,,t\,,t^2\,,t^3\},\\
&\{\frac{1}{z^3},\ \,\frac{1}{z^2},\,\ \frac{1}{z}\,\ ,
z\,\ , z^2\,\ ,z^3\,\ \}.
\end{align*}
We conclude  that $C_2$ must contain 
exactly six different monomials because the multiplicities
in the last line are all $1$. Analyzing the ``extreme" terms,
$$
C_2=\frac{1}{q^3\cdot 1}+ 1\cdot t^3 + 
\frac{1}{q^2}\cdot t + \frac{1}{q}\cdot t^2 + \ldots\, ,
$$
due to the lists 
in the first line. However this means that the
terms corresponding to $z$ and  $z^{-1}$ cannot be obtained
using the remaining $q$\~powers and $t$\~powers; indeed,
$1/q^3, t^3,$ $1/q^2, t^2$ and all $1$ have been already taken.

Thus, conditions $(i,ii,iii)$ result in the negativity of at least
one coefficient in $C_2$ and are sufficient to claim that
any such $\h_{3,2}(b\,;\,q,t,a)$ must have some negative 
coefficients.

If we allow changing the DAHA super-polynomials
by rational $q,t$\~factors, the positivity can be 
potentially saved at the expense of switching
to infinite $q,t$\~{\em series\,} instead of the polynomials.
It somehow matches the expectations of the specialists
in the Khovanov-Rozansky homology;
we are thankful Mikhail Khovanov for a discussion.

Evgeny Gorsky observed that our 
$H\!D_{3,2}(b=\om_1+\om_2)$
and the next one for the torus knot $\{4,3\}$ satisfy a 
strikingly simple positivity property. The expansion of 
$H\!D/(1-t)$, which is an infinite  series 
in terms of $t$, has only positive coefficients.
We are thankful for his participation. 
One can divide here by the coefficient of the maximal
power of $a$, which is a $q,t$\~monomial. Then the
result will be  a series in terms of non-negative powers of
$a,q^{-1},t\,$ with non-negative coefficients, infinite
only with respect to the powers of $t$. 
\medskip

{\em The case of \, $b=\om_1+\om_2$, \ knot$\,=\{4,3\}$.}

$$
H\!D_{4,3}(b=\om_1+\om_2;\,q,t,a)\ =
$$
\renewcommand{\baselinestretch}{0.5} 
{\small

\(
1+\frac{a^6 q^{17}}{t^2}+2 q t+2 q^2 t-q t^2+3 q^2 t^2+3 q^3 t^2
+2 q^4 t^2+q^5 t^2-2 q^2 t^3+3 q^3 t^3+5 q^4 t^3+4 q^5 t^3
+2 q^6 t^3-3 q^3 t^4+2 q^4 t^4+5 q^5 t^4+7 q^6 t^4+3 q^7 t^4
+q^8 t^4-4 q^4 t^5+q^5 t^5+5 q^6 t^5+9 q^7 t^5+4 q^8 t^5
+q^9 t^5-4 q^5 t^6+3 q^7 t^6+9 q^8 t^6+5 q^9 t^6
+q^{10} t^6-4 q^6 t^7+4 q^8 t^7+9 q^9 t^7+4 q^{10} t^7
+q^{11} t^7-4 q^7 t^8+3 q^9 t^8+9 q^{10} t^8
+3 q^{11} t^8-4 q^8 t^9+5 q^{10} t^9+7 q^{11} t^9
+2 q^{12} t^9-4 q^9 t^{10}+q^{10} t^{10}+5 q^{11} t^{10}
+4 q^{12} t^{10}+q^{13} t^{10}-4 q^{10} t^{11}+2 q^{11} t^{11}
+5 q^{12} t^{11}+2 q^{13} t^{11}-3 q^{11} t^{12}+3 q^{12} t^{12}
+3 q^{13} t^{12}-2 q^{12} t^{13}+3 q^{13} t^{13}
+2 q^{14} t^{13}-q^{13} t^{14}+2 q^{14} t^{14}+q^{15} t^{15}
\)

\noindent
\(
+a^1 \bigl(3 q^2+4 q^3+2 q^4
+q^5+\frac{q}{t}+\frac{q^2}{t}-q^2 t+4 q^3 t+8 q^4 t+7 q^5 t
+4 q^6 t+q^7 t-3 q^3 t^2+3 q^4 t^2+11 q^5 t^2+14 q^6 t^2+9 q^7 t^2
+3 q^8 t^2+q^9 t^2-5 q^4 t^3+q^5 t^3+11 q^6 t^3+20 q^7 t^3
+13 q^8 t^3+5 q^9 t^3+q^{10} t^3-6 q^5 t^4-2 q^6 t^4+9 q^7 t^4
+24 q^8 t^4+18 q^9 t^4+7 q^{10} t^4
+2 q^{11} t^4-7 q^6 t^5-4 q^7 t^5+7 q^8 t^5+25 q^9 t^5
+19 q^{10} t^5+7 q^{11} t^5+q^{12} t^5-7 q^7 t^6-5 q^8 t^6
+5 q^9 t^6+25 q^{10} t^6+18 q^{11} t^6+5 q^{12} t^6
+q^{13} t^6-7 q^8 t^7-5 q^9 t^7+7 q^{10} t^7+24 q^{11} t^7
+13 q^{12} t^7+3 q^{13} t^7-7 q^9 t^8-4 q^{10} t^8
+9 q^{11} t^8+20 q^{12} t^8+9 q^{13} t^8
+q^{14} t^8-7 q^{10} t^9-2 q^{11} t^9+11 q^{12} t^9
+14 q^{13} t^9+4 q^{14} t^9-6 q^{11} t^{10}+q^{12} t^{10}
+11 q^{13} t^{10}+7 q^{14} t^{10}+q^{15} t^{10}-5 q^{12} t^{11}
+3 q^{13} t^{11}+8 q^{14} t^{11}+2 q^{15} t^{11}-3 q^{13} t^{12}
+4 q^{14} t^{12}+4 q^{15} t^{12}-q^{14} t^{13}+3 q^{15} t^{13}
+q^{16} t^{13}+q^{16} t^{14}\bigr).
\)

\noindent
\(
+a^2 \bigl(q^4+7 q^5+8 q^6+7 q^7+3 q^8+q^9+\frac{q^3}{t^2}
+\frac{q^3}{t}+\frac{3 q^4}{t}+\frac{3 q^5}{t}+\frac{2 q^6}{t}
+\frac{q^7}{t}-q^4 t+8 q^6 t+15 q^7 t+14 q^8 t+8 q^9 t
+3 q^{10} t+q^{11} t-2 q^5 t^2-3 q^6 t^2+8 q^7 t^2+20 q^8 t^2
+22 q^9 t^2+12 q^{10} t^2+5 q^{11} t^2
+q^{12} t^2-3 q^6 t^3-6 q^7 t^3+4 q^8 t^3+22 q^9 t^3+27 q^{10} t^3
+15 q^{11} t^3+5 q^{12} t^3+q^{13} t^3-3 q^7 t^4-8 q^8 t^4
+2 q^9 t^4+22 q^{10} t^4+30 q^{11} t^4+15 q^{12} t^4+5 q^{13} t^4
+q^{14} t^4-3 q^8 t^5-9 q^9 t^5+22 q^{11} t^5+27 q^{12} t^5
+12 q^{13} t^5+3 q^{14} t^5-3 q^9 t^6-9 q^{10} t^6+2 q^{11} t^6
+22 q^{12} t^6+22 q^{13} t^6+8 q^{14} t^6
+q^{15} t^6-3 q^{10} t^7-8 q^{11} t^7+4 q^{12} t^7+20 q^{13} t^7
+14 q^{14} t^7+3 q^{15} t^7-3 q^{11} t^8-6 q^{12} t^8
+8 q^{13} t^8+15 q^{14} t^8+7 q^{15} t^8
+q^{16} t^8-3 q^{12} t^9-3 q^{13} t^9+8 q^{14} t^9+8 q^{15} t^9
+2 q^{16} t^9-2 q^{13} t^{10}+7 q^{15} t^{10}
+3 q^{16} t^{10}-q^{14} t^{11}+q^{15} t^{11}+3 q^{16} t^{11}
+q^{16} t^{12}+q^{17} t^{12}\bigr)
\)

\noindent
\(
+a^3 \bigl(-q^6+2 q^7+6 q^8+10 q^9+8 q^{10}+4 q^{11}+2 q^{12}
+\frac{q^5}{t^2}+\frac{q^6}{t^2}+\frac{q^7}{t^2}+\frac{q^8}{t^2}
+\frac{2 q^6}{t}+\frac{4 q^7}{t}+\frac{5 q^8}{t}+\frac{4 q^9}{t}
+\frac{2 q^{10}}{t}+\frac{q^{11}}{t}-2 q^7 t+q^8 t+7 q^9 t
+15 q^{10} t+13 q^{11} t+7 q^{12} t+3 q^{13} t
+q^{14} t-3 q^8 t^2-q^9 t^2+6 q^{10} t^2+18 q^{11} t^2
+15 q^{12} t^2+7 q^{13} t^2+3 q^{14} t^2-4 q^9 t^3-3 q^{10} t^3
+4 q^{11} t^3+20 q^{12} t^3+15 q^{13} t^3+7 q^{14} t^3
+2 q^{15} t^3-4 q^{10} t^4-3 q^{11} t^4+4 q^{12} t^4
+18 q^{13} t^4+13 q^{14} t^4+4 q^{15} t^4
+q^{16} t^4-4 q^{11} t^5-3 q^{12} t^5+6 q^{13} t^5+15 q^{14} t^5
+8 q^{15} t^5+2 q^{16} t^5-4 q^{12} t^6-q^{13} t^6+7 q^{14} t^6
+10 q^{15} t^6+4 q^{16} t^6-3 q^{13} t^7+q^{14} t^7+6 q^{15} t^7
+5 q^{16} t^7+q^{17} t^7-2 q^{14} t^8+2 q^{15} t^8+4 q^{16} t^8
+q^{17} t^8-q^{15} t^9+2 q^{16} t^9+q^{17} t^9+q^{17} t^{10}\bigr)
\)

\noindent
\(
+a^4 \bigl(-q^9+3 q^{11}
+6 q^{12}+4 q^{13}+3 q^{14}+q^{15}+\frac{q^8}{t^2}+\frac{q^9}{t^2}
+\frac{2 q^{10}}{t^2}+\frac{q^{11}}{t^2}+\frac{q^{12}}{t^2}
+\frac{q^9}{t}+\frac{3 q^{10}}{t}+\frac{4 q^{11}}{t}
+\frac{3 q^{12}}{t}+\frac{2 q^{13}}{t}
+\frac{q^{14}}{t}-q^{10} t-q^{11} t+4 q^{12} t+7 q^{13} t
+5 q^{14} t+3 q^{15} t+q^{16} t-q^{11} t^2-2 q^{12} t^2+3 q^{13} t^2
+7 q^{14} t^2+4 q^{15} t^2+2 q^{16} t^2-q^{12} t^3-2 q^{13} t^3
+4 q^{14} t^3
+6 q^{15} t^3+3 q^{16} t^3+q^{17} t^3-q^{13} t^4-q^{14} t^4
+3 q^{15} t^4+4 q^{16} t^4+q^{17} t^4-q^{14} t^5+3 q^{16} t^5
+2 q^{17} t^5-q^{15} t^6+q^{16} t^6+q^{17} t^6+q^{17} t^7\bigr)
\)

\noindent
\(
+a^5 \bigl(-q^{13}+q^{14}+q^{15}+q^{16}+q^{17}+\frac{q^{12}}{t^2}
+\frac{q^{13}}{t^2}+\frac{q^{14}}{t^2}+\frac{q^{15}}{t^2}
+\frac{q^{13}}{t}+\frac{q^{14}}{t}+\frac{q^{15}}{t}
+\frac{q^{16}}{t}-q^{14} t+q^{15} t+q^{16} t+q^{17} t-q^{15} t^2
+q^{16} t^2+q^{17} t^2+q^{17} t^3\bigr).
\)

}
\renewcommand{\baselinestretch}{1.2} 
\smallskip

\subsection{Further examples}
\medskip
Let us consider the knot $\{9,4\}$ and $b=\om_2$.
The corresponding super-polynomial becomes
exactly that at the end of Section
\ref{sect:rank-one} upon the following transformations.
First, one needs to apply the duality substitution
$q\mapsto t^{-1}$, $t\mapsto q^{-1}$
to $H\!D_{9,4}$. Second, the result must be
multiplied by $q^{48}t^{24}$. Finally, we
substitute $\,a\mapsto -t^2$ (notice the minus).
The super-formula is as follows:

$$
H\!D_{9,4}(b=\om_2;\,q,t,a)\ =
$$
\renewcommand{\baselinestretch}{0.5} 
{\small

\noindent
\(
a^0 \bigl(1+q t+q^2 t+q^3 t+q t^2+2 q^2 t^2+2 q^3 t^2+2 q^4 t^2
+q^5 t^2+q^6 t^2+q^2 t^3+3 q^3 t^3+4 q^4 t^3+4 q^5 t^3+3 q^6 t^3
+q^7 t^3+q^2 t^4+2 q^3 t^4+5 q^4 t^4+6 q^5 t^4+7 q^6 t^4
+5 q^7 t^4+3 q^8 t^4+q^3 t^5+3 q^4 t^5+7 q^5 t^5+9 q^6 t^5
+11 q^7 t^5+8 q^8 t^5+3 q^9 t^5+q^3 t^6+2 q^4 t^6+5 q^5 t^6
+10 q^6 t^6+13 q^7 t^6+15 q^8 t^6+11 q^9 t^6+4 q^{10} t^6
+q^4 t^7+3 q^5 t^7+7 q^6 t^7+13 q^7 t^7+18 q^8 t^7
+20 q^9 t^7+12 q^{10} t^7+3 q^{11} t^7+q^4 t^8+2 q^5 t^8
+5 q^6 t^8+10 q^7 t^8+18 q^8 t^8+23 q^9 t^8+24 q^{10} t^8
+11 q^{11} t^8+3 q^{12} t^8+q^5 t^9+3 q^6 t^9+7 q^7 t^9
+13 q^8 t^9+23 q^9 t^9+29 q^{10} t^9+25 q^{11} t^9
+9 q^{12} t^9+q^{13} t^9+q^5 t^{10}+2 q^6 t^{10}+5 q^7 t^{10}
+10 q^8 t^{10}+18 q^9 t^{10}+29 q^{10} t^{10}+34 q^{11} t^{10}
+23 q^{12} t^{10}+6 q^{13} t^{10}+q^{14} t^{10}+q^6 t^{11}
+3 q^7 t^{11}+7 q^8 t^{11}+13 q^9 t^{11}+23 q^{10} t^{11}
+35 q^{11} t^{11}+34 q^{12} t^{11}+18 q^{13} t^{11}
+3 q^{14} t^{11}+q^6 t^{12}+2 q^7 t^{12}+5 q^8 t^{12}
+10 q^9 t^{12}+18 q^{10} t^{12}+29 q^{11} t^{12}+41 q^{12} t^{12}
+31 q^{13} t^{12}+11 q^{14} t^{12}+q^{15} t^{12}+q^7 t^{13}
+3 q^8 t^{13}+7 q^9 t^{13}+13 q^{10} t^{13}+23 q^{11} t^{13}
+35 q^{12} t^{13}+41 q^{13} t^{13}+23 q^{14} t^{13}
+5 q^{15} t^{13}+q^7 t^{14}+2 q^8 t^{14}+5 q^9 t^{14}
+10 q^{10} t^{14}+18 q^{11} t^{14}+29 q^{12} t^{14}
+40 q^{13} t^{14}+37 q^{14} t^{14}+13 q^{15} t^{14}
+2 q^{16} t^{14}+q^8 t^{15}+3 q^9 t^{15}+7 q^{10} t^{15}
+13 q^{11} t^{15}+23 q^{12} t^{15}+35 q^{13} t^{15}
+40 q^{14} t^{15}+26 q^{15} t^{15}+5 q^{16} t^{15}+q^8 t^{16}
+2 q^9 t^{16}+5 q^{10} t^{16}+10 q^{11} t^{16}
+18 q^{12} t^{16}+29 q^{13} t^{16}+40 q^{14} t^{16}
+34 q^{15} t^{16}+14 q^{16} t^{16}+q^{17} t^{16}+q^9 t^{17}
+3 q^{10} t^{17}+7 q^{11} t^{17}+13 q^{12} t^{17}
+23 q^{13} t^{17}+34 q^{14} t^{17}+39 q^{15} t^{17}
+22 q^{16} t^{17}+4 q^{17} t^{17}+q^9 t^{18}+2 q^{10} t^{18}
+5 q^{11} t^{18}+10 q^{12} t^{18}+18 q^{13} t^{18}
+29 q^{14} t^{18}+38 q^{15} t^{18}+32 q^{16} t^{18}
+9 q^{17} t^{18}+q^{18} t^{18}+q^{10} t^{19}+3 q^{11} t^{19}
+7 q^{12} t^{19}+13 q^{13} t^{19}+23 q^{14} t^{19}
+33 q^{15} t^{19}+34 q^{16} t^{19}+18 q^{17} t^{19}
+2 q^{18} t^{19}+q^{10} t^{20}+2 q^{11} t^{20}+5 q^{12} t^{20}
+10 q^{13} t^{20}+18 q^{14} t^{20}+28 q^{15} t^{20}
+35 q^{16} t^{20}+24 q^{17} t^{20}+6 q^{18} t^{20}
+q^{11} t^{21}+3 q^{12} t^{21}+7 q^{13} t^{21}+13 q^{14} t^{21}
+23 q^{15} t^{21}+31 q^{16} t^{21}+29 q^{17} t^{21}
+10 q^{18} t^{21}+q^{19} t^{21}+q^{11} t^{22}+2 q^{12} t^{22}
+5 q^{13} t^{22}+10 q^{14} t^{22}+18 q^{15} t^{22}+27 q^{16} t^{22}
+30 q^{17} t^{22}+17 q^{18} t^{22}+2 q^{19} t^{22}+q^{12} t^{23}
+3 q^{13} t^{23}+7 q^{14} t^{23}+13 q^{15} t^{23}+22 q^{16} t^{23}
+28 q^{17} t^{23}+20 q^{18} t^{23}+5 q^{19} t^{23}+q^{12} t^{24}
+2 q^{13} t^{24}+5 q^{14} t^{24}+10 q^{15} t^{24}
+18 q^{16} t^{24}+25 q^{17} t^{24}+24 q^{18} t^{24}
+8 q^{19} t^{24}+q^{20} t^{24}+q^{13} t^{25}+3 q^{14} t^{25}
+7 q^{15} t^{25}+13 q^{16} t^{25}+21 q^{17} t^{25}
+23 q^{18} t^{25}+12 q^{19} t^{25}+q^{20} t^{25}+q^{13} t^{26}
+2 q^{14} t^{26}+5 q^{15} t^{26}+10 q^{16} t^{26}
+17 q^{17} t^{26}+22 q^{18} t^{26}+15 q^{19} t^{26}
+3 q^{20} t^{26}+q^{14} t^{27}+3 q^{15} t^{27}+7 q^{16} t^{27}
+13 q^{17} t^{27}+19 q^{18} t^{27}+17 q^{19} t^{27}
+4 q^{20} t^{27}+q^{14} t^{28}+2 q^{15} t^{28}+5 q^{16} t^{28}
+10 q^{17} t^{28}+16 q^{18} t^{28}+17 q^{19} t^{28}
+8 q^{20} t^{28}+q^{15} t^{29}+3 q^{16} t^{29}+7 q^{17} t^{29}
+12 q^{18} t^{29}+16 q^{19} t^{29}+9 q^{20} t^{29}+q^{21} t^{29}
+q^{15} t^{30}+2 q^{16} t^{30}+5 q^{17} t^{30}+10 q^{18} t^{30}
+14 q^{19} t^{30}+11 q^{20} t^{30}+2 q^{21} t^{30}+q^{16} t^{31}
+3 q^{17} t^{31}+7 q^{18} t^{31}+11 q^{19} t^{31}+11 q^{20} t^{31}
+3 q^{21} t^{31}+q^{16} t^{32}+2 q^{17} t^{32}+5 q^{18} t^{32}
+9 q^{19} t^{32}+11 q^{20} t^{32}+4 q^{21} t^{32}+q^{17} t^{33}
+3 q^{18} t^{33}+7 q^{19} t^{33}+9 q^{20} t^{33}+6 q^{21} t^{33}
+q^{17} t^{34}+2 q^{18} t^{34}+5 q^{19} t^{34}+8 q^{20} t^{34}
+6 q^{21} t^{34}+q^{22} t^{34}+q^{18} t^{35}+3 q^{19} t^{35}
+6 q^{20} t^{35}+6 q^{21} t^{35}+q^{22} t^{35}+q^{18} t^{36}
+2 q^{19} t^{36}+5 q^{20} t^{36}+6 q^{21} t^{36}+2 q^{22} t^{36}
+q^{19} t^{37}+3 q^{20} t^{37}+5 q^{21} t^{37}+2 q^{22} t^{37}
+q^{19} t^{38}+2 q^{20} t^{38}+4 q^{21} t^{38}+3 q^{22} t^{38}
+q^{20} t^{39}+3 q^{21} t^{39}+3 q^{22} t^{39}+q^{20} t^{40}
+2 q^{21} t^{40}+3 q^{22} t^{40}+q^{21} t^{41}+2 q^{22} t^{41}
+q^{23} t^{41}+q^{21} t^{42}+2 q^{22} t^{42}+q^{23} t^{42}
+q^{22} t^{43}+q^{23} t^{43}+q^{22} t^{44}+q^{23} t^{44}
+q^{23} t^{45}+q^{23} t^{46}+q^{24} t^{48}\bigr)
\)
\vfill

\noindent
\(
+ a^1 \bigl(q+2 q^2+3 q^3+3 q^4+2 q^5+q^6+\frac{q}{t}
+\frac{q^2}{t}+\frac{q^3}{t}+2 q^2 t+5 q^3 t+8 q^4 t
+8 q^5 t+6 q^6 t+3 q^7 t+q^8 t+q^2 t^2+4 q^3 t^2+9 q^4 t^2
+14 q^5 t^2+16 q^6 t^2+14 q^7 t^2+8 q^8 t^2+2 q^9 t^2
+2 q^3 t^3+7 q^4 t^3+16 q^5 t^3+24 q^6 t^3+29 q^7 t^3
+24 q^8 t^3+13 q^9 t^3+3 q^{10} t^3+q^3 t^4+4 q^4 t^4
+12 q^5 t^4+24 q^6 t^4+37 q^7 t^4+44 q^8 t^4+38 q^9 t^4
+19 q^{10} t^4+4 q^{11} t^4+2 q^4 t^5+7 q^5 t^5
+19 q^6 t^5+36 q^7 t^5+55 q^8 t^5+64 q^9 t^5+50 q^{10} t^5
+21 q^{11} t^5+4 q^{12} t^5+q^4 t^6+4 q^5 t^6+12 q^6 t^6
+28 q^7 t^6+51 q^8 t^6+75 q^9 t^6+84 q^{10} t^6
+58 q^{11} t^6+21 q^{12} t^6+3 q^{13} t^6+2 q^5 t^7
+7 q^6 t^7+19 q^7 t^7+40 q^8 t^7+71 q^9 t^7+100 q^{10} t^7
+101 q^{11} t^7+58 q^{12} t^7+17 q^{13} t^7+2 q^{14} t^7
+q^5 t^8+4 q^6 t^8+12 q^7 t^8+28 q^8 t^8+56 q^9 t^8
+93 q^{10} t^8+124 q^{11} t^8+107 q^{12} t^8+51 q^{13} t^8
+11 q^{14} t^8+q^{15} t^8+2 q^6 t^9+7 q^7 t^9+19 q^8 t^9
+40 q^9 t^9+76 q^{10} t^9+120 q^{11} t^9+142 q^{12} t^9
+102 q^{13} t^9+36 q^{14} t^9+5 q^{15} t^9+q^6 t^{10}
+4 q^7 t^{10}+12 q^8 t^{10}+28 q^9 t^{10}+56 q^{10} t^{10}
+99 q^{11} t^{10}+146 q^{12} t^{10}+147 q^{13} t^{10}
+82 q^{14} t^{10}+21 q^{15} t^{10}+2 q^{16} t^{10}
+2 q^7 t^{11}+7 q^8 t^{11}+19 q^9 t^{11}+40 q^{10} t^{11}
+76 q^{11} t^{11}+126 q^{12} t^{11}+164 q^{13} t^{11}
+134 q^{14} t^{11}+54 q^{15} t^{11}+9 q^{16} t^{11}
+q^7 t^{12}+4 q^8 t^{12}+12 q^9 t^{12}+28 q^{10} t^{12}
+56 q^{11} t^{12}+99 q^{12} t^{12}+152 q^{13} t^{12}
+167 q^{14} t^{12}+101 q^{15} t^{12}+27 q^{16} t^{12}
+2 q^{17} t^{12}+2 q^8 t^{13}+7 q^9 t^{13}+19 q^{10} t^{13}
+40 q^{11} t^{13}+76 q^{12} t^{13}+126 q^{13} t^{13}
+170 q^{14} t^{13}+146 q^{15} t^{13}+61 q^{16} t^{13}
+9 q^{17} t^{13}+q^8 t^{14}+4 q^9 t^{14}+12 q^{10} t^{14}
+28 q^{11} t^{14}+56 q^{12} t^{14}+99 q^{13} t^{14}
+151 q^{14} t^{14}+168 q^{15} t^{14}+103 q^{16} t^{14}
+26 q^{17} t^{14}+2 q^{18} t^{14}+2 q^9 t^{15}+7 q^{10} t^{15}
+19 q^{11} t^{15}+40 q^{12} t^{15}+76 q^{13} t^{15}
+125 q^{14} t^{15}+166 q^{15} t^{15}+139 q^{16} t^{15}
+53 q^{17} t^{15}+7 q^{18} t^{15}+q^9 t^{16}+4 q^{10} t^{16}
+12 q^{11} t^{16}+28 q^{12} t^{16}+56 q^{13} t^{16}
+99 q^{14} t^{16}+148 q^{15} t^{16}+157 q^{16} t^{16}
+87 q^{17} t^{16}+18 q^{18} t^{16}+q^{19} t^{16}+2 q^{10} t^{17}
+7 q^{11} t^{17}+19 q^{12} t^{17}+40 q^{13} t^{17}
+76 q^{14} t^{17}+123 q^{15} t^{17}+156 q^{16} t^{17}
+119 q^{17} t^{17}+37 q^{18} t^{17}+3 q^{19} t^{17}
+q^{10} t^{18}+4 q^{11} t^{18}+12 q^{12} t^{18}
+28 q^{13} t^{18}+56 q^{14} t^{18}+98 q^{15} t^{18}
+141 q^{16} t^{18}+137 q^{17} t^{18}+62 q^{18} t^{18}
+9 q^{19} t^{18}+2 q^{11} t^{19}+7 q^{12} t^{19}
+19 q^{13} t^{19}+40 q^{14} t^{19}+76 q^{15} t^{19}
+119 q^{16} t^{19}+139 q^{17} t^{19}+88 q^{18} t^{19}
+20 q^{19} t^{19}+q^{20} t^{19}+q^{11} t^{20}+4 q^{12} t^{20}
+12 q^{13} t^{20}+28 q^{14} t^{20}+56 q^{15} t^{20}
+96 q^{16} t^{20}+129 q^{17} t^{20}+106 q^{18} t^{20}
+35 q^{19} t^{20}+3 q^{20} t^{20}+2 q^{12} t^{21}
+7 q^{13} t^{21}+19 q^{14} t^{21}+40 q^{15} t^{21}
+75 q^{16} t^{21}+112 q^{17} t^{21}+114 q^{18} t^{21}
+53 q^{19} t^{21}+7 q^{20} t^{21}+q^{12} t^{22}+4 q^{13} t^{22}
+12 q^{14} t^{22}+28 q^{15} t^{22}+56 q^{16} t^{22}
+92 q^{17} t^{22}+111 q^{18} t^{22}+70 q^{19} t^{22}
+14 q^{20} t^{22}+2 q^{13} t^{23}+7 q^{14} t^{23}
+19 q^{15} t^{23}+40 q^{16} t^{23}+73 q^{17} t^{23}
+100 q^{18} t^{23}+81 q^{19} t^{23}+24 q^{20} t^{23}
+q^{21} t^{23}+q^{13} t^{24}+4 q^{14} t^{24}+12 q^{15} t^{24}
+28 q^{16} t^{24}+55 q^{17} t^{24}+85 q^{18} t^{24}
+85 q^{19} t^{24}+35 q^{20} t^{24}+3 q^{21} t^{24}
+2 q^{14} t^{25}+7 q^{15} t^{25}+19 q^{16} t^{25}
+40 q^{17} t^{25}+69 q^{18} t^{25}+82 q^{19} t^{25}
+46 q^{20} t^{25}+6 q^{21} t^{25}+q^{14} t^{26}
+4 q^{15} t^{26}+12 q^{16} t^{26}+28 q^{17} t^{26}
+53 q^{18} t^{26}+73 q^{19} t^{26}+54 q^{20} t^{26}
+11 q^{21} t^{26}+2 q^{15} t^{27}+7 q^{16} t^{27}
+19 q^{17} t^{27}+39 q^{18} t^{27}+62 q^{19} t^{27}
+57 q^{20} t^{27}+17 q^{21} t^{27}+q^{15} t^{28}
+4 q^{16} t^{28}+12 q^{17} t^{28}+28 q^{18} t^{28}
+49 q^{19} t^{28}+55 q^{20} t^{28}+24 q^{21} t^{28}
+q^{22} t^{28}+2 q^{16} t^{29}+7 q^{17} t^{29}
+19 q^{18} t^{29}+37 q^{19} t^{29}+50 q^{20} t^{29}
+29 q^{21} t^{29}+3 q^{22} t^{29}+q^{16} t^{30}
+4 q^{17} t^{30}+12 q^{18} t^{30}+27 q^{19} t^{30}
+42 q^{20} t^{30}+32 q^{21} t^{30}+5 q^{22} t^{30}
+2 q^{17} t^{31}+7 q^{18} t^{31}+19 q^{19} t^{31}
+33 q^{20} t^{31}+33 q^{21} t^{31}+8 q^{22} t^{31}
+q^{17} t^{32}+4 q^{18} t^{32}+12 q^{19} t^{32}
+25 q^{20} t^{32}+30 q^{21} t^{32}+11 q^{22} t^{32}
+2 q^{18} t^{33}+7 q^{19} t^{33}+18 q^{20} t^{33}
+26 q^{21} t^{33}+14 q^{22} t^{33}+q^{18} t^{34}
+4 q^{19} t^{34}+12 q^{20} t^{34}+21 q^{21} t^{34}
+15 q^{22} t^{34}+q^{23} t^{34}+2 q^{19} t^{35}
+7 q^{20} t^{35}+16 q^{21} t^{35}+15 q^{22} t^{35}
+2 q^{23} t^{35}+q^{19} t^{36}+4 q^{20} t^{36}
+11 q^{21} t^{36}+14 q^{22} t^{36}+3 q^{23} t^{36}
+2 q^{20} t^{37}+7 q^{21} t^{37}+12 q^{22} t^{37}
+4 q^{23} t^{37}+q^{20} t^{38}+4 q^{21} t^{38}
+9 q^{22} t^{38}+5 q^{23} t^{38}+2 q^{21} t^{39}
+6 q^{22} t^{39}+6 q^{23} t^{39}+q^{21} t^{40}
+4 q^{22} t^{40}+5 q^{23} t^{40}+2 q^{22} t^{41}
+4 q^{23} t^{41}+q^{24} t^{41}+q^{22} t^{42}
+3 q^{23} t^{42}+q^{24} t^{42}+2 q^{23} t^{43}
+q^{24} t^{43}+q^{23} t^{44}+q^{24} t^{44}
+q^{24} t^{45}+q^{24} t^{46}\bigr)
\)
\vfill

\noindent
\(
a^2 \bigl(2 q^3+5 q^4+10 q^5+12 q^6+12 q^7+7 q^8+3 q^9
+\frac{q^3}{t^2}+\frac{q^4}{t^2}+\frac{q^5}{t^2}+\frac{q^2}{t}
+\frac{2 q^3}{t}+\frac{4 q^4}{t}+\frac{4 q^5}{t}
+\frac{4 q^6}{t}+\frac{2 q^7}{t}+\frac{q^8}{t}+q^3 t
+5 q^4 t+12 q^5 t+21 q^6 t+26 q^7 t+25 q^8 t+17 q^9 t
+7 q^{10} t+q^{11} t+2 q^4 t^2+9 q^5 t^2+21 q^6 t^2+37 q^7 t^2
+47 q^8 t^2+47 q^9 t^2+29 q^{10} t^2+11 q^{11} t^2
+q^{12} t^2+q^4 t^3+5 q^5 t^3+17 q^6 t^3+36 q^7 t^3
+61 q^8 t^3+76 q^9 t^3+72 q^{10} t^3+42 q^{11} t^3
+14 q^{12} t^3+2 q^{13} t^3+2 q^5 t^4+9 q^6 t^4+27 q^7 t^4
+55 q^8 t^4+91 q^9 t^4+112 q^{10} t^4+99 q^{11} t^4
+50 q^{12} t^4+15 q^{13} t^4+q^{14} t^4+q^5 t^5+5 q^6 t^5
+17 q^7 t^5+43 q^8 t^5+83 q^9 t^5+130 q^{10} t^5+151 q^{11} t^5
+117 q^{12} t^5+52 q^{13} t^5+12 q^{14} t^5+q^{15} t^5+2 q^6 t^6
+9 q^7 t^6+27 q^8 t^6+63 q^9 t^6+116 q^{10} t^6+174 q^{11} t^6
+183 q^{12} t^6+123 q^{13} t^6+44 q^{14} t^6+8 q^{15} t^6+q^6 t^7
+5 q^7 t^7+17 q^8 t^7+43 q^9 t^7+92 q^{10} t^7+159 q^{11} t^7
+220 q^{12} t^7+201 q^{13} t^7+111 q^{14} t^7+31 q^{15} t^7
+4 q^{16} t^7+2 q^7 t^8
+9 q^8 t^8+27 q^9 t^8+63 q^{10} t^8+126 q^{11} t^8+206 q^{12} t^8
+254 q^{13} t^8+194 q^{14} t^8+84 q^{15} t^8+16 q^{16} t^8
+q^{17} t^8+q^7 t^9+5 q^8 t^9+17 q^9 t^9+43 q^{10} t^9
+92 q^{11} t^9+170 q^{12} t^9+254 q^{13} t^9+268 q^{14} t^9
+161 q^{15} t^9+51 q^{16} t^9+6 q^{17} t^9+2 q^8 t^{10}
+9 q^9 t^{10}+27 q^{10} t^{10}+63 q^{11} t^{10}+126 q^{12} t^{10}
+218 q^{13} t^{10}+288 q^{14} t^{10}+245 q^{15} t^{10}
+109 q^{16} t^{10}+23 q^{17} t^{10}+q^{18} t^{10}+q^8 t^{11}
+5 q^9 t^{11}+17 q^{10} t^{11}+43 q^{11} t^{11}
+92 q^{12} t^{11}+170 q^{13} t^{11}+267 q^{14} t^{11}
+294 q^{15} t^{11}+189 q^{16} t^{11}+57 q^{17} t^{11}
+7 q^{18} t^{11}+2 q^9 t^{12}+9 q^{10} t^{12}+27 q^{11} t^{12}
+63 q^{12} t^{12}+126 q^{13} t^{12}+217 q^{14} t^{12}
+298 q^{15} t^{12}+255 q^{16} t^{12}+114 q^{17} t^{12}
+20 q^{18} t^{12}+q^{19} t^{12}+q^9 t^{13}+5 q^{10} t^{13}
+17 q^{11} t^{13}+43 q^{12} t^{13}+92 q^{13} t^{13}
+170 q^{14} t^{13}+264 q^{15} t^{13}+294 q^{16} t^{13}
+179 q^{17} t^{13}+50 q^{18} t^{13}+4 q^{19} t^{13}
+2 q^{10} t^{14}+9 q^{11} t^{14}+27 q^{12} t^{14}
+63 q^{13} t^{14}+126 q^{14} t^{14}+216 q^{15} t^{14}
+286 q^{16} t^{14}+238 q^{17} t^{14}+91 q^{18} t^{14}
+14 q^{19} t^{14}+q^{10} t^{15}+5 q^{11} t^{15}
+17 q^{12} t^{15}+43 q^{13} t^{15}+92 q^{14} t^{15}
+169 q^{15} t^{15}+257 q^{16} t^{15}+267 q^{17} t^{15}
+145 q^{18} t^{15}+31 q^{19} t^{15}+2 q^{20} t^{15}
+2 q^{11} t^{16}+9 q^{12} t^{16}+27 q^{13} t^{16}
+63 q^{14} t^{16}+126 q^{15} t^{16}+211 q^{16} t^{16}
+266 q^{17} t^{16}+190 q^{18} t^{16}+60 q^{19} t^{16}
+5 q^{20} t^{16}+q^{11} t^{17}+5 q^{12} t^{17}
+17 q^{13} t^{17}+43 q^{14} t^{17}+92 q^{15} t^{17}
+167 q^{16} t^{17}+240 q^{17} t^{17}+224 q^{18} t^{17}
+95 q^{19} t^{17}+14 q^{20} t^{17}+2 q^{12} t^{18}
+9 q^{13} t^{18}+27 q^{14} t^{18}+63 q^{15} t^{18}
+125 q^{16} t^{18}+202 q^{17} t^{18}+227 q^{18} t^{18}
+133 q^{19} t^{18}+27 q^{20} t^{18}+q^{21} t^{18}
+q^{12} t^{19}+5 q^{13} t^{19}+17 q^{14} t^{19}
+43 q^{15} t^{19}+92 q^{16} t^{19}+162 q^{17} t^{19}
+215 q^{18} t^{19}+161 q^{19} t^{19}+49 q^{20} t^{19}
+3 q^{21} t^{19}+2 q^{13} t^{20}+9 q^{14} t^{20}
+27 q^{15} t^{20}+63 q^{16} t^{20}+123 q^{17} t^{20}
+185 q^{18} t^{20}+176 q^{19} t^{20}+71 q^{20} t^{20}
+8 q^{21} t^{20}+q^{13} t^{21}+5 q^{14} t^{21}
+17 q^{15} t^{21}+43 q^{16} t^{21}+91 q^{17} t^{21}
+153 q^{18} t^{21}+174 q^{19} t^{21}+97 q^{20} t^{21}
+15 q^{21} t^{21}+2 q^{14} t^{22}+9 q^{15} t^{22}
+27 q^{16} t^{22}+63 q^{17} t^{22}+118 q^{18} t^{22}
+159 q^{19} t^{22}+112 q^{20} t^{22}+27 q^{21} t^{22}
+q^{14} t^{23}+5 q^{15} t^{23}+17 q^{16} t^{23}
+43 q^{17} t^{23}+89 q^{18} t^{23}+136 q^{19} t^{23}
+123 q^{20} t^{23}+40 q^{21} t^{23}+2 q^{22} t^{23}
+2 q^{15} t^{24}+9 q^{16} t^{24}+27 q^{17} t^{24}
+62 q^{18} t^{24}+109 q^{19} t^{24}+119 q^{20} t^{24}
+54 q^{21} t^{24}+4 q^{22} t^{24}+q^{15} t^{25}
+5 q^{16} t^{25}+17 q^{17} t^{25}+43 q^{18} t^{25}
+84 q^{19} t^{25}+110 q^{20} t^{25}+65 q^{21} t^{25}
+9 q^{22} t^{25}+2 q^{16} t^{26}+9 q^{17} t^{26}
+27 q^{18} t^{26}+60 q^{19} t^{26}+92 q^{20} t^{26}
+72 q^{21} t^{26}+14 q^{22} t^{26}+q^{16} t^{27}
+5 q^{17} t^{27}+17 q^{18} t^{27}+42 q^{19} t^{27}
+75 q^{20} t^{27}+72 q^{21} t^{27}+22 q^{22} t^{27}
+2 q^{17} t^{28}+9 q^{18} t^{28}+27 q^{19} t^{28}
+55 q^{20} t^{28}+67 q^{21} t^{28}+27 q^{22} t^{28}
+q^{23} t^{28}+q^{17} t^{29}+5 q^{18} t^{29}+17 q^{19} t^{29}
+40 q^{20} t^{29}+58 q^{21} t^{29}+33 q^{22} t^{29}
+2 q^{23} t^{29}+2 q^{18} t^{30}+9 q^{19} t^{30}
+26 q^{20} t^{30}+46 q^{21} t^{30}+34 q^{22} t^{30}
+4 q^{23} t^{30}+q^{18} t^{31}+5 q^{19} t^{31}
+17 q^{20} t^{31}+35 q^{21} t^{31}+35 q^{22} t^{31}
+6 q^{23} t^{31}+2 q^{19} t^{32}+9 q^{20} t^{32}
+24 q^{21} t^{32}+30 q^{22} t^{32}+9 q^{23} t^{32}
+q^{19} t^{33}+5 q^{20} t^{33}+16 q^{21} t^{33}
+26 q^{22} t^{33}+11 q^{23} t^{33}+2 q^{20} t^{34}
+9 q^{21} t^{34}+19 q^{22} t^{34}+12 q^{23} t^{34}
+q^{20} t^{35}+5 q^{21} t^{35}+14 q^{22} t^{35}
+12 q^{23} t^{35}+q^{24} t^{35}+2 q^{21} t^{36}
+8 q^{22} t^{36}+11 q^{23} t^{36}+q^{24} t^{36}
+q^{21} t^{37}+5 q^{22} t^{37}+9 q^{23} t^{37}
+2 q^{24} t^{37}+2 q^{22} t^{38}+6 q^{23} t^{38}
+2 q^{24} t^{38}+q^{22} t^{39}+4 q^{23} t^{39}
+3 q^{24} t^{39}+2 q^{23} t^{40}+2 q^{24} t^{40}
+q^{23} t^{41}+2 q^{24} t^{41}+q^{24} t^{42}+q^{24} t^{43}\bigr)
\)
\vfill

\noindent
\(
+a^3 \bigl(2 q^5+8 q^6+15 q^7+22 q^8+24 q^9+18 q^{10}
+9 q^{11}+2 q^{12}+\frac{q^6}{t^3}+\frac{q^4}{t^2}
+\frac{2 q^5}{t^2}+\frac{3 q^6}{t^2}+\frac{3 q^7}{t^2}
+\frac{2 q^8}{t^2}+\frac{q^9}{t^2}+\frac{q^4}{t}
+\frac{3 q^5}{t}+\frac{6 q^6}{t}+\frac{9 q^7}{t}
+\frac{10 q^8}{t}+\frac{8 q^9}{t}+\frac{4 q^{10}}{t}
+\frac{q^{11}}{t}+q^5 t+5 q^6 t+15 q^7 t+29 q^8 t
+41 q^9 t+45 q^{10} t+34 q^{11} t+16 q^{12} t+4 q^{13} t
+2 q^6 t^2+10 q^7 t^2+27 q^8 t^2+50 q^9 t^2+71 q^{10} t^2
+73 q^{11} t^2+50 q^{12} t^2+21 q^{13} t^2+4 q^{14} t^2
+q^6 t^3+5 q^7 t^3+18 q^8 t^3+44 q^9 t^3+78 q^{10} t^3
+106 q^{11} t^3+103 q^{12} t^3+63 q^{13} t^3+23 q^{14} t^3
+4 q^{15} t^3+2 q^7 t^4+10 q^8 t^4+30 q^9 t^4+68 q^{10} t^4
+116 q^{11} t^4+148 q^{12} t^4+127 q^{13} t^4+68 q^{14} t^4
+19 q^{15} t^4+2 q^{16} t^4+q^7 t^5+5 q^8 t^5+18 q^9 t^5
+48 q^{10} t^5+99 q^{11} t^5+159 q^{12} t^5+184 q^{13} t^5
+138 q^{14} t^5+60 q^{15} t^5+13 q^{16} t^5+q^{17} t^5
+2 q^8 t^6+10 q^9 t^6+30 q^{10} t^6+72 q^{11} t^6+140 q^{12} t^6
+206 q^{13} t^6+208 q^{14} t^6+129 q^{15} t^6+43 q^{16} t^6
+6 q^{17} t^6+q^8 t^7+5 q^9 t^7+18 q^{10} t^7+48 q^{11} t^7
+104 q^{12} t^7+186 q^{13} t^7+245 q^{14} t^7+206 q^{15} t^7
+100 q^{16} t^7+24 q^{17} t^7+2 q^{18} t^7+2 q^9 t^8+10 q^{10} t^8
+30 q^{11} t^8+72 q^{12} t^8+145 q^{13} t^8+235 q^{14} t^8
+263 q^{15} t^8+176 q^{16} t^8+61 q^{17} t^8+9 q^{18} t^8
+q^9 t^9+5 q^{10} t^9+18 q^{11} t^9+48 q^{12} t^9+104 q^{13} t^9
+192 q^{14} t^9+273 q^{15} t^9+246 q^{16} t^9+121 q^{17} t^9
+28 q^{18} t^9+2 q^{19} t^9+2 q^{10} t^{10}+10 q^{11} t^{10}
+30 q^{12} t^{10}+72 q^{13} t^{10}+145 q^{14} t^{10}
+241 q^{15} t^{10}+283 q^{16} t^{10}+192 q^{17} t^{10}
+63 q^{18} t^{10}+8 q^{19} t^{10}+q^{10} t^{11}+5 q^{11} t^{11}
+18 q^{12} t^{11}+48 q^{13} t^{11}+104 q^{14} t^{11}
+192 q^{15} t^{11}+275 q^{16} t^{11}+249 q^{17} t^{11}
+114 q^{18} t^{11}+23 q^{19} t^{11}+q^{20} t^{11}+2 q^{11} t^{12}
+10 q^{12} t^{12}+30 q^{13} t^{12}+72 q^{14} t^{12}
+145 q^{15} t^{12}+238 q^{16} t^{12}+274 q^{17} t^{12}
+173 q^{18} t^{12}+49 q^{19} t^{12}+4 q^{20} t^{12}+q^{11} t^{13}
+5 q^{12} t^{13}+18 q^{13} t^{13}+48 q^{14} t^{13}
+104 q^{15} t^{13}+190 q^{16} t^{13}+263 q^{17} t^{13}
+219 q^{18} t^{13}+87 q^{19} t^{13}+12 q^{20} t^{13}
+2 q^{12} t^{14}+10 q^{13} t^{14}+30 q^{14} t^{14}
+72 q^{15} t^{14}+144 q^{16} t^{14}+229 q^{17} t^{14}
+242 q^{18} t^{14}+130 q^{19} t^{14}+27 q^{20} t^{14}
+q^{21} t^{14}+q^{12} t^{15}+5 q^{13} t^{15}+18 q^{14} t^{15}
+48 q^{15} t^{15}+104 q^{16} t^{15}+185 q^{17} t^{15}
+238 q^{18} t^{15}+168 q^{19} t^{15}+49 q^{20} t^{15}
+4 q^{21} t^{15}+2 q^{13} t^{16}+10 q^{14} t^{16}
+30 q^{15} t^{16}+72 q^{16} t^{16}+142 q^{17} t^{16}
+212 q^{18} t^{16}+192 q^{19} t^{16}+78 q^{20} t^{16}
+9 q^{21} t^{16}+q^{13} t^{17}+5 q^{14} t^{17}+18 q^{15} t^{17}
+48 q^{16} t^{17}+103 q^{17} t^{17}+175 q^{18} t^{17}
+196 q^{19} t^{17}+106 q^{20} t^{17}+19 q^{21} t^{17}
+2 q^{14} t^{18}+10 q^{15} t^{18}+30 q^{16} t^{18}
+72 q^{17} t^{18}+137 q^{18} t^{18}+183 q^{19} t^{18}
+129 q^{20} t^{18}+33 q^{21} t^{18}+q^{22} t^{18}+q^{14} t^{19}
+5 q^{15} t^{19}+18 q^{16} t^{19}+48 q^{17} t^{19}
+101 q^{18} t^{19}+157 q^{19} t^{19}+141 q^{20} t^{19}
+49 q^{21} t^{19}+3 q^{22} t^{19}+2 q^{15} t^{20}
+10 q^{16} t^{20}+30 q^{17} t^{20}+71 q^{18} t^{20}
+127 q^{19} t^{20}+140 q^{20} t^{20}+66 q^{21} t^{20}
+7 q^{22} t^{20}+q^{15} t^{21}+5 q^{16} t^{21}+18 q^{17} t^{21}
+48 q^{18} t^{21}+96 q^{19} t^{21}+128 q^{20} t^{21}
+80 q^{21} t^{21}+13 q^{22} t^{21}+2 q^{16} t^{22}
+10 q^{17} t^{22}+30 q^{18} t^{22}+69 q^{19} t^{22}
+109 q^{20} t^{22}+87 q^{21} t^{22}+21 q^{22} t^{22}
+q^{16} t^{23}+5 q^{17} t^{23}+18 q^{18} t^{23}
+47 q^{19} t^{23}+86 q^{20} t^{23}+87 q^{21} t^{23}
+29 q^{22} t^{23}+q^{23} t^{23}+2 q^{17} t^{24}
+10 q^{18} t^{24}+30 q^{19} t^{24}+64 q^{20} t^{24}
+81 q^{21} t^{24}+37 q^{22} t^{24}+2 q^{23} t^{24}
+q^{17} t^{25}+5 q^{18} t^{25}+18 q^{19} t^{25}
+45 q^{20} t^{25}+68 q^{21} t^{25}+42 q^{22} t^{25}
+4 q^{23} t^{25}+2 q^{18} t^{26}+10 q^{19} t^{26}
+29 q^{20} t^{26}+54 q^{21} t^{26}+44 q^{22} t^{26}
+7 q^{23} t^{26}+q^{18} t^{27}+5 q^{19} t^{27}
+18 q^{20} t^{27}+40 q^{21} t^{27}+42 q^{22} t^{27}
+10 q^{23} t^{27}+2 q^{19} t^{28}+10 q^{20} t^{28}
+27 q^{21} t^{28}+37 q^{22} t^{28}+13 q^{23} t^{28}
+q^{19} t^{29}+5 q^{20} t^{29}+17 q^{21} t^{29}
+30 q^{22} t^{29}+15 q^{23} t^{29}+2 q^{20} t^{30}
+10 q^{21} t^{30}+22 q^{22} t^{30}+16 q^{23} t^{30}
+q^{24} t^{30}+q^{20} t^{31}+5 q^{21} t^{31}+15 q^{22} t^{31}
+15 q^{23} t^{31}+q^{24} t^{31}+2 q^{21} t^{32}+9 q^{22} t^{32}
+13 q^{23} t^{32}+2 q^{24} t^{32}+q^{21} t^{33}+5 q^{22} t^{33}
+10 q^{23} t^{33}+3 q^{24} t^{33}+2 q^{22} t^{34}
+7 q^{23} t^{34}+3 q^{24} t^{34}+q^{22} t^{35}+4 q^{23} t^{35}
+3 q^{24} t^{35}+2 q^{23} t^{36}+3 q^{24} t^{36}+q^{23} t^{37}
+2 q^{24} t^{37}+q^{24} t^{38}+q^{24} t^{39}
\bigr)
\)
\vfill

\noindent
\(
+a^4 \bigl(q^7+5 q^8+12 q^9+21 q^{10}+24 q^{11}+21 q^{12}
+11 q^{13}+4 q^{14}+\frac{q^7}{t^3}+\frac{q^8}{t^3}
+\frac{q^9}{t^3}+\frac{q^6}{t^2}+\frac{2 q^7}{t^2}
+\frac{4 q^8}{t^2}+\frac{4 q^9}{t^2}+\frac{4 q^{10}}{t^2}
+\frac{2 q^{11}}{t^2}+\frac{q^{12}}{t^2}+\frac{2 q^7}{t}
+\frac{5 q^8}{t}+\frac{10 q^9}{t}+\frac{12 q^{10}}{t}
+\frac{11 q^{11}}{t}+\frac{6 q^{12}}{t}+\frac{2 q^{13}}{t}
+2 q^8 t+9 q^9 t+21 q^{10} t+36 q^{11} t+42 q^{12} t
+34 q^{13} t+17 q^{14} t+4 q^{15} t+q^8 t^2+5 q^9 t^2
+17 q^{10} t^2+36 q^{11} t^2+58 q^{12} t^2+61 q^{13} t^2
+45 q^{14} t^2+18 q^{15} t^2+4 q^{16} t^2+2 q^9 t^3
+9 q^{10} t^3+27 q^{11} t^3+54 q^{12} t^3+82 q^{13} t^3
+79 q^{14} t^3+50 q^{15} t^3+16 q^{16} t^3+2 q^{17} t^3
+q^9 t^4+5 q^{10} t^4+17 q^{11} t^4+43 q^{12} t^4+80 q^{13} t^4
+107 q^{14} t^4+90 q^{15} t^4+46 q^{16} t^4+10 q^{17} t^4
+q^{18} t^4+2 q^{10} t^5+9 q^{11} t^5+27 q^{12} t^5+62 q^{13} t^5
+106 q^{14} t^5+124 q^{15} t^5+87 q^{16} t^5+33 q^{17} t^5
+5 q^{18} t^5+q^{10} t^6+5 q^{11} t^6+17 q^{12} t^6+43 q^{13} t^6
+89 q^{14} t^6+132 q^{15} t^6+129 q^{16} t^6+68 q^{17} t^6
+18 q^{18} t^6+q^{19} t^6+2 q^{11} t^7+9 q^{12} t^7+27 q^{13} t^7
+62 q^{14} t^7+117 q^{15} t^7+147 q^{16} t^7+112 q^{17} t^7
+41 q^{18} t^7+7 q^{19} t^7+q^{11} t^8+5 q^{12} t^8+17 q^{13} t^8
+43 q^{14} t^8+89 q^{15} t^8+143 q^{16} t^8+144 q^{17} t^8
+78 q^{18} t^8+18 q^{19} t^8+q^{20} t^8+2 q^{12} t^9+9 q^{13} t^9
+27 q^{14} t^9+62 q^{15} t^9+116 q^{16} t^9+154 q^{17} t^9
+112 q^{18} t^9+40 q^{19} t^9+4 q^{20} t^9+q^{12} t^{10}
+5 q^{13} t^{10}+17 q^{14} t^{10}+43 q^{15} t^{10}
+89 q^{16} t^{10}+140 q^{17} t^{10}+140 q^{18} t^{10}
+67 q^{19} t^{10}+13 q^{20} t^{10}+2 q^{13} t^{11}
+9 q^{14} t^{11}+27 q^{15} t^{11}+62 q^{16} t^{11}
+115 q^{17} t^{11}+143 q^{18} t^{11}+96 q^{19} t^{11}
+26 q^{20} t^{11}+2 q^{21} t^{11}+q^{13} t^{12}+5 q^{14} t^{12}
+17 q^{15} t^{12}+43 q^{16} t^{12}+88 q^{17} t^{12}
+134 q^{18} t^{12}+116 q^{19} t^{12}+46 q^{20} t^{12}
+5 q^{21} t^{12}+2 q^{14} t^{13}+9 q^{15} t^{13}
+27 q^{16} t^{13}+62 q^{17} t^{13}+110 q^{18} t^{13}
+124 q^{19} t^{13}+65 q^{20} t^{13}+12 q^{21} t^{13}
+q^{14} t^{14}+5 q^{15} t^{14}+17 q^{16} t^{14}
+43 q^{17} t^{14}+86 q^{18} t^{14}+117 q^{19} t^{14}
+85 q^{20} t^{14}+21 q^{21} t^{14}+q^{22} t^{14}+2 q^{15} t^{15}
+9 q^{16} t^{15}+27 q^{17} t^{15}+61 q^{18} t^{15}
+101 q^{19} t^{15}+92 q^{20} t^{15}+34 q^{21} t^{15}
+2 q^{22} t^{15}+q^{15} t^{16}+5 q^{16} t^{16}+17 q^{17} t^{16}
+43 q^{18} t^{16}+81 q^{19} t^{16}+94 q^{20} t^{16}
+46 q^{21} t^{16}+6 q^{22} t^{16}+2 q^{16} t^{17}+9 q^{17} t^{17}
+27 q^{18} t^{17}+59 q^{19} t^{17}+84 q^{20} t^{17}
+56 q^{21} t^{17}+10 q^{22} t^{17}+q^{16} t^{18}+5 q^{17} t^{18}
+17 q^{18} t^{18}+42 q^{19} t^{18}+72 q^{20} t^{18}
+61 q^{21} t^{18}+17 q^{22} t^{18}+2 q^{17} t^{19}+9 q^{18} t^{19}
+27 q^{19} t^{19}+54 q^{20} t^{19}+61 q^{21} t^{19}
+22 q^{22} t^{19}+q^{23} t^{19}+q^{17} t^{20}+5 q^{18} t^{20}
+17 q^{19} t^{20}+40 q^{20} t^{20}+55 q^{21} t^{20}
+29 q^{22} t^{20}+2 q^{23} t^{20}+2 q^{18} t^{21}
+9 q^{19} t^{21}+26 q^{20} t^{21}+45 q^{21} t^{21}
+31 q^{22} t^{21}+4 q^{23} t^{21}+q^{18} t^{22}+5 q^{19} t^{22}
+17 q^{20} t^{22}+35 q^{21} t^{22}+33 q^{22} t^{22}
+6 q^{23} t^{22}+2 q^{19} t^{23}+9 q^{20} t^{23}+24 q^{21} t^{23}
+29 q^{22} t^{23}+9 q^{23} t^{23}+q^{19} t^{24}+5 q^{20} t^{24}
+16 q^{21} t^{24}+26 q^{22} t^{24}+11 q^{23} t^{24}
+2 q^{20} t^{25}+9 q^{21} t^{25}+19 q^{22} t^{25}
+12 q^{23} t^{25}+q^{20} t^{26}+5 q^{21} t^{26}+14 q^{22} t^{26}
+12 q^{23} t^{26}+q^{24} t^{26}+2 q^{21} t^{27}+8 q^{22} t^{27}
+11 q^{23} t^{27}+q^{24} t^{27}+q^{21} t^{28}+5 q^{22} t^{28}
+9 q^{23} t^{28}+2 q^{24} t^{28}+2 q^{22} t^{29}+6 q^{23} t^{29}
+2 q^{24} t^{29}+q^{22} t^{30}+4 q^{23} t^{30}+3 q^{24} t^{30}
+2 q^{23} t^{31}+2 q^{24} t^{31}+q^{23} t^{32}+2 q^{24} t^{32}
+q^{24} t^{33}+q^{24} t^{34}
\bigr)
\)
\vfill

\noindent
\(
+a^5 \bigl(q^{10}
+4 q^{11}+9 q^{12}+13 q^{13}+12 q^{14}+7 q^{15}+2 q^{16}
+\frac{q^9}{t^3}+\frac{q^{10}}{t^3}+\frac{q^{11}}{t^3}
+\frac{q^9}{t^2}+\frac{2 q^{10}}{t^2}+\frac{3 q^{11}}{t^2}
+\frac{3 q^{12}}{t^2}+\frac{2 q^{13}}{t^2}+\frac{q^{14}}{t^2}
+\frac{2 q^{10}}{t}+\frac{5 q^{11}}{t}+\frac{8 q^{12}}{t}
+\frac{7 q^{13}}{t}+\frac{4 q^{14}}{t}+\frac{q^{15}}{t}
+2 q^{11} t+7 q^{12} t+16 q^{13} t+21 q^{14} t+18 q^{15} t
+8 q^{16} t+q^{17} t+q^{11} t^2+4 q^{12} t^2+12 q^{13} t^2
+23 q^{14} t^2+28 q^{15} t^2+20 q^{16} t^2+7 q^{17} t^2
+q^{18} t^2+2 q^{12} t^3+7 q^{13} t^3+19 q^{14} t^3+32 q^{15} t^3
+34 q^{16} t^3+19 q^{17} t^3+4 q^{18} t^3+q^{12} t^4+4 q^{13} t^4
+12 q^{14} t^4+27 q^{15} t^4+40 q^{16} t^4+34 q^{17} t^4
+13 q^{18} t^4+2 q^{19} t^4+2 q^{13} t^5+7 q^{14} t^5
+19 q^{15} t^5+36 q^{16} t^5+44 q^{17} t^5+27 q^{18} t^5
+7 q^{19} t^5+q^{13} t^6+4 q^{14} t^6+12 q^{15} t^6+27 q^{16} t^6
+44 q^{17} t^6+39 q^{18} t^6+16 q^{19} t^6+2 q^{20} t^6
+2 q^{14} t^7+7 q^{15} t^7+19 q^{16} t^7+37 q^{17} t^7
+46 q^{18} t^7+27 q^{19} t^7+6 q^{20} t^7+q^{14} t^8+4 q^{15} t^8
+12 q^{16} t^8+27 q^{17} t^8+44 q^{18} t^8+37 q^{19} t^8
+13 q^{20} t^8+q^{21} t^8+2 q^{15} t^9+7 q^{16} t^9+19 q^{17} t^9
+36 q^{18} t^9+42 q^{19} t^9+22 q^{20} t^9+3 q^{21} t^9
+q^{15} t^{10}+4 q^{16} t^{10}+12 q^{17} t^{10}+27 q^{18} t^{10}
+40 q^{19} t^{10}+29 q^{20} t^{10}+7 q^{21} t^{10}+2 q^{16} t^{11}
+7 q^{17} t^{11}+19 q^{18} t^{11}+34 q^{19} t^{11}+34 q^{20} t^{11}
+12 q^{21} t^{11}+q^{22} t^{11}+q^{16} t^{12}+4 q^{17} t^{12}
+12 q^{18} t^{12}+26 q^{19} t^{12}+34 q^{20} t^{12}+17 q^{21} t^{12}
+2 q^{22} t^{12}+2 q^{17} t^{13}+7 q^{18} t^{13}+19 q^{19} t^{13}
+30 q^{20} t^{13}+22 q^{21} t^{13}+4 q^{22} t^{13}+q^{17} t^{14}
+4 q^{18} t^{14}+12 q^{19} t^{14}+24 q^{20} t^{14}+24 q^{21} t^{14}
+7 q^{22} t^{14}+2 q^{18} t^{15}+7 q^{19} t^{15}+18 q^{20} t^{15}
+23 q^{21} t^{15}+10 q^{22} t^{15}+q^{18} t^{16}+4 q^{19} t^{16}
+12 q^{20} t^{16}+20 q^{21} t^{16}+12 q^{22} t^{16}+q^{23} t^{16}
+2 q^{19} t^{17}+7 q^{20} t^{17}+16 q^{21} t^{17}+13 q^{22} t^{17}
+2 q^{23} t^{17}+q^{19} t^{18}+4 q^{20} t^{18}+11 q^{21} t^{18}
+13 q^{22} t^{18}+3 q^{23} t^{18}+2 q^{20} t^{19}+7 q^{21} t^{19}
+12 q^{22} t^{19}+4 q^{23} t^{19}+q^{20} t^{20}+4 q^{21} t^{20}
+9 q^{22} t^{20}+5 q^{23} t^{20}+2 q^{21} t^{21}+6 q^{22} t^{21}
+6 q^{23} t^{21}+q^{21} t^{22}+4 q^{22} t^{22}+5 q^{23} t^{22}
+2 q^{22} t^{23}+4 q^{23} t^{23}+q^{24} t^{23}+q^{22} t^{24}
+3 q^{23} t^{24}+q^{24} t^{24}+2 q^{23} t^{25}+q^{24} t^{25}
+q^{23} t^{26}+q^{24} t^{26}+q^{24} t^{27}+q^{24} t^{28}
\bigr)
\)
\vfill

\noindent
\(
+a^6 \bigl((q^{14}+3 q^{15}+3 q^{16}+2 q^{17}+\frac{q^{12}}{t^3}
+\frac{q^{13}}{t^2}+\frac{q^{14}}{t^2}+\frac{q^{15}}{t^2}
+\frac{q^{13}}{t}+\frac{2 q^{14}}{t}+\frac{2 q^{15}}{t}
+\frac{q^{16}}{t}+q^{14} t+2 q^{15} t+5 q^{16} t+4 q^{17} t
+q^{18} t+q^{15} t^2+3 q^{16} t^2+6 q^{17} t^2+3 q^{18} t^2
+q^{19} t^2+q^{15} t^3+2 q^{16} t^3+5 q^{17} t^3+7 q^{18} t^3
+2 q^{19} t^3+q^{16} t^4+3 q^{17} t^4+6 q^{18} t^4+5 q^{19} t^4
+q^{20} t^4+q^{16} t^5+2 q^{17} t^5+5 q^{18} t^5+6 q^{19} t^5
+3 q^{20} t^5+q^{17} t^6+3 q^{18} t^6+6 q^{19} t^6+4 q^{20} t^6
+q^{21} t^6+q^{17} t^7+2 q^{18} t^7+5 q^{19} t^7+6 q^{20} t^7
+q^{21} t^7+q^{18} t^8+3 q^{19} t^8+5 q^{20} t^8+3 q^{21} t^8
+q^{18} t^9+2 q^{19} t^9+5 q^{20} t^9+4 q^{21} t^9+q^{22} t^9
+q^{19} t^{10}+3 q^{20} t^{10}+4 q^{21} t^{10}+q^{22} t^{10}
+q^{19} t^{11}+2 q^{20} t^{11}+4 q^{21} t^{11}+2 q^{22} t^{11}
+q^{20} t^{12}+3 q^{21} t^{12}+2 q^{22} t^{12}+q^{20} t^{13}
+2 q^{21} t^{13}+3 q^{22} t^{13}+q^{21} t^{14}+2 q^{22} t^{14}
+q^{23} t^{14}+q^{21} t^{15}+2 q^{22} t^{15}+q^{23} t^{15}
+q^{22} t^{16}+q^{23} t^{16}+q^{22} t^{17}+q^{23} t^{17}
+q^{23} t^{18}+q^{23} t^{19}+q^{24} t^{21}\bigr).
\)

\renewcommand{\baselinestretch}{1.2} 
\smallskip
}

\medskip

It is a long formula, but
we hope that it can be of interest to the
experts; the torus knots
$\{3,2\}$ and $\{4,3\}$ are special in several ways.

We note that all formulas in this
paper are (finite) polynomials, due to the fact that we
use $P_b/P_b(t^\rho)$ instead of $P_b$. We do not focus
on it too much, but it is an important (actually,
surprising) part of Conjecture \ref{HOMFLY}. In the Khovanov-
Rozansky theory such choice corresponds to reduced 
KhR polynomials vs. unreduced ones.

The corresponding ``dimension"
(the sum of all coefficients) is $249^2$, where
$249$ is the dimension in the case of $b=\om_1$
for this knot. Moreover, 
$H\!D_{9,4}(\om_2\,;\, q,t=1,a)$ is the square of 
$H\!D_{9,4}(\om_1\,;\, q,t=1,a)$, which confirms
Conjecture \ref{CONJEVAL} 
on the evaluations at $t=1$. We note that 
$H\!D_{9,4}(\om_2\,;\, q=1,t,a)$
has no relation to that for $\om_1$ and is complicated;
by duality, it coincides with   
$H\!D_{9,4}(2\om_1\,;\, t^{-1},1,a)$ up to a power of $t$.

The case $t=1$ is of course very exceptional (the DAHA
becomes essentially the Weyl algebra); the coinvariants
we construct by formula (\ref{jones-d})
are identically $1$ in the tilde-normalization.
The involvement of $a$ makes Conjecture \ref{CONJEVAL}
non-trivial and challenging. 
 
\medskip

{\em The positivity for $2\om_3$ and trefoil.}
We will finish this section with a  
demonstration of the positivity claim from
Conjecture \ref{HOMFLY} (and other symmetries)
for the $3\times 2$\~rectangle ($b=2\om_3$).
Its duality with the $2\times 3$\~rectangle ($b=3\om_2$)
was also checked, so we will skip the formula for the 
corresponding super-polynomial. One has:

$$
H\!D_{3,2}(b=2\om_3;\,q,t,a)\ =
$$
\renewcommand{\baselinestretch}{0.5} 
{\small

\noindent
\(
1+\frac{a^6 q^{15}}{t^6}
+q^2 t+q^3 t+q^2 t^2+q^3 t^2+q^4 t^2+q^2 t^3+q^3 
t^3+q^4 t^3+q^5 t^3+2 q^4 t^4+2 q^5 t^4+q^6 t^4+q^4 t^5+2 q^5 t^5+2 
q^6 t^5+q^7 t^5+q^4 t^6+q^5 t^6+3 q^6 t^6+2 q^7 t^6+2 q^6 t^7+3 q^7 
t^7+q^8 t^7+q^6 t^8+2 q^7 t^8+2 q^8 t^8+q^6 t^9+q^7 t^9+2 q^8 t^9+2 
q^9 t^9+2 q^8 t^{10}+2 q^9 t^{10}+q^{10} t^{10}+q^8 t^{11}+2 q^9 
t^{11}+q^{10} t^{11}+q^8 t^{12}+q^9 t^{12}+q^{10} t^{12}+q^{10} 
t^{13}+q^{11} t^{13}+q^{10} t^{14}+q^{11} t^{14}+q^{10} 
t^{15}+q^{11} t^{15}+q^{12} t^{18}
\)

\noindent
\(
+a^1 \bigl(q^2+q^3+2 q^4+3 q^5+q^6+\frac{q^2}{t^2}+
\frac{q^3}{t^2}+\frac{q^2}{t}+\frac{q^3}{t}+\frac{q^4}{t}
+\frac{q^5}{t}+3 
q^4 t+5 q^5 t+3 q^6 t+q^7 t+2 q^4 t^2+4 q^5 t^2+5 q^6 t^2+4 q^7 
t^2+q^8 t^2+q^4 t^3+2 q^5 t^3+5 q^6 t^3+7 q^7 t^3+3 q^8 t^3+4 q^6 
t^4+8 q^7 t^4+6 q^8 t^4+2 q^9 t^4+2 q^6 t^5+5 q^7 t^5+7 q^8 t^5+5 
q^9 t^5+q^{10} t^5+q^6 t^6+2 q^7 t^6+6 q^8 t^6+8 q^9 t^6+3 q^{10} 
t^6+4 q^8 t^7+8 q^9 t^7+5 q^{10} t^7+q^{11} t^7+2 q^8 t^8+5 q^9 
t^8+5 q^{10} t^8+2 q^{11} t^8+q^8 t^9+2 q^9 t^9+4 q^{10} t^9+4 
q^{11} t^9+q^{12} t^9+3 q^{10} t^{10}+5 q^{11} t^{10}+2 q^{12} 
t^{10}+2 q^{10} t^{11}+4 q^{11} t^{11}+2 q^{12} t^{11}+q^{10} 
t^{12}+2 q^{11} t^{12}+q^{12} t^{12}+q^{12} t^{13}+q^{13} 
t^{13}+q^{12} t^{14}+q^{13} t^{14}+q^{12} t^{15}+q^{13} 
t^{15}\bigr) 
\)

\noindent
\(
+a^2 \bigl(q^5+3 q^6+8 q^7+6 q^8+2 q^9+\frac{q^5}{t^4}+
\frac{q^4}{t^3}+\frac{2 q^5}{t^3}+\frac{q^6}{t^3}
+\frac{q^4}{t^2}+\frac{3 q^5}{t^2}+\frac{2 
q^6}{t^2}+\frac{2 q^7}{t^2}+\frac{q^8}{t^2}+\frac{q^4}{t}+\frac{2 
q^5}{t}+\frac{3 q^6}{t}+\frac{5 q^7}{t}+\frac{3 q^8}{t}+2 q^6 t+7 
q^7 t+8 q^8 t+5 q^9 t+2 q^{10} t+q^6 t^2+4 q^7 t^2+7 q^8 t^2+10 q^9 
t^2+5 q^{10} t^2+q^{11} t^2+q^7 t^3+5 q^8 t^3+11 q^9 t^3+9 q^{10} 
t^3+2 q^{11} t^3+2 q^8 t^4+9 q^9 t^4+9 q^{10} t^4+5 q^{11} 
t^4+q^{12} t^4+q^8 t^5+4 q^9 t^5+7 q^{10} t^5+7 q^{11} t^5+3 q^{12} 
t^5+q^9 t^6+4 q^{10} t^6+9 q^{11} t^6+6 q^{12} t^6+q^{13} t^6+2 
q^{10} t^7+7 q^{11} t^7+6 q^{12} t^7+q^{13} t^7+q^{10} t^8+4 q^{11} 
t^8+4 q^{12} t^8+2 q^{13} t^8+q^{11} t^9+2 q^{12} t^9+2 q^{13} 
t^9+q^{14} t^9+q^{12} t^{10}+3 q^{13} t^{10}+q^{14} t^{10}+q^{12} 
t^{11}+2 q^{13} t^{11}+q^{14} t^{11}+q^{13} t^{12}\bigr)
\)

\noindent
\(
+a^3 \bigl(q^8+6 q^9+8 
q^{10}+4 q^{11}+q^{12}+\frac{q^7}{t^5}+\frac{q^8}{t^5}+\frac{2 
q^7}{t^4}+\frac{2 q^8}{t^4}+\frac{q^6}{t^3}+\frac{3 
q^7}{t^3}+\frac{3 q^8}{t^3}+\frac{2 
q^9}{t^3}+\frac{q^{10}}{t^3}+\frac{2 q^7}{t^2}+\frac{3 
q^8}{t^2}+\frac{4 q^9}{t^2}+\frac{4 
q^{10}}{t^2}+\frac{q^{11}}{t^2}+\frac{q^7}{t}+\frac{2 
q^8}{t}+\frac{6 q^9}{t}+\frac{7 q^{10}}{t}+\frac{2 q^{11}}{t}+3 q^9 
t+6 q^{10} t+6 q^{11} t+4 q^{12} t+q^{13} t+q^9 t^2+3 q^{10} t^2+7 
q^{11} t^2+7 q^{12} t^2+2 q^{13} t^2+q^{10} t^3+6 q^{11} t^3+8 
q^{12} t^3+3 q^{13} t^3+3 q^{11} t^4+5 q^{12} t^4+3 q^{13} 
t^4+q^{14} t^4+q^{11} t^5+2 q^{12} t^5+3 q^{13} t^5+2 q^{14} 
t^5+q^{12} t^6+3 q^{13} t^6+3 q^{14} t^6+q^{15} t^6+2 q^{13} t^7+2 
q^{14} t^7+q^{13} t^8+q^{14} t^8\bigr)
\)

\noindent
\(
+a^4 \bigl(q^{11}+3 q^{12}+2 
q^{13}+q^{14}+\frac{q^{10}}{t^6}+\frac{q^9}{t^5}+\frac{2 
q^{10}}{t^5}+\frac{q^{11}}{t^5}+\frac{q^9}{t^4}+\frac{3 
q^{10}}{t^4}+\frac{q^{11}}{t^4}+\frac{q^9}{t^3}+\frac{2 
q^{10}}{t^3}+\frac{2 q^{11}}{t^3}+\frac{2 
q^{12}}{t^3}+\frac{q^{13}}{t^3}+\frac{q^{10}}{t^2}+\frac{2 
q^{11}}{t^2}+\frac{4 q^{12}}{t^2}+\frac{2 q^{13}}{t^2}+\frac{2 
q^{11}}{t}+\frac{5 q^{12}}{t}+\frac{3 q^{13}}{t}+q^{12} t+2 q^{13} 
t+2 q^{14} t+q^{15} t+q^{13} t^2+3 q^{14} t^2+q^{15} t^2+q^{13} 
t^3+2 q^{14} t^3+q^{15} t^3+q^{14} t^4\bigr)
\)

\noindent
\(
a^5 \bigl(\frac{q^{12}}{t^6}+\frac{q^{13}}{t^6}+\frac{q^{12}}{t^5}+
\frac{q^{13}}{t^5}+\frac{q^{12}}{t^4}+\frac{q^{13}}{t^4}+
\frac{q^{14}}{t^3}+\frac{q^{15}}{t^3}+\frac{q^{14}}{t^2}+
\frac{q^{15}}{t^2}+\frac{q^{14}}{t}+\frac{q^{15}}{t}\bigr).
\)
}
\renewcommand{\baselinestretch}{1.2} 
\smallskip

The $t$\~evaluation is as follows:
\begin{align*}
&H\!D_{3,2}(b=2\om_3;\,q,t=1,a)\\
&=\ 
\bigr(1+q^2+a q^2+q^3+a q^3+q^4+a q^4+a q^5+a^2 q^5\bigl)^3,
\end{align*}
which matches the evaluation conjecture since
$H\!D_{3,2}(b=2\om_1;\,q,t,a)$
$$
=1 + a^2 q^5 + q^2 t + q^3 t + q^4 t^2 + a(q^2 + q^3 + q^4 t + q^5t);
$$
the cube of its evaluation at $t=1$ is exactly the quantity above.
The dimension is $3^6$, as it is supposed to be. The evaluation
at $q=1$ is the square of that for $\om_3$ (the $3$\~column).
\smallskip

\subsection {Reduced KhR polynomials}\label{sec:Khr-p}

We will provide several sample formulas for
DAHA-KhR polynomials for the smallest 
exceptional $n$. Hopefully these cases
can be verified now or in future within
the Khovanov and  Khovanov-Rozansky theory.

The first exceptional polynomial is always $1$
in the examples below; we will disregard it. 
It is actually a general feature of rectangle diagrams,
which is not present for more general diagrams. 

To simplify
using our formulas, {\em we switch in this section
standard parameters}
from (\ref{qtareli}),(\ref{qtarel}). In our
parameters, the KhR formulas would contain $\sqrt{q/t}$.

We note that there are some links of Khovanov
polynomials to generalized Verlinde algebras 
in \cite{AS} similar to what we use in this note,
for instance, the KhR substitution $a_{st}=q_{st}^n$
recalculated in terms of the DAHA parameters is present there.
Generally, the source of 
super-polynomials is not important for applying the 
procedure from (\ref{redprocedure}) or similar algorithms
(originated in \cite{DGR}).

This section has a clear practical meaning. It is
expected that only sufficiently small knots satisfy  
(\ref{nonzerod}) for all exceptional $n$. However, various 
applications, including those to modular forms and in physics, 
require only simple knots.
Also, Conjecture \ref{CONJKHR} provides potential examples
with non-trivial weights, which is beyond the current level of the
Khovanov-Rozansky homology and related categorification theory.  
And last but not least, the super-polynomials (where DAHA 
are on solid ground) can be 
checked now mainly against well-understood Khovanov ($n=2$)
polynomials, so $d_2$ and other differentials are really 
needed for establishing this and other connections.
\smallskip
    
{\em Uncolored examples.}
We will begin with the torus knot $\{9,4\}$. The stabilization 
of the DAHA-KR polynomials to the corresponding
super-polynomial for this knot is at $n=7$. Recall that
$\tilde{J\!D}^n$ and $K\!R^n$ are obtained
from the super-polynomial $H\!D$ 
upon the corresponding substitutions $a^2=-q^{2n}/t$ and $a^2=q^{2n}$ 
({\em in the standard parameters}) followed by further reductions 
of the monomials. The stabilization means that there will be 
no actual reduction of the number of monomials after these 
substitutions at such $n$. For the first non-stable $n=2,3$, 
one has:

$$
K\!R^2_{9,4}(b=\om_1;\,q,t)\ =
$$
\renewcommand{\baselinestretch}{0.5} 
{\small
\noindent
\(
1+q^4 t^2+q^6 t^3+q^6 t^4+q^{10} t^5+q^8 t^6+q^{12} t^7
+2 q^{12} t^8+2 q^{14} t^9+q^{14} t^{10}+2 q^{18} t^{11}
+q^{16} t^{12}+q^{18} t^{12}+q^{20} t^{12}+2 q^{20} t^{13}
+q^{20} t^{14}+q^{22} t^{15}
+q^{24} t^{15}+q^{24} t^{16}+q^{26} t^{16}+q^{26} t^{17},
\)
}
\renewcommand{\baselinestretch}{1.2} 

$$
K\!R^3_{9,4}(b=\om_1;\,q,t)\ =
$$
\renewcommand{\baselinestretch}{0.5} 
{\small
\noindent
\(
1+q^4 t^2+q^8 t^3+q^6 t^4+q^8 t^4+q^{10} t^5+q^{12} t^5+q^8 t^6
+q^{10} t^6+2 q^{14} t^7+q^{16} t^7+2 q^{12} t^8+q^{14} t^8
+q^{18} t^8+2 q^{16} t^9+2 q^{18} t^9+q^{14} t^{10}+2 q^{16} t^{10}
+q^{22} t^{10}+q^{18} t^{11}+3 q^{20} t^{11}+q^{22} t^{11}
+q^{16} t^{12}+2 q^{18} t^{12}+q^{20} t^{12}+q^{24} t^{12}
+q^{26} t^{12}+3 q^{22} t^{13}+3 q^{24} t^{13}+q^{20} t^{14}
+2 q^{22} t^{14}+q^{26} t^{14}+q^{28} t^{14}+q^{24} t^{15}
+4 q^{26} t^{15}+q^{28} t^{15}+2 q^{24} t^{16}+q^{26} t^{16}
+2 q^{30} t^{16}+q^{32} t^{16}+2 q^{28} t^{17}+3 q^{30} t^{17}
+q^{28} t^{18}+q^{30} t^{18}+2 q^{34} t^{18}+2 q^{32} t^{19}
+2 q^{34} t^{19}+q^{32} t^{20}+q^{36} t^{20}
+q^{38} t^{20}+q^{36} t^{21}+q^{38} t^{21}+q^{42} t^{22}.
\)
}
\renewcommand{\baselinestretch}{1.2} 
\smallskip

The torus knot $\{8,5\}$ will be next;
the stabilization to the
super-polynomial occurs at $n=9$:

$$
K\!R^2_{8,5}(b=\om_1;\,q,t)\ =
$$
\renewcommand{\baselinestretch}{0.5} 
{\small
\noindent
\(
1+q^4 t^2+q^6 t^3+q^6 t^4+q^{10} t^5+q^8 t^6+q^{12} t^7+q^{10} t^8
+q^{12} t^8+2 q^{14} t^9+2 q^{14} t^{10}+2 q^{16} t^{11}
+q^{18} t^{11}+2 q^{16} t^{12}+q^{20} t^{12}+3 q^{20} t^{13}
+2 q^{20} t^{14}+q^{22} t^{14}+3 q^{22} t^{15}+q^{22} t^{16}
+q^{26} t^{16}+2 q^{26} t^{17}+q^{26} t^{18}+q^{28} t^{18}
+q^{28} t^{19},
\)
}
\renewcommand{\baselinestretch}{1.2} 

$$
K\!R^3_{8,5}(b=\om_1;\,q,t)\ =
$$
\renewcommand{\baselinestretch}{0.5} 
{\small
\noindent
\(
1+q^4 t^2+q^8 t^3+q^6 t^4+q^8 t^4+q^{10} t^5+q^{12} t^5+q^8 t^6
+q^{10} t^6+2 q^{14} t^7+q^{16} t^7+q^{10} t^8+2 q^{12} t^8
+q^{18} t^8+2 q^{16} t^9+2 q^{18} t^9+2 q^{14} t^{10}
+2 q^{16} t^{10}+q^{22} t^{10}+3 q^{18} t^{11}+3 q^{20} t^{11}
+2 q^{16} t^{12}+3 q^{18} t^{12}+q^{24} t^{12}+q^{26} t^{12}
+q^{20} t^{13}+5 q^{22} t^{13}+2 q^{24} t^{13}+3 q^{20} t^{14}
+2 q^{22} t^{14}+2 q^{26} t^{14}+q^{28} t^{14}+5 q^{24} t^{15}
+5 q^{26} t^{15}+q^{22} t^{16}+4 q^{24} t^{16}+2 q^{28} t^{16}
+3 q^{30} t^{16}+q^{26} t^{17}+7 q^{28} t^{17}+2 q^{30} t^{17}
+2 q^{26} t^{18}+3 q^{28} t^{18}+4 q^{32} t^{18}+2 q^{34} t^{18}
+4 q^{30} t^{19}+6 q^{32} t^{19}+q^{36} t^{19}+2 q^{30} t^{20}
+2 q^{32} t^{20}+2 q^{34} t^{20}+4 q^{36} t^{20}+5 q^{34} t^{21}
+4 q^{36} t^{21}+q^{40} t^{21}+q^{34} t^{22}+q^{36} t^{22}
+4 q^{38} t^{22}+2 q^{40} t^{22}+3 q^{38} t^{23}+2 q^{40} t^{23}
+q^{42} t^{23}+q^{40} t^{24}+3 q^{42} t^{24}+q^{44} t^{24}
+q^{42} t^{25}+q^{44} t^{25}+q^{46} t^{25}+q^{46} t^{26}.
\)
}
\renewcommand{\baselinestretch}{1.2} 
\smallskip

The $n=2$ polynomials in the last two examples
coincide with the corresponding {\em reduced
and tilde-normalized} Khovanov polynomials, so there are
chances that the $n=3$ reductions will match the
corresponding Khovanov -Rozansky $SL(3)$\~polynomials 
(which are not known for these knots).
\smallskip

We note that the first stable $n$ 
can be significantly greater than $s$ or $r$:
$n=4$ for $\{7,3\}$ (see below),
$7$ for $\{9,4\}$, $9$ for $\{8,5\}$, $19$ for  $\{12,7\}$,
$26$ for $\{10,9\}$.
It was conjectured in \cite{DGR} that $n=2$ is already stable 
for $s=2$ and that the first stable $n$ 
for $\{r=3m+1,s=3\}$ is $n=m+2$; let us provide an
example.
\medskip

{\em The case of $\{7,3\}$.}
Here the DAHA super-polynomial is
$$
H\!D_{7,3}(b=\om_1;\,q,t,a)\ =
$$
\renewcommand{\baselinestretch}{0.5} 
{\small
\noindent
\(
1+q t+q^2 t+q^2 t^2+q^3 t^2+q^4 t^2+q^3 t^3+q^4 t^3
+q^4 t^4+q^5 t^4+q^5 t^5+q^6 t^6+a(q+q^2+q^2 t+2 q^3 t+q^4 t+q^3 t^2
+2 q^4 t^2+q^5 t^2+q^4 t^3+2 q^5 t^3+q^5 t^4+q^6 t^4+q^6 t^5)
+a^2 (q^3+q^4 t+q^5 t+q^5 t^2+q^6 t^3)
\)
}
\renewcommand{\baselinestretch}{1.2} 
\smallskip
using the DAHA parameters.

In terms of the {\em standard parameters}, it reads:
\smallskip

\renewcommand{\baselinestretch}{0.5} 
{\small
\noindent
\(
1+q^4 t^2+q^6 t^4+q^8 t^4+q^{10} t^6+q^{12} t^6
+q^{12} t^8+q^{14} t^8+q^{16} t^8+q^{18} t^{10}
+q^{20} t^{10}+q^{24} t^{12}+a^2 (q^2 t^3
+q^4 t^5+q^6 t^5+2 q^8 t^7+q^{10} t^7+q^{10} t^9
+2 q^{12} t^9+q^{14} t^9+q^{14} t^{11}+2 q^{16} t^{11}
+q^{18} t^{11}+q^{20} t^{13}+q^{22} t^{13})
+a^4 (q^6 t^8
+q^{10} t^{10}+q^{12} t^{12}+q^{14} t^{12}+q^{18} t^{14});
\)
}
\renewcommand{\baselinestretch}{1.2} 

it coincides with $a^{-12} q^{12}$ times the
expression conjectured for the torus knots $\{r=3m+1,s=3\}$
in formula (93) from \cite{DGR}.
Accordingly, our first stable value $n=4$ matches
the formula $n=m+2$ there formally deduced from (93); 
see right after formula (101).  These coincidences were 
also verified for $m=3,4$.

The corresponding Khovanov-Rozansky polynomials are
$$
K\!R^2_{7,3}(b=\om_1;\,q,t)\ =
$$
\renewcommand{\baselinestretch}{0.5} 
{\small
\noindent
\ \ \ \ \ \ \ \(
1+q^4 t^2+q^6 t^3+q^6 t^4+q^{10} t^5+q^{10} t^6+q^{12} t^7+
q^{12} t^8+q^{16} t^9,
\)
}
\renewcommand{\baselinestretch}{1.2} 

$$
K\!R^3_{7,3}(b=\om_1;\,q,t)\ =
$$
\renewcommand{\baselinestretch}{0.5} 
{\small
\noindent
\(
1+q^4 t^2+q^8 t^3+q^6 t^4+q^8 t^4+q^{10} t^5
+q^{12} t^5+q^{10} t^6+q^{12} t^6+2 q^{14} t^7
+q^{16} t^7+q^{12} t^8+q^{14} t^8+q^{18} t^8
+2 q^{18} t^9+q^{20} t^9+q^{18} t^{10}+q^{22} t^{10}
+2 q^{22} t^{11}+q^{24} t^{12}
+q^{26} t^{12}+q^{26} t^{13}+q^{28} t^{13}+q^{30} t^{14},
\)
}
\renewcommand{\baselinestretch}{1.2} 
where the first really coincides with the reduced
and tilde-normalized Khovanov polynomial, 
the second one is expected to be such renormalization
of the real polynomial KhR${}^{3}_{7,3}$.
\medskip

{\em Counterexamples.}
Let us provide the DAHA-Khovanov polynomial
for $\{12,7\}$ (the stabilization to the
super-polynomial is at $n=19$):

$$
K\!R^2_{12,7}(b=\om_1;\,q,t)\ =
$$
\renewcommand{\baselinestretch}{0.5} 
{\small
\noindent
\(
1+q^4 t^2+q^6 t^3+q^6 t^4+q^{10} t^5+q^8 t^6+q^{12} t^7+q^{10} t^8
+q^{12} t^8+2 q^{14} t^9+q^{12} t^{10}+q^{14} t^{10}+2 q^{16} t^{11}
+q^{18} t^{11}+q^{14} t^{12}+2 q^{16} t^{12}+q^{20} t^{12}
+2 q^{18} t^{13}+2 q^{20} t^{13}+3 q^{18} t^{14}+q^{22} t^{14}
+2 q^{20} t^{15}+3 q^{22} t^{15}+3 q^{20} t^{16}+q^{24} t^{16}
+q^{26} t^{16}+5 q^{24} t^{17}+2 q^{22} t^{18}+2 q^{24} t^{18}
+q^{26} t^{18}+q^{28} t^{18}+6 q^{26} t^{19}+2 q^{24} t^{20}
+3 q^{26} t^{20}+2 q^{30} t^{20}+6 q^{28} t^{21}+2 q^{30} t^{21}
+5 q^{28} t^{22}+3 q^{32} t^{22}+5 q^{30} t^{23}+5 q^{32} t^{23}
+4 q^{30} t^{24}+5 q^{34} t^{24}+9 q^{34} t^{25}+2 q^{32} t^{26}
+2 q^{34} t^{26}+4 q^{36} t^{26}+2 q^{38} t^{26}+8 q^{36} t^{27}
+q^{40} t^{27}+q^{34} t^{28}+4 q^{36} t^{28}+q^{38} t^{28}
+4 q^{40} t^{28}+7 q^{38} t^{29}+q^{40} t^{29}+q^{42} t^{29}
+4 q^{38} t^{30}+5 q^{42} t^{30}+4 q^{40} t^{31}+5 q^{42} t^{31}
+q^{44} t^{31}+2 q^{40} t^{32}+2 q^{42} t^{32}+6 q^{44} t^{32}
+q^{42} t^{33}+8 q^{44} t^{33}+q^{48} t^{33}+3 q^{44} t^{34}
+3 q^{46} t^{34}+4 q^{46} t^{35}+q^{50} t^{35}+q^{46} t^{36}
+2 q^{48} t^{36}+3 q^{50} t^{36}+2 q^{48} t^{37}+q^{50} t^{37}
+q^{52} t^{37}+q^{48} t^{38}+q^{50} t^{38}+q^{52} t^{38}
+3 q^{50} t^{39}+q^{50} t^{40}+q^{52} t^{40}+q^{52} t^{41}
+q^{54} t^{41}+q^{54} t^{42}+q^{56} t^{42}+q^{56} t^{43}.
\)
}
\renewcommand{\baselinestretch}{1.2} 
\smallskip

It is the only example we found so far where $(ii)$ from
Conjecture \ref{CONJKHR} fails, so $d_2$ here must
be of more complicated nature
than in (\ref{nonzerod}).

The actual reduced tilde-normalized Khovanov polynomial, 
calculated using 
http://katlas.org/wiki/KnotTheory, equals 
the sum of the DAHA polynomial 
$K\!R^2_{12,7}(b=\om_1;\,q,t)$ (posted above) and
the correction, which is:
\renewcommand{\baselinestretch}{0.5} 
{\small
\noindent
\begin{align*}
&q^{40} \bigl(t^{29}-t^{31}\bigr)+q^{42} \bigl(t^{31}-t^{33}\bigr)
+q^{48} \bigl(2 t^{34}+2 t^{35}-t^{36}-2 t^{37}-t^{38}\bigr)\\
&+q^{50} \bigl(t^{36}+3 t^{37}-3 t^{39}-t^{40}\bigr)
+q^{52} \bigl(t^{38}+t^{39}-t^{40}-t^{41}\bigr)\\
&+q^{54} \bigl(t^{39}
+t^{40}-t^{41}-t^{42}\bigr)
+q^{56} \bigl(t^{40}+t^{41}-t^{42}-t^{43}\bigr).
\end{align*}
}
\renewcommand{\baselinestretch}{1.2} 

This correction terms here do not change the dimension of
the polynomial $K\!R^2_{12,7}(b=\om_1;\,q,t)\ $ and the corresponding
$q$\~filtration; so the family of the DAHA-KhR polynomials
$\{K\!R^n\}$ still may satisfy claim $(iii)$ from this 
conjecture. 
Interestingly, the $t$\~powers in the 
real Khovanov polynomial are no greater than in
$K\!R^2_{12,7}$  (recall that in this section, 
$t$ is $t_{st}$, $q$ is $q_{st}$). Also, the 
correction from $K\!R^2$ to the real Khovanov polynomial in
this example affects only the greatest $t$\~terms, so certain 
stabilization of both constructions can be expected for 
sufficiently large torus knots (which matches the approach 
by Gorsky and Rasmussen).
\smallskip

The last uncolored example will be the torus knot $\{9,10\}$, where
the $a$\~degree of the super-polynomial is $8$. The
stabilization begins at $n=26$.  We did not find/obtain the 
formula for the corresponding real Khovanov polynomial, but it is 
not expected to coincide with the following DAHA-Khovanov 
polynomial:
$$
K\!R^2_{9,10}(b=\om_1;\,q,t)\ =
$$
\renewcommand{\baselinestretch}{0.5} 
{\small
\noindent
\(
1+q^4 t^2+q^6 t^3+q^6 t^4+q^{10} t^5+q^8 t^6+q^{12} t^7
+q^{10} t^8+q^{12} t^8+2 q^{14} t^9+q^{12} t^{10}
+q^{14} t^{10}+q^{20} t^{10}+2 q^{16} t^{11}+q^{18} t^{11}
+q^{14} t^{12}+2 q^{16} t^{12}+q^{22} t^{12}+2 q^{18} t^{13}
+2 q^{20} t^{13}+q^{16} t^{14}+2 q^{18} t^{14}+q^{24} t^{14}
+q^{26} t^{14}+2 q^{20} t^{15}+3 q^{22} t^{15}+q^{18} t^{16}
+3 q^{20} t^{16}+q^{26} t^{16}+q^{28} t^{16}+2 q^{22} t^{17}
+4 q^{24} t^{17}+3 q^{22} t^{18}+q^{28} t^{18}+2 q^{30} t^{18}
+6 q^{26} t^{19}+3 q^{24} t^{20}+q^{26} t^{20}+3 q^{32} t^{20}
+6 q^{28} t^{21}+q^{30} t^{21}+2 q^{26} t^{22}+2 q^{28} t^{22}
+4 q^{34} t^{22}+6 q^{30} t^{23}+2 q^{32} t^{23}+2 q^{28} t^{24}
+3 q^{30} t^{24}+4 q^{36} t^{24}+q^{38} t^{24}+5 q^{32} t^{25}
+4 q^{34} t^{25}+q^{40} t^{25}+q^{30} t^{26}+3 q^{32} t^{26}
+4 q^{38} t^{26}+2 q^{40} t^{26}+4 q^{34} t^{27}+5 q^{36} t^{27}
+q^{42} t^{27}+q^{32} t^{28}+3 q^{34} t^{28}+2 q^{40} t^{28}
+3 q^{42} t^{28}+2 q^{36} t^{29}+6 q^{38} t^{29}+q^{44} t^{29}
+3 q^{36} t^{30}+q^{42} t^{30}+3 q^{44} t^{30}+q^{38} t^{31}
+5 q^{40} t^{31}+q^{48} t^{31}+2 q^{38} t^{32}+q^{40} t^{32}
+2 q^{46} t^{32}+5 q^{42} t^{33}+q^{40} t^{34}+q^{42} t^{34}
+2 q^{48} t^{34}+3 q^{44} t^{35}+q^{46} t^{35}+q^{42} t^{36}
+q^{44} t^{36}+q^{50} t^{36}+2 q^{46} t^{37}+q^{48} t^{37}
+q^{46} t^{38}+q^{52} t^{38}
+q^{48} t^{39}+q^{50} t^{39}+q^{48} t^{40}+q^{52} t^{41}.
\)
}
\renewcommand{\baselinestretch}{1.2} 
\medskip

{\em Colored polynomials}.
Let as provide several examples of colored polynomials. 
Presumably they can be used in the ongoing research
toward colored Khovanov-Rozansky theory.

Recall that we use the
{\em standard parameters}
from (\ref{qtareli}),(\ref{qtarel}).
\smallskip

The only exceptional DAHA-KhR-polynomial for the
trefoil $\{3,2\}$ and $b=3\om_1$ is

$$
K\!R^2_{3,2}(b=3\om_1;\,q,t)\ =
$$
\renewcommand{\baselinestretch}{0.5} 
{\small

\noindent
\(
1+q^8 t^6+q^{10} t^7+q^{10} t^8+q^{12} t^9+q^{12} t^{10}
+q^{14} t^{11}+q^{16} t^{12}+q^{18} t^{13}+q^{18} t^{14}
+2 q^{20} t^{15}+q^{20} t^{16}+q^{22} t^{16}+2 q^{22} t^{17}
+q^{24} t^{18}+q^{28} t^{21}+q^{30} t^{22}
+q^{30} t^{23}+q^{32} t^{24}+q^{34} t^{26}+q^{36} t^{27}.
\)
}
\renewcommand{\baselinestretch}{1.2} 
Actually the first such polynomial is
$K\!R^1_{3,2}(b=3\om_1;\,q,t,a)=1$, but
we disregard it. 

\smallskip
The next case will be for $\{4,3\}$ and $2\om_2$ 
(the $2\times 2$\~square);
the stabilization begins at $n=15$. The first nontrivial exceptional
DAHA-KR polynomial is:
$$
K\!R^3_{4,3}(b=2\om_2;\,q,t)\ =
$$
\renewcommand{\baselinestretch}{0.5} 
{\small
\noindent
\(
1+q^6 t^4+q^{10} t^5+q^8 t^6+q^{12} t^7+q^{10} t^8+q^{12} t^8
+q^{14} t^9+q^{16} t^9+q^{12} t^{10}+q^{14} t^{10}+q^{16} t^{11}
+2 q^{18} t^{11}+2 q^{16} t^{12}+q^{18} t^{12}+q^{22} t^{12}
+3 q^{20} t^{13}+q^{22} t^{13}+q^{18} t^{14}+q^{20} t^{14}
+q^{24} t^{14}+3 q^{22} t^{15}+2 q^{24} t^{15}+q^{20} t^{16}
+2 q^{22} t^{16}+2 q^{26} t^{16}+q^{28} t^{16}+q^{24} t^{17}
+4 q^{26} t^{17}+2 q^{24} t^{18}+2 q^{30} t^{18}+4 q^{28} t^{19}
+2 q^{28} t^{20}+3 q^{32} t^{20}+q^{34} t^{20}+2 q^{30} t^{21}
+3 q^{32} t^{21}+q^{36} t^{21}+q^{30} t^{22}+2 q^{34} t^{22}
+2 q^{36} t^{22}+q^{32} t^{23}+3 q^{34} t^{23}+q^{38} t^{23}
+4 q^{38} t^{24}+q^{34} t^{25}+q^{42} t^{25}+q^{40} t^{26}
+q^{38} t^{27}+q^{44} t^{27}+q^{42} t^{28}+q^{40} t^{29}
+q^{48} t^{29}+2 q^{44} t^{30}
+q^{48} t^{31}+q^{50} t^{31}+q^{46} t^{32}+q^{52} t^{32}.
\)
}
\renewcommand{\baselinestretch}{1.2} 
\smallskip

The next example is the DAHA-Khovanov polynomial
for $\{4,3\}$ and $b=4\om_1$. Here the degree of $a$
is $8$ (as predicted by Conjecture \ref{CONJEVAL}) and
there is a lot of room for possible cancelations of
monomials apart from those between neighboring degrees
of $a$, the only type which is allowed in the construction
of $H\!D^{[n]}$. Nevertheless Conjecture \ref{CONJKHR},$(i)$
holds. The stabilization to the super-polynomial 
occurs at $n=7$. 

$$
K\!R^2_{3,4}(b=4\om_1;\,q,t)\ =
$$
\renewcommand{\baselinestretch}{0.5} 
{\small
\noindent
\(
1+q^{10} t^8+q^{12} t^9+q^{12} t^{10}+q^{14} t^{11}
+q^{14} t^{12}+q^{16} t^{13}+q^{16} t^{14}+q^{18} t^{15}
+q^{18} t^{16}+q^{22} t^{17}+q^{20} t^{18}+2 q^{24} t^{19}
+q^{22} t^{20}+q^{24} t^{20}+q^{26} t^{20}+3 q^{26} t^{21}
+q^{24} t^{22}+q^{26} t^{22}+q^{28} t^{22}+4 q^{28} t^{23}
+3 q^{28} t^{24}+2 q^{30} t^{24}+4 q^{30} t^{25}+3 q^{30} t^{26}
+q^{32} t^{26}+5 q^{32} t^{27}+4 q^{32} t^{28}+q^{34} t^{28}
+q^{36} t^{28}+5 q^{34} t^{29}+q^{36} t^{29}+3 q^{34} t^{30}
+2 q^{38} t^{30}+4 q^{36} t^{31}+2 q^{38} t^{31}+3 q^{36} t^{32}
+4 q^{40} t^{32}+2 q^{38} t^{33}+5 q^{40} t^{33}+q^{42} t^{33}
+2 q^{38} t^{34}+5 q^{42} t^{34}+8 q^{42} t^{35}+q^{44} t^{35}
+2 q^{40} t^{36}+3 q^{42} t^{36}+7 q^{44} t^{36}+12 q^{44} t^{37}
+q^{46} t^{37}+q^{42} t^{38}+4 q^{44} t^{38}+8 q^{46} t^{38}
+13 q^{46} t^{39}+q^{48} t^{39}+q^{44} t^{40}+5 q^{46} t^{40}
+8 q^{48} t^{40}+q^{50} t^{40}+13 q^{48} t^{41}+q^{52} t^{41}
+6 q^{48} t^{42}+7 q^{50} t^{42}+2 q^{52} t^{42}+12 q^{50} t^{43}
+2 q^{54} t^{43}+5 q^{50} t^{44}+3 q^{52} t^{44}+6 q^{54} t^{44}
+8 q^{52} t^{45}+3 q^{54} t^{45}+3 q^{56} t^{45}+4 q^{52} t^{46}
+10 q^{56} t^{46}+5 q^{54} t^{47}+9 q^{56} t^{47}+5 q^{58} t^{47}
+3 q^{54} t^{48}+16 q^{58} t^{48}+q^{60} t^{48}+14 q^{58} t^{49}
+6 q^{60} t^{49}+2 q^{56} t^{50}+2 q^{58} t^{50}+18 q^{60} t^{50}
+18 q^{60} t^{51}+7 q^{62} t^{51}+q^{58} t^{52}+5 q^{60} t^{52}
+21 q^{62} t^{52}+18 q^{62} t^{53}+6 q^{64} t^{53}+q^{66} t^{53}
+q^{60} t^{54}+4 q^{62} t^{54}+16 q^{64} t^{54}+16 q^{64} t^{55}
+3 q^{66} t^{55}+2 q^{68} t^{55}+5 q^{64} t^{56}+12 q^{66} t^{56}
+4 q^{68} t^{56}+q^{70} t^{56}+11 q^{66} t^{57}+6 q^{70} t^{57}
+3 q^{66} t^{58}+2 q^{68} t^{58}+11 q^{70} t^{58}+q^{72} t^{58}
+6 q^{68} t^{59}+4 q^{70} t^{59}+9 q^{72} t^{59}+2 q^{68} t^{60}
+18 q^{72} t^{60}+2 q^{74} t^{60}+2 q^{70} t^{61}+11 q^{72} t^{61}
+13 q^{74} t^{61}+q^{70} t^{62}+q^{72} t^{62}+22 q^{74} t^{62}
+3 q^{76} t^{62}+q^{72} t^{63}+14 q^{74} t^{63}+14 q^{76} t^{63}
+q^{72} t^{64}+2 q^{74} t^{64}+22 q^{76} t^{64}+3 q^{78} t^{64}
+13 q^{76} t^{65}+12 q^{78} t^{65}+2 q^{76} t^{66}+17 q^{78} t^{66}
+2 q^{80} t^{66}+9 q^{78} t^{67}+7 q^{80} t^{67}+q^{82} t^{67}
+q^{78} t^{68}+9 q^{80} t^{68}+q^{84} t^{68}+4 q^{80} t^{69}
+6 q^{84} t^{69}+q^{80} t^{70}+3 q^{82} t^{70}+3 q^{84} t^{70}
+2 q^{86} t^{70}+2 q^{82} t^{71}+10 q^{86} t^{71}+5 q^{84} t^{72}
+12 q^{86} t^{72}+4 q^{88} t^{72}+5 q^{84} t^{73}+2 q^{86} t^{73}
+12 q^{88} t^{73}+q^{84} t^{74}+3 q^{86} t^{74}+11 q^{88} t^{74}
+4 q^{90} t^{74}+4 q^{86} t^{75}+3 q^{88} t^{75}+12 q^{90} t^{75}
+q^{92} t^{75}+2 q^{88} t^{76}+10 q^{90} t^{76}+5 q^{92} t^{76}
+2 q^{88} t^{77}+q^{90} t^{77}+7 q^{92} t^{77}+q^{94} t^{77}
+q^{90} t^{78}+3 q^{92} t^{78}+q^{94} t^{78}+q^{90} t^{79}
+q^{92} t^{79}+2 q^{94} t^{79}+q^{92} t^{80}+q^{94} t^{80}
+q^{98} t^{80}+q^{96} t^{81}+3 q^{96} t^{82}+2 q^{100} t^{82}
+q^{96} t^{83}+5 q^{98} t^{83}+q^{100} t^{83}+5 q^{98} t^{84}
+3 q^{102} t^{84}+q^{98} t^{85}+5 q^{100} t^{85}+q^{102} t^{85}
+q^{104} t^{85}+4 q^{100} t^{86}+2 q^{104} t^{86}+4 q^{102} t^{87}
+2 q^{102} t^{88}+q^{106} t^{88}+3 q^{104} t^{89}+q^{110} t^{93}
+q^{110} t^{94}
+q^{112} t^{94}+q^{112} t^{95}+q^{114} t^{97}+q^{116} t^{98}.
\)
}
\renewcommand{\baselinestretch}{1.2} 
\smallskip

In the case of $4$\~column, dual to the one
considered above, the simplest DAHA-KR 
polynomial is:
$$
K\!R^5_{4,3}(b=\om_4;\,q,t)\ =
$$
\renewcommand{\baselinestretch}{0.5} 
{\small
\noindent
\(
1+q^4 t^2+q^{12} t^3+q^6 t^4+q^8 t^4+q^{14} t^5
+q^{16} t^5+q^{12} t^6+q^{18} t^7+q^{20} t^7+q^{26} t^8.
\)
}
\renewcommand{\baselinestretch}{1.2} 
\smallskip

It is much smaller than the one above; the stabilization
begins here only at $n=28$. 

\medskip

\setcounter{equation}{0}
\section{Kauffman polynomials etc.}\label{sec:Kauffman}

The purpose of this section is to extend the conjectures
above to the case of $C$-$D$ and, partially, to $B$.
The corresponding super-polynomials
will be called {\em hyper-polynomials} to distinguish them from
the super-polynomials. They depend
on 4 parameters and hardly have any positivity properties
even for minuscule weights.

We note that the duality is expected to hold in this case
and in the $C^\vee C$\~theory, which seems a natural
setting for the hyper-polynomials. It is worth mentioning 
that the hyper-polynomials are determined using infinitely
many relations with the DAHA-Jones polynomials; so if they
exist (we conjecture that it is always the case), their 
coefficients are completely rigid. 

This is different from
quite a few papers on super-polynomials (mainly physical),
where  only finitely many relations of evaluation type
(the differentials) are involved in the process
of calculation of the corresponding
super-polynomials. If there are no established restrictions
for the degree of $a$, the positivity 
can be always achieved by increasing this degree and
the construction becomes non-rigorous.
 
The $a$\~degree of DAHA
super-polynomials and hyper-polynomials are expected to
satisfy the evaluation Conjecture \ref{CONJEVAL} and 
its counterpart below. This alone is sufficient to see that 
the positivity cannot hold in certain examples; 
cf. Section \ref{sec:3hook}.
\smallskip

\subsection{Quasi-minuscule weights}\label{sec:quasim}
The definitions from Section \ref{sect:Aut} are quite sufficient
to manage the minuscule weights, which are
\begin{align*}
&\{\om_j,\,1\le j\le n\,\} \for A_n,\ \,\om_n \for B_n,\ \,
\om_1 \for C_n,\\
&\{\om_1,\om_n,\om_{n-1}\}
\for D_n (n\ge 3),\ \{\om_1,\om_6\} \for E_6,\ \om_7 \for E_7.
\end{align*} 

We note that $\om_j$ and $\om_{n-j}$ for $A_n$,
$\om_n$ and $\om_{n-1}$ for $D_n$ and
$\om_1$ and $\om_6$ for $E_6$ always give the same
invariants of torus knots due to the existence of the
automorphism of the corresponding
Dynkin diagram (preserving $\al_0$) 
transposing them.

For theoretical and practical calculations,
we use directly formulas (\ref{tauplus}),(\ref{taumin}) and
the following important observation.
The polynomials $P_b/P_b(q^{-\rho_k})$ in (\ref{jones-d})
can be replaced by  $E_b/E_b(q^{-\rho_k})$ for the
nonsymmetric Macdonald polynomial $E_b$ for $b\in P_+$.
It is because the $t$\~symmetrization of $E_b$ gives $P_b$
and the coinvariant we use is $t$\~symmetric. 
See \cite{C4} and references therein or \cite{C101},
formula (3.3.14). For minuscule $b=\om_j$, one
readily has: $E_b=X_b$. This makes the calculations
reasonably straightforward.

Let us explain how to manage the case 
of the quasi-minuscule weight $b=\vth$, which is
used in the examples considered below.
Recall that this weight is the maximal
{\em short} root in this paper (the twisted setting), 
namely,
\begin{align*}
&\vth=\om_1 \,(B_n, n\ge 2),\ \vth=\om_2 \for 
C_n (n\ge 3),\, D_n (n\ge 4),\, E_6,\\ 
&\vth=\om_1\, (E_7),\ \,\vth=\om_8\, (E_8),\ \,
\vth=\om_4\, (F_4),\ \, \vth=\om_1\, (G_2).
\end{align*}

We disregard the case of $A$ here and below in
this section. Let 
$$ 
t=t_{\sht}= q_{\sht}^{k_{\sht}}=q^{k_{\sht}},\ \, 
u=t_{\lng} =q^{\nu_{\lng}k_{\lng}}.
$$
In the $A,D,E$ cases, the roots are treated
as short and $u=1$. Recall that $\nu_{\sht}=1$.

The corresponding $E$\~polynomial is given by
the following uniform formula:
\begin{align}\label{Evth}
&E_{\vth}\ =\ X_{\vth}+\frac{q(1-t)}{1-q^{1+(\vth,\rho_k)}},\,
\where\\
&q^{(\vth,\rho_k)}\ =\ X_{\vth}(q^{\rho_k})\ =\ 
t^{(\rho_{\sht},\vth)}u^{(\rho_{\lng},\vth)/\nu_{\lng})},\notag\\
&E_{\vth}(q^{-\rho_k})\ =\  q^{-(\vth,\rho_k)}
\frac{(1-t q^{1+(\vth,\rho_k)})}{(1-q t q^{1+(\vth,\rho_k)})}.\notag
\end{align}
We note that $Y_{\vth}(E_{\vth})=
q^{-2} q^{-(\vth,\rho_k)} E_{\vth}$.
These formulas are a simple application of the technique 
of intertwiners.

Thus it suffices to calculate the action of $\tga$ on $X_{\vth}$
in this case. Let us give an example (which will be used below).
Using (\ref{tauplus}), (\ref{taumin}), we obtain
the following recurrence relations. Setting 
$\tilde{X}=q^{-1}X_\vth$ and $\ka=t^{1/2}-t^{-1/2}$, one has:

\begin{align}\label{tautheta}
q^{-1}\tau_+^a\tau_-^m(X_{\vth})=
\bigl(\tilde{X} Y_\vth 
-&\ka (\sum_{j=1}^a \tilde{X}^j)T_{s_\vth}\bigr)
(\tau_+^a\tau_-^{m-1})(X_{\vth})\\
+&\ka\bigr(t^{-1/2}\tilde{X}^a-\ka
(\sum_{j=1}^{a-1} \tilde{X}^j)\bigr).
\notag
\end{align}

The summations over empty sets of indices are naturally zero
in this formula; $m\ge 0$ and  $\tau_+^a(X_{\vth})=X_{\vth}$.
Since this formula will be used inside the coinvariant
(see (\ref{jones-d})), one can
replace $\tilde{X}$ and $T_{s_\vth}$ by their
evaluations:
$$
\{\,\tilde{X}\,\}=\tilde{X}(t^{-\rho_k})=q^{-1}q^{-(\vth,\rho_k)},\ 
\, \{\,T_{s_\vth}\,\}=t^{-1/2}q^{(\vth,\rho_k)}.
$$

We applied here the following general formula valid
for any $\al\in R_+$ such that either this root is short or
the set $\la(s_{al})$ from (\ref{xlambda})
contains no short roots:
\begin{align}\label{tstheta}
T_{s_{\al}}^{-1}X_{\al}^{-1}=X_{\al}T_{s_{\al}}+ 
(t_\al^{1/2}-t_\al^{-1/2}).
\end{align}
It can be checked directly using the 
description of $\la(s_\al)$. We need it only for
$\al=\vth$, where it readily follows from the
existence of the anti-involution $\vph$.
Indeed, 
$$
\vph(T_0)=\vph(Y_\vth T_{s_\vth}^{-1})=
T_{s_\vth}^{-1} X_{\vth}^{-1}.
$$
Thus the latter satisfies the same quadratic relation
as $T_0$, which is exactly relation (\ref{tstheta}) for
$\al=\vth$.
\smallskip

\subsection{Hyper-polynomials}
Using these and similar formulas, we
calculated quite a few examples, resulting
in the following extension of Conjecture 
\ref{HOMFLY} to the case of Kauffman polynomials. 

The root system will be $C_n$; the notations from
\cite{Bo} will be used. Let
$b=\sum_{i=1}^n c_i\ep_i$ with $c_i\in \Z_+$
be a dominant weight for all $C_m$
with $m\ge n\ge 2$.

\begin{conjecture}[Hyper-polynomials]\label{KAUF}
Given a knot $K=K_{r,s}$ for $r,s\in \Z$,
there exists a polynomial
$\h^{C}_{r,s}(b\,;\,q,t,u,a)$ with integral
coefficients  in terms of positive powers of 
$a,q,u$ and $t^{\pm 1}$ such that
for all $m\ge n$,
\begin{align}\label{hyper-b-c}
&\h^{C}_{r,s}(b\,;\,q,t,u,a=-t^{m-1})=
\tilde{J\!D}^{C_m}_{r,s}(b\,;\,q,t,u).
\end{align}
Moreover, for any $4\le m\ge n$,
\begin{align}\label{kauf}
&\h^{C}_{r,s}(b\,;\,q,t,u=1,a=-t^{m-1})\ =\ 
\tilde{J\!D}^{D_m}_{r,s}(b\,;\,q,t,u) \\
&\and \h^{C}_{r,s}(\om_1\,;\,q\mapsto \tilde{q}^2, 
t\mapsto \tilde{q}^2, u\mapsto 1,a\mapsto -\la /\tilde{q})
\label{kauf-d}
\end{align}
coincides with the reduced Kauffman polynomial 
$\tilde{\k}_{r,s}(\la,\tilde{q})$
in the standard notations upon the division
by the smallest power of $\la$ and
applying the tilde- normalization
(with respect to $\k(\la=0,\tilde{q})$).
\end{conjecture}
 
Using the substitution $a=-t^{m-1}$ here instead of 
$a=-t^{m+1}$ in the case of
$A_m$ is not very significant, but it removes 
the denominators (powers of $t$) in the formulas for 
some simple knots. Following the notations in this
section, $H\!D_{r,s}$ can be now denoted by $\h^A_{r,s}$,
but the above change of the substitutions for $a$
must be taken into consideration. See also (\ref{a-3-d}).

Due to (\ref{kauf-d}), we can define the 
{\em $D$\~hyper-polynomial} 
as the following specialization of the one for $C$:
\begin{align}\label{kauffm-d}
&\h^{D}_{r,s}(b\,;\,q,a)\,\equal\,
\h^{C}_{r,s}(b\,;\,q,t,u=1,a).
\end{align}

See, e.g., \cite{Kau} concerning the Kauffman polynomials.
In our setting, they appear from the
$D$\~hyper-polynomials when $\,t=q\,$ and $b=\om_1$.
Concerning using $\,\la,\tilde{q}\,$ in the
Kauffman polynomials, it matches the (reduced)
formulas from \cite{GW},\cite{GS}. For instance, 
(6.12) from \cite{GW} in the tilde- normalization
coincides with $\tilde{\k}_{3,2}(\la,\tilde{q}).$

More generally,
\begin{align}\label{GWsub}
\h^{D}_{r,s}(\om_1\,;\ q\mapsto -\tilde{q}^{\,2}\, \tilde{t},\ 
t\mapsto -\tilde{q}^{\,2}\, \tilde{t},\ a\mapsto 
\la\, \tilde{t}\,/\,\tilde{q}\,)
\end{align}
coincides with the reduced  
super-polynomial from (6.14) and that for the torus knot 
$\{5,2\}$ from Table 3 in \cite{GW}.
In these notations, the passage
to the Kauffman polynomials is $\,\tilde{t}=-1$.

For instance,
\begin{align}
\label{super-c-2-3}
&\h^{C}_{3,2}(\om_1\,;\,q,t,u,a)\ =\\ 
1+q t+a (q t-q u)&+a^2 \bigl(-q u+q^2 u-q^2 t u\bigr)
+a^3 \bigl(-q^2 t u+q^2 u^2\bigr),\notag\\
\label{super-d-2-3}
&\h^{D}_{3,2}(\om_1\,;\,q,t,u=1,a)\ =\\ 
1+q t+a (q t-q)&+a^2 \bigl(-q+q^2-q^2 t \bigr)
+a^3 \bigl(-q^2 t +q^2\bigr),\notag
\end{align}
and the latter results in $q^2 \f(3_1)=$
\renewcommand{\baselinestretch}{0.5} 
{\small
\noindent
$$
1+q^4 t^2 +\frac{\la t}{q} (q^4 t^2+q^2 t)+(\frac{\la t}{q})^2 
(q^2t+q^4t^2
+q^6 t^3) +(\frac{\la t}{q})^3 (q^6 t^3 +q^4 t^2)
$$
}
\renewcommand{\baselinestretch}{1.2} 

\noindent
from (6.14), \cite{GW} (we omit the bars). Let us provide
another example:
$$\h^{C}_{7,2}(\om_1\,;\,q,t,u,a)=$$

\renewcommand{\baselinestretch}{0.5} 
{\small
\noindent
\( 
1 + q t + q^2 t^2 + q^3 t^3 + 
a \bigl(q t + q^2 t^2 + q^3 t^3 - q u - q^2 t u - 
   q^3 t^2 u\bigr) + 
a^2 \bigl(-q u + q^2 u - 2 q^2 t u + q^3 t u - 2 q^3 t^2 u +
   q^4 t^2 u - q^4 t^3 u\bigr) + 
a^3 \bigl(-q^2 t u + q^3 t u - q^3 t^2 u + 
   q^4 t^2 u - q^4 t^3 u + q^2 u^2 - q^3 u^2 + q^3 t u^2 - 
q^4 t u^2 + q^4 t^2 u^2\bigr) + 
a^4 \bigl(-q^3 u^2 + q^4 u^2 + q^3 t u^2 - 2 q^4 t u^2 + 
   q^5 t u^2 + q^4 t^2 u^2 - q^5 t^2 u^2\bigr) + 
a^5 \bigl(-q^4 t u^2 + q^5 t u^2 - q^5 t^2 u^2 + q^4 u^3 - 
q^5 u^3 + q^5 t u^3\bigr) + 
a^6 \bigl(-q^5 u^3 + q^6 u^3 + q^5 t u^3 - q^6 t u^3\bigr) + 
a^7 \bigl(-q^6 t u^3 + q^6 u^4\bigr).
\)
}
\renewcommand{\baselinestretch}{1.2} 

Upon the substitution
$t\mapsto -t,q\mapsto -q$, $\h^{C}_{2n+1,2}(\om_1\,;\,q,t,u,a)$
are positive (in examples), but this fails for $\{4,3\}$.
Correspondingly, our $\h^{D}_{4,3}(\om_1\,;\,q,t,a)$ is very different
from that predicted in B.0.1 from \cite{GS} (symbol $a$ is used
there instead of $\la$), though both match the
Kauffman polynomial for the torus knot $\{4,3\}$ ($=8_{19}$).

The general formula (any $u$) in this case is as follows:
$$\h^{C}_{4,3}(\om_1\,;\,q,t,u,a)=$$

\renewcommand{\baselinestretch}{0.5} 
{\small
\noindent
\( 
1+q t+q^2 t+q^2 t^2+q^3 t^3
+a \bigl(-q u-q^2 u+q t+q^2 t-q^2 u t-q^3 u t+q^2 t^2
+q^3 t^2-q^3 u t^2+q^3 t^3\bigr)+a^2 \bigl(-q u+q^4 u
+q^3 u^2-2 q^2 u t-2 q^3 u t+q^3 t^2-2 q^3 u t^2-q^4 u t^3\bigr)
+a^3 \bigl(q^2 u^2+q^3 u^2-q^2 u t-q^3 u t
+q^3 u^2 t+q^4 u^2 t-q^3 u t^2-q^4 u t^2+q^4 u^2 t^2-q^4 u t^3\bigr)
+a^4 
\bigl(-q^4 u^3+q^3 u^2 t+q^4 u^2 t-q^4 u t^2+q^4 u^2 t^2\bigr).
\)
}
\renewcommand{\baselinestretch}{1.2} 

The $a$\~range here matches the conjecture below and
is smaller by $2$ than the range of $a$ 
in the super-polynomial suggested in \cite{GS}.

Moreover, there are negative terms in our formula. 
The best we can do to eliminate them algebraically is to use
the substitution $u\mapsto -u$, which makes all terms positive 
but one, namely $a^2q^4 u$. Certainly, the 
positivity of our hyper-polynomials can hold only for very 
simple knots.
\smallskip

\subsection{Duality}
Let us formulate counterparts of the conjectures on the 
DAHA super-polynomials in the $C$-$D$\~case.
We hope that it can be naturally
extended to the $C^\vee C$\~case (with $5$ parameters
instead of $q,t$). If it is true, then the corresponding
hyper-polynomials is expected to generalize those
for all classical root systems.

We identify the weights $b$ with the corresponding
Young diagrams and understand $b^{tr}$ as a
dominant weight for sufficiently large $C$\~systems
exactly in the same way as it was done in the $A$ case.

\begin{conjecture} [Hyper-duality]\label{HDUALIT}
(i) Upon multiplication by a proper monomial in
terms of $q,t,u$, the hyper-polynomial
$\h^{C}_{r,s}(b\,;\,q,t,u,a)$  coincides
with its dual counterpart
$$
\h^{C}_{r,s}(b^{tr}\,;\,t^{-1},q^{-1},u^{-1}t/q, -a q u).
$$
This substitution of the parameters is compatible with the 
specializations $u=-t$ and $u=-q^{-1}$, resulting in
$(at)\mapsto (at)$ and $a\mapsto a$ correspondingly.

(ii) The evaluation $\h^{C}_{r,s}(b\,;\,q,t=1,u,a)$
is the product of the
evaluations $\h^C_{r,s}(b_i;q,t=1,u,a)$ calculated
for $b_i=m_i\om_1$, where $m_i$ is the number of boxes in the
$i$-th row of the Young diagram for $b$. 
The maximal power of $a$ in  $\h^C_{r,s}$ is $\,(ord\times r)\,$
for the order of the Young diagram, where $1<s<r$.
\end{conjecture} 

We note that the symmetry from the conjecture
matches the symmetry 
$$ \tilde{q}\,\mapsto\, \tilde{q}^{-1},\ \la\,\mapsto\, -\la
$$ 
of the classical reduced Kauffman polynomials
$\tilde{\k}_{r,s}(\la,\tilde{q})$, where 
$u=1,t=q$ in our notations. Indeed, the switch of
the parameter in $(i)$ is compatible with making 
$u=1,t=q,b=\om_1$; it sends $\tilde{q}=q^{1/2}$ to 
$\tilde{q}^{-1}$ and
$$
\la=-a\sqrt{q}\ \mapsto\  -(-aq)\sqrt{q}^{\,-1}=-\la.
$$
\smallskip

The following is the simplest example of the colored
(non-minuscule) hyper- polynomial of type $C$:

$$\h^{C}_{3,2}(\om_2\,;\,q,t,u, a )=$$
\renewcommand{\baselinestretch}{0.5} 
{\small
\noindent
\( 
1+q t+q t^2+q^2 t^4+a \bigl(q+q t+q^2 t^2+q^2 t^3-q u
-\frac{q u}{t}-q^2 t u-q^2 t^2 u\bigr)+a^2 \bigl(q^2 t
-3 q^2 u-\frac{q u}{t^2}+\frac{q^2 u}{t^2}-\frac{q u}{t}
+\frac{q^3 u}{t}-2 q^2 t u+q^3 t u-q^3 t^2 u-q^3 t^3 u
+\frac{q^2 u^2}{t}\bigr)+a^3 \bigl(-q^2 u-\frac{q^2 u}{t^2}
-\frac{2 q^2 u}{t}+\frac{q^3 u}{t}-2 q^3 t u-q^3 t^2 u
+2 q^3 u^2+\frac{q^2 u^2}{t^3}+\frac{2 q^2 u^2}{t^2}
-\frac{q^3 u^2}{t^2}+\frac{q^2 u^2}{t}+q^3 t u^2\bigr)
+a^4 \bigl(-q^3 u-\frac{q^3 u}{t}+q^3 u^2-q^4 u^2
+\frac{q^2 u^2}{t^3}-\frac{q^3 u^2}{t^3}+\frac{q^3 u^2}{t^2}
+\frac{q^4 u^2}{t^2}+\frac{3 q^3 u^2}{t}-\frac{q^4 u^2}{t}
+q^4 t u^2-\frac{q^3 u^3}{t^3}-\frac{q^3 u^3}{t^2}\bigr)
+a^5 \bigl(q^4 u^2+\frac{q^3 u^2}{t^3}+\frac{q^3 u^2}{t^2}
-\frac{q^4 u^2}{t^2}-\frac{q^3 u^3}{t^4}-\frac{q^3 u^3}{t^3}
+\frac{q^4 u^3}{t^3}-\frac{q^4 u^3}{t}\bigr)
+a^6 \bigl(\frac{q^4 u^2}{t^2}-\frac{q^4 u^3}{t^3}
-\frac{q^5 u^3}{t^3}-\frac{q^4 u^3}{t^2}
+\frac{q^5 u^3}{t^2}+\frac{q^4 u^4}{t^4}\bigr).
\)
}
\renewcommand{\baselinestretch}{1.2} 

Its value at $t=1$ is 
$$
(1+q+a(q-q u)-a^2 q u+a^3 (-q^2 u +q^2 u^2))^2,
$$
the square of the evaluation at $t=1$ of

$$\h^{C}_{3,2}(\om_1\,;\,q,t,u, a )=$$
\renewcommand{\baselinestretch}{0.5} 
{\small
\noindent
$$
1+q t+a (q t-q u)+a^2 \bigl(-q u+q^2 u-q^2 t u\bigr)
+a^3 \bigl(-q^2 t u+q^2 u^2\bigr).
$$
}
\renewcommand{\baselinestretch}{1.2} 

Recall that when $u=1$,
we arrive at the reduced tilde-normalized
Kauffman super-polynomial for $\{3,2\}$
predicted in \cite{GW} upon the substitution 
from (\ref{GWsub}).

\smallskip

\subsection{The B-case}
The counterparts of $C$\~hyper-polynomials 
exist in the $B$ case, but have somewhat weaker
symmetries. We continue to assume that 
$b=\sum_{i=1}^n c_i\ep_i$ with $c_i\in \Z_+$
is a dominant weight in all $B_m$ starting
with $n\ge 2$.

\begin{conjecture}\label{KAUFB}
Given a knot $K=K_{r,s}$ for $r,s\in \Z$,
there exist a polynomial 
$\h^{B}_{r,s}(b\,;\,q,t,u,a)$ with integral
coefficients  in terms of positive powers of 
$a,q,t$ and $u^{\pm 1}$, such that
for all $m\ge n$,
\begin{align}\label{hyper-b}
&\h^{B}_{r,s}(b\,;\,q,t,u,a=-u^{m-1})=
\tilde{J\!D}^{B_m}_{r,s}(b\,;\,q,t,u).
\end{align}
Moreover, 
\begin{align}\label{kauf-b}
&\h^{B}_{r,s}(b\,;\,\sqrt{q},1,u,a)=
\h^{D}_{r,s}(b\,;\,q,t,a)=
\h^{C}_{r,s}(b\,;\,q,t,u=1,a)
\end{align}
and
\begin{align}\label{b-dualt} 
\h^{B}_{r,s}(b\,;\,q,t,u,a)\ =\ 
\h^{B}_{r,s}(b\,;\,q,t\mapsto u/(tq),u, a\mapsto a q t^2/u ).
\end{align}
\end{conjecture}

The reduction to the $D$\~hyper -polynomials and
the coincidence, 
\begin{align}\label{b-s-rel}
\h^{B}_{r,s}(b\,;\,q,1,u,a)\ =\ 
\h^{C}_{r,s}(b\,;\,q^2,t\mapsto u,1,a) 
\end{align}
are actually predictable. For instance,
it can be justified for $b=\om_1$. Notice
that $\om_1$ is minuscule for $C$, but becomes
quasi-minuscule for $B$, which requires more 
involved DAHA treatment.
It explains why the hyper-polynomials
of type $B$ for $\om_1$ are significantly longer
than those calculated for $C$.  

Let us comment on the isomorphisms among hyper-polynomials
in small ranks:
\begin{align*}
&\tilde{J\!D}^{B_2}_{r,s}(\om_1\,;\,q,t)\ =\ 
\tilde{J\!D}^{C_2}_{r,s}(\om_2,;\,q,t,u),\\ 
&\tilde{J\!D}^{B_2}_{r,s}(\om_2\,;\,q,t)\ =\ 
\tilde{J\!D}^{C_2}_{r,s}(\om_1,;\,q,t,u),
\end{align*}
which follows from the isomorphism of these
root systems. Similarly,
\begin{align*}
&\tilde{J\!D}^{D_3}_{r,s}(\om_1\,;\,q,t)\ =\ 
\tilde{J\!D}^{A_3}_{r,s}(\om_2,;\,q,t),\\ 
\end{align*}
though the DAHA super-polynomials $H\!D_{r,s}(\om_2)$ 
are very different 
from the hyper-polynomials $\h^D_{r,s}(\om_1)$.
The degree of the former is expected $2 (s-1)$,
when the degree of the latter is expected no greater 
than $r$.
For instance,
see (\ref{super3-2+2-0}) and (\ref{super-d-2-3}),
\begin{align}\label{a-3-d}
H\!D_{3,2}(\om_2)=& 1 + (a^2q^2)/t + qt + qt^2 + q^2t^4 
+ a(q + q/t + q^2t + q^2t^2),\notag\\
\h^D_{3,2}(\om_1)=& 
1 + qt + a(qt - q) + a^2(-q + q^2 - q^2t) + 
a^3(-q^2t + q^2)\,;\notag\\
&H\!D_{3,2}(\om_2\,;\,q,t,a=-t^4)\ = \ 
\h^D_{3,2}(\om_1\,;\,q,t,a=-t^2). 
\end{align}

We hope that the symmetries of the hyper-polynomials
of type $B,C$ discussed in this paper
will find a uniform interpretation in the
$C^\vee C$\~theory, as well as our prediction
for the maximal $a$\~degrees in type $C$ and other
relations.
\smallskip

The simplest example of the hyper-polynomial
of type $B$ is
$$\h^{B}_{3,2}(\om_1\,;\,q,t,u, a )=$$
\renewcommand{\baselinestretch}{0.5} 
{\small
\noindent
\( 
1+q^2 u+a \bigl(-q^2 t-q^3 t+q^3 t^2+q^2 u\bigr)
+a^2 \bigl(-q^2 t^2+q^4 t^2-q^4 t^2 u\bigr)
+a^3 \bigl(q^4 t^3+q^5 t^3-q^5 t^4-q^4 t^2 u\bigr).
\)
}
\renewcommand{\baselinestretch}{1.2} 

Its value at $t=1$  readily coincides
with 
$$
\h^{C}_{3,2}(\om_1\,;\,q^2,u,1, a )=
1 + q^2 u + aq^2(u - 1) +a^2q^2(q^2-1- q^2 u) + 
 a^3q^4(1-u),
$$
where
$\h^{C}_{3,2}(\om_1\,;\,q,t,u, a )$ 
is provided in (\ref{super-d-2-3}).

At level of formulas, the reduction from $B$
to $D$ occurs due to many reductions. For instance,

$$\h^{B}_{4,3}(\om_1\,;\,q,t,u, a )=$$
\renewcommand{\baselinestretch}{0.5} 
{\small
\noindent
\( 
 1+q^2 u+q^4 u+q^4 u^2+q^6 u^3+
 a \bigl(-q^2 t-q^3 t-q^4 t
-q^5 t+q^3 t^2+q^5 t^2+q^2 u+q^4 u-q^4 t u-q^5 t u-q^6 t u
-q^7 t u+q^5 t^2 u+q^7 t^2 u+q^4 u^2+q^6 u^2-q^6 t u^2
-q^7 t u^2+q^7 t^2 u^2+q^6 u^3\bigr)
+a^2 \bigl(-q^2 t^2
+q^5 t^2+2 q^6 t^2+q^7 t^2+q^8 t^2-q^5 t^3-q^6 t^3-q^7 t^3
-q^8 t^3+q^8 t^4-q^4 t u-q^5 t u-q^6 t u-q^7 t u-q^4 t^2 u
+q^5 t^2 u-q^6 t^2 u+2 q^7 t^2 u+q^8 t^2 u-q^7 t^3 u-q^8 t^3 u
+q^6 u^2-q^6 t u^2-q^7 t u^2-q^6 t^2 u^2+q^7 t^2 u^2
-q^8 t^2 u^3\bigr)
+a^3 \bigl(q^4 t^3+q^5 t^3+q^6 t^3
+q^7 t^3-q^8 t^3-q^9 t^3-q^5 t^4-q^7 t^4+q^8 t^4+q^9 t^4
-q^4 t^2 u-q^6 t^2 u+q^7 t^2 u+q^8 t^2 u+q^6 t^3 u+q^9 t^3 u
-q^7 t^4 u-q^9 t^4 u-q^6 t^2 u^2
-q^8 t^2 u^2+q^8 t^3 u^2+q^9 t^3 u^2-q^9 t^4 u^2-q^8 t^2 u^3\bigr)
+a^4 \bigl(-q^7 t^4-2 q^8 t^4
-q^9 t^4+q^7 t^5+q^8 t^5+q^9 t^5+q^{10} t^5-q^{10} t^6
+q^6 t^3 u+q^7 t^3 u+q^8 t^3 u+q^9 t^3 u-q^7 t^4 u-2 q^9 t^4 u
-q^{10} t^4 u+q^9 t^5 u+q^{10} t^5 u-q^8 t^2 u^2+q^8 t^3 u^2
+q^9 t^3 u^2-q^9 t^4 u^2\bigr)
+a^5 \bigl(q^{10} t^5
+q^{11} t^5-q^{10} t^6-q^{11} t^6-q^9 t^4 u-q^{10} t^4 u
+q^9 t^5 u+q^{10} t^5 u\bigr).
\)
}
\renewcommand{\baselinestretch}{1.2} 

Its reduction 
$t\mapsto 1, u\mapsto t, q^2\mapsto q$ 
from $B$ to $D$ results in the cancelation
of many terms and diminishing the
$a$\~degree:

$$\h^{D}_{4,3}(\om_1\,;\,q,t,a )=$$
\renewcommand{\baselinestretch}{0.5} 
{\small
\noindent
\( 
1+q t+q^2 t+q^2 t^2+q^3 t^3+a^4 \bigl(-q^4+q^3 t+q^4 t\bigr)
+a \bigl(-q-q^2+q t-q^3 t+q^2 t^2+q^3 t^3\bigr)
+a^2 \bigl(-q+q^3+q^4-2 q^2 t-2 q^3 t-q^3 t^2-q^4 t^3\bigr)
+a^3 \bigl(q^2+q^3-q^2 t+q^4 t-q^3 t^2-q^4 t^3\bigr).
\)
}
\renewcommand{\baselinestretch}{1.2}

\smallskip

\subsection{Exceptional root systems}\label{sec:escqg}
We will conclude this paper with the list of
the DAHA-Jones counterparts of QG-Jones invariants
for the simplest torus knots for 
the exceptional root systems. 
They are obtained from the DAHA-Jones polynomials
when $t=q$; see Conjecture \ref{MAINCONJ}.
Only minuscule weights (when present) and the
quasi-minuscule $\vth$ will be considered.
We use the approach and formulas described in Section 
\ref{sec:quasim}. 

In the case of trefoil, a direct identification 
of the DAHA formulas with the corresponding  
{\em reduced} quantum group invariants seems doable.
Based on extensive calculations with the classical
root systems and $G_2$, we expect that the coincidence
always holds.
The symbol $\tilde{\j^R}$ will be used in this section
for the DAHA counterparts of the reduced
tilde-normalized QG-Jones invariants. They are obtained
as specialization $t\mapsto q$ of the corresponding
$q,t$\~polynomials $J\!D$; correspondingly, for the 
root system $F_4$,
$r\mapsto q^2$.

We adjust the formulas listed below to
the setup of the program QuantumKnotInvariant[$\,q\,$]
from http://katlas.org/wiki/ upon changing the parameter
$q$ there to our $q^{-1/2}$. For instance, our
one reads:

$$\tilde{\j}^{G_2}_{4,3}(\om_1\,;\,q)=$$
\renewcommand{\baselinestretch}{0.5} 
{\small
\noindent
\( 
1+q^2+q^7-q^8+q^9+q^{12}-q^{13}+q^{14}-q^{15}
+q^{16}-q^{17}-q^{20}+q^{21}-2 q^{22}+q^{23}-q^{24}-q^{27}
+q^{28}-q^{29}+q^{30}-q^{31}+q^{35}-q^{36}+q^{37},
\)
}
\renewcommand{\baselinestretch}{1.2} 
\smallskip

where the program QuantumKnotInvariant[$\,q\,$] gives:

\smallskip
\renewcommand{\baselinestretch}{0.5} 
{\small
\noindent
\( 
\bigl(1-q+q^2-q^3+q^4-q^5+q^6\bigr) 
\bigl(1+q+q^2+q^3+q^4+q^5+q^6\bigr) \bigl(1-q^4+q^8\bigr)\times
\bigl(1-q^2+q^4-q^{12}+q^{14}-q^{16}+q^{18}-q^{20}-q^{26}
+q^{28}-2 q^{30}+q^{32}-q^{34}-q^{40}+q^{42}-q^{44}+q^{46}-q^{48}
+q^{50}+q^{56}-q^{58}+q^{60}+q^{70}+q^{74}\bigr)\times
\frac{1}{q^{144}}.
\)
}
\renewcommand{\baselinestretch}{1.2} 

Up to a proper power of $q$ (disregarded in the 
tilde-normalization),
the product of the first 3 factors is the $q$\~dimension of the
irreducible module with the weight $\om_1$ for $G_2$; 
when $q=1$, it is $7$. The longest factor coincides with our
expression upon $q\mapsto q^{-1/2}$ and the tilde-normalization.

\subsubsection{E6}
We will begin with  
$$\tilde{\j}^{E_6}_{3,2}(\om_1\,;\,q)=
1+q^2+q^5-q^{13}-q^{15}-q^{18}+q^{23},
$$
corresponding to
$$\tilde{J\!D}^{E_6}_{3,2}(\om_1\,;\,q,t)=
1 + q t + q t^4 + q^2 t^8 - q t^9 - q t^{12} - q^2 t^{13} - 
q^2 t^{16} + 
 q^2 t^{21}.
$$

Recall that the substitution to $\j$ from
$J\!D$ is $t\mapsto q$.

The next ones we provide here are:

$$\tilde{\j}^{E_6}_{4,3}(\om_1\,;\,q)=$$
\renewcommand{\baselinestretch}{0.5} 
{\small
\noindent
\( 
1+q^2+q^3+q^4+q^5+2 q^6+q^7+q^8+q^9+q^{12}-q^{13}
-q^{14}-q^{15}-2 q^{16}-2 q^{17}-q^{18}-3 q^{19}-2 q^{20}-q^{21}
-2 q^{22}-q^{23}-q^{26}+2 q^{27}+q^{29}+2 q^{30}+q^{31}+q^{32}
+q^{33}+q^{34}+q^{36}-q^{38}+q^{39}-q^{40}-q^{44}-q^{47}+q^{48},
\)
}
\renewcommand{\baselinestretch}{1.2} 
corresponding to
$$\tilde{J\!D}^{E_6}_{4,3}(\om_1\,;\,q,t)=$$
\renewcommand{\baselinestretch}{0.5} 
{\small
\noindent
\( 
1+q t+q^2 t+q^2 t^2+q^3 t^3+q t^4+q^2 t^4+q^2 t^5+q^3 t^5+q^3 t^6
+q^2 t^8+q^3 t^8+q^4 t^8-q t^9-q^2 t^9+q^3 t^9+q^4 t^9
-q^2 t^{10}-q^3 t^{10}+q^4 t^{10}-q^3 t^{11}-q t^{12}-q^2 t^{12}
+q^3 t^{12}+q^4 t^{12}-2 q^2 t^{13}-3 q^3 t^{13}+q^5 t^{13}
-2 q^3 t^{14}-q^4 t^{14}-q^4 t^{15}-q^2 t^{16}-2 q^3 t^{16}
+q^5 t^{16}-2 q^3 t^{17}-3 q^4 t^{17}+q^3 t^{18}-2 q^4 t^{18}
-q^5 t^{18}-q^3 t^{20}-2 q^4 t^{20}+q^2 t^{21}+2 q^3 t^{21}
-q^4 t^{21}-3 q^5 t^{21}+q^3 t^{22}+2 q^4 t^{22}+q^4 t^{23}
+q^3 t^{24}-q^4 t^{24}-2 q^5 t^{24}+q^6 t^{24}+q^3 t^{25}
+3 q^4 t^{25}-q^6 t^{25}+q^4 t^{26}+2 q^5 t^{26}
+q^4 t^{28}-q^6 t^{28}+q^4 t^{29}+3 q^5 t^{29}
-q^4 t^{30}+q^6 t^{30}+q^5 t^{32}-q^6 t^{32}-q^4 t^{33}
+2 q^6 t^{33}-q^5 t^{34}-q^6 t^{34}
+q^6 t^{36}-q^5 t^{37}-q^6 t^{38}-q^6 t^{41}+q^6 t^{42},
\)
}
\renewcommand{\baselinestretch}{1.2} 
$$\tilde{\j}^{E_6}_{3,2}(\om_2\,;\,q)=$$
\renewcommand{\baselinestretch}{0.5} 
{\small
\noindent
\( 
1+q^2+q^4+q^5+q^8-q^{12}-2 q^{14}-q^{15}-q^{16}
-q^{20}-q^{22}+q^{24}+3 q^{26}-q^{29}+q^{34}-q^{36}-q^{38}+q^{39},
\)
}
\renewcommand{\baselinestretch}{1.2} 
corresponding to

$$\tilde{J\!D}^{E_6}_{3,2}(\om_2\,;\,q,t)=$$
\renewcommand{\baselinestretch}{0.5} 
{\small
\noindent
\( 
1+q t+q t^3+q t^4+q^2 t^6+q^2 t^7-q t^8+q^2 t^8-q t^9
-q t^{11}+q^3 t^{11}-3 q^2 t^{12}-q^2 t^{13}-q^2 t^{14}
+q^3 t^{14}-q^2 t^{15}-q^3 t^{16}+q^2 t^{17}-q^3 t^{17}
-q^3 t^{18}+q^2 t^{19}-2 q^3 t^{19}+q^2 t^{20}+q^3 t^{21}
+q^4 t^{22}+2 q^3 t^{23}-q^4 t^{23}+q^3 t^{24}-q^4 t^{25}
+q^4 t^{27}-q^3 t^{28}+q^4 t^{30}-q^4 t^{32}-q^5 t^{33}+q^5 t^{34},
\)
}
\renewcommand{\baselinestretch}{1.2} 

$$\tilde{\j}^{E_6}_{4,3}(\om_2\,;\,q)=$$
\renewcommand{\baselinestretch}{0.5} 
{\small
\noindent
\( 
1+q^2+q^3+2 q^4+2 q^5+3 q^6+2 q^7+3 q^8+2 q^9
+3 q^{10}+2 q^{11}+2 q^{12}+q^{13}-q^{14}-q^{15}-3 q^{16}
-4 q^{17}-4 q^{18}-5 q^{19}-6 q^{20}-6 q^{21}-8 q^{22}
-6 q^{23}-5 q^{24}-3 q^{25}-2 q^{26}-q^{27}+q^{28}+2 q^{29}
+2 q^{30}+5 q^{31}+5 q^{32}+7 q^{33}+6 q^{34}+5 q^{35}
+5 q^{36}+4 q^{37}+3 q^{38}+2 q^{39}-q^{42}-3 q^{43}-2 q^{44}
-3 q^{45}-2 q^{46}-2 q^{47}-4 q^{48}-2 q^{49}-q^{50}-q^{51}
+q^{52}+2 q^{54}+q^{58}+q^{59}+q^{60}-q^{64}-q^{66}+q^{67},
\)
}
\renewcommand{\baselinestretch}{1.2} 

corresponding to
$$\tilde{J\!D}^{E_6}_{4,3}(\om_2\,;\,q,t)=$$
\renewcommand{\baselinestretch}{0.5} 
{\small
\noindent
\( 
1+q t+q^2 t+q^2 t^2+q t^3+q^2 t^3+q^3 t^3+q t^4+2 q^2 t^4+q^3 t^4
+q^2 t^5+2 q^3 t^5+q^2 t^6+2 q^3 t^6+q^4 t^6+q^2 t^7+3 q^3 t^7
+2 q^4 t^7-q t^8+2 q^3 t^8+3 q^4 t^8-q t^9-2 q^2 t^9+q^3 t^9
+3 q^4 t^9-q^2 t^{10}-q^3 t^{10}+4 q^4 t^{10}+2 q^5 t^{10}
-q t^{11}-q^2 t^{11}-q^3 t^{11}+3 q^4 t^{11}+3 q^5 t^{11}
-4 q^2 t^{12}-5 q^3 t^{12}+q^4 t^{12}+3 q^5 t^{12}+q^6 t^{12}
-q^2 t^{13}-6 q^3 t^{13}-3 q^4 t^{13}+4 q^5 t^{13}-q^2 t^{14}
-3 q^3 t^{14}-4 q^4 t^{14}+4 q^5 t^{14}+2 q^6 t^{14}-q^2 t^{15}
-5 q^3 t^{15}-7 q^4 t^{15}+2 q^5 t^{15}+2 q^6 t^{15}-3 q^3 t^{16}
-10 q^4 t^{16}-3 q^5 t^{16}+2 q^6 t^{16}+q^2 t^{17}-7 q^4 t^{17}
-5 q^5 t^{17}+3 q^6 t^{17}+q^7 t^{17}+q^3 t^{18}-6 q^4 t^{18}
-9 q^5 t^{18}+4 q^6 t^{18}+q^7 t^{18}+q^2 t^{19}-4 q^4 t^{19}
-13 q^5 t^{19}-q^6 t^{19}+q^7 t^{19}+q^2 t^{20}+4 q^3 t^{20}
-11 q^5 t^{20}-3 q^6 t^{20}+q^7 t^{20}+3 q^3 t^{21}+7 q^4 t^{21}
-7 q^5 t^{21}-6 q^6 t^{21}+2 q^7 t^{21}+q^3 t^{22}+4 q^4 t^{22}
-3 q^5 t^{22}-9 q^6 t^{22}+q^7 t^{22}+q^8 t^{22}+3 q^3 t^{23}
+6 q^4 t^{23}-13 q^6 t^{23}-q^7 t^{23}+q^3 t^{24}+8 q^4 t^{24}
+11 q^5 t^{24}-10 q^6 t^{24}-2 q^7 t^{24}+3 q^4 t^{25}
+10 q^5 t^{25}-4 q^6 t^{25}-3 q^7 t^{25}+q^8 t^{25}+2 q^4 t^{26}
+10 q^5 t^{26}-2 q^6 t^{26}-7 q^7 t^{26}+2 q^4 t^{27}+11 q^5 t^{27}
+9 q^6 t^{27}-8 q^7 t^{27}-q^8 t^{27}-q^3 t^{28}-q^4 t^{28}
+7 q^5 t^{28}+13 q^6 t^{28}-4 q^7 t^{28}+q^8 t^{28}-3 q^4 t^{29}
+q^5 t^{29}+13 q^6 t^{29}-5 q^7 t^{29}-q^8 t^{29}-q^5 t^{30}
+12 q^6 t^{30}+q^7 t^{30}-3 q^8 t^{30}-2 q^4 t^{31}-q^5 t^{31}
+11 q^6 t^{31}+8 q^7 t^{31}-2 q^8 t^{31}-2 q^4 t^{32}-7 q^5 t^{32}
+2 q^6 t^{32}+11 q^7 t^{32}-2 q^8 t^{32}-4 q^5 t^{33}-5 q^6 t^{33}
+8 q^7 t^{33}-2 q^8 t^{33}+q^9 t^{33}-2 q^5 t^{34}-3 q^6 t^{34}
+12 q^7 t^{34}-2 q^8 t^{34}-q^9 t^{34}-3 q^5 t^{35}-9 q^6 t^{35}
+7 q^7 t^{35}+5 q^8 t^{35}-q^5 t^{36}-8 q^6 t^{36}-3 q^7 t^{36}
+3 q^8 t^{36}+q^5 t^{37}-4 q^6 t^{37}-3 q^7 t^{37}+5 q^8 t^{37}
-q^9 t^{37}-2 q^6 t^{38}-7 q^7 t^{38}+7 q^8 t^{38}-2 q^6 t^{39}
-10 q^7 t^{39}+2 q^8 t^{39}-q^9 t^{39}+q^5 t^{40}+2 q^6 t^{40}
-7 q^7 t^{40}-2 q^8 t^{40}+3 q^6 t^{41}-q^7 t^{41}-q^8 t^{41}
+q^9 t^{41}-2 q^7 t^{42}-4 q^8 t^{42}+2 q^9 t^{42}+q^6 t^{43}
+q^7 t^{43}-8 q^8 t^{43}+q^9 t^{43}+q^6 t^{44}+5 q^7 t^{44}
-2 q^8 t^{44}+2 q^7 t^{45}-q^8 t^{45}+2 q^9 t^{45}+q^7 t^{46}
-q^8 t^{46}-3 q^9 t^{46}+q^7 t^{47}+3 q^8 t^{47}-2 q^9 t^{47}
-q^{10} t^{47}+3 q^8 t^{48}+q^{10} t^{48}-q^7 t^{49}+q^8 t^{49}
-q^9 t^{49}+q^{10} t^{49}+q^8 t^{50}-q^9 t^{50}-q^{10} t^{50}
+q^8 t^{51}+3 q^9 t^{51}-q^8 t^{52}+q^9 t^{52}+q^{10} t^{52}
-q^8 t^{53}-q^9 t^{53}-q^{10} t^{53}+q^9 t^{54}-q^{10} t^{54}
+q^{10} t^{55}-q^{11} t^{55}-q^9 t^{56}+q^{11} t^{56}.
\)
}
\renewcommand{\baselinestretch}{1.2}

\subsubsection{E7} 
$$\tilde{\j}^{E_7}_{3,2}(\om_7\,;\,q)=$$
\renewcommand{\baselinestretch}{0.5} 
{\small
\noindent
\( 
1+q^2+q^6+q^{10}-q^{11}+q^{12}-q^{15}+q^{16}-q^{17}-q^{19}
+q^{20}-2 q^{21}-2 q^{25}+q^{26}-q^{29}+2 q^{30}-q^{31}
+q^{34}-q^{35}+q^{36}-q^{39}+q^{40}+q^{44}-q^{45},
\)
}
\renewcommand{\baselinestretch}{1.2} 

corresponding to
$$\tilde{J\!D}^{E_7}_{3,2}(\om_7\,;\,q,t)=$$
\renewcommand{\baselinestretch}{0.5} 
{\small
\noindent
\( 
1+q t+q t^5+q t^9-q t^{10}+q^2 t^{10}-q t^{14}+q^2 t^{14}
-q^2 t^{15}-q t^{18}+q^2 t^{18}-2 q^2 t^{19}-2 q^2 t^{23}
+q^2 t^{24}-q^2 t^{27}+q^3 t^{27}+q^2 t^{28}-q^3 t^{28}
+q^2 t^{32}-q^3 t^{32}+q^3 t^{33}-q^3 t^{36}+q^3 t^{37}
+q^3 t^{41}-q^3 t^{42},
\)
}
\renewcommand{\baselinestretch}{1.2} 

$$\tilde{\j}^{E_7}_{4,3}(\om_7\,;\,q)=$$
\renewcommand{\baselinestretch}{0.5} 
{\small
\noindent
\( 
1+q^2+q^3+q^4+2 q^6+q^7+q^8+q^9+2 q^{10}+q^{12}+q^{13}+2 q^{14}
-q^{15}+q^{16}-2 q^{19}+q^{20}-2 q^{21}-2 q^{22}-2 q^{23}
-4 q^{25}-2 q^{26}-2 q^{27}-2 q^{28}-4 q^{29}-2 q^{31}-2 q^{32}
-2 q^{33}+q^{34}-2 q^{35}+q^{37}+q^{38}-q^{39}+2 q^{40}+2 q^{41}
+q^{42}+q^{43}+3 q^{44}+q^{45}+q^{46}+2 q^{47}+2 q^{48}+q^{50}
+q^{51}+q^{54}-q^{56}-q^{59}-q^{60}-q^{62}-q^{63}-q^{66}+q^{71},
\)
}
\renewcommand{\baselinestretch}{1.2} 

corresponding to
$$\tilde{J\!D}^{E_7}_{4,3}(\om_7\,;\,q,t)=$$
\renewcommand{\baselinestretch}{0.5} 
{\small
\noindent
\( 
1+q t+q^2 t+q^2 t^2+q^3 t^3+q t^5+q^2 t^5+q^2 t^6+q^3 t^6+q^3 t^7
+q t^9+q^2 t^9-q t^{10}+q^2 t^{10}+2 q^3 t^{10}+q^4 t^{10}
-q^2 t^{11}+q^3 t^{11}+q^4 t^{11}-q^3 t^{12}+q^4 t^{12}-q t^{14}
+2 q^3 t^{14}+q^4 t^{14}-2 q^2 t^{15}-q^3 t^{15}+q^4 t^{15}
-2 q^3 t^{16}+q^4 t^{16}+q^5 t^{16}-q^4 t^{17}-q t^{18}+q^3 t^{18}
+q^4 t^{18}-3 q^2 t^{19}-3 q^3 t^{19}+q^4 t^{19}+q^5 t^{19}
-3 q^3 t^{20}-q^4 t^{20}+q^5 t^{20}-3 q^4 t^{21}-2 q^2 t^{23}
-3 q^3 t^{23}+q^5 t^{23}+q^2 t^{24}-2 q^3 t^{24}-4 q^4 t^{24}
+q^3 t^{25}-3 q^4 t^{25}-q^5 t^{25}+q^4 t^{26}-q^5 t^{26}
-q^2 t^{27}-q^3 t^{27}+q^5 t^{27}+q^6 t^{27}+q^2 t^{28}-q^3 t^{28}
-5 q^4 t^{28}-2 q^5 t^{28}+2 q^3 t^{29}-q^4 t^{29}-2 q^5 t^{29}
+2 q^4 t^{30}-2 q^5 t^{30}+q^2 t^{32}-3 q^4 t^{32}-2 q^5 t^{32}
+3 q^3 t^{33}+3 q^4 t^{33}-2 q^5 t^{33}+3 q^4 t^{34}-q^5 t^{34}
-q^6 t^{34}+2 q^5 t^{35}-q^4 t^{36}-q^5 t^{36}+2 q^3 t^{37}
+3 q^4 t^{37}-q^5 t^{37}-q^6 t^{37}+2 q^4 t^{38}+2 q^5 t^{38}
-q^6 t^{38}+3 q^5 t^{39}+q^3 t^{41}+2 q^4 t^{41}-q^6 t^{41}
-q^3 t^{42}+q^4 t^{42}+4 q^5 t^{42}-q^4 t^{43}+2 q^5 t^{43}
+q^6 t^{43}-q^5 t^{44}+q^6 t^{44}-q^6 t^{45}+3 q^5 t^{46}
+q^6 t^{46}-q^4 t^{47}-q^5 t^{47}+q^6 t^{47}-q^5 t^{48}
+q^6 t^{48}+q^5 t^{50}+q^6 t^{50}
-q^4 t^{51}-2 q^5 t^{51}+q^6 t^{51}-q^5 t^{52}-q^6 t^{53}
-q^5 t^{55}-q^6 t^{56}-q^6 t^{57}-q^6 t^{60}+q^6 t^{65},
\)
}
\renewcommand{\baselinestretch}{1.2} 

$$\tilde{\j}^{E_7}_{3,2}(\om_1\,;\,q)=$$
\renewcommand{\baselinestretch}{0.5} 
{\small
\noindent
\( 
1+q^2+q^5+q^7+q^{10}+q^{12}-q^{13}+q^{14}-q^{15}-2 q^{18}
+q^{19}-2 q^{20}-q^{22}-q^{23}+q^{24}-2 q^{25}+q^{26}
-q^{27}-q^{30}+q^{31}-2 q^{32}+2 q^{33}-q^{34}+q^{35}+q^{36}
-q^{37}+3 q^{38}-q^{39}+q^{40}-q^{42}+q^{43}-q^{44}
+q^{45}-q^{46}+q^{50}-q^{51}+q^{52}-q^{53}-q^{56}+q^{57},
\)
}
\renewcommand{\baselinestretch}{1.2} 

corresponding to
$$\tilde{J\!D}^{E_7}_{3,2}(\om_1\,;\,q,t)=$$
\renewcommand{\baselinestretch}{0.5} 
{\small
\noindent
\( 
1+q t+q t^4+q t^6+q^2 t^8+q^2 t^{10}-q t^{12}+q^2 t^{12}-q t^{14}
-q^2 t^{16}-q t^{17}+q^2 t^{17}+q^3 t^{17}-3 q^2 t^{18}-q^2 t^{20}
-q^2 t^{21}+q^3 t^{21}-q^3 t^{22}-q^2 t^{23}+q^3 t^{23}-q^3 t^{24}
-q^3 t^{25}+q^2 t^{26}-q^3 t^{27}+q^2 t^{29}-2 q^3 t^{29}+q^3 t^{30}
+q^2 t^{31}-q^3 t^{31}+q^3 t^{32}+q^3 t^{33}-q^3 t^{34}
+q^4 t^{34}+2 q^3 t^{35}-q^4 t^{35}+q^3 t^{37}
-q^4 t^{38}+q^4 t^{39}-q^4 t^{40}+q^4 t^{41}-q^3 t^{43}
+q^4 t^{46}-q^4 t^{47}+q^4 t^{48}-q^4 t^{49}-q^5 t^{51}+q^5 t^{52},
\)
}
\renewcommand{\baselinestretch}{1.2}

$$\tilde{\j}^{E_7}_{4,3}(\om_1\,;\,q)=$$
\renewcommand{\baselinestretch}{0.5} 
{\small
\noindent
\( 
1+q^2+q^3+q^4+q^5+2 q^6+2 q^7+2 q^8+2 q^9+2 q^{10}
+2 q^{11}+3 q^{12}+2 q^{13}+3 q^{14}+q^{15}+2 q^{16}
+q^{17}+q^{18}-q^{20}-q^{21}-2 q^{22}-3 q^{23}-4 q^{24}
-4 q^{25}-5 q^{26}-5 q^{27}-5 q^{28}-7 q^{29}-5 q^{30}
-6 q^{31}-5 q^{32}-5 q^{33}-5 q^{34}-2 q^{35}-4 q^{36}
-2 q^{38}+q^{39}+2 q^{40}+2 q^{41}+5 q^{42}+2 q^{43}
+6 q^{44}+5 q^{45}+5 q^{46}+7 q^{47}+3 q^{48}+8 q^{49}
+3 q^{50}+6 q^{51}+3 q^{52}+2 q^{53}+4 q^{54}-q^{55}
+3 q^{56}-2 q^{57}-q^{58}-4 q^{60}-5 q^{62}-q^{63}
-3 q^{64}-3 q^{65}-q^{66}-5 q^{67}-3 q^{69}-q^{71}
-2 q^{72}+2 q^{73}-2 q^{74}+2 q^{75}+3 q^{78}-q^{79}
+2 q^{80}-q^{81}+q^{82}+q^{83}+q^{85}-q^{86}
+q^{87}-q^{91}+q^{92}-q^{93}-q^{96}+q^{97},
\)
}
\renewcommand{\baselinestretch}{1.2} 
\smallskip

which corresponds to
$$\tilde{J\!D}^{E_7}_{4,3}(\om_1\,;\,q,t)=$$

\renewcommand{\baselinestretch}{0.5} 
{\small
\noindent
\( 
1+q t+q^2 t+q^2 t^2+q^3 t^3+q t^4+q^2 t^4+q^2 t^5+q^3 t^5+q t^6
+q^2 t^6+q^3 t^6+q^2 t^7+q^3 t^7+q^2 t^8+2 q^3 t^8+q^4 t^8
+q^3 t^9+q^4 t^9+q^2 t^{10}+2 q^3 t^{10}+2 q^4 t^{10}+q^3 t^{11}
+q^4 t^{11}-q t^{12}+2 q^3 t^{12}+3 q^4 t^{12}-q^2 t^{13}
+2 q^4 t^{13}+q^5 t^{13}-q t^{14}-q^2 t^{14}+3 q^4 t^{14}
+q^5 t^{14}-q^2 t^{15}-q^3 t^{15}+q^4 t^{15}+q^5 t^{15}
-q^2 t^{16}-2 q^3 t^{16}+2 q^4 t^{16}+2 q^5 t^{16}-q t^{17}
+q^4 t^{17}+3 q^5 t^{17}-4 q^2 t^{18}-4 q^3 t^{18}+3 q^5 t^{18}
+q^6 t^{18}-4 q^3 t^{19}-2 q^4 t^{19}+2 q^5 t^{19}-q^2 t^{20}
-3 q^3 t^{20}-5 q^4 t^{20}+2 q^5 t^{20}+q^6 t^{20}-q^2 t^{21}
-2 q^3 t^{21}-q^4 t^{21}+2 q^5 t^{21}+q^6 t^{21}-4 q^3 t^{22}
-7 q^4 t^{22}+2 q^6 t^{22}-q^2 t^{23}-q^3 t^{23}-3 q^4 t^{23}
+q^6 t^{23}-3 q^3 t^{24}-7 q^4 t^{24}-4 q^5 t^{24}+3 q^6 t^{24}
-q^3 t^{25}-5 q^4 t^{25}-3 q^5 t^{25}+2 q^6 t^{25}+q^7 t^{25}
+q^2 t^{26}+q^3 t^{26}-5 q^4 t^{26}-8 q^5 t^{26}+q^6 t^{26}
-2 q^4 t^{27}-5 q^5 t^{27}+2 q^6 t^{27}+q^7 t^{27}+q^3 t^{28}
-2 q^4 t^{28}-8 q^5 t^{28}-2 q^6 t^{28}+q^2 t^{29}-3 q^4 t^{29}
-7 q^5 t^{29}+q^6 t^{29}+2 q^7 t^{29}+3 q^3 t^{30}+2 q^4 t^{30}
-7 q^5 t^{30}-6 q^6 t^{30}+q^2 t^{31}+q^3 t^{31}+q^4 t^{31}
-6 q^5 t^{31}-3 q^6 t^{31}+2 q^7 t^{31}+3 q^3 t^{32}+5 q^4 t^{32}
-q^5 t^{32}-8 q^6 t^{32}-q^7 t^{32}+q^3 t^{33}+4 q^4 t^{33}
-2 q^5 t^{33}-5 q^6 t^{33}+2 q^7 t^{33}+3 q^4 t^{34}+4 q^5 t^{34}
-6 q^6 t^{34}-q^7 t^{34}+q^8 t^{34}+3 q^3 t^{35}+3 q^4 t^{35}
-q^5 t^{35}-9 q^6 t^{35}+6 q^4 t^{36}+9 q^5 t^{36}-4 q^6 t^{36}
-3 q^7 t^{36}+q^3 t^{37}+3 q^4 t^{37}+6 q^5 t^{37}-5 q^6 t^{37}
-3 q^7 t^{37}+2 q^4 t^{38}+7 q^5 t^{38}+2 q^6 t^{38}-q^7 t^{38}
+q^8 t^{38}+3 q^4 t^{39}+8 q^5 t^{39}-q^6 t^{39}-8 q^7 t^{39}
-q^8 t^{39}-q^4 t^{40}+5 q^5 t^{40}+7 q^6 t^{40}+q^7 t^{40}
+q^8 t^{40}+2 q^4 t^{41}+6 q^5 t^{41}+4 q^6 t^{41}-10 q^7 t^{41}
-q^8 t^{41}+4 q^5 t^{42}+9 q^6 t^{42}+q^7 t^{42}+q^8 t^{42}
-q^3 t^{43}-q^4 t^{43}+5 q^5 t^{43}+11 q^6 t^{43}-5 q^7 t^{43}
-2 q^8 t^{43}-2 q^4 t^{44}-2 q^5 t^{44}+7 q^6 t^{44}+3 q^7 t^{44}
+q^8 t^{44}-q^4 t^{45}+8 q^6 t^{45}-3 q^8 t^{45}-q^5 t^{46}
+4 q^6 t^{46}+2 q^7 t^{46}-2 q^4 t^{47}-2 q^5 t^{47}+7 q^6 t^{47}
+8 q^7 t^{47}-2 q^8 t^{47}-2 q^5 t^{48}+3 q^7 t^{48}-2 q^8 t^{48}
-2 q^4 t^{49}-4 q^5 t^{49}+q^6 t^{49}+12 q^7 t^{49}+q^8 t^{49}
-3 q^5 t^{50}-4 q^6 t^{50}+3 q^7 t^{50}-3 q^8 t^{50}-3 q^5 t^{51}
-6 q^6 t^{51}+6 q^7 t^{51}+4 q^8 t^{51}+q^9 t^{51}-q^6 t^{52}
+3 q^7 t^{52}-6 q^8 t^{52}-q^9 t^{52}-3 q^5 t^{53}-7 q^6 t^{53}
+5 q^7 t^{53}+8 q^8 t^{53}-3 q^6 t^{54}-q^7 t^{54}-q^8 t^{54}
-q^5 t^{55}-6 q^6 t^{55}-3 q^7 t^{55}+4 q^8 t^{55}-2 q^6 t^{56}
-2 q^7 t^{56}+4 q^8 t^{56}-q^9 t^{56}+q^5 t^{57}-3 q^6 t^{57}
-8 q^7 t^{57}+q^9 t^{57}+q^6 t^{58}-q^7 t^{58}+5 q^8 t^{58}
-q^9 t^{58}-q^6 t^{59}-7 q^7 t^{59}-2 q^8 t^{59}-q^9 t^{59}
-q^6 t^{60}-4 q^7 t^{60}+5 q^8 t^{60}+2 q^9 t^{60}+q^5 t^{61}
+2 q^6 t^{61}-5 q^7 t^{61}-7 q^8 t^{61}-2 q^9 t^{61}+q^6 t^{62}
+2 q^8 t^{62}+3 q^9 t^{62}+2 q^6 t^{63}-4 q^8 t^{63}-3 q^9 t^{63}
-q^7 t^{64}-2 q^8 t^{64}+4 q^9 t^{64}+q^6 t^{65}+3 q^7 t^{65}
-2 q^8 t^{65}-2 q^9 t^{65}-q^7 t^{66}-4 q^8 t^{66}+2 q^9 t^{66}
+q^6 t^{67}+4 q^7 t^{67}-q^8 t^{67}+q^7 t^{68}-2 q^8 t^{68}
-3 q^9 t^{68}+2 q^7 t^{69}+3 q^8 t^{69}+5 q^9 t^{69}-3 q^8 t^{70}
-5 q^9 t^{70}+q^7 t^{71}+4 q^8 t^{71}+3 q^9 t^{71}-q^8 t^{72}
-4 q^9 t^{72}-q^{10} t^{72}+3 q^8 t^{73}+2 q^9 t^{73}
+2 q^{10} t^{73}+q^8 t^{74}-2 q^9 t^{74}-2 q^{10} t^{74}
-q^7 t^{75}+q^8 t^{75}+2 q^9 t^{75}+2 q^{10} t^{75}-q^9 t^{76}
-q^{10} t^{76}-q^9 t^{77}+q^8 t^{78}+2 q^9 t^{78}-q^8 t^{79}
+q^9 t^{80}+q^{10} t^{80}-q^8 t^{81}-q^9 t^{81}-2 q^{10} t^{81}
+q^9 t^{82}+2 q^{10} t^{82}-q^9 t^{83}-2 q^{10} t^{83}
+q^9 t^{84}+q^{10} t^{84}-q^9 t^{85}-q^{11} t^{85}+q^{11} t^{86}.
\)
}
\renewcommand{\baselinestretch}{1.2}

\subsubsection{E8}
Here there are no minuscule weights. The simplest
DAHA-generated Jones polynomials are:

$$\tilde{\j}^{E_8}_{3,2}(\om_8\,;\,q)=$$
\renewcommand{\baselinestretch}{0.5} 
{\small
\noindent
\( 
1+q^2+q^7+q^{11}+q^{14}+q^{18}-q^{21}+q^{22}-q^{25}
-q^{28}-q^{30}+q^{31}-2 q^{32}-q^{36}-q^{37}+q^{38}
-q^{39}-q^{41}+q^{42}-q^{43}-q^{44}+q^{46}-q^{48}+q^{51}
-2 q^{52}+q^{53}+q^{55}-q^{56}+q^{57}+q^{58}-q^{61}
+3 q^{62}-q^{63}+q^{66}-q^{68}+q^{69}-q^{72}
+q^{73}-q^{76}+q^{82}-q^{83}+q^{86}-q^{87}-q^{92}+q^{93},
\)
}
\renewcommand{\baselinestretch}{1.2} 

which corresponds to
$$\tilde{J\!D}^{E_8}_{4,3}(\om_8\,;\,q,t)=$$

\renewcommand{\baselinestretch}{0.5} 
{\small
\noindent
\( 
1+q t+q t^6+q t^{10}+q^2 t^{12}+q^2 t^{16}-q t^{20}
+q^2 t^{20}-q t^{24}-q^2 t^{26}-q t^{29}+q^2 t^{29}
+q^3 t^{29}-3 q^2 t^{30}-q^2 t^{34}-q^2 t^{35}+q^3 t^{35}
-q^3 t^{36}-q^2 t^{39}+q^3 t^{39}-q^3 t^{40}-q^3 t^{41}
+q^2 t^{44}-q^3 t^{45}+q^2 t^{49}-2 q^3 t^{49}+q^3 t^{50}
+q^2 t^{53}-q^3 t^{53}+q^3 t^{54}+q^3 t^{55}-q^3 t^{58}
+q^4 t^{58}+2 q^3 t^{59}-q^4 t^{59}+q^3 t^{63}-q^4 t^{64}
+q^4 t^{65}-q^4 t^{68}+q^4 t^{69}-q^3 t^{73}
+q^4 t^{78}-q^4 t^{79}+q^4 t^{82}-q^4 t^{83}-q^5 t^{87}+q^5 t^{88},
\)
}
\renewcommand{\baselinestretch}{1.2}

$$\tilde{\j}^{E_8}_{4,3}(\om_8\,;\,q)=$$
\renewcommand{\baselinestretch}{0.5} 
{\small
\noindent
\( 
1+q^2+q^3+q^4+q^6+q^7+q^8+q^9+q^{10}+2 q^{11}+q^{12}+q^{13}
+2 q^{14}+2 q^{15}+2 q^{16}+q^{17}+2 q^{18}+2 q^{19}+2 q^{20}
+q^{21}+2 q^{22}+q^{23}+2 q^{24}+2 q^{26}+q^{27}-q^{29}+q^{31}
-2 q^{32}-q^{33}-2 q^{34}-q^{35}-3 q^{36}-3 q^{37}-2 q^{38}
-3 q^{39}-4 q^{40}-5 q^{41}-2 q^{42}-4 q^{43}-5 q^{44}-6 q^{45}
-3 q^{46}-4 q^{47}-4 q^{48}-5 q^{49}-3 q^{50}-3 q^{51}-6 q^{52}
-3 q^{53}-3 q^{54}-5 q^{56}-q^{57}-q^{58}+q^{59}-2 q^{60}
+3 q^{62}+2 q^{64}+q^{65}+6 q^{66}+2 q^{67}+2 q^{68}+3 q^{69}
+6 q^{70}+5 q^{71}+2 q^{72}+5 q^{73}+4 q^{74}+6 q^{75}+q^{76}
+6 q^{77}+4 q^{78}+5 q^{79}+4 q^{81}+5 q^{82}+2 q^{83}+q^{84}
+5 q^{86}-q^{87}+q^{88}-q^{89}+2 q^{90}-2 q^{91}-3 q^{92}-q^{94}
-5 q^{96}-2 q^{98}-q^{99}-5 q^{100}-q^{101}-q^{102}-3 q^{103}
-2 q^{104}-2 q^{105}+q^{106}-4 q^{107}-2 q^{108}-2 q^{109}
+2 q^{110}-3 q^{111}-q^{112}-2 q^{116}+2 q^{117}+q^{119}
-2 q^{120}+2 q^{121}+q^{122}+3 q^{126}-q^{127}+q^{128}+2 q^{130}
-q^{131}+q^{133}+q^{134}-q^{136}+q^{137}
-q^{140}+q^{141}-q^{147}+q^{150}-q^{151}-q^{156}+q^{157}.
\)
}
\renewcommand{\baselinestretch}{1.2} 
\smallskip

The corresponding $J\!D$\~polynomial is

$$\tilde{J\!D}^{E_8}_{4,3}(\om_8\,;\,q,t)=$$

\renewcommand{\baselinestretch}{0.5} 
{\small
\noindent
\( 
1+q t+q^2 t+q^2 t^2+q^3 t^3
+q t^6+q^2 t^6+q^2 t^7+q^3 t^7+q^3 t^8+q t^{10}+q^2 t^{10}
+q^2 t^{11}+q^3 t^{11}+q^2 t^{12}+2 q^3 t^{12}+q^4 t^{12}
+q^3 t^{13}+q^4 t^{13}+q^4 t^{14}+q^2 t^{16}+2 q^3 t^{16}
+q^4 t^{16}+q^3 t^{17}+q^4 t^{17}+q^3 t^{18}+2 q^4 t^{18}
+q^4 t^{19}+q^5 t^{19}-q t^{20}+q^3 t^{20}+q^4 t^{20}-q^2 t^{21}
+q^4 t^{21}+3 q^4 t^{22}+q^5 t^{22}+q^4 t^{23}+q^5 t^{23}
-q t^{24}-q^2 t^{24}+q^4 t^{24}+q^5 t^{24}-q^2 t^{25}
-q^3 t^{25}+q^5 t^{25}-q^2 t^{26}-2 q^3 t^{26}
+q^4 t^{26}+q^5 t^{26}-q^3 t^{27}+q^5 t^{27}+2 q^5 t^{28}-q t^{29}
+q^3 t^{29}+q^4 t^{29}+2 q^5 t^{29}-4 q^2 t^{30}-4 q^3 t^{30}
+2 q^5 t^{30}+q^6 t^{30}-4 q^3 t^{31}-2 q^4 t^{31}+q^5 t^{31}
-q^3 t^{32}-4 q^4 t^{32}+q^5 t^{32}+q^6 t^{32}-q^4 t^{33}
-q^2 t^{34}-2 q^3 t^{34}-q^4 t^{34}+q^5 t^{34}+q^6 t^{34}
-q^2 t^{35}-2 q^3 t^{35}+2 q^5 t^{35}+q^6 t^{35}-4 q^3 t^{36}
-7 q^4 t^{36}-q^5 t^{36}+2 q^6 t^{36}-4 q^4 t^{37}-2 q^5 t^{37}
-q^4 t^{38}-2 q^5 t^{38}+q^6 t^{38}-q^2 t^{39}-q^3 t^{39}
+q^4 t^{39}+q^5 t^{39}+q^6 t^{39}-3 q^3 t^{40}-6 q^4 t^{40}
-2 q^5 t^{40}+2 q^6 t^{40}-q^3 t^{41}-5 q^4 t^{41}-2 q^5 t^{41}
+2 q^6 t^{41}+q^7 t^{41}-4 q^4 t^{42}-7 q^5 t^{42}-4 q^5 t^{43}
-q^6 t^{43}+q^2 t^{44}+q^3 t^{44}-q^4 t^{44}-2 q^5 t^{44}
-2 q^4 t^{45}-q^5 t^{45}+3 q^6 t^{45}+q^7 t^{45}-2 q^4 t^{46}
-7 q^5 t^{46}-2 q^6 t^{46}-q^4 t^{47}-5 q^5 t^{47}+q^7 t^{47}
+q^3 t^{48}-4 q^5 t^{48}-4 q^6 t^{48}+q^2 t^{49}-2 q^4 t^{49}
-2 q^5 t^{49}+q^7 t^{49}+3 q^3 t^{50}+2 q^4 t^{50}-3 q^5 t^{50}
-3 q^6 t^{50}+2 q^4 t^{51}-4 q^5 t^{51}-q^6 t^{51}+2 q^7 t^{51}
-2 q^5 t^{52}-6 q^6 t^{52}-q^7 t^{52}+q^2 t^{53}+q^3 t^{53}
-q^4 t^{53}-3 q^5 t^{53}-3 q^6 t^{53}+q^7 t^{53}+3 q^3 t^{54}
+5 q^4 t^{54}+q^5 t^{54}-5 q^6 t^{54}-q^7 t^{54}+q^3 t^{55}
+4 q^4 t^{55}-q^5 t^{55}-3 q^6 t^{55}+q^7 t^{55}+3 q^4 t^{56}
+4 q^5 t^{56}-3 q^6 t^{56}-q^7 t^{56}+2 q^5 t^{57}-2 q^6 t^{57}
+q^7 t^{57}-2 q^6 t^{58}-q^7 t^{58}+q^8 t^{58}+3 q^3 t^{59}
+3 q^4 t^{59}-3 q^5 t^{59}-8 q^6 t^{59}-q^7 t^{59}+6 q^4 t^{60}
+9 q^5 t^{60}-q^6 t^{60}-2 q^7 t^{60}+q^4 t^{61}+6 q^5 t^{61}
-q^6 t^{61}-2 q^7 t^{61}+3 q^5 t^{62}+2 q^6 t^{62}-q^7 t^{62}
+q^3 t^{63}+2 q^4 t^{63}-3 q^6 t^{63}-2 q^7 t^{63}+2 q^4 t^{64}
+4 q^5 t^{64}+q^8 t^{64}+3 q^4 t^{65}+8 q^5 t^{65}-q^6 t^{65}
-8 q^7 t^{65}-q^8 t^{65}+6 q^5 t^{66}+9 q^6 t^{66}+q^7 t^{66}
+q^5 t^{67}+4 q^6 t^{67}-2 q^7 t^{67}-q^4 t^{68}-q^5 t^{68}
+q^6 t^{68}+2 q^7 t^{68}+q^8 t^{68}+2 q^4 t^{69}+5 q^5 t^{69}
-8 q^7 t^{69}-q^8 t^{69}+4 q^5 t^{70}+6 q^6 t^{70}
+q^8 t^{70}+3 q^5 t^{71}+8 q^6 t^{71}-3 q^7 t^{71}
-2 q^8 t^{71}+6 q^6 t^{72}+5 q^7 t^{72}-q^3 t^{73}
-q^4 t^{73}+2 q^5 t^{73}+4 q^6 t^{73}-q^7 t^{73}-2 q^4 t^{74}
-2 q^5 t^{74}+q^6 t^{74}-q^7 t^{74}+q^8 t^{74}+6 q^6 t^{75}
-3 q^8 t^{75}+4 q^6 t^{76}+4 q^7 t^{76}-q^4 t^{77}
+4 q^6 t^{77}+3 q^7 t^{77}-q^8 t^{77}-q^5 t^{78}
+q^8 t^{78}-2 q^4 t^{79}-2 q^5 t^{79}+4 q^6 t^{79}
+6 q^7 t^{79}-q^8 t^{79}-3 q^5 t^{80}-2 q^6 t^{80}
+q^7 t^{80}-2 q^8 t^{80}+5 q^7 t^{81}+q^5 t^{82}
+2 q^6 t^{82}+q^7 t^{82}-q^8 t^{82}-2 q^4 t^{83}
-4 q^5 t^{83}+q^6 t^{83}+9 q^7 t^{83}+2 q^8 t^{83}
-3 q^5 t^{84}-4 q^6 t^{84}+q^7 t^{84}-3 q^8 t^{84}
-2 q^5 t^{85}-5 q^6 t^{85}+4 q^7 t^{85}+3 q^8 t^{85}
-3 q^6 t^{86}-2 q^7 t^{86}-q^8 t^{86}-q^5 t^{87}-q^6 t^{87}
+q^7 t^{87}+2 q^8 t^{87}+q^9 t^{87}+2 q^6 t^{88}+5 q^7 t^{88}
-5 q^8 t^{88}-q^9 t^{88}-3 q^5 t^{89}-7 q^6 t^{89}+3 q^7 t^{89}
+7 q^8 t^{89}-4 q^6 t^{90}-4 q^7 t^{90}-2 q^6 t^{91}-2 q^7 t^{91}
+2 q^8 t^{91}+q^6 t^{92}+q^7 t^{92}-2 q^8 t^{92}-q^5 t^{93}
-4 q^6 t^{93}-q^7 t^{93}+2 q^8 t^{93}-2 q^6 t^{94}
+5 q^8 t^{94}-q^9 t^{94}-3 q^6 t^{95}-7 q^7 t^{95}
+q^8 t^{95}+q^9 t^{95}-3 q^7 t^{96}-q^8 t^{96}
+q^5 t^{97}-3 q^7 t^{97}-3 q^8 t^{97}+q^6 t^{98}
+2 q^7 t^{98}+6 q^8 t^{98}-q^9 t^{98}-q^6 t^{99}
-5 q^7 t^{99}-q^8 t^{99}-q^9 t^{99}-2 q^7 t^{100}
+2 q^8 t^{100}+2 q^9 t^{100}-3 q^7 t^{101}-4 q^8 t^{101}
-q^6 t^{102}-2 q^7 t^{102}+3 q^8 t^{102}+q^5 t^{103}
+2 q^6 t^{103}-2 q^7 t^{103}-4 q^8 t^{103}-2 q^9 t^{103}
+q^6 t^{104}+q^7 t^{104}+2 q^8 t^{104}+3 q^9 t^{104}
-q^7 t^{105}-4 q^8 t^{105}-q^9 t^{105}-q^7 t^{106}
+q^9 t^{106}+2 q^6 t^{107}+q^7 t^{107}-q^8 t^{107}
-3 q^9 t^{107}-q^7 t^{108}-2 q^8 t^{108}+3 q^9 t^{108}
+q^6 t^{109}+3 q^7 t^{109}-2 q^8 t^{109}-q^9 t^{109}
+q^7 t^{110}+q^8 t^{110}+q^9 t^{110}+q^8 t^{111}
-q^9 t^{111}-2 q^7 t^{112}-5 q^8 t^{112}+q^9 t^{112}
+q^6 t^{113}+4 q^7 t^{113}+q^9 t^{113}+q^7 t^{114}
-q^8 t^{114}-2 q^9 t^{114}+q^7 t^{115}+2 q^8 t^{115}
-q^8 t^{116}-q^9 t^{116}+q^7 t^{117}+q^8 t^{117}
+5 q^9 t^{117}-3 q^8 t^{118}-5 q^9 t^{118}
+q^7 t^{119}+4 q^8 t^{119}+q^9 t^{119}-q^9 t^{120}
+q^8 t^{121}+3 q^9 t^{121}-q^8 t^{122}-3 q^9 t^{122}
-q^{10} t^{122}+2 q^8 t^{123}+q^9 t^{123}
+2 q^{10} t^{123}-2 q^9 t^{124}-q^{10} t^{124}
+q^8 t^{125}+2 q^9 t^{125}+q^8 t^{126}-q^{10} t^{126}
-q^7 t^{127}+2 q^{10} t^{127}-q^9 t^{128}-q^{10} t^{128}
+q^9 t^{129}-2 q^9 t^{131}+q^8 t^{132}+2 q^9 t^{132}
-q^8 t^{133}+q^9 t^{136}+q^{10} t^{136}
-q^8 t^{137}-q^9 t^{137}-2 q^{10} t^{137}
+q^9 t^{138}+q^{10} t^{138}-q^9 t^{139}
+q^{10} t^{140}-2 q^{10} t^{141}+q^9 t^{142}
+q^{10} t^{142}-q^9 t^{143}-q^{11} t^{145}+q^{11} t^{146}.
\)
}
\renewcommand{\baselinestretch}{1.2}

\subsubsection{F4} 
The quasi-minuscule weight is $\om_4$ in this case.
For the trefoil, the formula is:

$$\tilde{\j}^{F_4}_{3,2}(\om_4\,;\,q)=$$
\renewcommand{\baselinestretch}{0.5} 
{\small
\noindent
\( 
1+q^2-q^4+q^5+q^8-2 q^9+2 q^{10}-q^{12}+2 q^{13}-2 q^{14}
+q^{15}+q^{16}-2 q^{17}+2 q^{18}-2 q^{19}-q^{20}+2 q^{21}
-3 q^{22}+2 q^{23}-q^{24}-2 q^{25}+3 q^{26}-2 q^{27}+q^{29}
-2 q^{30}+2 q^{31}+q^{34}-2 q^{35}+q^{36}-q^{38}+q^{39}.
\)
}
\renewcommand{\baselinestretch}{1.2} 

It corresponds to the $J\!D$\~polynomial
$$\tilde{J\!D}^{F_4}_{3,2}(\om_4\,;\,q,t,r)=$$

\renewcommand{\baselinestretch}{0.5} 
{\small
\noindent
\( 
1+q t+q r t^2+q r^2 t^2-q t^3+q r^3 t^3-q r t^4
+q^2 r t^4-q r^2 t^4+q^2 r^2 t^4
+q^2 r^3 t^4-q^2 r t^5-q^2 r^2 t^5-q r^3 t^5+q^2 r^3 t^5
+q^3 r^3 t^5+q^2 r^4 t^5+q^2 r^5 t^5-q^2 r^2 t^6
-3 q^2 r^3 t^6-q^2 r^4 t^6+q^2 r t^7+q^2 r^2 t^7
-2 q^2 r^4 t^7+q^3 r^4 t^7-2 q^2 r^5 t^7
+q^3 r^5 t^7+2 q^2 r^3 t^8-q^3 r^3 t^8
-q^3 r^4 t^8-q^3 r^5 t^8-q^2 r^6 t^8+q^3 r^6 t^8
+q^3 r^3 t^9+q^2 r^4 t^9-q^3 r^4 t^9+q^2 r^5 t^9
-2 q^3 r^5 t^9-2 q^3 r^6 t^9-q^3 r^7 t^9
+2 q^3 r^4 t^{10}+2 q^3 r^5 t^{10}-q^3 r^6 t^{10}
+q^4 r^6 t^{10}-q^3 r^7 t^{10}-q^3 r^8 t^{10}
-q^3 r^3 t^{11}+q^3 r^5 t^{11}+3 q^3 r^6 t^{11}
-q^4 r^6 t^{11}+q^3 r^7 t^{11}-q^3 r^4 t^{12}
-q^3 r^5 t^{12}+q^3 r^7 t^{12}-q^4 r^7 t^{12}
+q^3 r^8 t^{12}-q^4 r^8 t^{12}-q^3 r^6 t^{13}
+q^4 r^6 t^{13}+q^4 r^7 t^{13}
+q^4 r^8 t^{13}-q^4 r^9 t^{13}-q^4 r^6 t^{14}
+q^4 r^8 t^{14}+q^4 r^9 t^{14}
+q^4 r^{10} t^{14}-q^4 r^7 t^{15}
-q^4 r^8 t^{15}-q^5 r^9 t^{15}
+q^4 r^6 t^{16}-q^4 r^9 t^{16}+q^5 r^9 t^{16}.
\)
}
\renewcommand{\baselinestretch}{1.2} 

Recall that 
$\tilde{\j}^{F_4}_{r,s}(\om_4\,;\,q)\ =\ 
\tilde{J\!D}^{F_4}_{r,s}(\om_4\,;\,q,q,q^2).$
\medskip

The last DAHA generated QG-Jones polynomial
and the corresponding $J\!D$\~polynomial
we provide will be:

$$\tilde{\j}^{F_4}_{4,3}(\om_4\,;\,q)=$$
\renewcommand{\baselinestretch}{0.5} 
{\small
\noindent
\( 
1+q^2+q^3+q^6+q^8+q^{10}+q^{13}+q^{14}-q^{15}+q^{16}+q^{19}
-2 q^{20}-2 q^{23}+q^{24}-q^{25}-2 q^{26}+q^{27}-3 q^{28}
-3 q^{31}+q^{32}-q^{33}-2 q^{34}+3 q^{35}-3 q^{36}+q^{38}
-2 q^{39}+2 q^{40}-q^{42}+3 q^{43}-2 q^{44}+q^{45}+2 q^{46}
-2 q^{47}+2 q^{48}-q^{50}+3 q^{51}-q^{52}+2 q^{54}-2 q^{55}
+2 q^{56}-q^{58}+q^{59}-q^{60}-2 q^{63}+q^{64}-q^{66}+q^{67},
\)
}
\renewcommand{\baselinestretch}{1.2} 

corresponding to 

$$\tilde{J\!D}^{F_4}_{4,3}(\om_4\,;\,q,t,r)=$$

\renewcommand{\baselinestretch}{0.5} 
{\small
\noindent
\( 
1+q t+q^2 t+q^2 t^2+q r t^2+q^2 r t^2+q r^2 t^2+q^2 r^2 t^2
-q t^3-q^2 t^3+q^3 t^3+q^2 r t^3+q^3 r t^3+q^2 r^2 t^3
+q^3 r^2 t^3+q r^3 t^3+q^2 r^3 t^3-q^2 t^4-q^3 t^4-q r t^4
+q^3 r t^4+q^4 r t^4-q r^2 t^4+2 q^3 r^2 t^4+q^4 r^2 t^4
+2 q^2 r^3 t^4+3 q^3 r^3 t^4+q^4 r^3 t^4+q^3 r^4 t^4-q^3 t^5
-2 q^2 r t^5-2 q^3 r t^5-2 q^2 r^2 t^5-2 q^3 r^2 t^5
-q r^3 t^5+3 q^3 r^3 t^5+2 q^4 r^3 t^5+q^5 r^3 t^5+q^2 r^4 t^5
+2 q^3 r^4 t^5+q^4 r^4 t^5+q^2 r^5 t^5+2 q^3 r^5 t^5+q^4 r^5 t^5
+q^3 t^6-2 q^3 r t^6-q^2 r^2 t^6-3 q^3 r^2 t^6-4 q^2 r^3 t^6
-5 q^3 r^3 t^6+2 q^4 r^3 t^6+q^5 r^3 t^6+q^6 r^3 t^6-q^2 r^4 t^6
+2 q^4 r^4 t^6+q^5 r^4 t^6+q^3 r^5 t^6+2 q^4 r^5 t^6+q^5 r^5 t^6
+q^3 r^6 t^6+q^2 r t^7+q^3 r t^7-q^4 r t^7+q^2 r^2 t^7
-2 q^4 r^2 t^7-5 q^3 r^3 t^7-4 q^4 r^3 t^7-2 q^2 r^4 t^7
-3 q^3 r^4 t^7+q^4 r^4 t^7+2 q^5 r^4 t^7+q^6 r^4 t^7
-2 q^2 r^5 t^7-2 q^3 r^5 t^7+3 q^4 r^5 t^7+3 q^5 r^5 t^7
+q^6 r^5 t^7+q^3 r^6 t^7+3 q^4 r^6 t^7+2 q^5 r^6 t^7+q^4 r^7 t^7
+q^5 r^7 t^7+q^3 r t^8+q^4 r t^8+q^3 r^2 t^8+2 q^2 r^3 t^8
+q^3 r^3 t^8-6 q^4 r^3 t^8-3 q^5 r^3 t^8-4 q^3 r^4 t^8
-8 q^4 r^4 t^8-q^5 r^4 t^8+q^6 r^4 t^8-4 q^3 r^5 t^8-7 q^4 r^5 t^8
-q^5 r^5 t^8+q^6 r^5 t^8-q^2 r^6 t^8-q^3 r^6 t^8+q^4 r^6 t^8
+3 q^5 r^6 t^8+q^6 r^6 t^8+q^4 r^7 t^8+q^5 r^7 t^8+q^4 r^8 t^8
+q^5 r^8 t^8+q^4 r t^9+q^3 r^2 t^9+2 q^4 r^2 t^9+4 q^3 r^3 t^9
+5 q^4 r^3 t^9-q^6 r^3 t^9+q^2 r^4 t^9+q^3 r^4 t^9-4 q^4 r^4 t^9
-3 q^5 r^4 t^9-q^6 r^4 t^9+q^7 r^4 t^9+q^2 r^5 t^9-q^3 r^5 t^9
-8 q^4 r^5 t^9-4 q^5 r^5 t^9+q^7 r^5 t^9-4 q^3 r^6 t^9
-8 q^4 r^6 t^9-q^5 r^6 t^9+4 q^6 r^6 t^9+q^7 r^6 t^9-q^3 r^7 t^9
-3 q^4 r^7 t^9+2 q^6 r^7 t^9+q^5 r^8 t^9+q^6 r^8 t^9-q^4 r t^{10}
+4 q^4 r^3 t^{10}+q^5 r^3 t^{10}+q^6 r^3 t^{10}+4 q^3 r^4 t^{10}
+6 q^4 r^4 t^{10}-3 q^5 r^4 t^{10}-2 q^6 r^4 t^{10}-q^7 r^4 t^{10}
+4 q^3 r^5 t^{10}+4 q^4 r^5 t^{10}-5 q^5 r^5 t^{10}
-4 q^6 r^5 t^{10}-q^7 r^5 t^{10}-6 q^4 r^6 t^{10}-8 q^5 r^6 t^{10}
-2 q^6 r^6 t^{10}+q^7 r^6 t^{10}+q^8 r^6 t^{10}-q^3 r^7 t^{10}
-4 q^4 r^7 t^{10}-2 q^5 r^7 t^{10}+q^7 r^7 t^{10}-q^3 r^8 t^{10}
-3 q^4 r^8 t^{10}+3 q^6 r^8 t^{10}+q^7 r^8 t^{10}+q^5 r^9 t^{10}
+q^6 r^9 t^{10}+q^6 r^{10} t^{10}-q^3 r^3 t^{11}-2 q^4 r^3 t^{11}
+q^5 r^3 t^{11}+q^6 r^3 t^{11}+5 q^4 r^4 t^{11}+4 q^5 r^4 t^{11}
+q^6 r^4 t^{11}+q^3 r^5 t^{11}+6 q^4 r^5 t^{11}+2 q^5 r^5 t^{11}
+4 q^3 r^6 t^{11}+5 q^4 r^6 t^{11}-8 q^5 r^6 t^{11}-8 q^6 r^6 t^{11}
+q^3 r^7 t^{11}-q^4 r^7 t^{11}-9 q^5 r^7 t^{11}-5 q^6 r^7 t^{11}
+q^7 r^7 t^{11}-2 q^4 r^8 t^{11}-7 q^5 r^8 t^{11}-4 q^6 r^8 t^{11}
+q^7 r^8 t^{11}-q^4 r^9 t^{11}-q^5 r^9 t^{11}-q^4 r^3 t^{12}
-q^5 r^3 t^{12}-q^6 r^3 t^{12}-q^3 r^4 t^{12}-q^4 r^4 t^{12}
+4 q^5 r^4 t^{12}+2 q^6 r^4 t^{12}-q^3 r^5 t^{12}+q^4 r^5 t^{12}
+8 q^5 r^5 t^{12}+2 q^6 r^5 t^{12}-q^7 r^5 t^{12}+8 q^4 r^6 t^{12}
+11 q^5 r^6 t^{12}-4 q^7 r^6 t^{12}+q^3 r^7 t^{12}+3 q^4 r^7 t^{12}
-5 q^6 r^7 t^{12}-q^7 r^7 t^{12}+q^8 r^7 t^{12}+q^3 r^8 t^{12}
+q^4 r^8 t^{12}-6 q^5 r^8 t^{12}-7 q^6 r^8 t^{12}+q^7 r^8 t^{12}
+q^8 r^8 t^{12}-q^4 r^9 t^{12}-5 q^5 r^9 t^{12}-3 q^6 r^9 t^{12}
+3 q^7 r^9 t^{12}-2 q^5 r^{10} t^{12}-2 q^6 r^{10} t^{12}
+q^7 r^{10} t^{12}-q^5 r^3 t^{13}-2 q^4 r^4 t^{13}-3 q^5 r^4 t^{13}
+q^7 r^4 t^{13}-2 q^4 r^5 t^{13}-q^5 r^5 t^{13}+q^6 r^5 t^{13}
+2 q^7 r^5 t^{13}-q^3 r^6 t^{13}+10 q^5 r^6 t^{13}+7 q^6 r^6 t^{13}
-q^8 r^6 t^{13}+4 q^4 r^7 t^{13}+11 q^5 r^7 t^{13}-6 q^7 r^7 t^{13}
-q^8 r^7 t^{13}+4 q^4 r^8 t^{13}+9 q^5 r^8 t^{13}-2 q^6 r^8 t^{13}
-9 q^7 r^8 t^{13}-q^8 r^8 t^{13}-q^5 r^9 t^{13}-6 q^6 r^9 t^{13}
-3 q^7 r^9 t^{13}+q^8 r^9 t^{13}-q^5 r^{10} t^{13}
-2 q^6 r^{10} t^{13}-q^7 r^{10} t^{13}-q^5 r^{11} t^{13}
-q^6 r^{11} t^{13}+q^7 r^{11} t^{13}+q^5 r^3 t^{14}
-2 q^5 r^4 t^{14}-q^6 r^4 t^{14}-q^7 r^4 t^{14}-q^4 r^5 t^{14}
-3 q^5 r^5 t^{14}-4 q^4 r^6 t^{14}-6 q^5 r^6 t^{14}
+4 q^6 r^6 t^{14}+6 q^7 r^6 t^{14}+q^8 r^6 t^{14}-q^4 r^7 t^{14}
+5 q^5 r^7 t^{14}+8 q^6 r^7 t^{14}+4 q^7 r^7 t^{14}+q^4 r^8 t^{14}
+8 q^5 r^8 t^{14}+7 q^6 r^8 t^{14}+2 q^4 r^9 t^{14}
+6 q^5 r^9 t^{14}-2 q^6 r^9 t^{14}-8 q^7 r^9 t^{14}
+q^4 r^{10} t^{14}+2 q^5 r^{10} t^{14}-2 q^6 r^{10} t^{14}
-6 q^7 r^{10} t^{14}-2 q^6 r^{11} t^{14}-3 q^7 r^{11} t^{14}
-q^7 r^{12} t^{14}+q^5 r^4 t^{15}-q^7 r^4 t^{15}-q^6 r^5 t^{15}
-2 q^7 r^5 t^{15}-5 q^5 r^6 t^{15}-5 q^6 r^6 t^{15}
-3 q^7 r^6 t^{15}-2 q^4 r^7 t^{15}-4 q^5 r^7 t^{15}
+5 q^6 r^7 t^{15}+5 q^7 r^7 t^{15}-q^8 r^7 t^{15}-2 q^4 r^8 t^{15}
-2 q^5 r^8 t^{15}+10 q^6 r^8 t^{15}+9 q^7 r^8 t^{15}
-2 q^8 r^8 t^{15}+4 q^5 r^9 t^{15}+11 q^6 r^9 t^{15}
+2 q^7 r^9 t^{15}-3 q^8 r^9 t^{15}+q^9 r^9 t^{15}
+2 q^5 r^{10} t^{15}+3 q^6 r^{10} t^{15}+q^5 r^{11} t^{15}
-q^6 r^{11} t^{15}-4 q^7 r^{11} t^{15}+q^8 r^{11} t^{15}
-q^6 r^{12} t^{15}-q^7 r^{12} t^{15}-q^7 r^{13} t^{15}
+q^7 r^4 t^{16}+q^7 r^5 t^{16}+q^4 r^6 t^{16}+q^5 r^6 t^{16}
-3 q^6 r^6 t^{16}-3 q^7 r^6 t^{16}-4 q^5 r^7 t^{16}
-6 q^6 r^7 t^{16}-q^7 r^7 t^{16}+2 q^8 r^7 t^{16}
-4 q^5 r^8 t^{16}-4 q^6 r^8 t^{16}+q^7 r^8 t^{16}
+2 q^8 r^8 t^{16}-q^4 r^9 t^{16}-q^5 r^9 t^{16}
+8 q^6 r^9 t^{16}+11 q^7 r^9 t^{16}-4 q^8 r^9 t^{16}
-q^9 r^9 t^{16}+q^5 r^{10} t^{16}+8 q^6 r^{10} t^{16}
+5 q^7 r^{10} t^{16}-5 q^8 r^{10} t^{16}+q^5 r^{11} t^{16}
+6 q^6 r^{11} t^{16}+3 q^7 r^{11} t^{16}-5 q^8 r^{11} t^{16}
+q^6 r^{12} t^{16}-q^8 r^{12} t^{16}+q^5 r^6 t^{17}
+2 q^6 r^6 t^{17}+q^7 r^6 t^{17}-q^8 r^6 t^{17}-4 q^6 r^7 t^{17}
-4 q^7 r^7 t^{17}-q^5 r^8 t^{17}-7 q^6 r^8 t^{17}-3 q^7 r^8 t^{17}
+3 q^8 r^8 t^{17}-4 q^5 r^9 t^{17}-8 q^6 r^9 t^{17}
+3 q^7 r^9 t^{17}+9 q^8 r^9 t^{17}-q^5 r^{10} t^{17}
-q^6 r^{10} t^{17}+6 q^7 r^{10} t^{17}+3 q^8 r^{10} t^{17}
+3 q^6 r^{11} t^{17}+6 q^7 r^{11} t^{17}-2 q^8 r^{11} t^{17}
+2 q^6 r^{12} t^{17}+q^7 r^{12} t^{17}-3 q^8 r^{12} t^{17}
+q^6 r^{13} t^{17}+q^7 r^{13} t^{17}-2 q^8 r^{13} t^{17}
+q^6 r^6 t^{18}+q^7 r^6 t^{18}+q^8 r^6 t^{18}+q^5 r^7 t^{18}
+2 q^6 r^7 t^{18}-q^7 r^7 t^{18}-3 q^8 r^7 t^{18}+q^5 r^8 t^{18}
+q^6 r^8 t^{18}-3 q^7 r^8 t^{18}-5 q^8 r^8 t^{18}-6 q^6 r^9 t^{18}
-8 q^7 r^9 t^{18}-q^5 r^{10} t^{18}-4 q^6 r^{10} t^{18}
+q^7 r^{10} t^{18}+6 q^8 r^{10} t^{18}-q^9 r^{10} t^{18}
-q^5 r^{11} t^{18}-3 q^6 r^{11} t^{18}+4 q^7 r^{11} t^{18}
+11 q^8 r^{11} t^{18}-q^9 r^{11} t^{18}+4 q^7 r^{12} t^{18}
+3 q^8 r^{12} t^{18}+q^7 r^{13} t^{18}+2 q^8 r^{13} t^{18}
-q^8 r^{14} t^{18}-q^6 r^6 t^{19}+q^8 r^6 t^{19}+q^6 r^7 t^{19}
+2 q^7 r^7 t^{19}+2 q^8 r^7 t^{19}+q^6 r^8 t^{19}+q^7 r^8 t^{19}
+q^8 r^8 t^{19}+q^5 r^9 t^{19}+2 q^6 r^9 t^{19}-5 q^7 r^9 t^{19}
-8 q^8 r^9 t^{19}+q^9 r^9 t^{19}-2 q^6 r^{10} t^{19}
-7 q^7 r^{10} t^{19}-3 q^8 r^{10} t^{19}-2 q^6 r^{11} t^{19}
-5 q^7 r^{11} t^{19}-q^9 r^{11} t^{19}-q^6 r^{12} t^{19}
+2 q^7 r^{12} t^{19}+5 q^8 r^{12} t^{19}-3 q^9 r^{12} t^{19}
+2 q^7 r^{13} t^{19}+4 q^8 r^{13} t^{19}-q^9 r^{13} t^{19}
+q^7 r^{14} t^{19}+2 q^8 r^{14} t^{19}+q^8 r^{15} t^{19}
-q^8 r^6 t^{20}+q^8 r^7 t^{20}+q^7 r^8 t^{20}+2 q^8 r^8 t^{20}
+2 q^6 r^9 t^{20}+4 q^7 r^9 t^{20}+2 q^8 r^9 t^{20}-q^9 r^9 t^{20}
-2 q^7 r^{10} t^{20}-4 q^8 r^{10} t^{20}+3 q^9 r^{10} t^{20}
-5 q^7 r^{11} t^{20}-6 q^8 r^{11} t^{20}+5 q^9 r^{11} t^{20}
-q^6 r^{12} t^{20}-5 q^7 r^{12} t^{20}-q^8 r^{12} t^{20}
+3 q^9 r^{12} t^{20}-2 q^7 r^{13} t^{20}+q^8 r^{14} t^{20}
-2 q^9 r^{14} t^{20}-q^8 r^7 t^{21}-q^8 r^8 t^{21}
+2 q^7 r^9 t^{21}+3 q^8 r^9 t^{21}-q^9 r^9 t^{21}
+q^6 r^{10} t^{21}+2 q^7 r^{10} t^{21}-q^8 r^{10} t^{21}
-3 q^9 r^{10} t^{21}+q^6 r^{11} t^{21}+2 q^7 r^{11} t^{21}
-2 q^8 r^{11} t^{21}-3 q^9 r^{11} t^{21}-2 q^7 r^{12} t^{21}
-5 q^8 r^{12} t^{21}+4 q^9 r^{12} t^{21}-q^7 r^{13} t^{21}
-q^8 r^{13} t^{21}+4 q^9 r^{13} t^{21}-q^7 r^{14} t^{21}
+4 q^9 r^{14} t^{21}+q^9 r^{15} t^{21}-q^7 r^9 t^{22}
+2 q^9 r^9 t^{22}+q^7 r^{10} t^{22}+3 q^8 r^{10} t^{22}
+q^7 r^{11} t^{22}+2 q^8 r^{11} t^{22}-2 q^9 r^{11} t^{22}
+q^7 r^{12} t^{22}-3 q^8 r^{12} t^{22}-6 q^9 r^{12} t^{22}
-3 q^8 r^{13} t^{22}-3 q^9 r^{13} t^{22}-q^{10} r^{13} t^{22}
-2 q^8 r^{14} t^{22}-q^{10} r^{14} t^{22}+q^9 r^{16} t^{22}
-q^8 r^9 t^{23}-q^9 r^9 t^{23}+q^8 r^{10} t^{23}+2 q^9 r^{10} t^{23}
+2 q^8 r^{11} t^{23}+3 q^9 r^{11} t^{23}+q^7 r^{12} t^{23}
+4 q^8 r^{12} t^{23}+q^{10} r^{12} t^{23}-q^9 r^{13} t^{23}
+2 q^{10} r^{13} t^{23}-q^8 r^{14} t^{23}-3 q^9 r^{14} t^{23}
+2 q^{10} r^{14} t^{23}-q^8 r^{15} t^{23}-q^9 r^{15} t^{23}
-q^{10} r^{15} t^{23}-q^9 r^{16} t^{23}-q^9 r^9 t^{24}
-q^8 r^{10} t^{24}-q^9 r^{10} t^{24}-q^8 r^{11} t^{24}
-q^9 r^{11} t^{24}+q^8 r^{12} t^{24}+3 q^9 r^{12} t^{24}
-2 q^{10} r^{12} t^{24}+q^8 r^{13} t^{24}-q^{10} r^{13} t^{24}
+q^8 r^{14} t^{24}-q^9 r^{15} t^{24}+2 q^{10} r^{15} t^{24}
+q^{10} r^{16} t^{24}+q^9 r^9 t^{25}+q^{10} r^{12} t^{25}
+q^9 r^{13} t^{25}-q^{10} r^{13} t^{25}+q^9 r^{14} t^{25}
-2 q^{10} r^{14} t^{25}-q^{10} r^{15} t^{25}-q^{11} r^{15} t^{25}
-q^{10} r^{16} t^{25}-q^9 r^{12} t^{26}+q^{10} r^{12} t^{26}
+q^9 r^{13} t^{26}+q^{10} r^{13} t^{26}+q^9 r^{14} t^{26}
+q^{10} r^{14} t^{26}+q^9 r^{15} t^{26}-q^{10} r^{15} t^{26}
+q^{11} r^{15} t^{26}-q^{10} r^{12} t^{27}
-q^9 r^{13} t^{27}-q^9 r^{14} t^{27}+q^{10} r^{15} t^{27}.
\)
\vfill\eject
}
\renewcommand{\baselinestretch}{1.2}


In conclusion, we want to comment on explicit numerical
formulas (sometimes very long) in this paper. They are 
needed to justify our conjectures
and can potentially help to discover/confirm further properties 
of the super and hyper- polynomials, which have extremely 
rich symmetries. They can also provide support for
ongoing research toward the Khovanov-Rozansky polynomials, 
including incorporating colors into this theory and its  
$C^\vee C$\~generalization. 
The IMRN version of this paper is equivalent to this preprint,
but contains a reduced number of numerical formulas.

\bibliographystyle{unsrt}

\end{document}